\documentclass[10pt]{amsart}
\pdfoutput=1
\usepackage{ mathrsfs }
\usepackage{mathtools}
\usepackage{amsmath, amsthm, amssymb,slashed,stmaryrd}
\usepackage{ifpdf}
\usepackage[pdftex]{graphicx}
\usepackage{tikz}
\usetikzlibrary{decorations.markings}
\newcommand*{\halfway}{0.5*\pgfdecoratedpathlength+.5*3pt}
\usetikzlibrary{matrix,arrows,calc}
\usepackage[all]{xy}
\UseRawInputEncoding

\usepackage[pdftex,plainpages=false,hypertexnames=false,pdfpagelabels]{hyperref}
\usepackage{amsthm}

\theoremstyle{definition}
 \newtheorem{thm}{Theorem}[section]
    
 \newtheorem{cor}[thm]{Corollary}
 \newtheorem{lem}[thm]{Lemma}
  \newtheorem{hyp}[thm]{Hypothesis}

 \newtheorem{prop}[thm]{Proposition}
 
 \newtheorem{defn}[thm]{Definition}
 \newtheorem{notation}[thm]{Notation}
 \newtheorem{ex}[thm]{Example}
  
 \newtheorem*{thm*}{Theorem}
 \theoremstyle{remark}
 \newtheorem{rmk}[thm]{Remark}

\def\beq{\begin{eqnarray}}
\def\eeq{\end{eqnarray}}
 \newcommand{\bp}{\begin{proof}[Proof]}
 \newcommand{\ep}{\end{proof}}

\DeclareSymbolFont{bbold}{U}{bbold}{m}{n}
\DeclareSymbolFontAlphabet{\mathbbold}{bbold}
\def\one{\mathbbold{1}}

\def\Pf{{\sf{Pf}}}
\def\inv{{\rm inv}}
\def\u{{\sf u}}
\def\s{{\sf s}}
\def\t{{\sf t}}
\def\m{{\sf c}}
\def\c{{\sf c}}
\def\sL{{\sf L}}
\def\sR{{\sf R}}
\def\sI{{\sf I}}
\def\sS{{\sf S}}

\def\proj{{\rm proj}}
\def\act{{\rm act}}

\def\fl{{\sf fl}}
\def\Fl{{\sf Fl}}

\def\M{{\mathbb{M}}}
\def\V{{\mathcal{V}}}
\def\SMfld{{\sf SMfld}}
\def\Mfld{{\sf Mfld}}
\def\inn{{\rm in}}
\def\Path{{\sf Path}}
\def\sP{{\sf sPath}}

\def\Or{{\sf Or}}
\def\orr{{\sf or}}
\def\rr{{\sf r}}
\def\RR{{\sf R}}
\def\RRR{{\mathscr R}}

\def\out{{\rm out}}
\def\sTr{{\sf sTr}}
\def\Tr{{\sf Tr}}
\def\deg{{\rm deg}}

\def\KO{{\rm KO}}
\def\Ind{{\rm Ind}}
\def\K{{\rm K}}
\def\Mod{{\sf Mod}}

\def\Ch{{\rm Ch}}

\def\dKO{\widehat{\rm KO}{}}
\def\Cl{{\rm Cl}}
\def\cCl{\C{\rm l}}

\def\twist{\mathscr{T}}

\def\Ob{{\rm Ob}}
\def\Mor{{\rm Mor}}
\def\cl{{\rm cl}}

\def\H{{\rm H}}
\def\HH{{\mathbb H}}
\def\Spin{{\rm Spin}}

\def\SO{{\rm SO}}

\def\Map{{\sf Map}}
\def\RG{{\sf RG}}
\def\rg{{\sf rg}}

\def\Vect{{\sf Vect}}

\def\Fun{{\sf Fun}}

\def\EFT{ \hbox{-{\sf EFT}}}
\def\eft{\hbox{-{\sf eft}}}
\def\Bord{\hbox{-{\sf Bord}}}
\def\EBord{\hbox{-}{\sf EBord}}
\def\ebord{\hbox{-{\sf ebord}}}
\def\TMF{{\rm TMF}}

\def\pt{{\rm pt}}
\def\spt{{\rm spt}}
\def\ev{{\rm ev}}

\def\eva{{\rm act}}
\def\odd{{\rm odd}}

\def\bS{{\mathbb{S}}}
\def\A{{\mathbb{A}}}

\def\R{{\mathbb{R}}}
\def\red{{\rm{red}}}

\def\TA{{\sf {Alg}}}
\def\TV{{\sf {Vect}}}

\def\id{{{\rm id}}}
\def\C{{\mathbb{C}}}

\def\Z{{\mathbb{Z}}}

\def\End{{\rm End}}

\def\Rep{{\rm Rep}}

\newcommand{\op}{{\sf{op}}}   
\newcommand{\sq}{/\!\!/}

\def\twocommute{\ensuremath{\rotatebox[origin=c]{30}{$\Rightarrow$}}}

\newcommand\nc{\newcommand}
\nc\mf\mathfrak
\nc\mc\mathcal
\nc\mb\mathbb

\begin{document}

\title{The families analytic index for $1|1$-dimensional Euclidean field theories}
\author{Daniel Berwick-Evans}

\begin{abstract}
We construct a KO-valued families index for a class of $1|1$-dimensional Euclidean field theories. This realizes a conjectured cocycle map in the Stolz--Teichner program. We further show that a bundle of spin manifolds leads to a family of partially-defined $1|1$-Euclidean field theories, yielding a cocycle refinement of the families analytic index. The methods are chosen with the goal of generalizing to $2|1$-dimensional field theories, where analogous structures are expected to provide an analytic index valued in topological modular forms. 
\end{abstract}
\date{\today}

\maketitle 
%
%\setcounter{tocdepth}{1}
%\tableofcontents

\section{Introduction and statement of results}

Stolz and Teichner have outlined constructions of categories of $1|1$- and $2|1$-dimensional Euclidean field theories over a smooth manifold $M$ indexed by a degree $n\in \Z$~\cite{ST04,HST,ST11}. 
%(see~\S\ref{sec:ebord} below for a review)
%\beq
%d|1\EFT(M)=\Fun^{\otimes}(d|1\EBord(M),\TV),\qquad d=1,2 \label{eq:defnFT}
%\eeq
%as symmetric monoidal functors from certain bordism categories $d|1\EBord(M)$ to a category of topological $\C$-vector spaces $\TV$. 
%%There are various subtle points in this definition, e.g., the flavor of topological vector space, as well as certain unitary conditions; see~\S\ref{sec:STFT} below. 
They conjecture that these categories determine geometric cocycle models for real $\K$-theory and topological modular forms~\cite[\S1.5-1.6]{ST11},
%\beq
%1|1\EFT^n(M)\dashrightarrow \KO^n(M),\qquad 2|1\EFT^n(M)\to \TMF^n(M)
%\eeq
\beq
&&\begin{tikzpicture}[baseline=(basepoint)];
\node (A) at (0,0) {$\left\{\begin{array}{c} d|1\hbox{-}{\rm dimensional} \\ {\rm Euclidean\ field\ theories} \\ {\rm of \ degree} \ n \ { \rm over} \ M\end{array} \right\}$};
\node (B) at (7,0) {$\left\{\begin{array}{ll} \KO^{n}(M) & d=1\\ \TMF^{n}(M) & d=2.\end{array}\right.$};
%\node (C) at (6,0) {$\left\{\begin{array}{c} 2|1\hbox{-}{\rm Euclidean\ degree\ } n \\ {\rm field\ theories \ over} \ M\end{array} \right\}$};
\draw[->>,dashed] (A) to node [above] {\small cocycle}  (B);
%\draw[->>,dashed] (C) to node [above] {\small cocycle} (D);
\path (0,-.05) coordinate (basepoint);
\end{tikzpicture}\label{eq:conjecture}
%&&1|1\EFT(M)\stackrel{{\rm cocycle}}{\dashrightarrow}  \K(M)  \qquad  2|1\EFT(M) \stackrel{{\rm cocycle}}{\dashrightarrow} \TMF(M)
\eeq
%where $n\in \Z$ parameterizes the \emph{degree} of a field theory, specializing to~\eqref{eq:defnFT} when~$n=0$. 
It is further expected that these cocycle models inherit an analytic pushforward from quantization of field theories. The ultimate goal is to use these pushforwards to prove a TMF-generalization of the Atiyah--Singer index theorem. 
%So far, no candidate cocycle maps for either K-theory or TMF have been proposed, nor has a construction of the pushforwards been suggested. 
%This is partly due to the fact that a complete definition of $d|1\EFT^n(M)$ has yet to appear. 
%To date, no candidate cocycle maps for either K-theory of TMF have been proposed. 
%We caution that a fully extended (i.e., 2-categorical) refinement of $2|1\EBord(M)$ is required to have a chance of recovering $\TMF$; see~\cite[Conjecture~1.17]{ST11}. 
%However, with the existing definitions from~\cite{ST11} one might hope to construct a cocycle map valued in an approximation to $\TMF$, just as 1-extended field theories approximate 2-extended ones. 

The feeling has always been that a proper understanding of the relationship between $1|1$-Euclidean field theories and the index theorem in K-theory will provide significant insight into the $2|1$-dimensional case~\cite[\S6]{Segal_Elliptic}, \cite[\S3]{ST04}. With this in mind, the goal of this paper is to formulate and prove a version of~\eqref{eq:conjecture} when $d=1$ that communicates with the families analytic index. We use a variation on the notion of degree~$n$ field theories over~$M$, denoted $1|1\eft^n(M)$, that incorporates reflection positivity and reality data as well as a condition demanding the existence of smooth cutoff theories; see~\S\ref{sec:framework} below. 
%In fact, our main result realizes cocycles for \emph{differential} KO-theory in the sense of Hopkins--Singer~\cite{hopsing}.
% that incorporates Fei Han's character map for $1|1$-Euclidean field theories, explained below.
% valued in closed differential forms~\cite{Han}. 

\begin{thm} \label{thm:cocycle} 
%There is a full subcategory $1|1\EFT^n_\eff(M)\subset 1|1\EFT^n(M)$ of field theories for which 
An object of $1|1\eft^n(M)$ determines a self-adjoint, Clifford linear superconnection $\A$ on a real vector bundle $V\to M$ (not necessarily finite-rank). The families index construction for self-adjoint superconnections determines a cocycle map 
\beq\label{eq:cocycleK}
1|1\eft^n(M)\xrightarrow{{\rm Index}} \dKO^n(M)
\eeq
%\beq
%&&\begin{tikzpicture}[baseline=(basepoint)];
%%\node (AA) at (-4,0) {$1|1\EFT_\eff^n(M)$};
%\node (A) at (0,0) {$1|1\EFT^n(M)$};
%\node (B) at (4,0) {$\dKO^n(M)$};
%%\draw[->] (AA) to node [above] {\small restrict}  (A);
%\draw[->] (A) to node [above] {\small Index}  (B);
%\path (0,-.05) coordinate (basepoint);
%\end{tikzpicture}
%%&&1|1\EFT(M)\stackrel{{\rm cocycle}}{\dashrightarrow}  \K(M)  \qquad  2|1\EFT(M) \stackrel{{\rm cocycle}}{\dashrightarrow} \TMF(M)
%\eeq
valued in the differential $\KO$-theory of $M$. If the renormalization group flow sends an object $E\in 1|1\eft^n(M)$ to the zero theory, then its underlying class in $\KO^n(M)$ is zero. 
\end{thm}

The categories $1|1\eft^n(M)$ also accommodate a refinement of the KO-valued families analytic index for a bundle of spin manifolds.
% for a bundle of spin manifolds. 

%To connect Theorem~\ref{thm:cocycle} with the families index theorem, we construct objects in $1|1\eft^n(M)$ from families of spin manifolds.

%This makes the category $1|1\eft^n_\eff(M)$ easier to study and examples easier to construct. 

\begin{thm} \label{thm:index}
The Bismut superconnection of a family of Riemannian spin manifolds $X\to M$ determines an object $\alpha(X)\in 1|1\eft^{-n}(M)$ whose image under~\eqref{eq:cocycleK} is the families differential analytic index of $X\to M$. 
\end{thm}

With the knowledge that the existing definitions of supersymmetric Euclidean field theory are unlikely to be the final word (see Remark~\ref{rmk:extenddown} and~\ref{rmk:badcat}), the cocycle map~\eqref{eq:cocycleK} is designed to be as robust as possible. 
%categories $1|1\eft^n(M)$ are designed to receive a map from a broad class of possible notions of $1|1$-Euclidean field theory.
Indeed, the categories $1|1\eft^n(M)$ and Theorem~\ref{thm:cocycle} come from evaluating a field theory on specific families of $1|1$-Euclidean bordisms, and the structures relevant to index theory are inherited from the geometry of these bordisms. The upshot is that any reasonable category of $1|1$-Euclidean field theories over~$M$ has a subcategory with a cocycle map~\eqref{eq:cocycleK}; see Remark~\ref{rmk:generalcocycle}. Perhaps more importantly, the categories $1|1\eft^n(M)$ have a generalization to the $d=2$ version of the conjecture~\eqref{eq:conjecture}, and the techniques behind Theorems~\ref{thm:cocycle} and~\ref{thm:index} lay the groundwork for this~$d=2$ sequel~\cite{DBEtorsion}.
% that extends the results of \cite{DBEChern}. 

%Corollary~\ref{cor:STcocycle} and Remark~\ref{rmk:concordance} explain how Theorems~\ref{thm:cocycle} and~\ref{thm:index} compare with Stolz and Teichner's conjectures relating K-theory and $1|1$-Euclidean field theories~\cite{ST11}. 

\subsection{The main idea} 

In the physics literature, a $1|1$-dimensional field theory is more commonly called \emph{$\mathcal{N}=1$ supersymmetric quantum mechanics.} The basic datum is a $\Z/2$-graded inner product space $\mathcal{H}=\mathcal{H}^\ev\oplus\mathcal{H}^\odd$ and an odd operator~$Q\colon \mathcal{H}^{\ev/\odd}\to \mathcal{H}^{\odd/\ev}$~\cite[\S1]{susymorse}. This determines a representation of the super semigroup 
\beq\label{eq:superrep}
\R^{1|1}_{\ge 0}\to \End(\mathcal{H}),\qquad (t,\theta)\mapsto e^{-tQ^2+\theta Q},
\eeq
that enhances the usual time-evolution semigroup in Wick-rotated quantum mechanics. Depending on the type of theory desired, one demands additional structure and property on~\eqref{eq:superrep}. For example, \emph{compact} theories require that~\eqref{eq:superrep} takes values in trace-class operators, and \emph{reflection-positive} theories require $Q$ to be self-adjoint. The main example motivating~\eqref{eq:conjecture} takes $Q=\slashed{D}$ to be the Dirac operator for a compact spin manifold~$X$~\cite[Example 3.2.9]{ST04}. This theory is both compact and reflection positive. 

The super semigroup $\R^{1|1}_{\ge 0}$ has a geometric interpretation as a moduli space of super intervals~\cite[\S3]{ST04} for which multiplication in the super semigroup corresponds to gluing of super intervals. Super intervals comprise a subcategory of the $1|1$-Euclidean bordism category~\cite[\S6]{HST},
\beq\label{eq:includesuperint}
&&\left\{\!\!\begin{array}{c}{\rm super} \\ {\rm intervals}\end{array}\!\!\right\}\xhookrightarrow{\iota} \left\{\!\!\begin{array}{c} 1|1\hbox{-}{\rm Euclidean} \\ {\rm bordisms}\end{array}\!\! \right\} \implies \left\{\!\!\begin{array}{c} 1|1\hbox{-}{\rm Euclidean} \\ {\rm field\ theories}\end{array}\!\!\right\}\xrightarrow{\iota^*} \Rep(\R^{1|1}_{\ge 0})
\eeq
and $1|1$-Euclidean field theories are defined as representations of this bordism category. In particular, a functor out of the $1|1$-Euclidean bordism category restricts along $\iota$ to a representation of $\R^{1|1}_{\ge 0}$. The representations~\eqref{eq:superrep} that come from this restriction have additional properties, e.g., compactness is automatic~\cite[\S2]{STTraces}. The initial definition of $1|1$-Euclidean field theory from~\cite[\S3]{ST04} simply equips representations of $\R^{1|1}_{\ge 0}$ with these expected additional properties. The more recent definition from~\cite{ST11} instead gives a precise notion of the $1|1$-Euclidean bordism category, and defines a $1|1$-Euclidean field theory as a representation of this bordism category. One upshot of the newer definition is a natural generalization to field theories over a smooth manifold, gotten by equipping $1|1$-Euclidean bordisms with a map to a (background) manifold~$M$, see Definition~\ref{defn:ogFT} below. 
%~\cite[Proposition 3.2.6]{ST04}
%On the other hand, reflection positivity is an additional structure on a representation of the bordism category~\cite{FreedHopkins}.  

Super intervals with a map to $M$ determine the category of \emph{superpaths} in $M$. This leads to a families-generalization of~\eqref{eq:includesuperint}
\beq\label{restrict1}
&&\resizebox{.9\textwidth}{!}{$
\left\{\!\!\!\begin{array}{c} {\rm super \ paths} \\ {\rm in} \ M \end{array}\!\!\!\right\}\xhookrightarrow{\iota} \left\{\!\!\!\begin{array}{c} 1|1\hbox{-}{\rm Euclidean} \\ {\rm bordisms\ over} \ M\end{array} \!\!\!\right\} \implies \left\{\!\!\!\begin{array}{c} 1|1\hbox{-}{\rm Euclidean\ field} \\ {\rm theories\ over} \ M \end{array}\!\!\!\right\}\xrightarrow{\iota^*} \Rep\left(\!\!\!\begin{array}{c} {\rm super \ paths} \\ {\rm in} \ M \end{array}\!\!\!\right),$}
\eeq
where a functor out of the $1|1$-Euclidean bordism category over $M$ determines a representation of the category of superpaths in $M$ by restriction. Examples of representations of the category of superpaths arise from \emph{super parallel transport} along a superconnection~\cite{florin}; a folklore theorem extends these representations to $1|1$-Euclidean field theories over $M$~\cite[\S1.3]{ST11}. 

The basic mechanism behind Theorem~\ref{thm:cocycle} is an inverse to this super parallel transport construction: a $1|1$-Euclidean field theory over $M$ determines a superconnection via restriction to \emph{nearly constant superpaths} in $M$, denoted $\sP_0(M)$. These are superpaths whose maps to $M$ factor through the projection $\R^{1|1}\to \R^{0|1}$. 
%This property implies that $\sP_0(M)$ is a \emph{super Lie category}, meaning its objects and morphisms are finite-dimensional supermanifolds. 

\begin{prop}[Proposition~\ref{prop:superconn}]\label{lem:mainA}
There is a functor 
$$
\Rep(\sP_0(M))\to \Vect^{\A}(M)
$$
from (not necessarily finite-rank) representations of $\sP_0(M)$ to the groupoid of (not necessarily finite-rank) super vector bundles with superconnection on~$M$. 
\end{prop}

%We emphasize that the super vector bundles underlying objects in $\Vect^{\A}(M)$ need not have finite rank. 

This leads to a natural candidate for the cocycle map~\eqref{eq:conjecture}
%Given a superconnection, one may attempt a families index construction (e.g., as in~\cite[Chapter~9]{BGV}). The composition 
\beq\label{eq:FTtosConn}
&&\left\{\!\!\!\begin{array}{c} 1|1\hbox{-}{\rm Euclidean\ field} \\ {\rm theories\ over} \ M \end{array}\!\!\!\right\}\xrightarrow{\iota^*} \Rep(\sP_0(M))\xrightarrow{{\rm Prop~\ref{lem:mainA}}}\Vect^{\A}(M)\stackrel{{\rm Index}}{\dashrightarrow} \KO^0(M)
\eeq
where the last arrow is the families index construction applied to a superconnection, e.g., see~\cite[Chapter~9]{BGV}. 

For an arbitrary superconnection, the index construction may not be valid---one needs certain additional properties. Some of these properties are automatic when a representation of $\sP_0(M)$ comes from the restriction of a $1|1$-Euclidean field theory. For example, the operators $e^{-t\A^2}$ are self-adjoint and trace class for $t>0$. However, these conditions alone do not guarantee the existence of an index bundle. The precise conditions on the superconnection can be complicated to state and depend sensitively on the analytical framework, e.g., compare Lott and Gorokhovsky's superconnections on Hilbert bundles~\cite{LottGorokhovsky}. 

From the field theory point of view, the additional property required to construct~\eqref{eq:FTtosConn} is quite straightforward: one needs the existence of \emph{smooth energy cutoffs}. Roughly, this means that a $1|1$-Euclidean field theory over~$M$ determines a compatible $\R_{>0}$-family of theories that discard states of energy greater than $\lambda\in \R_{>0}$. The states with energy~$<\lambda$ are used to construct the families index. The existence of smooth energy cutoffs (and hence the existence of a families index) is a \emph{property} of a field theory, not additional data. Furthermore, this property makes sense for Euclidean field theories in any super dimension, in particular dimension~$2|1$. The existence of smooth cutoffs is analogous to the ellipticity property for differential operators, see Definition~\ref{defn:effective1} and Remark~\ref{rmk:elliptic}. The categories $1|1\eft^n(M)$ are roughly defined as the value of a field a $1|1$-Euclidean field theory on nearly constant superpaths (see~\S\ref{sec:framework}) where in additional we require a smooth cutoff property; the following provides an index-theoretic characterization. 
% are defined in terms of the values of  satisfying certain additional properties described in~\S\ref{sec:framework}. 
% these bordisms are chosen so that the category encodes precisely the field-theoretic data required to construct the families index map~\eqref{eq:FTtosConn}.  
%In terms of a superconnection $\A$ extracted via~\eqref{eq:FTtosConn}, the condition of smooth energy cutoffs is basically the same as the condition that smooth index bundles exist, see Definition~\ref{defn:superconncutoff}. 
%Objects of the category $1|1\eft^n(M)$ are determined by superconnections with this smooth cutoff property. 
%Proposition~\ref{prop:geodata} summarizes the resulting bridge between $1|1$-Euclidean field theories and index theory.

\begin{prop}[Proposition~\ref{prop:dataofeft2}] \label{prop:dataofeft}
There is an equivalence of categories
\beq\label{eq:dataofeft}
{\rm sPar}\colon \Vect^{\A}_{\rm ind}(M) \xrightarrow{\sim} 1|1\eft^0(M)
\eeq
from the category of real, self-adjoint superconnections admitting smooth index bundles to the category~$1|1\eft^0(M)$. 
\end{prop}

In short, $\Vect^{\A}_{\rm ind}(M)\subset \Vect^{\A}(M)$ is the subcategory of superconnections on which the (differential) index map is defined. The functor ${\rm sPar}$ is totally explicit, e.g., 
\beq\label{eq:bigspar}
e^{-t\A^2+\theta\A}\colon \R^{1|1}_{\ge 0} \to \End(\Omega^\bullet(M;V_\C))
\eeq
determines the underlying representation of $\sP_0(M)$ under~\eqref{eq:dataofeft}. The formula~\eqref{eq:bigspar} also coincides with super parallel transport along a superconnection $\A$ on a finite rank vector bundle; e.g., see \cite[Theorem 1]{florin2}, \cite[Proposition 6.7]{TwistAugusto} or \cite[\S2]{DanFei}. This motivates the notation ${\rm sPar}$ for the functor~\eqref{eq:dataofeft}, as it provides a partial construction of super parallel transport for more general superconnections; see Remark~\ref{rmk:SParinfinite} below. Proposition~\ref{prop:dataofeft} allows us to construct the cocycle map in Theorem~\ref{thm:cocycle} as the composition~\eqref{eq:FTtosConn}. 
%In Proposition~\ref{prop:geodata}, we characterize the properties of a superconnection inherited from restriction~\eqref{eq:FTtosConn} of a $1|1$-Euclidean field theory over $M$. Theorem~\ref{thm:cocycle} then follows from~\eqref{eq:FTtosConn} after demanding the existence of smooth energy cutoffs; this guarantees the validity of the index bundle construction, producing a map to (differential) KO-theory. 

\begin{rmk}\label{rmk:SParinfinite}
Dumitrescu's super parallel transport~\cite{florin} extends~\eqref{eq:dataofeft} to a functor
$$
\Vect^\A_{{\rm fin}}(M)\to 1|1\EFT^0(M),
$$ 
from superconnection on finite-rank super vector bundles; compare Corollary~\ref{prop:florin1}. An infinite-rank generalization of super parallel transport would provide a complete dictionary between $1|1$-Euclidean field theories and superconnections.
%A superconnection $\A$ is \emph{$\theta$-summable} if the trace $\Tr(\exp(-\theta \A^2))$ exists for all $\theta\in \R_{>0}$, e.g., see \cite[\S2.1]{LottGorokhovsky} (of course, the occurrence of $\theta$ in this definition has nothing to do with our prior use of~$\theta$ as an odd coordinate). For a self-adjoint, $\theta$-summable
However, in our present analytical setting (following~\cite{ST11})) super parallel transport for infinite rank bundles need not be defined owing to the failure of existence and uniqueness results for differential equations in topological vector spaces, e.g., see \cite[Example~5.6.1]{Hamilton} for counterexamples in Fr\'echet spaces. From the physics point of view, one expects field theories to be determined by infinitesimal data. 
% than the definitions involving general topological vector spaces in~\cite{ST11}.
%it would appear that one must demand more analytical control to construct super parallel transport in general. 
We note that some of the prior definitions of $1|1$-Euclidean field theory leverage existence and uniqueness results in Hilbert spaces~\cite[\S6]{HST}. In a different direction, there are existence and uniqueness results for tame Fr\'echet spaces and tame maps. Adopting a more Wilsonian approach, one could define $1|1$-Euclidean field theories as an $\Ind$-completion of finite-rank objects with compatibilities via the energy filtration; this also forces existence and uniqueness. When attempting to weigh these options, the main concern is that the correct analytical context for the $d=2$ case in conjecture~\eqref{eq:conjecture} is not yet clear. As such, we leave these choices of analysis flexible to avoid over-constraining the possible definitions in this $d=2$ generalization.
% putting just enough existence and uniqueness requirements into the definition of $1|1\eft^n(M)$.
%However, it seems plausible that a bit of additional analytical control would avoid these problems. Optimism comes from the fact that super semigroup representations on Hilbert spaces are equivalent to the data of an infinitesimal generator~\cite[Proposition 39]{HST}, which proves the desired result for super parallel transport on a Hilbert bundle $V\to \pt$. Unbounded superconnections on general Hilbert bundles are studied in~\cite{LottGorokhovsky}. 
\end{rmk}

% structure for the candidate cocycle map to work, as well as its enhancements to non-zero degree and differential KO-theory. 
\subsection{Technical overview}\label{sec:framework}

The categories $1|1\eft^n(M)$ in Theorems~\ref{thm:cocycle} and~\ref{thm:index} come from distilling the information of a degree~$n$ Euclidean field theory down to the essential pieces required to formulate the index map~\eqref{eq:cocycleK}. We work in the framework for field theories developed by Stolz and Teichner in~\cite{ST11}, and begin
%wherein supersymmetric Euclidean field theories are defined as functors out of a \emph{super Euclidean bordism category}. 
%We review their definitions and the associated notions of (twisted) field theories in~\S\ref{sec:ebord}.
%One technical contribution of this paper is to give a first complete definition of 
with a brief recollection of their main definition; see~\S\ref{sec:ebord} for a more detailed overview.

\begin{defn}[{\cite[\S4-5]{ST11}}]\label{defn:ogFT} For $n\in \Z$ the \emph{degree~$n$ twist} is a $\otimes$-invertible functor 
\beq\label{eq:maintwist}
 \twist^n\colon 1|1\EBord(\pt)\to \TA
\eeq
where $\TA$ is a Morita category of (topological) algebras, bimodules, and bimodule maps. Pulling back along the canonical functor $1|1\EBord(M)\to 1|1\EBord(\pt)$ we obtain a twist for the $1|1$-Euclidean bordism category over $M$ that we will also denote by $\twist^n$. 
%We note that this degree~$n$ twist is an essential part of the $1|1$-dimensional version of conjecture~\eqref{eq:conjecture}. 
A \emph{degree~$n$ field theory over $M$} is a natural transformation
\beq
&&\begin{tikzpicture}[baseline=(basepoint)];
\node (A) at (0,0) {$1|1\EBord(M)$};
\node (B) at (5,0) {$\TA.$};
\node (C) at (2.5,0) {$E \Downarrow$};
\draw[->,bend left=15] (A) to node [above] {$\one$} (B);
\draw[->,bend right=15] (A) to node [below] {$\twist^n$} (B);
\path (0,0) coordinate (basepoint);
\end{tikzpicture}\label{eq:twistedEFT}
\eeq
\end{defn} 

\begin{rmk}
Some details of the above definition have not appeared in the literature; we refer to \S\ref{sec:degreendefns} for further discussion that culminates in Hypothesis~\ref{hyp:twist}. 
\end{rmk}

\begin{defn}[Definitions~\ref{defn:RP0},~\ref{defn:Real0}, and~\ref{defn:RPFT}, assuming Hypothesis~\ref{hyp:twist}] \label{defn:STfieldtheories} Define a category $1|1\EFT^n(M)$ whose objects are real, reflection positive degree~$n$ field theories~\eqref{eq:twistedEFT}. Morphisms in $1|1\EFT^n(M)$ are isomorphisms between internal natural transformations $E\simeq E'$ that are compatible with the real and reflection positive structures. 
\end{defn}

\begin{rmk}\label{rmk:KTate}
Reflection positivity is perhaps unnecessary in the $d=1$ version of conjecture~\eqref{eq:conjecture}: unoriented, real field theories (with a different positivity condition) recover KO-theory~\cite{HST}, see Remark~\ref{rmk:RPFT}. However, unoriented and real structures do not generalize to the $2|1$-dimensional Euclidean bordism category, and the Miller character $\TMF\to \KO(\!(q)\!)$ necessitates a map from $2|1$-Euclidean field theories to real K-theory if conjecture~\eqref{eq:conjecture} is to hold. Reflection positivity is one way to obtain real structures on Hilbert spaces of states in dimension~$2|1$; see~\cite[3.2.2]{GPPV} and~\cite{DBEtorsion}.
\end{rmk}

In~\S\ref{sec:superEuclideanpath}-\S\ref{sec:iota}, we construct a category $\sP_0(M)$ of nearly constant super Euclidean paths in~$M$ with an internal functor $\iota$ 
\beq
\label{eq:twistedrestriction}
&&\begin{tikzpicture}[baseline=(basepoint)];
\node (AA) at (-3,0) {$\sP_0(M)$};
\node (A) at (0,0) {$1|1\EBord(M)$};
\node (B) at (4,0) {$\TA.$};
\node (C) at (2.1,0) {$\Downarrow$};
\draw[->] (AA) to node [above] {$\iota$} (A);
\draw[->,bend left=10] (A) to node [above] {$\one$} (B);
\draw[->,bend right=10] (A) to node [below] {$\twist^n$} (B);
\path (0,0) coordinate (basepoint);
\end{tikzpicture}
\eeq 
When applied to degree zero theories, restriction along $\iota$ makes~\eqref{restrict1} precise. Together with Proposition~\ref{lem:mainA}, we obtain the composition~\eqref{eq:FTtosConn} that extract a superconnection from a $1|1$-Euclidean field theory. 
%This superconnection inherits additional structure and property from a further factorization of $\iota$,
%\beq\label{eq:twistedrestriction2}
%\sP_0(M)\hookrightarrow 1|1\ebord(M)\to 1|1\EBord(M)
%\eeq
%through a category $1|1\ebord(M)$ of \emph{small bordisms}, which is built from nearly constant superpaths with different choices of source and target data; see the pictures~\eqref{diag:supergenerators}. 
%The objects and morphisms of $\sP_0(M)$ are (finite-dimensional) supermanifolds, and so the restriction of a degree~$n$ theory along~\eqref{eq:twistedrestriction} provides data that can be analyzed using standard techniques in super geometry. 

%\begin{rmk}
%One can describe the composite functors in~\eqref{eq:twistedrestriction} without delving into the details of the intermediate bordism category $1|1\EBord(M)$. The definition of the categories $1|1\eft^n(M)$ are phrased in this way, i.e., in terms of the finite-dimensional super geometry of $\sP_0(M)$ rather than the category $1|1\EBord(M)$ internal to symmetric monoidal stacks from~\cite{ST11}. This allows the cocycle map~\eqref{eq:cocycleK} to apply to future definitions of field theory that deviate from~\cite{ST11}; see Remark~\ref{rmk:robusto}.
%\end{rmk}

\begin{prop}\label{prop1}\label{prop:geodata}
The restriction of a degree~$n$ Euclidean field theory~\eqref{eq:twistedEFT} along the functor~\eqref{eq:twistedrestriction} determines:
\begin{enumerate}
\item vector bundles $V_+,V_-\to M$ with fiberwise $\cCl_{+n}$- and $\cCl_{-n}$-actions and Clifford linear superconnections $\A_+$ and $\A_-$, respectively (Proposition~\ref{prop:superconn});
\item a pairing between section-valued differential forms (Proposition~\ref{eq:proppairing})
\beq\label{eq:mainhermitian}
&&\langle-,-\rangle\colon \cCl_{n}\otimes_{\cCl_{-n}\otimes \cCl_n} \Big(\Omega^\bullet(M;V_-)\otimes \Omega^\bullet(M;V_+)\Big)\to \Omega^\bullet(M)
\eeq
for which the superconnections satisfy $\langle \A_-x,y\rangle=(-1)^{|x|}\langle x,\A_+y\rangle$; 
\item differential forms $Z,Z_\ell\in \Omega^\bullet(M;C^\infty(\R_{>0}))$ satisfying 
$$
dZ=0, \ \ \partial_\ell Z=dZ_\ell, \qquad \deg(Z)=n \mod 2,\ \ \deg(Z_\ell)=n-1\mod 2
$$
where $\ell$ is the standard coordinate on $\R_{>0}$ (Proposition~\ref{prop:trace}).
\end{enumerate}
The renormalization group flow on degree~$n$ field theories induces the Getzler rescaling of superconnections (Proposition~\ref{prop:Getzler}): renormalization flow by $\mu\in \R_{>0}$ restricts to
\beq
&&\A_\pm\mapsto \mu (\A_\pm)_0+(\A_\pm)_1+\mu^{-1}(\A_\pm)_2+\mu^{-2}(\A_\pm)_3 + \dots\label{eq:RGsconn}
\eeq
Reflection and real structures restrict along~\eqref{eq:twistedrestriction} as follows. 
\begin{enumerate}
\item A reflection structure determines isomorphisms $V_\pm \simeq \overline{V}_\mp$ of bundles of Clifford modules relative to the $*$-superalgebra structure $\cCl_{\pm n} \to \overline{\cCl}_{\mp n}$ defined in~\eqref{eq:superstarCl}. Using this data, the pairing~\eqref{eq:mainhermitian} determines graded hermitian pairings on~$V_+$ and~$V_-$ for which the superconnections $\A_+$ and $\A_-$ are graded self-adjoint. The positivity condition requires the hermitian pairing to be positive (Proposition~\ref{prop:selfadjoint}).
\item A real structure determines isomorphisms $V_\pm\simeq \overline{V}_\pm$ of bundles of modules relative to the standard real structure on $\cCl_n$. The superconnections are compatible with this real structure on~$V_\pm$. For a reflection positive theory with real structure, the hermitian inner product is also real (Proposition~\ref{prop:realA}). 
\end{enumerate}
\end{prop}

As described before Proposition~\ref{prop:dataofeft}, the upshot of Proposition~\ref{prop:geodata} is that a degree~$n$ Euclidean field theory~\eqref{eq:twistedEFT} with real and reflection positive structures determines a superconnection for which a $\KO^n$-valued families index construction \emph{might} apply. The existence of this index bundle can be a rather subtle analytical question. The main condition one requires is that a certain direct sum of eigenspaces for eigenvalues less than a fixed $\lambda\in \R_{>0}$ form a smooth vector bundle over (an open submanifold) of $M$; see Definition~\ref{defn:superconncutoff}. Physically, this subbundle of eigenspaces is the subspace of states with energy less than~$\lambda$. Rephrasing in physical language, the existence of the index bundle is the same as demanding the existence of a smooth cutoff theory in the sense of Wilsonian effective field theory, e.g., see~\cite[\S1.3]{costbook} for a mathematical exposition of the key ideas. The following definition captures the precise pieces of a $1|1$-Euclidean field theory required to extract a KO-valued families index.
%Demanding this property from an object of $1|1\EFT^n(M)$ leads to the following.

%\begin{defn}[Sketch of Definition~\ref{defn:effective}] Let $1|1\eft^n(M)$ be the groupoid whose objects are degree~$n$ representations of $1|1\ebord(M)$ (see Definition~\ref{defn:smallbordisms}) with a reflection structure and a real structure. These data are required to satisfy a positivity property and the property of admitting 
%%admiting satisfy the property that the associated self-adjoint superconnection from Proposition~\ref{prop1} admits 
%smooth energy cutoffs.
%\end{defn} 

 \begin{defn}\label{defn:effective1}\label{defn:effective}
The groupoid $1|1\eft^n(M)$ has objects degree~$n$ representations of~$\sP_0(M)$ (Definition~\ref{defn:degreensP}) that are reflection positive (Definition~\ref{defn:rssP}), trace class (Definition~\ref{defn:tracestru}), real (Definition~\ref{defn:rrsP}), and admit cutoffs (Definition~\ref{def:hypothesis}). Morphisms are isomorphisms between degree~$n$ representations of~$\sP_0(M)$ that are compatible with the reflection structure, real structure, and trace. 
\end{defn}

\begin{cor}\label{cor:STcocycle}
Theorem~\ref{thm:cocycle} constructs a cocycle map for a subcategory of Stolz and Teichner's degree~$n$ theories,
\beq\nonumber
&&\{E\in 1|1\EFT^n(M)\mid E \ {\rm admits\ energy\ cutoffs}\}\xrightarrow{{\rm cocycle}} \dKO^n(M).
\eeq

\end{cor}
\bp
The definition of $2|1\eft^n(M)$ is set up so that an object $E\in 1|1\EFT^n(M)$ determines the \emph{data} of an object in $1|1\eft^n(M)$ by restriction to specific families of super Euclidean bordisms; it remains to check that these data satisfy the required properties. The positivity property for an object in $1|1\eft^n(M)$ follows from reflection positivity for $1|1\EFT^n(M)$. Field theories satisfying the cutoff property provide a subcategory of field theories, and we have a span
\beq\label{eq:STcomparison}
&&1|1\EFT^n(M)\hookleftarrow \{E\in 1|1\EFT^n(M)\mid E \ {\rm admits\ energy\ cutoffs}\}\xrightarrow{{\rm restrict}} 1|1\eft^n(M). 
\eeq
Postcomposing the map on the right with~\eqref{eq:cocycleK} provides the claimed map to $\dKO^n(M)$. 
\ep

Theorem~\ref{thm:index} follows from turning Proposition~\ref{prop:geodata} around: the Bismut superconnection determines a representation of $\sP_0(M)$ from the formula~\eqref{eq:bigspar}, and this representation naturally comes with the additional data and properties required to determine an object in $1|1\eft^n(M)$. In this example, the existence of cutoffs relies on Bismut's families index theorem~\cite{Bismutindex}. Despite only being part of the data of a field theory (see Remark~\ref{rmk:smallfieldtheory}), $\alpha(X)\in 1|1\eft^{-n}(M)$ contains all the relevant data for the families differential index. 

We conclude with several remarks about the framework developed in this paper. 

\begin{rmk}\label{rmk:generalcocycle}
One can formulate the span~\eqref{eq:STcomparison} for any definition of degree~$n$, $1|1$-Euclidean field theories defined using a bordism category that contains the families of small bordisms (see Definition~\ref{defn:smallbordisms}) where the degree twist satisfies Hypothesis~\ref{hyp:twist}. Hence, such field theories admit a subcategory with a cocycle map valued in $\dKO^n(M)$ by composing with~\eqref{eq:cocycleK}. For finite-rank examples, the cutoff condition is automatically satisfied, see Corollary~\ref{prop:florin1} which follows from~\cite[\S9.1]{BGV}. Hence, the cutoff condition in~\eqref{eq:STcomparison} is only interesting for infinite-rank examples. 
\end{rmk}
% bundle and what type of superconnection one considers. E.g., general sheaves rather than vector bundles, Hilbert bundles, etc.

\begin{rmk}\label{rmk:elliptic}
Cutoffs make sense for a broad class of field theories. Indeed, in Witten's study of the $\mathcal{N}=(0,1)$ supersymmetric sigma model~\cite[\S1]{Witten_Dirac}, he suggested that cutoffs be used to analyze the field theories in the $d=2$ version of the conjectures~\eqref{eq:conjecture}. 
%We emphasize that the existence of a cutoff is a \emph{property} of a field theory rather than additional data. 
%Hence, field theories with cutoffs form a subcategory of the category of all field theories, as in~\eqref{eq:STcomparison}. 
For $1|1$-dimensional field theories, the property of having smooth energy cutoffs is analogous to the property of a differential operator being elliptic; in particular, cutoffs require eigenfunctions of the degree zero part of $\A^2$ to determine smooth sections over~$M$.
\end{rmk}

\begin{rmk}\label{rmk:smallfieldtheory}
Assuming a generators and relations description of the category $1|1\EBord(M)$ (following \cite[6.7]{HST} and \cite[3.2]{ST11}), the data of an object of $1|1\eft^n(M)$ is the data of a degree~$n$ field theory defined on a subcategory of \emph{small bordisms in $M$}
$$
1|1\ebord(M)\hookrightarrow 1|1\EBord(M)
$$
generated by the $S$-families of bordisms in Definition~\ref{defn:smallbordisms}. 
%whose map to~$M$ factors through~$S\times \R^{0|1}$. 
\end{rmk}

\begin{rmk} \label{rmk:concordance}
By results of Hopkins--Hovey~\cite{HopkinsHovey}, Theorems~\ref{thm:cocycle}  and~\ref{thm:index} yield a surjection 
\beq\label{eq:KOsurjection}
1|1\eft^n(M)\twoheadrightarrow \KO^n(M),
\eeq 
verifying a version of the $d=1$ case of~\eqref{eq:conjecture}.
A stronger form of the conjectures~\eqref{eq:conjecture} asserts that the map~\eqref{eq:KOsurjection} induces an isomorphism on concordance classes, so that the corresponding cohomology theories are a complete invariant of field theories~\cite[\S1.5]{ST11}. It is easy to see that~\eqref{eq:KOsurjection} factors through concordance classes, but we do not expect an isomorphism in our framework used below for reasons we explain presently. 
%that the isomorphism class of a vector bundle on which a superconnection is defined is a concordance invariant of the corresponding field theory. 
In view of the factorization~\eqref{eq:FTtosConn} of the cocycle map, the stronger form of conjecture~\eqref{eq:conjecture} in particular requires a category of superconnections whose concordance classes are KO-theory. But in existing definitions the isomorphism class of the underlying vector bundle is a concordance invariant. In particular, the rank of a vector bundle is a concordance invariant, but is not an invariant of the underlying $\KO$-class. We see two possible ways to avoid this problem. The first option is to define field theories on (infinite rank) Hilbert bundles, using that such bundles are always trivializable~\cite{Kuiper}; this resonates with previous definitions from \cite{ST04} and \cite{HST}, and likely makes contact with Lott and Gorokhovsky's model for differential K-theory~\cite{LottGorokhovsky}. The second option follows the suggestion in~\cite[Remark~3.16]{ST11}, and requires a framework for superconnections defined on (not necessarily locally free) sheaves of topological vector spaces, e.g., see~\cite{PeterU} for finite-rank results in this direction. It is not clear if either of these approaches generalize to the desired statement in dimension~$2|1$. 
%Although more technical, this stronger conjecture appears quite reasonable from the physical point of view. In brief, the missing concordances between superconnections (in the standard definition) connect invertible superconnections to the zero bundle; morally, the Getzler rescaling~\eqref{eq:RGsconn} furnishes such a concordance. By Proposition~\ref{prop:geodata}, these concordances should come from the renormalization group flow. 
%Theorem~\ref{thm:cocycle} shows that this requires a concordance connecting a field theory with invertible superconnection to the zero theory, where the concordance itself is implemented by the renormalization group flow. Given the factorization of~\eqref{eq:cocycleK} through superconnections, this requires the existence of an associated concordance in the category of super vector bundles. , and the main constructions in index theory do not require it. 

% for a closely related result in the finite-rank case. 
% RG-trivial theories should be sent to zero, and the question is whether or not there is an honest family of field theories implementing this equivalence.
\end{rmk}

%\subsection{Overview of previous work}
%
%The super bordism category in \cite[Definition 3.2.2]{ST04} is not spelled out; roughly, field theories in this treatment are taken to be self-adjoint, trace-class representation of the non-unital super semigroup $\R^{1|1}_{>0}$ (via \cite[Proposition 3.2.6]{ST04}), and it is indicated how this super semigroup is related to a category of super intervals. A definition of the $1|1$-bordism category is sketched in \cite{HST}, and a similarl relationship with KO-theory is described (passing through super semigroup representations), but the proofs of certain key technical lemmas were omited (e.g., generators and relations of the bordism category) and an ad-hoc positivity condition is imposed (which differs from the standard notion of reflection positivity). In \cite{ST11}, the bordism category is defined in greater detail, but no map to K-theory has been constructed with this definition, and the degree~$n$ twist hasn't been constructed. Furthermore, the definition of the bordism category there is certainly not final, as it doesn't easily generali
\subsection*{Outline}

The motivation from physics is explained in~\S\ref{sec:SUSYQMmot}. In~\S\ref{sec:superEuclideanpath} we define categories of superpaths in a smooth manifold. In~\S\ref{sec:iota} we show that superpaths give examples of bordisms in Stolz and Teichner's $1|1$-Euclidean bordism category from~\cite{ST11}. In~\S\ref{sec:degreendefns} we sketch a construction of the degree~$n$ twist and the values of the resulting degree~$n$ twisted field theories. The main technical results are proved in~\S\ref{sec:repsP}, where we restrict these values of field theories to the nearly constant superpaths to prove Proposition~\ref{prop1}. In~\S\ref{sec:index}, the relationship between cutoff theories and the families index construction leads us to the definition of the groupoids~$1|1\eft^n(M)$. Theorems~\ref{thm:cocycle} and~\ref{thm:index} then follow from standard constructions in index theory. In~\S\ref{eq:superlinear} we review some background super algebra and super geometry,~\S\ref{sec:ebord} summarizes the main definitions and constructions from~\cite{ST11} used in this paper, and~\S\ref{sec:KOindex} reviews the families index in differential $\KO$-theory from~\cite{DBEBis}. 

%We have attempted to make the paper as self-contained as possible: the definitions behind Theorems~\ref{thm:cocycle} and~\ref{thm:index} are given in full below, and (although deeply inspired by their ideas) the central definitions are independent from~\cite{ST11}. Indeed, 
Readers without detailed knowledge of the Stolz--Teichner program may prefer to skim sections~\S\ref{sec:iota} and~\S\ref{sec:degreendefns}; these sections explain how our results fit with~\cite{ST11}, but are logically independent from the main theorems. 

This paper is part of a series~\cite{DBEChern,DBEBis,DBEtorsion} that solves a 1-categorical version of the conjectures~\eqref{eq:conjecture}. The zero categorical version is the main result of~\cite{DBEChern}. In the 1-categorical enhancement, the target of~\eqref{eq:conjecture} requires a new model for (differential) $\KO$-theory; this is not directly related to field theories and is handled in~\cite{DBEBis}, reviewed in~\S\ref{sec:KOindex} below. The present paper connects $1|1$-Euclidean field theories with these differential KO-cocycles via an index map. In~\cite{DBEtorsion} we generalize this index to $2|1$-dimensional field theories.

\subsection*{Acknowledgements} It is a pleasure to thank Kevin Costello, Theo Johnson-Freyd, Matthias Ludewig, Stephan Stolz, Peter Teichner, and Arnav Tripathy for stimulating conversations on this material. This work was supported by the National Science Foundation under grant number DMS-2205835.

\setcounter{tocdepth}{1}
\tableofcontents

\newpage

\section{Motivation from supersymmetric quantum mechanics} \label{sec:SUSYQMmot}

In this section we overview the arguments in supersymmetric quantum mechanics that lead to Theorems~\ref{thm:cocycle} and~\ref{thm:index} in the case that $M=\pt$. The discussion is intentionally informal with precise definitions delayed until later sections. 

\subsection{Quantum mechanics on $X$}\label{sec:QM}
To build intuition, we start by explaining ordinary quantum mechanics without supersymmetry. 

Let $X$ be a closed, oriented Riemannian manifold. Quantum mechanics on~$X$ has as its space of states $\mathcal{H}=C^\infty(X;\C)$ endowed with the hermitian inner product inherited from integration with respect to the Riemannian volume form on $X$. The Laplacian~$H=\Delta$ on~$X$ is the Hamiltonian of the theory, generating time-evolution. In the Wick-rotated theory, this gives a representation of the semigroup~$\R_{\ge 0}$,
\beq
&&\R_{\ge 0}\to \End(\mathcal{H}),\qquad t\mapsto e^{-tH},\label{eq:semigroup}
\eeq
where $t\in \R_{\ge 0}$ is identified with a time parameter. This representation has important additional properties: $e^{-tH}$ is self-adjoint, and for $t>0$ the trace is finite 
\beq
Z:=\Tr(e^{-tH})\in C^\infty(\R_{>0})\label{eq:partitionfunctioneasy}
\eeq
since $X$ was assumed to be compact. The function of $t$ given by~\eqref{eq:partitionfunctioneasy} is called the \emph{partition function}. A general quantum mechanical system is given by a semigroup representation~\eqref{eq:semigroup} satisfying a self-adjoint and trace-class property. 
%extends to the Hilbert completion of $V$ where it takes values in trace class operators. 

\begin{rmk}
Often one takes the Hilbert completion of the Fr\'echet space $\mathcal{H}=C^\infty(X;\C)$, resulting in a Hilbert space of states and an unbounded Hamiltonian $H=\Delta$. For the definition of field theory from~\cite{ST11} the domain of the operator $\Delta$ (as a Fr\'echet space) turns out to be more fundamental. By contrast, the Hilbert completion is the basic object of study for definitions of field theories in~\cite{ST04} and \cite{HST}. 
%This flexibility in Stolz and Teichner's definitions is an important feature, as it allows one to consider various analytic models for K-theory in their framework. However, 
It is unclear which point of view generalizes more readily to the $d=2$ conjecture in~\eqref{eq:conjecture}. In this paper we will largely follow~\cite{ST11}, though we will ignore these functional analytic distinctions for the remainder of this section. 
\end{rmk}

%\begin{rmk} Smooth functions $C^\infty(X;\C)$ naturally form a nuclear Fr\'echet space, and the tensor product in~\eqref{eq:semigroup} is the projective tensor product of nuclear Fr\'echet spaces. We refer to~\cite[Appendix~2]{costbook} for an excellent overview of nuclear Fr\'echet spaces. 
%%As is customary in quantum mechanics, one can complete $V$ with respect to the hermitian inner product to obtain a Hilbert space. However, (particular in later examples) the Fr\'echet space $V$ turns out to be the more fundamental object. 
%\end{rmk} 
%%$$
%%\langle f,g\rangle=\int_X \overline{f}g \ d{\rm vol}. 
%$$
In the formalism of~\cite{ST04,HST,ST11}, a quantum mechanical system determines a \emph{1-dimensional Euclidean field theory}, defined as a representation of a 1-dimensional oriented Riemannian bordism category $1\EBord$. The objects of $1\EBord$ are finite disjoint unions of $\bullet=\pt^+$ and $\circ=\pt^-$, the point with its two choices of orientation. The morphisms are disjoint unions of 1-dimensional Riemannian bordisms:
\beq
&&\begin{tikzpicture}[baseline=(basepoint)];
%\draw [fill] (0,0) circle [radius=0.1] node [black,above=4] {$\pt^+$};
%\draw (0,-1.5) circle [radius=0.1] node [black,above=4] {$\pt^-$};
\draw [fill]  (.5,0) circle [radius=0.1] node [black,left=4] {$\pt^+$};
\draw  [fill] (1.5,0) circle [radius=0.1] node [black,right=4] {$\pt^+$};
\draw[decoration={markings,
        mark=at position \halfway with \arrow{<}},
        postaction=decorate,thick] (.5,0) to [out=0,in=0] (1.5,0);
\node (Iplus) at (1,.3) {$\sI_t^+$};

\draw (.5,-1.5) circle [radius=0.1] node [black,left=4] {$\pt^-$};
\draw (1.5,-1.5) circle [radius=0.1] node [black,right=4] {$\pt^-$};
\draw[decoration={markings,
        mark=at position \halfway with \arrow{>}},
        postaction=decorate,thick] (.5,-1.5) to [out=0,in=0] (1.4,-1.5);
\node (Iminus) at (1,-1.2) {$\sI_t^-$};

\draw (4,-1.5) circle [radius=0.1] node [black,left=4] {$\pt^-$};
\draw [fill] (4,0) circle [radius=0.1] node [black,left=4] {$\pt^+$};
\draw[decoration={markings,
        mark=at position \halfway with \arrow{>}},
        postaction=decorate,thick] (4,-1.5) to [out=35,in=-35] (4,0);
\node (Rplus) at (3.75,-.75) {$\sR_t^-$};

\draw [fill] (6,-1.5) circle [radius=0.1] node [black,left=4] {$\pt^+$};
\draw  (6,0) circle [radius=0.1] node [black,left=4] {$\pt^-$};
\draw[decoration={markings,
        mark=at position \halfway with \arrow{>}},
        postaction=decorate,thick]  (6,0) to [in=35,out=-35] (6,-1.5) ;
\node (Rminus) at (5.75,-.75) {$\sR_t^+$};

\draw  (8.5,-1.5) circle [radius=0.1] node [black,left=4] {$\pt^-$};
\draw [fill] (8.5,0) circle [radius=0.1] node [black,left=4] {$\pt^+$};
\draw[decoration={markings,
        mark=at position \halfway with \arrow{>}},
        postaction=decorate,thick]  (8.5,0) to [in=145,out=225] (8.5,-1.5);
\node (Lplus) at (7.75,-.75) {$\sL_t^+$};

\draw [fill] (10.5,-1.5) circle [radius=0.1] node [black,left=4] {$\pt^+$};
\draw  (10.5,0) circle [radius=0.1] node [black,left=4] {$\pt^-$};
\draw[decoration={markings,
        mark=at position \halfway with \arrow{<}},
        postaction=decorate,thick]  (10.5,0) to [in=145,out=225] (10.5,-1.5) ;
\node (Lmius) at (9.75,-.75) {$\sL_t^-$};

\draw[decorate,thick] (12,-.75) circle [radius=.35] node [left=7] {$\sS_t$};

%\draw[decoration={markings,
%        mark=at position \halfway with \arrow{>}},
%        postaction=decorate,thick] (12,-1.25) circle [radius=.35] node [left=7] {$S_t^-$};

%\draw[decoration={markings,
%        mark=at position \halfway with \arrow{<}},
%        postaction=decorate,thick] (12,-.2) circle [radius=.35] node [left=7] {$S_t^+$};

\path (0,-.75) coordinate (basepoint);
\end{tikzpicture}\label{diag:generators}
%&&1|1\EFT(M)\stackrel{{\rm cocycle}}{\dashrightarrow}  \K(M)  \qquad  2|1\EFT(M) \stackrel{{\rm cocycle}}{\dashrightarrow} \TMF(M)
\eeq
%The 1-manifolds $S_t^\pm$ are the circles of circumference~$t$, viewed as bordisms from the empty set to itself. 
The 1-manifold $\sS_t$ is the circle of circumference~$t$, viewed as a bordism from the empty set to itself. 
The 1-manifolds underlying the bordisms $\sI_t^\pm$, $\sR_t^\pm$ and~$\sL_t^\pm$ are all intervals of length~$t$; the distinctions between them come from their source and target data. The notation in~\eqref{diag:generators} follows~\cite{HST,ST11}, where $\sL_t^\pm$ and $\sR_t^\pm$ are the ``left-" and ``right-elbows" of length~$t$, respectively, where the $+$ or $-$ decoration refers to the orientation, indicated pictorially in~\eqref{diag:generators} by an arrow. We also follow the convention from~\cite{HST,ST11} to read these pictures of bordisms from right to left. For example, $\sR_t^+$ is a bordism from the empty set to $\pt^-\coprod \pt^+$, while $\sL_t^+$ is a bordism from $\pt^+\coprod \pt^-$ to the empty set. 
%This convention makes it easy to translate into formulas like~\eqref{eq:relations} below where the direction of composition is also from right to left. 
%The 1-manifold underlying $S_t^\pm$ is the circle of circumference~$t$. 
Bordisms compose by gluing along matching source and target data. The disjoint union of 0- and 1-manifolds endows the category $1\EBord$ with a symmetric monoidal structure. 
%We refer to~\S\ref{sec:ex1EB} below for details.

One can compute relations associated with compositions of bordisms. For example, there are isometries of Riemannian 1-manifolds: 
\beq
&&\sI_s^\pm\circ \sI_t^\pm \simeq \sI_{s+t}^\pm,  \qquad \sL_0^\pm\circ \sigma^\pm \simeq \sL_0^\mp,\qquad  \sL_r^\pm \circ (\sI_t^\pm \coprod \sI_s^\mp)\simeq \sL_{r+s+t}^\pm,
%\qquad R_{t}^{\pm} \circ_{\pt^\mp} L_0^\mp= I_{t}^\pm 
%\quad L_t^\pm\circ R_s^\mp\simeq S_{t+s}^\pm.
\label{eq:relations}
\eeq
where in the middle identity $\sigma^\pm \colon \pt^\mp \coprod \pt^\pm \stackrel{\sim}{\to} \pt^\pm \coprod \pt^\mp$ is the symmetry for the monoidal structure. 
The first identity is a geometric counterpart to the semigroup $\R_{\ge 0}$ from~\eqref{eq:semigroup}: the lengths of intervals add. 
%One can show~\cite[Theorem~58]{HST} that all morphisms are generated by~$R_t^\pm$ and~$L_0^\pm$, where the latter is the length zero specialization of $L_t^\pm$. These generators satisfy the relations (using notation from~\cite[Theorem~58]{HST})
%\beq
%R_s^\pm\circ_{L_0^\pm} R_t^\pm\simeq R_{s+t}^\pm,\qquad \sigma^\pm \circ R_t^\pm\simeq R_{t}^\mp, \qquad L_0^\pm\circ \sigma^\pm \simeq L_0^\mp\label{eq:genrelate}
%\eeq
%The semigroup of Riemannian intervals is recovered from these generators via $I_t^\pm\simeq R_t^\mp \circ_{\pt^\pm} L_0^\pm$. 
%The presentation~\eqref{eq:genrelate} follows from the same argument in in the unoriented case from~\cite[Theorem~58]{HST}. 
Define the category $\Fun^\otimes (1\EBord,\Vect)$ whose objects are symmetric monoidal functors
\beq
E\colon 1\EBord\to \Vect\label{eq:repof1EB0}
\eeq
and morphisms are natural isomorphisms of functors. Above, $\Vect$ denotes the category of super (i.e., $\Z/2$-graded) nuclear Fr\'echet spaces over $\C$ endowed with the $\Z/2$-graded projective tensor product. Physically, the $\Z/2$-grading $\mathcal{H}=\mathcal{H}^\ev\oplus \mathcal{H}^\odd$ keeps track of the even subspace~$\mathcal{H}^\ev$ as \emph{bosonic} states and the odd subspace $\mathcal{H}^\odd$ as \emph{fermionic} states. This grading will become more important when we add supersymmetry. 
%Restriction to the bordisms $I_t^+$ in~\eqref{eq:relations} determines a functor
%\beq
%&&\begin{tikzpicture}[baseline=(basepoint)];
%\node (A) at (0,0) {$\Fun^\otimes (1\EBord,\Vect)$};
%\node (B) at (5,0) {$\Rep(\R_{\ge 0}),$};
%\node (C) at (8,0) {$E\mapsto E(I_t^+)$,};
%\draw[->] (A) to node [below] {\small restrict} (B);
%%\draw[->,dashed, bend right=15] (B) to (A);
%\path (0,-.05) coordinate (basepoint);
%\end{tikzpicture}\label{eq:restrict1EB0}
%%&&1|1\EFT(M)\stackrel{{\rm cocycle}}{\dashrightarrow}  \K(M)  \qquad  2|1\EFT(M) \stackrel{{\rm cocycle}}{\dashrightarrow} \TMF(M)
%\eeq
%with values in representations of the semigroup $\R_{\ge 0}$. Lifting a semigroup representation to a field theory is both data and property, namely, values on the remaining bordisms in~\eqref{diag:generators} and relations amongst these bordisms~\eqref{eq:relations}. 

For the data of a functor~\eqref{eq:repof1EB0} to match the standard data of a quantum mechanical system~\eqref{eq:semigroup}, 
%Hence,~\eqref{eq:repof1EB0} specifies values
%\beq
%\begin{array}{c}
%E(\pt^\pm)=V_\pm,\quad E(I_t^\pm)\colon V_\pm\to V_\pm, \quad  E(L_t^\pm)\colon V_\pm\otimes V_\mp \to \C,\\ E(R_t^\pm) \colon \C\to V_\pm\otimes V_\mp,\quad E(S^\pm_t)\in \C\end{array} \label{eq:FTdata}
%\eeq
%where relations amongst compositions of bordisms impose conditions.
%%In the context of topological field theory, these relations are dualizability conditions involving the given data~\cite{BaezDolancobord,Lurie_cob}. In the non-topological setting, this dualizability argument breaks down essentially because $E(I_t^\pm)$ need not be the identity linear map. Instead, 
%For example, the maps $E(I_t^\pm)$ determine semigroup representations on the vector spaces~$V_\pm$. The representations must be adjoint to each other under the pairings determined by $E(L_0^\pm )$. Finally, for $t>0$ the linear maps $E(I_t^\pm)$ must possess the trace property, with the trace being~$E(S_t^\pm)\in \C$, the value on the circle of circumference~$t$. 
%For a representation~\eqref{eq:repof1EB0} to make contract with what a physicist would recognize as a quantum mechanical theory, 
one must demand an additional structure on the functors~\eqref{eq:repof1EB0} called \emph{reflection positivity}, e.g., see~\cite[\S3]{FreedHopkins} or~\cite[page~18]{KontsevichSegal}. This begins with the $\Z/2$-action on $1\EBord$ generated by orientation reversal of 1-manifolds; in terms of~\eqref{diag:generators}, this is 
\beq
&&\pt^\pm\mapsto \pt^\mp,\qquad \sI_t^\pm\mapsto \sI_t^\mp,\qquad \sR_t^\pm\mapsto \sR_t^\mp,\qquad \sL_t^\pm\mapsto \sL_t^\mp. \label{eq:orientationreverse}
\eeq
Physically, this action reverses the arrow of time. There is also a $\Z/2$-action on vector spaces generated by complex conjugation. A \emph{reflection structure} for the functor~\eqref{eq:repof1EB0} is $\Z/2$-equivariant data with respect to the actions on the source and target (e.g., see~\cite[\S{B}]{FreedHopkins} or Definition~\ref{defn:equivariancedata} below), which includes the data of an isomorphism
$$
\mathcal{H}:=E(\pt^+)\stackrel{\sim}{\to} \overline{E(\pt^-)}.
%,\qquad E(\pt^-)=V_-\stackrel{\sim}{\to} \overline{V}_+=E(\pt^+). 
$$
A representation~\eqref{eq:repof1EB0} with a reflection structure is \emph{reflection positive} if the hermitian pairing
$$
\langle-,-\rangle\colon \overline{\mathcal{H}}\otimes \mathcal{H}\simeq E(\pt^-)\otimes E(\pt^+)\xrightarrow{E(\sL_0^-)} \C
$$
is positive. A \emph{1-dimensional Euclidean field theory} is a reflection positive representation~\eqref{eq:repof1EB0}. A property demanded by reflection positivity is that $E(\sI_t^\pm)$ be a self-adjoint representation of the semigroup $\R_{\ge 0}$ on $\mathcal{H}$. 

%\begin{rmk} Reflection positivity data can be viewed as a categorical generalization of Atiyah's Real vector bundles. Indeed, recall that a Real vector bundle on a space $Z$ with involution $\tau\colon Z\to Z$ is a vector bundle $V\to Z$ with a conjugate linear involution $\tilde{\tau}\colon V\to V$ covering $\tau\colon Z\to Z$. A representation of the category $1\EBord$ generalizes the notion of vector bundle; orientation reversal gives an involution on $1\EBord$, and reflection positivity data lifts this to a conjugate linear involution on the representation. 
%% has an involution, and we demand representations be $\Z/2$-equivariant in the \emph{conjugate linear} sense. 
%\end{rmk}

From the discussion above, a 1-dimensional Euclidean field theory is determined by its restriction along~\eqref{eq:orientationreverse} as the data of a self-adjoint semigroup representation $(t\mapsto E(\sI_t^+))$ on the vector space~$\mathcal{H}=E(\pt^+)$ endowed with the hermitian inner product $\langle-,-\rangle=E(L_0^-)$. Furthermore, when restricted to~$\R_{>0}\subset \R_{\ge 0}$, this semigroup representation is trace class and the partition function
$$
Z=\Tr(E(\sI_t^+))=E(\sS_t),\qquad t\in \R_{>0}
$$
can be viewed as a function on the moduli of metrized circles with circumference $t>0$. Quantum mechanics on a closed, oriented Riemannian manifold~$X$ provides an example of a field theory with values
$$
E(\pt^+)=C^\infty(X;\C),\quad E(\sI_t^+)=e^{-t\Delta}, \quad E(\sL_0^-)(s_1,s_2)= \int_X \overline{s_1}s_2\ d{\rm vol }.
%\quad E(I_t^-)=\overline{e^{-t\Delta}}\colon \overline{V}\to \overline{V}
$$

%By adding additional structures (e.g., supersymmetry) one can hope that the resulting space of field theories will have an interesting topology, and hence admit interesting invariants. 

\begin{rmk}\label{rmk:collars} We comment briefly on an important technical detail ignored above. Various quantities in a 1-dimensional Euclidean field theory ought to depend smoothly on geometric parameters, e.g., the semigroup $t\mapsto E(\sI_t^+)$ and the partition function $Z=E(\sS_t)$ should depend smoothly on~$t$. In Stolz and Teichner's formalism~\cite{HST,ST11} this is addressed by defining field theories as fibered functors, assigning smoothly varying quantities to smoothly varying families of bordisms;
%More precisely, $1\EBord$ and $\Vect$ are defined as categories fibered over smooth manifolds, and the functors~\eqref{eq:repof1EB0} are fibered functors; 
see~\S\ref{sec:ebord} below. 
\end{rmk}

\subsection{Supersymmetric quantum mechanics} \label{sec:SUSYQM}

Let $X$ be a closed, even-dimensional Riemannian spin manifold, and $\slashed{S}\to X$ be the $\Z/2$-graded spinor bundle with Dirac operator~$\slashed{D}$. Supersymmetric quantum mechanics on~$X$ studies the space of states $\mathcal{H}=\Gamma(X;\slashed{S})$ as a $\Z/2$-graded vector space with hermitian inner product determined by pairing and integrating spinors. The Dirac Laplacian $H=\slashed{D}{}^2$ is the Hamiltonian for this theory, leading to an odd square root~$Q=\slashed{D}$ of the Hamiltonian
\beq
\frac{1}{2}[Q,Q]=\frac{1}{2}(QQ+QQ)=Q^2=\slashed{D}^2=H,\label{eq:Lie11}
\eeq
where $[-,-]$ is the graded commutator. The operator $Q$ defines a \emph{supersymmetry}, e.g., see~\cite[\S1]{susymorse}. A first consequence of supersymmetry is that the partition function
\beq\label{eq:partition function2}
&&Z=\sTr(e^{-tH})=\sTr(e^{-tQ^2})=\dim(\ker(Q))_+-\dim(\ker(Q))_-=\Ind(Q)
\eeq
is independent of $t$ and equal to the super dimension of the kernel of $Q$, i.e., the index of the Dirac operator~$Q=\slashed{D}$. This follows from the McKean--Singer formula~\cite{McKeanSinger}; see also~\cite[Theorem~3.50]{BGV}. It is crucial in this argument that the partition function~\eqref{eq:partition function2} uses the \emph{supertrace}, defined as
$$
\sTr(e^{-tH}):=\Tr((-1)^{\sf F}\circ e^{-tH})
$$
where $\Tr$ is the ordinary trace, and $(-1)^{\sf F}$ acts by $+1$ on the even subspace of $\mathcal{H}$ and $-1$ on the odd subspace. 

The supersymmetry $Q$ allows one to extend the usual time-evolution~\eqref{eq:semigroup} to a representation of a \emph{super semigroup},
\beq
\R^{1|1}_{\ge 0}\to \End(\mathcal{H})\qquad (t,\theta)\mapsto e^{-tH+\theta Q},\label{eq:N1super}
\eeq
where $(t,\theta)$ are coordinates on the supermanifold with boundary $\R^{1|1}_{\ge 0}$, defined by restricting the structure sheaf of~$\R^{1|1}$ along the inclusion~$\R_{\ge 0}\subset \R$ (see Example~\ref{ex:Roneonege}). The semigroup multiplication on~$\R^{1|1}_{\ge 0}$ is the restriction of multiplication in the super Lie group defined as follows. 

\begin{defn}
Define the super Lie group $\R^{1|1}$ with multiplication 
\beq
(t,\theta)\cdot (s,\eta)=(t+s+\theta\eta,\theta+\eta)\qquad (t,\theta),(s,\eta)\in \R^{1|1}(S), \label{eq:superEuc}
\eeq
where the above formula is understood in terms of the functor of points of~$\R^{1|1}$. 
\end{defn}

\begin{rmk}\label{rmk:SUSYmorse}
A closely related (and perhaps better known) example of supersymmetric quantum mechanics takes as its space of states~$\mathcal{H}=\Omega^\bullet(X)$, differential forms on an oriented Riemannian manifold~\cite[\S2]{susymorse}. The Hamiltonian is then defined to be the Hodge Laplacian $H=dd^*+d^*d$, which has a pair of commuting odd square roots, 
\beq
Q_1=d+d^*, \quad Q_2=i(d-d^*), \quad [Q_1,Q_2]=0, \quad Q_1^2=Q_2^2=H.\label{eq:Nis2}
\eeq
This is the data of $\mathcal{N}=2$ supersymmetry, rather than the $\mathcal{N}=1$ example above. With our intended connections to index theory, $\mathcal{N}=1$ supersymmetry is more fundamental. Of course, forgetting one of the $Q_i$ in~\eqref{eq:Nis2} extracts an $\mathcal{N}=1$ theory. 
\end{rmk}

%with multiplication
%\beq
%(t,\theta)\cdot (s,\eta)=(t+s+\theta\eta,\theta+\eta)\qquad (t,\theta),(s,\eta)\in \R^{1|1} \label{eq:superEuc0}
%\eeq
%In particular, restriction along the canonical inclusion $\R_{\ge 0}\subset \R^{1|1}_{\ge 0}$ is a homomorphism. 
A general supersymmetric quantum mechanical system is given by a supersemigroup representation~\eqref{eq:N1super} satisfying a self-adjoint and trace class condition. In the formalism of~\cite{ST04,HST,ST11}, such data determine a \emph{$1|1$-dimensional Euclidean field theory}, as we review below.

One can repackage supersymmetric quantum mechanics on $X$ as a representation of a $1|1$-Euclidean bordism category $1|1\EBord$, arising as a generalization of $1\EBord$. The objects of $1|1\EBord$ are finite disjoint unions of $\bullet=\spt^+$ and $\circ=\spt^-$, the \emph{super point} with its two choices of orientation. The morphisms can be visualized as
\beq
&&\begin{tikzpicture}[baseline=(basepoint)];
%\draw [fill] (0,0) circle [radius=0.1] node [black,above=4] {$\pt^+$};
%\draw (0,-1.5) circle [radius=0.1] node [black,above=4] {$\pt^-$};
\draw [fill]  (.5,0) circle [radius=0.1] node [black,left=4] {$\spt^+$};
\draw  [fill] (1.5,0) circle [radius=0.1] node [black,right=4] {$\spt^+$};
\draw [thick] (.5,0) to [out=0,in=0] (1.5,0);
\node (Iplus) at (1,.3) {$\sI_{t,\theta}^+$};

\draw (.5,-1.5) circle [radius=0.1] node [black,left=4] {$\spt^-$};
\draw (1.5,-1.5) circle [radius=0.1] node [black,right=4] {$\spt^-$};
\draw [thick] (.5,-1.5) to [out=0,in=0] (1.4,-1.5);
\node (Iminus) at (1,-1.2) {$\sI_{t,\theta}^-$};

\draw (4,-1.5) circle [radius=0.1] node [black,left=4] {$\spt^-$};
\draw [fill] (4,0) circle [radius=0.1] node [black,left=4] {$\spt^+$};
\draw [thick] (4,-1.5) to [out=35,in=-35] (4,0);
\node (Rplus) at (3.75,-.75) {$\sR_{t,\theta}^-$};

\draw [fill] (6,-1.5) circle [radius=0.1] node [black,left=4] {$\spt^+$};
\draw  (6,0) circle [radius=0.1] node [black,left=4] {$\spt^-$};
\draw [thick]  (6,0) to [in=35,out=-35] (6,-1.5) ;
\node (Rminus) at (5.75,-.75) {$\sR_{t,\theta}^+$};

\draw  (8.5,-1.5) circle [radius=0.1] node [black,left=4] {$\spt^-$};
\draw [fill] (8.5,0) circle [radius=0.1] node [black,left=4] {$\spt^+$};
\draw [thick]  (8.5,0) to [in=145,out=225] (8.5,-1.5);
\node (Lplus) at (7.75,-.75) {$\sL_{t,\theta}^+$};

\draw [fill] (10.5,-1.5) circle [radius=0.1] node [black,left=4] {$\spt^+$};
\draw  (10.5,0) circle [radius=0.1] node [black,left=4] {$\spt^-$};
\draw [thick]  (10.5,0) to [in=145,out=225] (10.5,-1.5) ;
\node (Lmius) at (9.75,-.75) {$\sL_{t,\theta}^-$};

\draw [thick] (12,-1.25) circle [radius=.35] node [left=7] {$\sS_{t,\theta}^-$};
\draw [thick] (12,-.2) circle [radius=.35] node [left=7] {$\sS_{t,\theta}^+$};

%\node (C) at (6,0) {$2|1\EFT^n(M)$};
%\node (D) at (9,0) {$\TMF^{-n}(M)$};
%\draw[->>,dashed] (A) to node [above] {\small cocycle}  (B);
%\draw[->>,dashed] (C) to node [above] {\small cocycle} (D);
\path (0,-.75) coordinate (basepoint);
\end{tikzpicture}\label{diag:supergenerators}
%&&1|1\EFT(M)\stackrel{{\rm cocycle}}{\dashrightarrow}  \K(M)  \qquad  2|1\EFT(M) \stackrel{{\rm cocycle}}{\dashrightarrow} \TMF(M)
\eeq
much the same as~\eqref{diag:generators}, but where the parameter~$t$ measuring the Riemannian length is replaced by a pair $(t,\theta)\in \R^{1|1}$ measuring a ``super" length. This super length has the usual even parameter $t$, but also an odd parameter~$\theta$; these parameters determine (super) families of bordisms as we explain in~\S\ref{eq:superpath}. It turns out to be convenient\footnote{The isomorphism $\R^{1|1}_{\ge 0}\simeq (\R^{1|1}_{\le 0})^\op$ equivalently allows one to view the bordisms $\sI^-_{t,\theta},\sR^-_{t,\theta},\sL^-_{t,\theta}$ and $\sS^-_{t,\theta}$ with $(t,\theta)$ in~$\R^{1|1}_{\ge 0}$, but then compositions involving these bordisms use the opposite super semigroup structure. This is a non-issue for ordinary Euclidean bordisms, since $\R_{\ge 0}$ and $\R_{\le 0}$ are isomorphic semigroups.} for $\sI^+_{t,\theta},\sR^+_{t,\theta},\sL^+_{t,\theta}$ and $\sS^+_{t,\theta}$ to have super length parameters $(t,\theta)$ in~$\R^{1|1}_{\ge 0}$, whereas $\sI^-_{t,\theta},\sR^-_{t,\theta},\sL^-_{t,\theta}$ and $\sS^-_{t,\theta}$ have $(t,\theta)$ in~$\R^{1|1}_{\le 0}$. The supermanifolds underlying the bordisms $\sI_{t,\theta}^\pm$, $\sR_{t,\theta}^\pm$ and~$\sL_{t,\theta}^\pm$ are again all the same, being super intervals. We similarly find relations associated with compositions of bordisms, e.g.,
\beq
&&\sI_{s,\eta}^\pm\circ \sI_{t,\theta}^\pm \simeq \sI_{s+t+\eta\theta,\eta+\theta}^\pm,\quad \sL_0^\pm\circ \sigma^\pm\simeq \sL_0^\mp,\quad  \sL_0^\pm \circ (\sI_{t,\theta}^\pm \coprod \sI_{s,\eta}^\mp)\simeq \sL_{(s,\eta)^{-1}(t,\theta)}^\pm \label{eq:11generatoreasy}
%\qquad (I_t^+\coprod I_s^-) \circ R_r^+=R_{r+s+t}^+,\qquad L_t^\pm\circ R_s^\mp=S_{t+s}^\pm.\label{eq:relations}
\eeq
where in the middle identity, $\sigma^\pm \colon \spt^\mp \coprod \spt^\pm \stackrel{\sim}{\to} \spt^\pm \coprod \spt^\mp$ is again the symmetry for the monoidal structure. 
The first identity is the geometric counterpart to multiplication in the super semigroups $\R_{\ge 0}^{1|1}$ from~\eqref{eq:N1super}. Note that the gluing of super intervals encodes a non-commutative super group; this is expected  as~\eqref{eq:Lie11} is not an abelian super Lie algebra. 

With this in place, we may consider the category $\Fun^\otimes (1|1\EBord,\Vect)$ whose objects are symmetric monoidal functors
\beq
E\colon 1|1\EBord\to \Vect\label{eq:repof1EB}.
\eeq
%In this case, restriction determines a functor
%\beq
%&&\begin{tikzpicture}[baseline=(basepoint)];
%\node (A) at (0,0) {$\Fun^\otimes (1|1\EBord,\Vect)$};
%\node (B) at (5,0) {$\Rep(\R_{\ge 0}^{1|1}),$};
%\node (C) at (8,0) {$E\mapsto E(I_{t,\theta}^+)$};
%\draw[->] (A) to node [below] {\small restrict} (B);
%\draw[->,dashed, bend right=15] (B) to (A);
%\path (0,-.05) coordinate (basepoint);
%\end{tikzpicture}\label{eq:restrict11EB0}
%%&&1|1\EFT(M)\stackrel{{\rm cocycle}}{\dashrightarrow}  \K(M)  \qquad  2|1\EFT(M) \stackrel{{\rm cocycle}}{\dashrightarrow} \TMF(M)
%\eeq
%with values in representations of the super semigroup $\R_{\ge 0}^{1|1}$. As before, lifting a super semigroup representation to a functor in~\eqref{eq:restrict11EB0} is additional data and property; 
For the data of a functor~\eqref{eq:repof1EB} to properly match with the physicists notion of a supersymmetric quantum mechanical theory, we again must impose a version of reflection positivity; this is defined analogously to the case without supersymmetry. We define a \emph{$1|1$-dimensional Euclidean field theory} to be a reflection-positive functor~\eqref{eq:repof1EB}.
One consequence is that a $1|1$-dimensional Euclidean field theory determines a self-adjoint representation of the super semigroup $\R_{\ge 0}^{1|1}$ on a $\Z/2$-graded topological vector space~$\mathcal{H}$ with hermitian inner product $\langle-,-\rangle$.  This data satisfies a trace property on $\R_{>0}^{1|1}\subset \R^{1|1}_{\ge 0}$, where the trace is the value of the functor~\eqref{eq:repof1EB} on a super circle (which turns out to be a constant function of~$(t,\theta)$). Supersymmetric quantum mechanics determines a reflection positive functor~\eqref{eq:repof1EB} with values
\beq
&&\begin{array}{c}
E(\spt^+)=\Gamma(X;\slashed{S}), \quad E(\sI_{t,\theta}^+)=e^{-t\slashed{D}+\theta\slashed{D}},\quad E(\sL_0^-)(s_1,s_2)=\int_X ( s_1,s_2) d{\rm vol}.
\end{array}\label{eq:supergensandrelations}
%\colon V\to V, \quad E(I_{t,\theta}^-)=\overline{e^{-t\slashed{D}+\theta\slashed{D}}}\colon \overline{V}\to \overline{V}
\eeq

\begin{rmk}
Restricting~\eqref{eq:N1super} along the homomorphism $\R_{\ge 0}\subset \R^{1|1}_{\ge 0}$ of semigroups forgets the supersymmetry, extracting an ordinary quantum mechanical theory. Similarly, any $1|1$-dimensional Euclidean field theory determines a 1-dimensional Euclidean field theory; see~\cite[Equation~4.14]{ST11}. 
%. This comes from restriction along a functor $1\EBord\to 1|1\EBord$
\end{rmk}

\begin{rmk}
Supersymmetric quantum mechanics can be understood entirely from the vantage of super Lie algebra representations~\eqref{eq:Lie11} or super semigroup representations~\eqref{eq:N1super}. Our primary motivation for pursuing a description in terms of the $1|1$-Euclidean bordism category is that there is an evident $2|1$-Euclidean generalization, illuminating a candidate definition for $2|1$-Euclidean field theory in the conjecture~\eqref{eq:conjecture}. 
\end{rmk}

\subsection{The renormalization group flow, low-energy cutoffs, and K-theory}\label{sec:RGcutoffmot}

Given $\mu\in \R_{>0}$, there is a functor $1\EBord\to 1\EBord$ that dilates the metric on Euclidean bordisms by the factor $\mu$. In terms of the elementary bordisms~\eqref{diag:generators}, this has the effect
\beq
&&\pt^\pm\mapsto \pt^\pm,\qquad \sI_t^\pm\mapsto \sI_{\mu t}^\pm,\qquad \sR_t^\pm\mapsto \sR_{\mu t}^\pm,\qquad \sL_t^\pm\mapsto \sL_{\mu t}^\pm, \quad \sS_t\mapsto \sS_{\mu t}.\nonumber
\eeq
%The functors $\RG^\mu$ compose as $\RG^{\mu_1}\circ\RG^{\mu_1}\simeq \RG^{\mu_1\mu_2}$, giving an $\R_{>0}$-action on $1\EBord(M)$. 
Precomposing a functor~\eqref{eq:repof1EB0} with this action defines the \emph{renormalization group flow} acting on 1-dimensional Euclidean field theories. A field theory is invariant under the renormalization group (RG) flow if and only if $E(\sI_t^+)=\id_{\mathcal{H}}$ for all $t$. The trace class condition on $E(\sI_t^+)$ implies that RG-invariant field theories have a finite-dimensional space of states~$\mathcal{H}$, and hence determine a 1-dimensional topological field theory. 

Similarly, there are functors $1|1\EBord\to 1|1\EBord$ that dilate the super lengths as $(t,\theta)\mapsto (\mu^2t,\mu\theta)$. On the elementary bordisms~\eqref{diag:supergenerators} this gives
\beq
\pt^\pm\mapsto \pt^\pm,\qquad \sI_{t,\theta}^\pm\mapsto \sI_{\mu^2 t,\mu\theta}^\pm,\qquad \sR_{t,\theta}^\pm\mapsto \sR_{\mu^2 t,\mu\theta}^\pm,\qquad \sL_{t,\theta}^\pm\mapsto \sL_{\mu^2 t,\mu\theta}^\pm, \quad \sS_{t,\theta}^\pm\mapsto \sS_{\mu^2 t,\mu\theta}^{\pm}. \nonumber
\eeq
The nontrivial action on the odd part of the super length~$\theta$ is necessary to make the dilation action compatible with composition (and hence a functor).  Precomposing~\eqref{eq:repof1EB} with the above functor defines the renormalization group flow on $1|1$-dimensional super Euclidean field theories. In the presence of supersymmetry, the partition function is invariant under the renormalization group flow: it is a constant function of the geometric parameters $(t,\theta)$. 
%In particular, the partition function is an RG-invariant quantity. 

Given a quantum mechanical system of the form~\eqref{eq:semigroup} or~\eqref{eq:N1super}, let $\mathcal{H}_\lambda\subset \mathcal{H}$ denote the $\lambda$-eigenspace of the self-adjoint Hamiltonian $H$ (so $\lambda\in \R$). Then $\mathcal{H}_\lambda$ is the space of \emph{states with energy~$\lambda$}. Define the vector space of states with energy $<\lambda$,
\beq
\mathcal{H}^{<\lambda}:=\bigoplus_{\nu<\lambda} \mathcal{H}_\nu\label{eq:lambdacut}
\eeq
with hermitian inner product given by the restriction of the inner product on $\mathcal{H}$ to the subspaces $\mathcal{H}^{<\lambda}\subset \mathcal{H}$. The trace class condition on $e^{-tH}$ implies that eigenspaces of $H$ are finite-dimensional, the spectrum of $H$ is discrete, and $\mathcal{H}_\lambda=\{0\}$ for $\lambda\ll 0$. Hence, $\mathcal{H}^{<\lambda}$ is finite-dimensional. Define \emph{$\lambda$-cutoff theory} as the quantum mechanical system given by the finite-dimensional vector space $\mathcal{H}^{<\lambda}$ equipped with the identity time-evolution operator. The $\lambda$-cutoff theory is therefore invariant under the renormalization group flow, i.e., a 1-dimensional topological field theory. 

\begin{rmk}
One could equally equip $\mathcal{H}^{<\lambda}$ with the restricted Hamiltonian $H|_{\mathcal{H}^{<\lambda}}$ when defining the cutoff theory. 
For our purposes, it turns out to be more convenient for the cutoff theory to be RG-invariant. 
% that taking the identity Hamiltonian communicates slightly better with index theory.
\end{rmk}

Since $\mathcal{H}$ is $\Z/2$-graded and~$H$ is an even map, $\mathcal{H}_\lambda$ is $\Z/2$-graded for each~$\lambda\in \R$. Given a fixed choice of $\lambda$ and 1- or $1|1$-dimensional Euclidean field theory~$E$, consider the assignment
\beq
E\mapsto [\mathcal{H}^{<\lambda}]=[(\mathcal{H}^{<\lambda})^\ev]-[(\mathcal{H}^{<\lambda})^\odd]\in \K^0(\pt).\label{eq:maptoKO}
\eeq
In the absence of supersymmetry, the map~\eqref{eq:maptoKO} generally depends on~$\lambda$. However, when $H$ has an odd square root~$Q$, then the assignment is independent of~$\lambda$: $Q$ gives isomorphisms 
$$
Q\colon \mathcal{H}_\lambda^{\ev/\odd}\stackrel{\sim}{\to} \mathcal{H}_\lambda^{\odd/\ev}, \qquad \lambda\ne 0
$$
and so 
$$
[\mathcal{H}^{<\lambda}]=[\mathcal{H}_0]\in \K^0(\pt)\simeq \Z
$$ 
for all $\lambda$. This integer-valued invariant in supersymmetric quantum mechanics is the \emph{Witten index}, originally studied in~\cite{susymorse}. The Witten index is a deformation invariant of a $1|1$-dimensional Euclidean field theory. 

\begin{ex} 
For the theory with $Q=\slashed{D}$ the Dirac operator on an even dimensional spin manifold, the Witten index agrees with the analytic index $\Ind(\slashed{D})=[\mathcal{H}_0]\in \K^0(\pt)$ of the Dirac operator. In the example discussed in Remark~\ref{rmk:SUSYmorse}, the Witten index is the Euler characteristic of~$X$. 
\end{ex} 

By contrast, 1-dimensional Euclidean field theories without supersymmetry form a contractible space~\cite[Theorem 6.35]{HST} and so have no interesting deformation invariants.

\begin{rmk}
We note that the McKean--Singer formula~\eqref{eq:partition function2} equates the Witten index and the partition function as elements of~$\Z$. In families, this identification is modified: the partition function is a differential form which agrees with the Chern character of the Witten index. \end{rmk}
\subsection{Symmetries and real structures}\label{eq:CliffordlinearFT} 
Given a quantum mechanical system in the form~\eqref{eq:semigroup} or~\eqref{eq:N1super}, a super algebra $A$ acting on $\mathcal{H}$ is an algebra of \emph{symmetries} if it commutes with the action of $H$ or~$Q$. Reflection positivity for a theory with $A$-symmetry again endows~$\mathcal{H}$ with a hermitian inner product, and additionally requires the data of a $*$-structure on~$A$ and the property that the $A$-action on $\mathcal{H}$ be a $*$-action. 

It is often convenient to organize field theories according to their symmetries: for a given super algebra $A$ one can consider the collection of theories with $A$-symmetry. Stolz and Teichner's degree~$n$ field theories are supersymmetric quantum mechanical systems~\eqref{eq:N1super} with symmetry specified by the the $n$th (complex) Clifford algebra
\beq
\cCl_n:=\left\{\begin{array}{ll} \langle f_1,\dots f_n \mid [f_j,f_k]=-2\delta_{jk}\rangle & n \ge 0 \\ \langle e_1,\dots e_n \mid [e_j,e_k]=2\delta_{jk}\rangle & n<0.\end{array}\right.
\label{eq:Clifford}
\eeq
where $\{e_1,\dots, e_n\}$ or $\{f_1,\dots,f_n\}$ denotes the standard basis of $\R^n$. These Clifford algebras carry a standard $*$-structure; see~\S\ref{sec:superstar}. 

One can also ask for a \emph{real structure} on a quantum mechanical system, meaning a grading-preserving $\C$-antilinear involution of $\mathcal{H}$ commuting with the Hamiltonian. When considering theories with $A$-symmetry, one should demand a compatible real structure on~$A$. In the case of the Clifford algebras, this is equivalent to declaring the generators in~\eqref{eq:Clifford} to be real so that the real structure is determined by complex conjugation. 
%
%\beq
%\Cl_n:=\left\{\begin{array}{ll} \langle e_1,\dots, e_n \mid e_je_k+e_ke_j=-2\delta_{jk}\rangle & n \ge 0 \\ \langle e_1,\dots, e_n \mid e_je_k+e_ke_j=2\delta_{jk}\rangle & n<0\end{array}\right.
%\label{eq:rClifford}
%\eeq
Using the Atiyah--Bott--Shapiro description of $\KO^n(\pt)$ in terms of Clifford modules~\cite{ABS}, such real, degree~$n$ theories determine classes in $\KO^n(\pt)$ using the construction~\eqref{eq:maptoKO}. 
% the previous discussion connecting supersymmetric quantum mechanics with $\K^0(\pt)$ can be generalized. Indeed, a
%We briefly indicate the geometry behind this generalization. 

Examples of real, degree~$n$ theories arise naturally in spin geometry. 
%The suitable generalization of supersymmetric quantum mechanics on a spin manifold comes by way of the Clifford linear Dirac operator. 
Let $X$ be an $n$-dimensional spin manifold with principal spin bundle~$P$, i.e., $P\to X$ is a $\Spin(n)$-principal bundle that double covers of the oriented frame bundle of~$TX$. Define the spinor bundle as
\beq
\bS:=P \times_{\Spin(n)} \Cl_n=P \times_{\Spin(n)} (\Cl_n^+\oplus \Cl_n^-)\label{eq:Cliffspin}
\eeq
with $\Z/2$-grading inherited from the $\Z/2$-grading on the real Clifford algebra $\Cl(\R^n)=\Cl_n$. This version of the spinor bundle carries a right $\Cl_n$-action, which we will identify with a left $\Cl_{-n}$-action (see Example~\ref{ex:leftright}). The odd, self-adjoint, $\Cl_{-n}$-linear Dirac operator~$\slashed{D}$ acts on global sections $\Gamma(X;\bS)$. The same constructions as before give a quantum mechanical system with additional structures: the space of states $\mathcal{H}=\Gamma(X;\slashed{S})$ is a real, self-adjoint $\Cl_{-n}$-module and the operator $Q=\slashed{D}$ commutes with the Clifford action.
% super semigroup representation~\eqref{eq:N1super} is $\Cl_n$-linear. 
 The energy cutoff~\eqref{eq:lambdacut} gives a finite-dimensional Clifford module. Via the Atiyah--Bott--Shapiro construction, the assignment~\eqref{eq:maptoKO} is the \emph{Clifford index} of the Dirac operator, taking values in $\KO^{-n}(\pt)$. This fundamental geometric example motivates Stolz and Teichner's definitions of degree~$n$ field theories in dimension~$1|1$, see~\cite[Example~2.3.3]{ST04}. 

\begin{rmk}
For an inner product space $W$, the Clifford algebra $\Cl(W)$ can be understood as the quantization of a free (fermionic) particle moving in the odd vector space $\Pi W$, e.g., see~\cite[page~19]{Freed5}. Free fermions in higher dimensions illuminate a path towards a general framework for degree~$n$ field theories. In particular, this suggests a $2|1$-dimensional generalization for the second conjecture in~\eqref{eq:conjecture}. 
\end{rmk}

\section{Super Euclidean paths in a manifold}\label{sec:superEuclideanpath}
\subsection{Warm-up: The category of (non-super) Euclidean paths in a manifold} \label{eq:ordpath}
Let~$\A^1$ be the affine real line, regarded as an oriented Riemannian manifold. The isometry group of $\A^1$ is the Lie group~$\R$, acting on $\A^1$ by translation. We will describe Euclidean paths modeled on~$\A^1$ with an eye towards the generalization to $1|1$-Euclidean paths. 

The Riemannian manifold $\A^1$, isometry group $\R$, and a choice of connected codimension~1 submanifold $\A^0\subset \A^1$ determines the 1-dimensional Euclidean geometry $(\R,\A^1,\A^0)$ in the sense of Stolz and Teichner~\cite[\S2.5]{ST11}; see~\S\ref{sec:rigidgeometry} for a review. Below, for $s\in \R$ let $T_s\colon \A^1\to \A^1$ denote translation by~$s$.

%We begin by describing a Lie category of constant (non-super) paths in a manifold. 

\begin{defn}
Let $M$ be a smooth manifold. A \emph{positively oriented Euclidean path} in~$M$ is a triple $(\phi,t_\inn,t_\out)$ for a smooth map $\phi\colon \A^1\to M$ and points~$t_\inn,t_\out\in \A^1$ with $t_\out-t_\inn\in \R_{\ge 0}$. Similarly, a \emph{negatively oriented Euclidean path} in $M$ is a triple $(\phi,t_\inn,t_\out)$ with $t_\out-t_\inn\in \R_{\le 0}$. A \emph{Euclidean path} is a positively or negatively oriented Euclidean path. The point~$\phi(t_\inn) \in M$ is the \emph{incoming} boundary and~$\phi(t_\out)\in M$ is the \emph{outgoing} boundary of a Euclidean path. 
%positively oriented path is $t\in \R_{\ge 0}$ and for a negative oriented path is $-t\in \R_{\ge 0}$. 
An \emph{isometry} $(\phi,t_\inn,t_\out)\stackrel{\sim}{\to} (\phi',t_\inn',t_\out')$ between a pair of Euclidean paths is the data of $s \in \R$ such that $\phi'=\phi\circ T_s^{-1}$, $T_s(t_\inn)=t_\inn'$, and $T_s(t_\out)=t_\out'$. 
% where $T_s\colon \A\to \A$ is translation by~$s$. 
\end{defn}

Euclidean paths and isometries form a groupoid that is equivalent to a discrete groupoid (i.e., a groupoid with only identity morphisms) whose objects are
\beq\label{eq:bigEpath}	
\R_{\ge 0}\times \Map(\R,M)\coprod \R_{\le 0}\times \Map(\R,M),
\eeq
where we regard $\Map(\R,M)$ as a presheaf on the site of manifolds,
$$
\Map(\R,M)\colon {\sf Mfld}^\op\to {\sf Set},\qquad S\mapsto C^\infty(S\times \R,M).
$$ 
%To verify that~\eqref{eq:bigEpath} is indeed equivalent to the groupoid of Euclidean paths, 
Indeed, for a given Euclidean path $(\phi,t_\inn,t_\out)$, the basepoint $t_\inn\in \A^1$ gives an identification $\A^1\simeq \R^1$, so that $\phi\colon \R\to M$ and $t_\out$ is determined by a point in $\R_{\ge 0}\coprod \R_{\le 0}$. 
%Given a positively oriented path $(\phi,t_\inn,t_\out)$, we have the identification $\A^1\simeq \R^1$ from the choice of basepoint $t_\inn\in \A^1$. Then $t=t_\out \in \R_{\ge 0}\subset \R$ can be identified with the length of the path. Similarly, given a negatively oriented path we have the identification $\A^1\simeq \R^1$ for the basepoint $t_\inn\in \A^1$, and~$t=t_\out\in \R_{\le 0}\subset \R$ is the negative length. 

The 1-dimensional Euclidean geometry has some additional structures that lead to additional structures on Euclidean paths. First, the basepoint $\A^0\subset \A^1$ provides an isomorphism $\A^1\simeq \R^1$ (that may be different than the identification using $t_\inn$). This gives a preferred orientation-reversing map of~$\A^1$,
\beq\label{Eq:easyorr}
\A^1\simeq \R^1\xrightarrow{\orr} \R^1\simeq \A^1, \qquad \orr(t)=-t, \ t\in \R. 
\eeq
%A different choice of basepoint $\A^0\hookrightarrow \A^1$ gives an orientation reversing automorphism 
%We emphasize that we regard $\A^1$ as a \emph{oriented} Riemannian manifold so that~\eqref{Eq:easyorr} is not an isometry. 
%Changing the basepoint $\A^0\subset \A^1$ result in a different identification $\A^1\simeq \R^1$, and the resulting orientation-reversing map differs from~\eqref{Eq:easyorr} by an isometry of $\A^1$. 
Next, for any $\mu\in \R_{>0}$ there is a diffeomorphism
\beq\label{Eq:easyrg}
\A^1\simeq \R^1\xrightarrow{\rg_\mu} \R^1\simeq \A^1,\qquad \rg_\mu(t)=\mu\cdot t
\eeq
that dilates the metric on $\A^1$. The notation $\rg_\mu$ is in reference to the connection with the renormalization group. We emphasize that~\eqref{Eq:easyorr} and~\eqref{Eq:easyrg} are not isometries in the 1-dimensional Euclidean geometry $(\R,\A^1,\A^0)$.

The structures~\eqref{Eq:easyorr} and~\eqref{Eq:easyrg} provide the following additional structures on Euclidean paths. The \emph{orientation reversal} of a Euclidean path $(\phi,t_\inn,t_\out)$ applies the orientation-reversing isometry~\eqref{Eq:easyorr}, giving an assignment
$$
(\phi,t_\inn,t_\out)\mapsto (\phi\circ\orr^{-1},\orr(t_\inn),\orr(t_\out)).
$$
For $\mu\in \R_{>0}$ the \emph{$\mu$-dilation} of a Euclidean path $(\phi,t_\inn,t_\out)$ applies the diffeomorphism~\eqref{Eq:easyrg} to a Euclidean path so that 
$$
(\phi,t_\inn,t_\out)\mapsto(\phi\circ\rg_\mu^{-1},\rg_\mu(t_\inn),\rg_\mu(t_\out)). 
$$
In terms of the description~\eqref{eq:bigEpath} of Euclidean paths, orientation reversal and $\mu$-dilation are 
\beq
\R_{\ge 0}\times \Map(\R,M)\coprod \R_{\le 0}\times \Map(\R,M)&\to& \R_{\ge 0}\times \Map(\R,M)\coprod \R_{\le 0}\times \Map(\R,M)\nonumber\\
(t,\phi)&\mapsto&(-t,\phi\circ \orr^{-1})\label{eq:Estructures}\\
(t,\phi)&\mapsto&(\mu\cdot t,\phi\circ \rg_\mu^{-1}).\nonumber
\eeq
In particular, orientation reversal exchanges positively and negatively oriented Euclidean paths. The incoming and outgoing boundary $\phi(t_\inn),\phi(t_\out)\in M$ are preserved by orientation reversal and $\mu$-dilation. 

Next, we consider the extent to which Euclidean paths can be viewed as morphisms in a category. There are maps 
$$
\inn,\out\colon \R_{\ge 0}\times \Map(\R,M)\coprod \R_{\le 0}\times \Map(\R,M)\rightrightarrows M\coprod M
$$
sending a Euclidean path to its incoming or outgoing boundary. We would like the above to determine the source and target maps in a category in which composition comes from concatenation of paths. However, this runs into a problem: concatenation is only piecewise smooth. There are different ways to resolve this problem. One option is to consider collared paths in~$M$, which is the route followed by Stolz and Teichner~\cite{ST11}; see~\S\ref{sec:STBord}. Another option is to encode cutting laws rather than gluing laws as in~\cite{GradyPavlov,LudewigStoffel,BEPTFT}. The resulting bordism categories are not equivalent, and it isn't clear which one is ``correct" when generalizing the quantum mechanical systems discussed in~\S\ref{sec:QM}. However, there is a subspace of~\eqref{eq:bigEpath} where the difficulties of concatenation disappear: the \emph{constant} Euclidean paths in~$M$. 
%It turns out that this subspace and its super geometric enhancement already capture much of the information relevant for index theory. 

%A cruder option is to restrict attention to a class of paths for which concatenation is always smooth
%Of course, the category of constant paths is not very interesting geometrically. However, its super generalization turns out to 
% that its representation theory affords a families generalization of the quantum mechanical theories reviewed in~\S\ref{sec:QM}, and the generalization to superpaths is sufficient for our purposes.
% (and closely related to~\cite[\S3.1]{ST04}). 

\begin{defn}\label{defn:Path} Define the (Lie) \emph{category of constant paths in $M$}, denoted $\Path_0(M)$, as having objects and morphisms the smooth manifolds
$$
\Ob(\Path_0(M))=M\coprod M,\qquad \Mor(\Path_0(M)) =\R_{\ge 0}\times M\coprod \R_{\le 0}\times M.
$$
%where a point in the first component $(t,x)\in \R_{\ge 0} \times M\subset \Mor(\Path_0(M)) $ determines a positively oriented path of length $t$ mapped constantly to $x$, while the same data in the second component determines a negatively oriented superpath. 
The source and target maps in $\Path_0(M)$ are given by the projections.
% where the source and target of a positively (respectively, negatively) oriented constant path is a positively (respectively, negatively) oriented point. 
 The unit map is determined by inclusion at $0\in \R$, and composition is determined by addition in $\R$. 
\end{defn}

Under the inclusion of the constant maps
$$
\Mor(\Path_0(M)) =\R_{\ge 0}\times M\coprod \R_{\le 0}\times M\subset \R_{\ge 0}\times \Map(\R,M)\coprod \R_{\le 0}\times \Map(\R,M)
$$
concatenation of paths corresponds to composition in the category $\Path_0(M)$. The assignments~\eqref{eq:Estructures} furthermore extend to functors
\beq\label{eq:lookfunctors}
\Or\colon \Path_0(M)\to \Path_0(M),\quad \RG_\mu \colon \Path_0(M)\to \Path_0(M). 
\eeq

Next we explain the relationship between Euclidean paths and Euclidean loops. There is a $\Z$-action on~\eqref{eq:bigEpath} given by $(t,\phi)\mapsto (t,n_t\cdot \phi)$ where $n_t\cdot \phi\colon \R\to M$ is the map
$$
(n_t\cdot \phi)(s)=\phi(s+nt),\qquad n\in \Z. 
$$
For paths of strictly positively length, the $\Z$-fixed subspace (as a presheaf) is the \emph{Euclidean loop space} of~$M$
$$
\mathcal{L}(M):=(\R_{> 0}\times \Map(\R,M)\coprod \R_{< 0}\times \Map(\R,M))^\Z\subset \R_{\ge 0}\times \Map(\R,M)\coprod \R_{\le 0}\times \Map(\R,M)
$$ 
whose points $(t,\phi)$ are loops in $M$ of circumference $t>0$. Furthermore, $\mathcal{L}(M)$ has an $S^1$-action by rotation of loops. Any constant path of positive length determines a loop, giving a span
\beq\label{eq:loopspace1}
\Mor(\Path_0(M))\hookleftarrow \R_{>0} \times M\coprod \R_{<0}\times M \hookrightarrow \mathcal{L}(M),
\eeq
and the $S^1$-action on $\mathcal{L}(M)$ restricts to the trivial action on $\Mor(\Path_0(M))$.

In summary, although constant Euclidean paths lack the geometric richness of the full path space, they provide two technical simplifications: (1) constant paths form a finite-dimensional smooth manifold, and (2) concatenation of constant paths always results in a (smooth) constant path. When $M=\pt$, the category $\Path_0(M)$ recovers the semigroup of intervals~\eqref{eq:semigroup}, and hence the categories $\Path_0(M)$ (and their representations) offer one possible families-generalization of the quantum mechanical systems studied previously. With this in mind, we purse a generalization of the above; the goal of the remainder of the section is to generalize Definition~\ref{defn:Path} to a category of (nearly) constant super Euclidean paths in~$M$.
% with functors~\eqref{eq:lookfunctors} and a map to the super loop space of $M$ generalizing~\eqref{eq:loopspace1}

\subsection{$1|1$-dimensional Euclidean geometry}\label{sec:11Eucgeo}
Recall from~\eqref{eq:superEuc} the super Lie group $\R^{1|1}$, and let $\A^{1|1}$ denote the (affine) supermanifold underlying the super Lie group $\R^{1|1}$. This notation emphasizes that~$\R^{1|1}$ has a distinguished basepoint (the identity $0\in \R^{1|1}$) whereas~$\A^{1|1}$ does not. Similarly, let $\A^{0|1}$ denote the supermanifold underlying the super Lie group~$\R^{0|1}$. Let $\A^{0|1}\subset \A^{1|1}$ be the embedding modeled on the standard\footnote{We caution that $\R^{0|1}\subset \R^{1|1}$ is not a homomorphism for the group structure~\eqref{eq:superEuc} on $\R^{1|1}$.} inclusion $\R^{0|1}\subset \R^{1|1}$. We again refer to~\S\ref{sec:rigidgeometry} for a review of rigid super geometries.

\begin{defn}[{\cite[Example 6.16]{HST} and \cite[\S4.2]{ST11}}] \label{ex:modelgeo} Form the semidirect product $\R^{1|1}\rtimes \Z/2$ where $\Z/2=\{\pm1\}$ acts on $\R^{1|1}$ by $(t,\theta)\mapsto (t,\pm \theta)$. The 
\emph{$1|1$-Euclidean geometry} is specified by the triple
$$
(G,\M,\M^c)=(\R^{1|1}\rtimes \Z/2,\A^{1|1},\A^{0|1})
$$ 
for the evident action of $\R^{1|1}\rtimes \Z/2$ on $\A^{1|1}$ and the inclusion $\A^{0|1}\subset \A^{1|1}$ described above.
% that after a choice of basepoint on $\A^{0|1}$ is the standard inclusion $\R^{0|1}\subset \R^{1|1}$. 
\end{defn}

In the notation of Definition~\ref{defn:rigid}, the stack $(\R^{1|1}\rtimes \Z/2,\A^{1|1})$-$\SMfld$ of $1|1$-Euclidean manifolds has as objects submersions $Y\to S$ whose fibers are $1|1$-dimensional manifolds endowed with a super Euclidean structure. The stack of $1|1$-Euclidean pairs has as objects fiber bundles $Y\to S$ of $1|1$-manifolds endowed with super Euclidean structure together with a family of embedded $0|1$-dimensional supermanifolds $Y^c\subset Y$ over $S$. 

\begin{ex} \label{ex:Eucpair}
Our main example of an $S$-family of $1|1$-Euclidean pairs depends on $(t,\theta)\in \R^{1|1}(S)$ and is given by 
\beq\label{eq:Eucpairmain}
Y^c=S\times \R^{0|1}\xhookrightarrow{i_{t,\theta}}S\times \R^{1|1}=Y
\eeq
where
\beq\label{eq:itthetadef}
i_{t,\theta}\colon S\times \R^{0|1}\subset S\times \R^{1|1}\xrightarrow{T_{t,\theta}} S\times \R^{1|1}
\eeq
using the standard inclusion and the super translation map $T_{t,\theta}$ defined as the composition
\beq
&&S \times \R^{1|1}\xrightarrow{\id_S\times (t,\theta)\times \id_{\R^{1|1}}} S\times \R^{1|1} \times \R^{1|1}\xrightarrow{\id_S\times{\rm act}} S\times \R^{1|1},\label{eq;translation}
\eeq
where the second arrow is the action of $\R^{1|1}$ on itself. For any $(t,\theta)$ and $(s,\eta)$ the associated $1|1$-Euclidean pairs are isomorphic via the isometry
$$
S\times \R^{1|1}\xrightarrow{T_{(s,\eta)\cdot (t,\theta)^{-1}}} S\times \R^{1|1}.
$$
In particular, all pairs~\eqref{eq:Eucpairmain} are isomorphic to the standard one $S\times \R^{0|1}\subset S\times \R^{1|1}$. 
%
%We describe the main examples of $S$-families of $1|1$-Euclidean pairs in this paper. First we need some notation. Given $(t,\theta,\pm 1)\in (\R^{1|1}\rtimes \Z/2)(S)$, let $T_{t,\theta,\pm 1}\colon S\times \A^{1|1}\to S\times \A^{1|1}$ denote the map over $S$ defined by 
%
% If a choice of identification $\A^{1|1}\simeq \R^{1|1}$ has been fixed, we use the same notation for the corresponding map $T_{t,\theta,\pm 1}\colon S\times \R^{1|1}\to S\times \R^{1|1}$. For an $S$-point of the subgroup $\R^{1|1}<\R^{1|1}\rtimes \Z/2$, we use the notation~$T_{t,\theta}$. Now for any $(t,\theta,\pm1)\in (\R^{1|1}\rtimes \Z/2)(S)$, 
 \end{ex}
%We use the simpler notation $T_{s,\eta}$ in the special case that $(s,\eta,\pm 1)=(s,\eta,+1)$. 

The stack of $1|1$-dimensional super Euclidean manifolds has some extra structures generalizing orientation reversal~\eqref{Eq:easyorr} and dilation~\eqref{Eq:easyrg} of ordinary Euclidean 1-manifolds. 

\begin{defn} \label{defn:orientationbord}\emph{Orientation reversal} is the automorphism of the $1|1$-Euclidean geometry determined by the homomorphism
\beq
%\R^1\to \R^1,&&\qquad t\mapsto -t, \qquad t\in \R^1(S)\\
\orr \colon \R^{1|1}\to \R^{1|1},&&\qquad (t,\theta)\mapsto (-t,-i\theta),\qquad (t,\theta)\in \R^{1|1}(S).\label{eq:supertimereverse}
\eeq
We observe that when restricted to $\R\subset \R^{1|1}$~\eqref{eq:supertimereverse} is the standard orientation reversing map~\eqref{Eq:easyorr} on $\R$.
\end{defn}

\begin{rmk}
The map~\eqref{eq:supertimereverse} is called the \emph{pin generator} in~\cite[Example 6.16]{HST}.
\end{rmk}

\begin{defn} \label{defn:rg}For $\mu\in \R_{>0}$, the \emph{renormalization group action} is the automorphism of the $1|1$-Euclidean geometry determined by the homomorphism 
\beq
%\R^1\to \R^1,&&\qquad t\mapsto -t, \qquad t\in \R^1(S)\\
\rg_\mu \colon \R^{1|1}\to \R^{1|1},&&\qquad (t,\theta)\mapsto (\mu^2t,\mu\theta),\qquad (t,\theta)\in \R^{1|1}(S).\label{eq:superRG}
\eeq
\end{defn}

\begin{defn}[{\cite[Example~2.39 and \S4.1]{ST11}}] \label{defn:spinflip}The \emph{spin flip} is the automorphism of the $1|1$-Euclidean geometry determined by the action of $-1\in \{\pm 1\}\simeq \Z/2<\R^{1|1}\rtimes \Z/2$
\beq
\fl\colon \R^{1|1}\to \R^{1|1},\qquad (t,\theta)\mapsto (t,-\theta),\qquad (t,\theta)\in \R^{1|1}(S).\label{eq:spinflipfor11}
\eeq
Since $\Z/2<\R^{1|1}\rtimes \Z/2$ is in the center, the spin flip commutes with the action of the isometry group. 
\end{defn}

The $1|1$-Euclidean geometry also has a real structure in the sense of Definition~\ref{defn:realFT}. 

%\begin{rmk}
%In the notation of Definition~\ref{defn:11EBspt}, the automorphisms of the object $(\spt^\pm,\phi)\in 1|1\EBord(\pt)(S)$ (for the unique map $\phi\colon S\times \R^{1|1}\to \pt$) is isomorphic to $\Z/2(S)$ generated by the spin flip; this follows the same argument as \cite[Lemma~86]{HST}
%\end{rmk}

\begin{lem}[{\cite[Example 6.20]{HST}}]\label{rmk:real2}
The super Lie group $\R^{1|1}$ has a real structure
\beq\label{eq:real11}
\rr\colon \R^{1|1}\xrightarrow{\sim}  \overline{\R}^{1|1},\qquad \overline{\rr} \circ \rr=\id_{\R^{1|1}}.
\eeq
\end{lem}
\bp
Following Example~\ref{ex:real}, first choose an identification $\R^{1|1}\simeq \Pi (\underline{\R}_ \C)$ where $\underline{\R}$ is the trivial real line bundle on~$\R$ and $\underline{\R}_ \C$ is its complexification. This endows the supermanifold underlying $\R^{1|1}$ with a real structure. To promote this to a real structure on the super Lie group, let $\theta\in C^\infty(\R^{1|1})^\odd\simeq \Gamma(\R,\Lambda^1 (\underline{\R}_\C))$ denote the function associated with a (real) trivialization of $\underline{\R}$. With this choice, the group multiplication~\eqref{eq:superEuc} is real. 
%The orientation reversing automorphism~\eqref{eq:supertimereverse} does not preserve this real structure. 
\ep

\begin{defn} \label{defn:realbordisms}Define a real structure on the $1|1$-dimensional super Euclidean model geometry from the real structure on the super Lie group $\R^{1|1}\rtimes \Z/2$ from Lemma~\ref{rmk:real2}, using that $\Z/2$ preserves this real structure on $\R^{1|1}$. 
\end{defn}

\begin{ex}\label{ex:functorsuperpair}\label{ex:realstructureonSP}
By Lemma~\ref{lem:GMaddstructure1}, the maps \eqref{eq:supertimereverse}, \eqref{eq:superRG} and \eqref{eq:spinflipfor11} determine endofunctors of the stack of $1|1$-Euclidean manifolds,
$$
\Or,\RG_\mu,\Fl\colon (\R^{1|1}\rtimes \Z/2,\A^{1|1})\hbox{-}\SMfld\to (\R^{1|1}\rtimes \Z/2,\A^{1|1})\hbox{-}\SMfld,
$$
that send an $S$-family of $1|1$-Euclidean manifolds to another $S$-family. Applied to the super Euclidean pair~\eqref{eq:Eucpairmain}, we find
\beq
\Or(S\times \R^{0|1}\xhookrightarrow{i_{t,\theta}}S\times \R^{1|1})&=&(S\times \R^{0|1}\xhookrightarrow{i_{-t,-i\theta}}S\times \R^{1|1})\nonumber\\
\RG_\mu(S\times \R^{0|1}\xhookrightarrow{i_{t,\theta}}S\times \R^{1|1})&=&(S\times \R^{0|1}\xhookrightarrow{i_{\mu^2t,\mu\theta}}S\times \R^{1|1})\nonumber\\
\Fl(S\times \R^{0|1}\xhookrightarrow{i_{t,\theta}}S\times \R^{1|1})&=&(S\times \R^{0|1}\xhookrightarrow{i_{t,-\theta}}S\times \R^{1|1})\nonumber
\eeq
Similarly, the map~\eqref{eq:real11} determines a map of stacks 
$$
\RR\colon (\R^{1|1}\rtimes \Z/2,\A^{1|1})\hbox{-}\SMfld\to (\R^{1|1}\rtimes \Z/2,\A^{1|1})\hbox{-}\SMfld,\quad (Y^c,Y)\mapsto (\overline{Y}{}^c,\overline{Y})
$$
covering complex conjugation of supermanifolds, that sends an $S$-family of $1|1$-Euclidean manifolds to an $\overline{S}$-family. There are natural isomorphisms of functors
\beq\label{eq:Orrelation}
\Or^2\simeq \Fl,\quad \RG_\mu\circ \RG_{\mu'}\simeq \RG_{\mu\mu'}, \quad \RR \circ \Or\simeq \overline{\Fl}\circ \overline{\Or} \circ \RR, \quad \overline{\RR}\circ \RR\simeq \id
\eeq
coming from the equalities of maps
$$
\orr^2=\fl, \quad \rg_\mu\circ \rg_{\mu'}=\rg_{\mu\mu'}, \quad \overline{\orr}\circ \rr=\overline{\fl}\circ \rr\circ \orr, \quad \overline{\rr} \circ \rr=\id_{\R^{1|1}}. 
$$
\end{ex} 
% an automorphisms of the stacks of (super) Euclidean manifolds (e.g., following the discussion in~\cite[\S2.6]{ST11}). This then determines the internal functors~\eqref{eq:timereversal11EB} on the associated bordism categories. 

\subsection{Super Euclidean paths in a smooth manifold}\label{eq:superpath}

%Suppose we are given a 
%given a $1|1$-dimensional super Euclidean pair of the form
%$$
%S\times \A^{0|1}\coprod S\times \A^{0|1}\xhookrightarrow{t_\inn\coprod t_\out} S\times \A^{1|1}
%$$
Given a triple $(\phi,(t,\theta),(s,\eta))$ where $\phi\colon S\times \R^{1|1}\to M$ is a map and $(t,\theta),(s,\eta)\in \R^{1|1}(S)$, consider
\beq\label{diag:atriple}
\begin{tikzpicture}[baseline=(basepoint)];
\node (A) at (0,0) {$S\times \R^{0|1}$};
\node (B) at (2.5,0) {$S\times \R^{1|1}$};
\node (C) at (5,0) {$S\times \R^{0|1}$};
\node (D) at (2.5,-1.25) {$M$};
\draw[->,right hook-latex] (A) to node [above] {$i_{t,\theta}$} (B);
\draw[->,left hook-latex] (C) to node [above] {$i_{s,\eta}$} (B);
\draw[->] (B) to node [right] {$\phi$} (D);
\path (0,-.75) coordinate (basepoint);
\end{tikzpicture}
\eeq
for the maps defined in~\eqref{eq:itthetadef}. 

%is an inclusion, and $t_\out=T_{t,\theta}\circ t_\inn$ is this inclusion followed by a translation.

\begin{defn}
Let $M$ be a smooth manifold and $S$ a supermanifold. An $S$-family of \emph{positively oriented super Euclidean paths} in~$M$ is a triple $(\phi,(t_\out,\theta_\out),(t_\inn,\theta_\inn))$ determining a diagram~\eqref{diag:atriple} with the property 
\beq
%\label{eq:moreRG}
\nonumber
 (t_\out,\theta_\out)\cdot (t_\inn,\theta_\inn)^{-1} \in \R^{1|1}_{\ge 0}(S)\subset \R^{1|1}(S).
\eeq
%for $(t,\theta)$ from Lemma~\ref{lem:uniquetranslate}. 
Similarly, an $S$-family of \emph{negatively oriented super Euclidean paths} in~$M$ is a triple where 
\beq
%\label{eq:moreRG}
\nonumber
 (t_\out,\theta_\out)\cdot (t_\inn,\theta_\inn)^{-1} \in \R^{1|1}_{\le 0}(S)\subset \R^{1|1}(S).
\eeq
An $S$-family of \emph{super Euclidean paths} in~$M$ is a triple $(\phi,(t_\out,\theta_\out),(t_\inn,\theta_\inn))$ that over each component of $S$ is a family of positively or negatively oriented superpaths. The \emph{incoming} and \emph{outgoing} boundaries of $(\phi,(t,\theta),(s,\eta))$ are the compositions 
$$
\phi_\out \colon S\times \R^{0|1}\xhookrightarrow{i_{(t_\out,\theta_\out)}} S\times \R^{1|1}\xrightarrow{\phi} M,\quad \phi_\inn \colon S\times \R^{0|1}\xhookrightarrow{i_{(t_\inn,\theta_\inn)}} S\times \R^{1|1}\xrightarrow{\phi} M.
$$
An \emph{isometry} $(\phi,(t_\out,\theta_\out),(t_\inn,\theta_\inn))\stackrel{\sim}{\to}(\phi',(t_\out',\theta_\out'),(t_\inn',\theta_\inn'))$ between a pair of $S$-families of superpaths is the data of $(s,\eta) \in \R^{1|1}(S)$ such that 
$$
\phi'=\phi\circ T_{s,\eta}, \quad  T_{s,\eta}\circ \phi_\inn=\phi_\inn', \quad T_{s,\eta}\circ \phi_\out=\phi_\out'
$$
where $T_{s,\eta}$ is the super translation by $(s,\eta)\in \R^{1|1}(S)$ defined in~\eqref{eq;translation}.
\end{defn}

Super Euclidean paths in $M$ and their isometries form a stack: for each $S$ we have a groupoid of $S$-families of super Euclidean paths over~$S$, and for a map $f\colon S'\to S$ of supermanifolds we obtain a functor between groupoids by pulling back the triples $(\phi,(t_\out,\theta_\out),(t_\inn,\theta_\inn))$ along~$f$. Analogously to~\eqref{eq:bigEpath}, we have the following.

\begin{lem}\label{lem:bigsEpath}
The groupoid of super Euclidean paths over $S$ is equivalent to the discrete groupoid given by the value of the sheaf
\beq\label{eq:bigsEpath}
\R^{1|1}_{\ge 0}\times \Map(\R^{1|1},M)\coprod \R^{1|1}_{\le 0}\times \Map(\R^{1|1},M)
\eeq
on the supermanifold $S$, and hence the stack of super Euclidean paths in $M$ is equivalent to the sheaf~\eqref{eq:bigsEpath}. 
\end{lem}
\bp
%We must show that an $S$-family of nearly constant superpaths in $M$ is the same data as a map from $S$ to~\eqref{eq:fdrepresentable}. 
Given an $S$-family of positively (respectively, negatively) oriented super Euclidean paths $(\phi,(t_\out,\theta_\out),(t_\inn,\theta_\inn))$, there is a unique isomorphism to a superpath with $(t_\inn,\theta_\inn)=(0,0)\in \R^{1|1}(S)$. With this identification fixed, the remaining data is
\beq
&&\phi\colon S\times \R^{1|1}\to M,\qquad (t_\out,\theta_\out)=(t,\theta) \in (\R^{1|1}_{\ge 0}\coprod \R^{1|1}_{\le 0})(S). \label{eq:protocat}
\eeq
The result follows. 
\ep

\subsection{The category of nearly constant superpaths}

Concatenation of $S$-families of superpaths in $M$ runs into essentially the same issues as in the case of ordinary paths: to get a (smooth) superpath one requires agreement in a neighborhood of the gluing. Just as before, the problem goes away when one only considers suitably constant maps to $M$.

\begin{defn} An $S$-family of \emph{nearly constant super Euclidean paths in $M$} is the data of a super Euclidean path $(\phi,(t_\out,\theta_\out),(t_\inn,\theta_\inn))$ for which $\phi$ is determined by a map $\phi_0$, 
\beq
\phi\colon S\times \R^{1|1}\xrightarrow{q} S\times \R^{1|1}/\R\simeq S\times \R^{0|1}\xrightarrow{\phi_0} M,\label{eq:constantsuper}
\eeq
where the first arrow is the quotient of $S\times \R^{1|1}$ by the $\R$-action on $\R^{1|1}$. 
\end{defn}

%Hence,
%$$
%C^\infty(\Map(\R^{0|1},M))\simeq C^\infty(\Pi TM)\simeq \Omega^\bullet(M).
%$$

\begin{lem}\label{lem:st} The groupoid of $S$-families of nearly constant super Euclidean path in $M$ is equivalent to the discrete groupoid given by the value of the sheaf
%Up to unique isomorphism, a family of nearly constant oriented superpaths in~$M$ is determined by a map to the finite-dimensional supermanifold,
\beq
\R^{1|1}_{\ge 0}\times \Map(\R^{0|1},M)\coprod \R^{1|1}_{\le 0}\times \Map(\R^{0|1},M),\label{eq:fdrepresentable}
\eeq
on the supermanifold $S$. 
%where the two components correspond to the positively and negatively oriented nearly constant super Euclidean paths. 
%i.e., 
\end{lem}
\bp
The follows the same argument as the proof of Lemma~\ref{lem:bigsEpath}, but the map $\phi$ in~\eqref{eq:protocat} is required to be of the form~\eqref{eq:constantsuper} and hence determined by $\phi_0\colon S\times \R^{0|1}\to M$. 
\ep

\begin{lem}\label{lem:concat} Suppose that $\gamma=(\phi,(t_\out,\theta_\out),(t_\inn,\theta_\inn))$ and $\gamma'=(\phi',(t_\out',\theta_\out'),(t_\inn',\theta_\inn'))$ are $S$-families of positively (respectively, negatively) oriented nearly constant superpaths in $M$ and that $\phi_\out=\phi_\inn'\colon S\times \R^{0|1}\to M$. Then there is a well-defined concatenation $\gamma*\gamma'$ as a positively (respectively, negatively) oriented nearly constant superpath in $M$ that is unique up to unique isomorphism.
\end{lem}
\bp 
It suffices to work in the universal case~\eqref{eq:fdrepresentable}, for which the concatenation of positively oriented superpaths is given by
\beq
&&(\R^{1|1}_{\ge 0}\times \Map(\R^{0|1},M))\times_{\Map(\R^{0|1},M)} (\R^{1|1}_{\ge 0}\times \Map(\R^{0|1},M))\nonumber\\
&&\simeq \R^{1|1}_{\ge 0} \times \R^{1|1}_{\ge 0} \times \Map(\R^{0|1},M) \stackrel{m\times \id}{\longrightarrow} \R^{1|1}_{\ge 0}\times \Map(\R^{0|1},M)\label{eq:compose}
\eeq
where $m$ is the restriction of multiplication on $\R^{1|1}$. The concatenation map in the negatively oriented case is completely analogous.
% with one important subtlety: the inversion map~\eqref{eq:oppositesemigroup} determines an isomorphism to the \emph{opposite} semigroup
%\beq\label{eq:oppositesemi}
%\R^{1|1}_{\ge 0}\simeq (\R^{1|1}_{\le 0})^\op,\qquad (t,\theta)\mapsto (-t,-\theta).
%\eeq
%This follows from the fact that inversion is a homomorphism $(-)^{-1}\colon \R^{1|1}\to (\R^{1|1})^\op$ to the opposite super Lie group. Hence, concatenation of negatively oriented superpaths uses multiplication in the opposite semigroup in~\eqref{eq:compose}. 
% where one restricts the opposite group structure on~$\R^{1|1}$ to $\R^{1|1}_{\ge 0}$ to define the composition.
% because of the inversion for the incoming boundary in~\eqref{eq:negative}. 
%In the above, we use the isomorphism $\Map(\R^{0|1},M)\simeq \Pi TM_\C$ to shorten the length of the formula. 
\ep

Super Lie categories are reviewed in Definition~\ref{defn:superLiecat}.

\begin{defn}\label{defn:sP} The super Lie category of \emph{nearly constant superpaths in $M$}, denoted $\sP_0(M)$, has objects and morphisms,
$$
\Ob(\sP_0(M))=\Map(\R^{0|1},M)\coprod \Map(\R^{0|1},M),$$
$$
 \Mor(\sP_0(M)) =\R^{1|1}_{\ge 0}\times \Map(\R^{0|1},M)\coprod \R_{\le 0}^{1|1}\times \Map(\R^{0|1},M),
$$
where source, target, and unit data are given by the incoming and outgoing boundaries of a superpath and the identity superpath (using Lemma~\ref{lem:st}) while composition is the concatenation of superpaths from Lemma~\ref{lem:concat}.
% $\Map(\R^{0|1},M)\times \R_{\ge 0}^{1|1}$ parameterizes the positively oriented nearly constant superpaths in $M$, and $\Map(\R^{0|1},M)\times \R_{\le 0}^{1|1}$ parameterizes the negatively oriented nearly constant superpaths in $M$. The source is given by the projections, whereas the target is the map
%$$
%\R^{1|1}\times \Map(\R^{0|1},M)\to \Map(\R^{0|1},M)
%$$
%gotten from precomposition with the action of $\R^{1|1}$ on $\R^{0|1}$. The unit map is inclusion at $0\in \R_{\ge 0}^{1|1}$ or $0\in \R_{\le 0}^{1|1}$, and composition is given by multiplication in $\R_{\ge 0}^{1|1}$ and $\R_{\le 0}^{1|1}$. 
%We regard $\sP_0(M)$ as a \emph{super Lie category}, meaning its objects and morphisms are supermanifolds, and the structure maps (source, target, unit, composition) are maps of supermanifolds. 
\end{defn}

%\begin{defn}
%A \emph{functor} $F\colon {\sf C}\to \mathcal{D}$ between super Lie categories is a pair of maps of supermanifolds $F_0\colon \Ob({\sf C})\to \Ob(\mathcal{D})$ and $F_1\colon \Mor({\sf C})\to \Mor(\mathcal{D})$ satisfying the axioms of a functor. A \emph{natural transformation} $\eta\colon F\Rightarrow G$ between functors is a map of supermanifolds $\eta\colon \Ob({\sf C})\to \Mor(\mathcal{D})$ satisfying the axioms of a natural transformations. 
%\end{defn}
The above super Lie category can be understood as a subcategory of the following super Lie groupoid. 

\begin{lem}\label{lem:action}
%The super Lie category $\sP_0(M)$ is a subcategory of an action groupoid: 
There is an essentially surjective and faithful (but not full) functor
\beq
\sP_0(M)\hookrightarrow \left(\Map(\R^{0|1},M)\coprod \Map(\R^{0|1},M)\right)\sq \R^{1|1},\label{eq:itsafunctor}
\eeq
to the quotient super Lie groupoid for the $\R^{1|1}$-action factoring through the homomorphism $\R^{1|1}\to \R^{0|1}$ (see~\eqref{eq:evaluationaction}), with the functor~\eqref{eq:itsafunctor} determined by the standard inclusions $\R^{1|1}_{\ge 0}, \R^{1|1}_{\le 0}\subset \R^{1|1}$. 
%In particular, the source and target maps of $\sP_0(M)$ are given by restrictions along the inclusion 
%$$
%\R^{1|1}_{\ge 0}\times \Map(\R^{0|1},M)\coprod (\R^{1|1}_{\ge 0})^\op\times \Map(\R^{0|1},M) \hookrightarrow \R^{1|1}\times \Map(\R^{0|1},M)\coprod \R^{1|1}\times \Map(\R^{0|1},M)
%$$
%of the maps
%$$
%\s=p\coprod \eva,\qquad \t=\eva\coprod p
%$$
%where $\eva$ is defined in~\eqref{eq:evaluationaction} and $p$ is the projection. 
%Orientation reversal on $\Path_0(M)$ is the restriction of the automorphisms of $\Map(\R^{0|1},M)\sq \R^{1|1}\coprod \Map(\R^{0|1},M)\sq \R^{1|1}$ given by exchanging components and applying the automorphism~\eqref{eq:supertimereverse} of $\R^{1|1}$. 
\end{lem}
\bp
By definition, the source map on an $S$-point of $\sP_0(M)$ is
$$
(\phi,(t_\out,\theta_\out),(t_\inn,\theta_\inn))\mapsto \phi_\inn,\qquad \phi_\inn \colon S\times \R^{0|1} \xhookrightarrow{i_0} S\times \R^{1|1}\xrightarrow{\phi} M,
$$
with $\phi$ as in~\eqref{eq:constantsuper}. Since $q\circ i_0=\id_{S\times \R^{0|1}}$, we find $\phi_\inn=\phi_0$ in the notation of~\eqref{eq:constantsuper}. The source map is therefore the projection $\s=p\colon \Mor(\sP_0(M))\to \Ob(\sP_0(M))$. The target map on an $S$-point of $\sP_0(M)$ is 
\beq\label{eq:positive}
%&&\begin{array}{lllc} t_\inn&:=&i_0\colon S\times \R^{0|1}\hookrightarrow S\times \R^{1|1} \\
(\phi,(t_\out,\theta_\out),(t_\inn,\theta_\inn))\mapsto \phi_\out,\qquad \phi_\out \colon S\times \R^{0|1} \xhookrightarrow{T_{t,\theta}\circ i_0} S\times \R^{1|1}\xrightarrow{\phi} M.
%\end{array}\qquad {\rm (positively\ oriented\ super\ paths)}
\eeq
%\beq\label{eq:negative}
%&&\begin{array}{lllc} t_\out&:=&i_0\colon S\times \R^{0|1}\hookrightarrow S\times \R^{1|1} \\
%t_\inn&:=& T_{-t,-\theta}\circ i_0\colon S\times \R^{0|1} \hookrightarrow S\times \R^{1|1}
%\end{array}
%%\qquad {\rm (negative\ oriented\ super\ paths)}
%\eeq
%where $i_0=t_\inn$ denotes the map over~$S$ determined by the canonical inclusion $\R^{0|1}\subset \R^{1|1}$. 
Again using that  $q\circ i_0=\id_{S\times \R^{0|1}}$, the target map is the restriction of the (left) action 
\beq
\eva \colon \R^{1|1}\times \Map(\R^{0|1},M)\to \Map(\R^{0|1},M)\label{eq:evaluationaction}
\eeq
coming from the precomposition with the $\R^{1|1}$-action on $\R^{0|1}$, 
$$
\R^{1|1}\times \R^{0|1}\to \R^{0|1},\qquad (t,\theta)\cdot \eta=(\theta+\eta), \quad (t,\theta)\in \R^{1|1}(S), \ \eta\in \R^{0|1}(S).
$$
The lemma is proved.
\ep

We recall the isomorphism of supermanifolds $\Map(\R^{0|1},M)\simeq \Pi TM$ leading to the isomorphism between functions and differential forms
\beq\label{eq:oddtangentiso}
C^\infty(\Map(\R^{0|1},M))\simeq C^\infty(\Pi TM_\C)\simeq \Omega^\bullet(M). 
\eeq
Together with Lemma~\ref{lem:action}, this will allow us to give explicit formulas for the structure maps in the super Lie category $\sP_0(M)$, presented in terms of differential forms on~$M$. 

\begin{rmk}
Using~\eqref{eq:oddtangentiso}, functions on the quotient groupoid $\Map(\R^{0|1},M)\sq \R^{0|1}$ give a supergeometric interpretation of the de~Rham complex of~$M$, e.g., see~\cite[\S3]{HKST}. Lemma~\ref{lem:action} shows that the category of superpaths is closely related to this groupoid. 
% $\Map(\R^{0|1},M)\sq \R^{0|1}$. 
\end{rmk}

The isomorphism~\eqref{eq:oddtangentiso} implies
\beq\label{eq:oddtangentiso1}\begin{array}{c} 
C^\infty(\R^{1|1}_{\ge 0}\times \Map(\R^{0|1},M))\simeq \Omega^\bullet(M;C^\infty(\R^{1|1}_{\ge 0})),\\ 
C^\infty(\R^{1|1}_{\le 0}\times \Map(\R^{0|1},M))\simeq \Omega^\bullet(M;C^\infty(\R^{1|1}_{\le 0})).
\end{array}
\eeq

\begin{lem}\label{lem:stmap}
Using~\eqref{eq:oddtangentiso} and~\eqref{eq:oddtangentiso1}, the source and target maps $\s,\t\colon \Mor(\sP_0(M))\rightrightarrows \Ob(\sP_0(M))$ are determined by the super algebra maps
$$
\Omega^\bullet(M)\oplus \Omega^\bullet(M)\to \Omega^\bullet(M;C^\infty(\R_{\ge 0}^{1|1}))\oplus \Omega^\bullet(M;C^\infty(\R_{\le 0}^{1|1}))
$$
$$
\s^*(\alpha,\beta)=(\alpha,\beta),\qquad \t^*(\alpha,\beta)=(\alpha-\theta d\alpha,\beta-\theta d\beta),\quad \alpha,\beta\in \Omega^\bullet(M),
$$
where we use the canonical inclusions 
$$
\Omega^\bullet(M)\hookrightarrow \Omega^\bullet(M;C^\infty(\R_{\ge 0}^{1|1})),\qquad \Omega^\bullet(M)\hookrightarrow \Omega^\bullet(M;C^\infty(\R_{\le 0}^{1|1}))
$$ 
along the constant functions on $\R^{1|1}_{\ge 0}$ or $\R^{1|1}_{\le 0}$
\end{lem}
\bp 
Standard results (e.g.,~\cite[Lemma~3.8]{HKST}) compute the $\R^{0|1}$-action on $\Map(\R^{0|1},M)$ in terms of a map on differential forms:
% standard computations, e.g., from~\cite[\S3]{HKST}. 
%In terms of the functor of points, the action~\eqref{eq:evaluationaction} is given by
%\beq\label{eq:ev11}
%&&\eva\colon \R^{1|1}\times \Map(\R^{0|1},M)\to \Map(\R^{0|1},M), \quad (t,\theta,x,\psi)\mapsto (x-\theta\psi,\psi), \\
%\nonumber && (t,\theta)\in \R^{1|1}(S), \ (x,\psi)\in \Pi TM_\C(S). 
%\eeq
%This action is equivalent data to the super algebra homomorphism
\beq
\Omega^\bullet(M)\simeq C^\infty(\Map(\R^{0|1},M))&\to& C^\infty(\R^{0|1}\times \Map(\R^{0|1},M))\simeq \Omega^\bullet(M;C^\infty(\R^{0|1}))\nonumber\\
\alpha&\mapsto& \alpha-\theta  d\alpha,\quad \alpha\in \Omega^\bullet(M), \ C^\infty(\R^{0|1})\simeq C^\infty(\R)[\theta]\label{eq:actdeRham}
\eeq
where $d$ is the de~Rham differential on $\Omega^\bullet(M)$. From Lemma~\ref{lem:action}, the source and target maps in $\sP_0(M)$ are therefore determined by restrictions of~\eqref{eq:actdeRham}. 
% (following conventions in classical supersymmetric field theory, e.g., see~\cite[Page~44]{Freed5}). 
\ep

\begin{rmk}
The formula~\eqref{eq:actdeRham} differs from~\cite[Lemma~3.8]{HKST} by a sign; their formula is for the right action of $\R^{0|1}$ on $\Pi TM_\C$, whereas~\eqref{eq:actdeRham} is for the left action.
\end{rmk}

%\bp
%Define the functor~\eqref{eq:itsafunctor} as the identity map on objects, and take the inclusion on morphisms
%\beq
%\resizebox{\textwidth}{!}{$
%\begin{array}{rcl}
%\R^{1|1}_{\ge 0}\times \Map(\R^{0|1},M) \coprod \R^{1|1}_{\le 0} \times \Map(\R^{0|1},M) 
%%&\stackrel{\sim}{\to}& \R^{1|1}_{\ge 0} \times \Map(\R^{0|1},M)\coprod \R^{1|1}_{\le 0}\times \Map(\R^{0|1},M) \\
%&\hookrightarrow& \R^{1|1}\times \Map(\R^{0|1},M)\coprod \R^{1|1} \times \Map(\R^{0|1},M)
%\end{array}$}\nonumber
%\eeq
%determined by the standard inclusions $\R^{1|1}_{\ge 0}, \R^{1|1}_{\le 0}\subset \R^{1|1}$. 
%% on the second component of $\Mor(\sP_0(M))$
%It remains to verify compatibility with source, target, unit and composition. The compatibility with source and target comes from examining the maps~\eqref{eq:positive},~\eqref{eq:negative} and~\eqref{eq:evaluationaction}. Compatibility with units is clear, and compatibility with composition follows from the fact that~\eqref{eq:compose} is determined by the restriction of multiplication in $\R^{1|1}$, and similarly restriction of the opposite multiplication for negatively oriented superpaths. 
%\ep

\subsection{Automorphisms of the category of superpaths}

\begin{lem}\label{lem:structureonsP}
The maps \eqref{eq:supertimereverse}, \eqref{eq:superRG}, \eqref{eq:spinflipfor11} and~\eqref{eq:real11} determine smooth functors 
\beq
\Or\colon \sP_0(M)\to \sP_0(M),&& (t,\theta,\phi)\mapsto (-t,-i\theta,\phi\circ \orr^{-1}) \label{eq:formOr}\\ 
\RG\colon \R_{>0}\times \sP_0(M)\to \sP_0(M), && (\mu,t,\theta,\phi)\mapsto (\mu^2 t,\mu\theta, \phi\circ \rg_\mu^{-1}) \label{eq:formRG}\\ 
\Fl\colon \sP_0(M)\to \sP_0(M),&& (t,\theta,\phi)\mapsto (t,-\theta,\phi\circ \fl^{-1})\label{eq:formFl} \\ 
{\sf R} \colon \sP_0(M)\to \overline{\sP_0(M)}, && (t,\theta,\phi)\mapsto (\overline{t},\overline{\theta},\overline{\phi}\circ \rr) \label{eq:formRR}
\eeq
given by the indicated maps on $S$-points of morphisms. Furthermore, there are equalities of smooth functors
\beq\label{Eq:Orrelation}
\Or\circ \Or=\Fl,\qquad \RR \circ \Or= \overline{\Fl}\circ \overline{\Or} \circ \RR
\eeq
and $\RG$ gives the structure of a strict $\R_{>0}$-action on $\sP_0(M)$ in the sense that the diagram of super Lie categories strictly commutes
\beq
%\label{eq:moreRG}
\begin{tikzpicture}[baseline=(basepoint)];
\node (A) at (0,0) {$\R_{>0}\times \R_{>0}\times \sP_0(M)$};
\node (B) at (6,0) {$\R_{>0}\times \sP_0(M)$};
\node (C) at (0,-1.5) {$\R_{>0}\times \sP_0(M)$};
\node (D) at (6,-1.5) {$\sP_0(M)$};
\draw[->] (A) to node [above] {$\id_{\R_{>0}}\times \RG$} (B);
\draw[->] (A) to node [left] {$m\times \id_{\sP_0(M)}$} (C);
\draw[->] (C) to node [below] {$\RG$} (D);
\draw[->] (B) to node [right] {$\RG$} (D);
\path (0,-.75) coordinate (basepoint);
\end{tikzpicture}\label{eq:RGaction}
\eeq
where $\R_{>0}$ is regarded as a discrete super Lie category and $m\colon \R_{>0}\times \R_{>0}\to \R_{>0}$ is multiplication. 
\end{lem}
\bp 
The formulas~\eqref{eq:formOr}-\eqref{eq:formFl} follow from Example~\ref{ex:functorsuperpair}. The equalities~\eqref{Eq:Orrelation} and commutativity~\eqref{eq:RGaction} follow from~\eqref{eq:Orrelation}. The real structure~\eqref{eq:formRR} is inherited from the real structure on $1|1$-Euclidean supermanifolds in Example~\ref{ex:realstructureonSP}, together with the canonical real structure on $M$, giving the $\bar{S}$-point of $\Map(\R^{0|1},M)$
$$
\overline{S}\times \R^{1|1}\xrightarrow{\id_{\overline{S}}\times \rr} \overline{S\times \R^{1|1}}\xrightarrow{\overline{\phi}} \overline{M}\simeq M
$$
determined by the image of $\phi\colon S\times \R^{0|1}\to M$ under the conjugation functor on supermanifolds. 
%The action on morphisms is determined by the diagram
%\beq
%%\label{eq:moreRG}
%\begin{tikzpicture}[baseline=(basepoint)];
%\node (A) at (0,0) {$S\times \R^{0|1}$};
%\node (B) at (3,0) {$S\times \R^{1|1}$};
%\node (C) at (6,0) {$S\times \R^{0|1}$};
%\node (D) at (3,-1.25) {$M$};
%\node (E) at (3,1.5) {$S\times \R^{1|1}$};
%\draw[->,right hook-latex] (A) to node [below] {$i_{(t_\out,\theta_\out)}$} (B);
%\draw[->,left hook-latex] (C) to node [below] {$i_{(t_\inn,\theta_\inn)}$} (B);
%\draw[->] (B) to node [right] {$\phi$} (D);
%\draw[->] (B) to node [right] {$\orr,\rg_\mu,\fl$} (E);
%\path (0,-.75) coordinate (basepoint);
%\end{tikzpicture}\nonumber
%\eeq
%
\ep

We use the notation
$$
\RG_\mu\colon \sP_0(M)\to \sP_0(M)
$$
to denote the restriction of the functor $\RG$ along $\mu\in \R_{>0}$. 

\begin{cor} The functors~\eqref{eq:formOr}-\eqref{eq:formFl} are determined by the maps on morphisms induced by the super algebra maps
$$
\Omega^\bullet(M;C^\infty(\R^{1|1}_{\ge 0}))\oplus \Omega^\bullet(M;C^\infty(\R^{1|1}_{\le 0}))\to \Omega^\bullet(M;C^\infty(\R^{1|1}_{\ge 0}))\oplus \Omega^\bullet(M;C^\infty(\R^{1|1}_{\le 0}))
$$
\beq
\Or^*(\alpha(t,\theta),\beta(t,\theta))&=&(i^{|\beta|}\beta(-t,-i\theta),i^{|\alpha|}\alpha(-t,-i\theta)), \quad i=\sqrt{-1}\\
\RG_\mu^*(\alpha(t,\theta),\beta(t,\theta))&=&(\mu^{|\alpha|}\alpha(\mu^2t,\mu\theta),\mu^{|\beta|}\beta(\mu^2t,\mu\theta)\label{eq:RGonsP}\\
\Fl^*(\alpha(t,\theta),\beta(t,\theta))&=&((-1)^{|\alpha|}\alpha(t,-\theta),(-1)^{|\beta|}\beta(t,-\theta))
\eeq
for $\alpha,\beta\in \Omega^\bullet(M;C^\infty(\R^{1|1}_{\ge 0}))$, where $|\alpha|,|\beta|\in \Z$ denote the $\Z$-grading of $\alpha$ and $\beta$. 
\end{cor}
\bp
This follows from~\eqref{eq:formOr}-\eqref{eq:formFl} and standard results (e.g.,~\cite[Lemma~3.8]{HKST}) that compute the $\R^{0|1}\rtimes \C^\times$-action on $\Map(\R^{0|1},M)$. The above statements follow by specializing to the subgroups $\Z/4$, $\R_{>0}$ and $\Z/2$ of $\C^\times$. 
\ep

\subsection{Superpaths as super loops}

\begin{defn}\label{defn:supercircle}
Given an $S$-point $(\ell,\lambda)\in \R^{1|1}_{>0}(S)$, define the family of $1|1$-dimensional \emph{super Euclidean circles} as the quotient
\beq
S^{1|1}_{\ell,\lambda}:=(S\times \R^{1|1})/\Z\label{eq:supercircle}
\eeq
for the left $\Z$-action on $S\times \R^{1|1}$ given by
%composition
%\beq
%\Z\times S\times \R^{1|1}\stackrel{(\ell,\lambda)}{\hookrightarrow} \E^{1|1}\times S\times \R^{1|1}\to S\times \R^{1|1}\label{eq:circleact}
%\eeq
%where the second arrow is the $\E^{1|1}$-action on $\R^{1|1}$, while the first is the inclusion induced by the homomorphism over $S$ determined by the image of the generator,
%$$
%\Z\times S\to \E^{1|1}\times S, \quad \{1\}\times S\simeq S\stackrel{(\ell,\lambda)}{\to} \R^{1|1}_{>0}\times S\hookrightarrow \E^{1|1}\times S
%$$
%Explicitly~\eqref{eq:circleact} is determined by the 
\beq
n\cdot (t,\theta)= (t+n\ell+n\lambda\theta,n\lambda+\theta),\quad n\in \Z(S), (t,\theta)\in \R^{1|1}(S).\label{eq:Zact}
\eeq
%Equivalently this is the restriction of the left $\E^{1|1}$-action on $S\times \R^{1|1}$ to the $S$-family of subgroups $\Z\times S\subset \E^{1|1}\times S$ with generator $\{1\}\times S\simeq S\stackrel{(\ell,\lambda)}{\hookrightarrow} \R^{1|1}_{>0}\times S\subset \E^{1|1}\times S$. Define the \emph{standard super Euclidean circle} denoted $S^{1|1}=S^{1|1}_{1,0}=\R^{1|1}/\Z$ as the quotient by the action for the standard inclusion $\Z\subset \R\subset \E^{1|1}$. 
\end{defn}

Let $S^{1|1}=\R^{1|1}/\Z$ denote the quotient by the standard inclusion $\Z\subset \R\subset \R^{1|1}$. For any $(\ell,\lambda)$ there is an isomorphism 
\beq\label{eq:torusstnd}
S^{1|1}_{\ell,\lambda}\xrightarrow{\sim} S\times S^{1|1}
\eeq
determined by the $\Z$-equivariant map
\beq
S\times \R^{1|1}&\to& S\times \R^{1|1}\label{eq:Zequiv}\\
(\ell,\lambda,t,\theta)&\mapsto& (\ell,\lambda,t(\ell+\lambda\theta),\theta+t\lambda),\qquad (\ell,\lambda)\in \R^{1|1}_{>0}(S),\ (t,\theta)\in \R^{1|1}(S).\nonumber
\eeq
\begin{rmk}
The $S$-family of subgroups $S\times \Z\hookrightarrow S\times \R^{1|1}$ generated by $(\ell,\lambda)\in \R^{1|1}_{>0}(S)$ is normal if and only if $\lambda=0$. Hence, the standard super circle $S^{1|1}=\R^{1|1}/\Z$ inherits a group structure from $\R^{1|1}$, but a generic $S$-family of super Euclidean circles $S^{1|1}_{\ell,\lambda}$ does not. 
\end{rmk}

\begin{defn}\label{defn:sEucloop} Define the \emph{super Euclidean loop space} as the generalized supermanifold (see Example~\ref{ex:mappingsheaf})
$$
\mathcal{L}^{1|1}(M):=\R^{1|1}_{>0}\times \Map(S^{1|1},M).
$$ 
\end{defn}

We identify an $S$-point of $\mathcal{L}^{1|1}(M)$ with a map $S^{1|1}_{\ell,\lambda}\to M$ given by the composition 
\beq
S^{1|1}_{\ell,\lambda}\simeq  S\times S^{1|1}\to M\label{eq:superloop}
\eeq
using the isomorphism~\eqref{eq:torusstnd}. There is a left action of $\R^{1|1}\rtimes \Z/2$ on $\mathcal{L}^{1|1}(M)$ determined by the diagram
\beq
\begin{tikzpicture}[baseline=(basepoint)];
\node (A) at (0,0) {$S^{1|1}_{\ell,\lambda}$};
\node (B) at (3,0) {$S\times S^{1|1}$}; 
\node (C) at (0,-1.2) {$S^{1|1}_{\ell',\lambda'}$};
\node (D) at (3,-1.2) {$S\times S^{1|1}$};
\node (E) at (5,-.6) {$M$};
\draw[->] (A) to node [above] {$\simeq$} (B);
\draw[->] (A) to node [left] {$f$} (C);
\draw[->] (C) to node [above] {$\simeq$} (D);
\draw[->] (B) to node [above] {$\phi$} (E);
\draw[->,dashed] (D) to node [below] {$\phi'$} (E);
\path (0,-.6) coordinate (basepoint);
\end{tikzpicture}\label{eq:E11action}
\eeq
%\beq
%T^{2|1}_{\Lambda'}\simeq T^{2|1}_{\Lambda}\simeq S\times T^{2|1}\to M,
%\eeq
where the horizontal arrows are the isomorphism~\eqref{eq:torusstnd} and the super Euclidean isometry $f$ is inherited from the action of $\R^{1|1}\rtimes \Z/2$ on $S\times \R^{1|1}$. The arrow $\phi'$ is uniquely determined by these isomorphisms and the input map $\phi$. Hence, given $(\ell,\lambda,\phi)\in \R^{1|1}_{>0}(S)\times \Map(S^{1|1},M)(S)$ and an $S$-point of $\R^{1|1}\rtimes \Z/2$, the action on $\mathcal{L}^{1|1}(M)$ outputs $(\ell',\lambda',\phi')$ as in~\eqref{eq:E11action}. 
%Define a left action of $\Euc_{1|1}=\E^{1|1}\rtimes \Z/2$ on $\mathcal{L}^{1|1}(M)$ that sends an $S$-point~\eqref{eq:superloop} to the composition
%\beq
%S^{1|1}_{\ell',\lambda'}\simeq S^{1|1}_{\ell,\lambda}\simeq  S\times S^{1|1}\to M\label{eq:E11action}
%\eeq
%where the leftmost isomorphism is from Lemma~\ref{lem:supercircleE11}, so that 

There is an evident $S^1$-action on $\mathcal{L}^{1|1}(M)$ coming from the precomposition action of $S^1=\R/\Z<\R^{1|1}/\Z$ on $\Map(S^{1|1},M)$. Since the quotient is given by~$S^{1|1}/S^1\simeq \R^{0|1}$, the $S^1$-fixed points are
\beq
\mathcal{L}^{1|1}_0(M):=\R^{1|1}_{>0}\times\Map(\R^{0|1},M)\subset \R^{1|1}_{>0}\times \Map(S^{1|1},M)=\mathcal{L}^{1|1}(M).\label{eq:constantloops}
\eeq
We identify an $S$-point of $\mathcal{L}^{1|1}_0(M)$ with a map $S^{1|1}_{\ell,\lambda}\to M$ that factors as
\beq
S^{1|1}_{\ell,\lambda}\simeq S\times S^{1|1}= S\times \R^{1|1}/\Z\stackrel{p}{\to} S\times \R^{0|1}\to M\label{eq:R01proj}
\eeq
where the map $p$ is induced by the projection $\R^{1|1}\to \R^{0|1}$. The $\R^{1|1}\rtimes \Z/2$-action in~\eqref{eq:E11action} preserves this factorization condition, and is computed explicitly in~\cite[\S2.3]{DBEChern}.
% we give an explicit formula in Lemma~\ref{lem:E11action} below. Hence, the inclusion~\eqref{eq:constantloops} is $\Euc_{1|1}$-equivariant. 

The equalizer of the source and target maps in $\sP_0(M)$ is
$$
\R_{\ge 0}\times \Map(\R^{0|1},M)\coprod \R_{\le 0}\times \Map(\R^{0|1},M) \subset  \Mor(\sP_0(M))\rightrightarrows\Ob(\sP_0(M)).
$$
Geometrically, this equalizer parameterizes constant superpaths with the same source and target. We obtain nearly constant super loops by restricting to a further subspace. 

\begin{lem} \label{lem:spanofinclude}
There is a span of inclusions
$$
\Mor(\sP_0(M))\hookleftarrow \R_{>0}\times \Map(\R^{0|1},M)\coprod \R_{<0}\times \Map(\R^{0|1},M)\hookrightarrow \mathcal{L}_0^{1|1}(M)
$$
that views a constant superpath with the same source and target and nonzero super length as a constant super loop. 
\end{lem}

\section{Super Euclidean paths as bordisms}\label{sec:iota}

The main goal of this section is to construct an internal functor
\beq\label{eq:maininclude}
\iota\colon \sP_0(M)\to 1|1\EBord(M)
\eeq
from the super Lie category of constant superpaths in~$M$ to Stolz and Teichner's $1|1$-Euclidean bordism category over~$M$. We further show that~\eqref{eq:maininclude} is compatible with additional structures from Lemma~\ref{lem:structureonsP}: the spin flip, orientation reversal, renormalization group action, and real structure. The arguments are all phrased within the framework of~\cite{ST11} (see~\S\ref{sec:ebord} for a review). However,~\eqref{eq:maininclude} is a totally geometric construction that we expect to hold in any reasonable definition of the $1|1$-Euclidean bordism category.

\subsection{The $1|1$-Euclidean bordism category}

\begin{defn}\label{defn:11EB}
Let $1|1\EBord(M)$ denote the geometric bordism category associated with the super Euclidean geometry from Example~\ref{ex:modelgeo}. We call $1|1\EBord(M)$ the \emph{$1|1$-Euclidean bordism category over $M$}.
\end{defn}

We refer to~\S\ref{sec:STBord} for Stolz and Teichner's definition of a geometric bordism category associated to a rigid geometry. To briefly summarize, $1|1\EBord(M)$ is a category internal to the category of symmetric monoidal stacks with involution determined by the spin flip (see Definition~\ref{defn:spinflip} and Example~\ref{ex:flips}). Hence, part of the data of $1|1\EBord(M)$ are maps of symmetric monoidal stacks with involution
% In particular, $\Ob(1|1\EBord(M))$ and $\Mor(1|1\EBord(M))$ are symmetric monoidal stacks with involution, and 
$$
\s,\t\colon \Mor(1|1\EBord(M))\rightrightarrows \Ob(1|1\EBord(M)),\quad \u\colon \Ob(1|1\EBord(M))\to \Mor(1|1\EBord(M))
$$
$$
\c\colon \Mor(1|1\EBord(M))^{[2]}\to \Mor(1|1\EBord(M))
$$
giving the source, target, unit and composition maps, where we use the notation\footnote{By Lemma~\ref{lem:isofibration}, the strict fibered product agrees with the weak fibered product.}
$$
\Mor(1|1\EBord(M))^{[2]}=\Mor(1|1\EBord(M))\times_{\Ob(1|1\EBord(M))} \Mor(1|1\EBord(M)).
$$
for the fibered product relative to the source and target maps, i.e., $\Mor(1|1\EBord(M))^{[2]}$ is the stack that classifies pairs of composable bordisms. The full definition of $1|1\EBord(M)$ includes further associator and unitor data satisfying certain coherence properties; see~\S\ref{sec:internalcategories}. 

Structures in the $1|1$-Euclidean geometry discussed in~\S\ref{sec:11Eucgeo} equip the internal category $1|1\EBord(M)$ with certain endofunctors. 

\begin{lem} \label{lem:extrastrEB}
The internal category $1|1\EBord(M)$ admits internal endofunctors:
\begin{enumerate}
\item orientation reversal, $\Or\colon 1|1\EBord(M)\to 1|1\EBord(M)$,
\item the renormalization group action, $\RG_\mu\colon 1|1\EBord(M)\to 1|1\EBord(M)$, and
\item a real structure, $\RR\colon 1|1\EBord(M)\to 1|1\EBord(M)$,
\end{enumerate}
%The spin flip automorphism is an inner automorphism, whereas 
where $\Or$, and $\RG_\mu$ are maps of categories internal to symmetric monoidal stacks with involution, whereas the functor $\RR$ is an internal functor that covers complex conjugation on supermanifolds. 
\end{lem} 
\bp From Definition~\ref{defn:orientationbord} and Example~\ref{ex:orientationGM}, we obtain an orientation reversal functor
\beq
&&\Or\colon 1|1\EBord(M)\to 1|1\EBord(M).\label{eq:timereversal11EB}
%\qquad \Or \colon 1\EBord(M)\to 1\EBord(M).
\eeq
Using Definition~\ref{defn:rg} and Example~\ref{ex:RGgeneral}, we obtain an internal functor
\beq
&&\RG_\mu\colon 1|1\EBord(M)\to 1|1\EBord(M).\label{eq:RG11EB}
%\qquad \Or \colon 1\EBord(M)\to 1\EBord(M).
\eeq
Finally, Definition~\ref{defn:realbordisms} determines an internal functor (see Example~\ref{ex:realGM}) 
\beq
&&\RR\colon 1|1\EBord(M)\to 1|1\EBord(M) \label{eq:realstru11EB}
\eeq
whose values are given by
$$
(S\leftarrow Y\to M)\stackrel{\RR_0}{\mapsto} (\overline{S}\leftarrow \overline{Y}\to \overline{M}\simeq M) \qquad (S\leftarrow \Sigma\to M)\stackrel{\RR_1}{\mapsto} (\overline{S} \leftarrow \overline{\Sigma}\to \overline{M}\simeq M)
$$ 
where $\overline{Y}\to \overline{S}$ and $\overline{\Sigma}\to \overline{S}$ have the canonical structure of an $\overline{S}$-family of super Euclidean manifolds. We emphasize that the family parameter changes under the functor $\RR$ to the conjugate supermanifold. 
% and hence define objects and morphisms in the internal category $1|1\EBord(M)$. 
\ep

A map of smooth manifolds $M\to N$ determines an internal functor
$$
1|1\EBord(M)\to 1|1\EBord(N)
$$
see Lemma~\ref{lem:naturalityofGBord}. In particular, the $1|1$-Euclidean bordism category over $M$ is always an internal category over $1|1\EBord(\pt)$, the $1|1$-Euclidean bordism category over a point.

\subsection{Families of bordisms in $1|1\EBord(M)$}\label{sec:superpathbordism}

Definition~\ref{defn:GBord} reviews the construction of the symmetric monoidal stacks with involution $\Ob(1|1\EBord(M))$ and $\Mor(1|1\EBord(M))$. In this subsection, we use the same notation for the underlying stack, i.e., forgetting the monoidal structure and spin flip. Throughout, the family $S\times \R^{1|1} \to S$ will always have the canonical structure of a $1|1$-Euclidean family. When $M=\pt$, the definitions below specialize to the pictures~\eqref{diag:supergenerators}, and we use analogous notation to impart this intuition. 

\begin{defn}[Compare {\cite[Remark 6.40]{HST}}] \label{defn:11EBspt}
Define maps of stacks
\beq\label{eq:sptM}
\spt^+_M,\spt^-_M\colon \Map(\R^{1|1},M)\to \Ob(1|1\EBord(M))
\eeq
as follows. Given an $S$-point $\phi\colon S\times \R^{1|1}\to M$ of the source, define objects $(\spt^+,\phi),(\spt^-,\phi)$ in the groupoid $\Ob(1|1\EBord(M))(S)$ by
\beq\label{eq:collardataI}
\begin{array}{c}
Y^c=S\times \R^{0|1} \subset S\times \R^{1|1}=Y\stackrel{\phi}{\to} M,\label{eq:stdinclude} \\
Y\setminus Y^c \simeq S\times \R^{1|1}_{>0}\coprod S\times \R^{1|1}_{<0} =\left\{\begin{array}{ll} Y^+\coprod Y^- & (\spt^+,\phi) \\ Y^-\coprod Y^+ & (\spt^-,\phi)\end{array}\right.
\end{array}
\eeq
where $Y^c \hookrightarrow Y$ (the inclusion of the \emph{core}) is determined by the standard inclusion $0\colon \R^{0|1}\hookrightarrow \R^{1|1}$. As indicated in the second line, the $+$ or $-$ decoration on $\spt^\pm_M$ reflects the choice of the partition of the collar. By pulling back families, the above is natural in~$S$ and hence leads to the claimed maps~\eqref{eq:sptM}. 
%The universal family of the objects $(\spt^\pm,\phi)$ is given by the morphisms of stacks,
%$$
%(\spt^+,\eval),(\spt^-,\eval)\colon \Map(\R^{1|1},M)\to \Ob(1|1\EBord(M))
%$$
%where $
%\eval\colon \R^{1|1}\times \Map(\R^{1|1},M)\to M$ is the evaluation map.
%where the inclusion $Y^c \hookrightarrow Y$ is determined by the standard inclusion $\R^{0|1}\hookrightarrow \R^{1|1}$. As indicated in the second line, this determines a positively or negatively oriented super point depending on the partition of the collar~. 
\end{defn}

\begin{rmk} \label{rmk:germs1}
A pair of $S$-points of $\Map(\R^{1|1},M)$ determine isomorphic objects under $\spt^\pm_M$ if the the maps $S\times \R^{1|1}\to M$ agree on a half-open neighborhood of the core $S\times \R^{0|1}\hookrightarrow 
S\times \R^{1|1}$; see Definition~\ref{defn:GBord}. Hence, the map of stacks~\eqref{eq:sptM} factors through the quotient of $\Map(\R^{1|1},M)$ that parameterizes half-open germs of superpaths in~$M$.
\end{rmk} 

The following gives a map from super Euclidean paths in $M$ to bordisms.

%\begin{lem} \label{lem:superspinflip}
%Automorphisms of the object $\spt^\pm\in \Ob(1|1\EBord(\pt))(S)$ are in bijection with the set of maps $S\to \Z/2$, where the generator of $\Z/2$ acts on $Y=S\times \R^{1|1}$ in Definition~\ref{defn:11EBspt} by the spin flip. For $(\spt^\pm,\phi)\in \Ob(1|1\EBord(M))(S)$, this spin flip supplies an isomorphism $(\spt^\pm,\phi)\stackrel{\sim}{\to} (\spt^\pm,\phi\circ (\id_S\times \fl))$ where $\fl \colon \R^{1|1}\to \R^{1|1}$ is~\eqref{eq:spinflipfor11}. 
%% (i.e., the grading involution). 
%\end{lem} 
%\bp
%This follows from the same argument as~\cite[Lemma~86]{HST}: an automorphism of $\spt^\pm\in \Ob(1|1\EBord(\pt))(S)$ comes from an $S$-family of isometries of $Y=S\times \R^{1|1}$ preserving $S\times \R^{0|1}=Y^c\subset Y$. The only such family of isometries is given by a map from $S$ to $\Z/2<\R^{1|1}\rtimes \Z/2$, acting on $\R^{1|1}$ by the spin flip. For bordisms over $M$, this spin flip acts by an isometry, but it changes the map to $M$ as stated (using that $\fl^{-1}=\fl$). 
%\ep

\begin{defn}[Compare {\cite[Definition 6.41]{HST}}]\label{defn:11EBspath}
Define maps of stacks
\beq\label{eq:IM}
\begin{array}{c} \sI^+_M\colon \R^{1|1}_{\ge 0}\times \Map(\R^{1|1},M)\to \Mor(1|1\EBord(M))\\ 
\sI^-_M\colon \R^{1|1}_{\le 0}\times \Map(\R^{1|1},M)\to \Mor(1|1\EBord(M))\end{array}
\eeq
as follows. Given $\Phi\colon S\times \R^{1|1}\to M$ and $(t,\theta)\colon S\to \R^{1|1}$, consider
\beq
S\times \R^{0|1}\coprod S\times \R^{0|1} \xrightarrow{i_0 \coprod i_{t,\theta}} S\times \R^{1|1}=\Sigma \xrightarrow{\Phi} M,\label{eq:morphism11}
\eeq
%\beq
%&&S\times \R^{0|1}\coprod S\times \R^{0|1} \xhookrightarrow{i_0 \coprod i_0} S\times \R^{1|1}\coprod S\times \R^{1|1} \xrightarrow{\id\coprod T_{t,\theta}} S\times \R^{1|1}=\Sigma \xrightarrow{\Phi} M,\label{eq:morphism11}
%\eeq
using the notation from~\eqref{eq:itthetadef}.
To promote~\eqref{eq:morphism11} to an object in the groupoid $\Mor(1|1\EBord(M))(S)$, we require a choice of source and target data. Take
\beq\label{eq:morphisms11EB}
&&
(Y_{\inn}^c, Y_{\out}^c)=\left\{\begin{array}{lll} \Big((\spt^+,\Phi), (\spt^+,\Phi\circ T_{t,\theta}^{-1})\Big) &\implies (I_{t,\theta}^+,\Phi), & (t,\theta)\in \R^{1|1}_{\ge 0}(S)\\ 
\Big((\spt^-,\Phi), (\spt^-,\phi\circ T_{t,\theta}^{-1})\Big) &\implies (I_{t,\theta}^-,\Phi) & (t,\theta)\in \R^{1|1}_{\le 0}(S) \end{array}\right.
\eeq
i.e., the ordering in the coproduct in~\eqref{eq:morphism11} lists the source followed by the target $1|1$-Euclidean pair. The resulting options and the restriction on $(t,\theta)$ in~\eqref{eq:morphisms11EB} follow from requiring the existence of collars for the source and target. Any pair of choices of collar provide isomorphic functors~\eqref{eq:IM}; see Remark~\ref{rmk:collarexists}. 
\end{defn}

\begin{rmk} In parallel to Remark~\ref{rmk:germs1}, $S$-points of $\R^{1|1}_{\ge 0}\times \Map(\R^{1|1},M)$ or $\R^{1|1}_{\le 0}\times \Map(\R^{1|1},M)$ determine isomorphic objects in $\Mor(1|1\EBord(M))$ under $\sI^\pm_M$ if the associated maps $S\times \R^{1|1}\to M$ agree on a neighborhood of the core of the bordism, and so the family pulls back from the associated quotient. 
\end{rmk}
Changing the source and target data give variations on~\eqref{eq:IM}.

\begin{defn}[Compare {\cite[Definition 6.41]{HST}}]\label{defn:11EBspath2}
Define maps of stacks
\beq\label{eq:LM}
\begin{array}{c} \sL^+_M\colon \R^{1|1}_{\ge 0}\times \Map(\R^{1|1},M)\to \Mor(1|1\EBord(M))\\ 
\sL^-_M\colon \R^{1|1}_{\le 0}\times \Map(\R^{1|1},M)\to \Mor(1|1\EBord(M)),\\
\sR^+_M\colon \R^{1|1}_{> 0}\times \Map(\R^{1|1},M)\to \Mor(1|1\EBord(M))\\ 
\sR^-_M\colon \R^{1|1}_{< 0}\times \Map(\R^{1|1},M)\to \Mor(1|1\EBord(M))\end{array}
\eeq
as follows. Again for $\Phi\colon S\times \R^{1|1}\to M$ and $(t,\theta)\colon S\to \R^{1|1}$ consider~\eqref{eq:morphism11}, but choose source and target data
\beq\label{eq:spath11EB}
&&(Y_{\inn}^c, Y_{\out}^c)=\left\{\begin{array}{lll} \Big(\big((\spt^+,\Phi),(\spt^-,\Phi\circ T^{-1}_{t,\theta})\big),\emptyset\Big) &\implies (L_{t,\theta}^+,\Phi) & (t,\theta)\in \R^{1|1}_{\ge 0}(S) \\
\Big(\big((\spt^-,\Phi),(\spt^+,\Phi\circ T^{-1}_{t,\theta})\big),\emptyset\Big) & \implies(L_{t,\theta}^-,\Phi) & (t,\theta)\in \R^{1|1}_{\le 0}(S)  \\ 
\Big(\emptyset,\big((\spt^-,\Phi),(\spt^+,\Phi\circ T^{-1}_{t,\theta})\big)\Big)& \implies(R_{t,\theta}^+,\Phi) & (t,\theta)\in \R^{1|1}_{>0}(S) \\
\Big(\emptyset,\big((\spt^+,\Phi),(\spt^-,\Phi\circ T^{-1}_{t,\theta})\big)\Big) & \implies(R_{t,\theta}^-,\Phi) & (t,\theta)\in \R^{1|1}_{<0}(S)
\end{array}\right.
\eeq
where the restriction on $(t,\theta)$ is necessary for the existence of collars for the source and target. Any pair of choices provide isomorphic functors~\eqref{eq:LM}; see Remark~\ref{rmk:collarexists}. 
\end{defn}

\begin{rmk}\label{rmk:collarexists}
The requirements on~$(t,\theta)\in \R^{1|1}(S)$ in~\eqref{eq:morphisms11EB} and~\eqref{eq:spath11EB} come from demanding the existence of the collars $W_{\inn/\out}^\pm$ with properties demanded by Definition~\ref{defn:GBord}; any pair of choices of $W_{\inn/\out}^\pm$ lead to isomorphic objects in $\Mor(1|1\EBord(M))$ by intersecting the choices of collar. 
%Hereafter, we will assume a specific choice of collars to be fixed once and for all in~\eqref{eq:morphisms11EB}. 
%When $M=\pt$, there are universal families of the bordisms in~\eqref{eq:morphisms11EB} parameterized by $S=\R^{1|1}_{\le 0}$, $\R^{1|1}_{\ge 0}$, $\R^{1|1}_{>0}$ and $\R^{1|1}_{<0}$. 
We emphasize that the families $(R_{t,\theta}^\pm,\Phi)$ in~\eqref{eq:morphisms11EB} do not have a limit in $1|1\EBord(M)$ as $(t,\theta)\to 0$; compare \cite[Theorem 6.23]{HST}.
\end{rmk}

\begin{lem} The restriction of the spin flip on $\Ob(1|1\EBord(M))$ or $\Mor(1|1\EBord(M))$ along $\spt^\pm_M$, $\sI^\pm_M$, $\sL^\pm_M$ or $\sR^\pm_M$ is determined by the maps
\beq
\Map(\R^{1|1},M)&\to&\Map(\R^{1|1},M),\quad \phi\mapsto \phi\circ \fl^{-1}\nonumber\\
\R^{1|1}_{\ge 0}\times \Map(\R^{1|1},M)&\to&\R^{1|1}_{\ge 0}\times \Map(\R^{1|1},M),\quad (t,\theta,\phi)\mapsto (t,-\theta,\phi\circ \fl^{-1})\nonumber\\
\R^{1|1}_{\le 0}\times \Map(\R^{1|1},M)&\to&\R^{1|1}_{\le 0}\times \Map(\R^{1|1},M),\quad (t,\theta,\phi)\mapsto (t,-\theta,\phi\circ \fl^{-1})\nonumber
\eeq
\end{lem}
\bp
This follows from Definition~\ref{defn:spinflip} of the spin flip. 
\ep
%Recall from~\eqref{eq;translation} that $T_{t,\theta}\colon S\times \R^{1|1}\to S\times \R^{1|1}$ denotes translation by $(t,\theta)\in \R^{1|1}(S)$. Changing the standard inclusion $S\times \R^{0|1}{\hookrightarrow} S\times \R^{1|1} $ in~\eqref{eq:stdinclude} by postcomposition with an $S$-family of super translations $T_{t,\theta}$ results in objects of $\Ob(1|1\EBord(M))(S)$ isomorphic to $(\spt^+,\phi)$ and~$(\spt^-,\phi)$. Using such isomorphisms, 
\subsection{Composition of bordisms} 
%Define maps 
%\beq
%\proj,\ \act\colon \R^{1|1} \times \Map(\R^{1|1},M)&\to& \Map(\R^{1|1},M)\nonumber\\ 
%\proj(t,\theta,\Phi)&=&\Phi\label{eq:actproj}\\ 
%\act(t,\theta,\Phi)&=&\Phi\circ T_{t,\theta}^{-1}.\nonumber
%\eeq
%for $(t,\theta)\in \R^{1|1}(S)$ and $\Phi\colon S\times \R^{1|1}\to M$. 
%Restrictions of the above provide source and target data for the families of bordisms~\eqref{eq:IM} and~\eqref{eq:LM}. We start with the families~$\sI^\pm_M$. 
The source and target of the families of bordisms $\sI^\pm_M$, $\sL^\pm_M$ and $\sR^\pm_M$ are specified in their definitions. This source and target data can equivalently be phrased in terms of 2-commuting diagrams of stacks. For example, we have
\beq
&&\begin{tikzpicture}[baseline=(basepoint)];
\node (A) at (0,0) {$\R^{1|1}_{\ge 0} \times \Map(\R^{1|1},M)$};
\node (B) at (7,0) {$\Mor(1|1\EBord(M))$};
\node (C) at (0,-1.5) {$\Map(\R^{1|1},M)\times \Map(\R^{1|1},M)$};
\node (D) at (7,-1.5) {$\Ob(1|1\EBord(M))\times \Ob(1|1\EBord(M))$};
\node (E) at (3.25,-.75) {$\twocommute$};
\draw[->] (A) to node [above] {$\sI^+_M$} (B);
\draw[->] (A) to node [left] {$(\proj_+, \act_+)$} (C);
\draw[->] (C) to node [below] {$\spt^+_M\times \spt^+_M$} (D);
\draw[->] (B) to node [right] {$\s\times \t$} (D);
\path (0,-.75) coordinate (basepoint);
\end{tikzpicture} \label{eq:cyl2comm}
\eeq
where $\proj$ and $\act$ are the projection and action maps
\beq
\proj,\ \act\colon \R^{1|1} \times \Map(\R^{1|1},M)&\to& \Map(\R^{1|1},M)\nonumber\\ 
\proj(t,\theta,\Phi)&=&\Phi\label{eq:actproj}\\ 
\act(t,\theta,\Phi)&=&\Phi\circ T_{t,\theta}^{-1}.\nonumber
\eeq
for $(t,\theta)\in \R^{1|1}(S)$ and $\Phi\colon S\times \R^{1|1}\to M$, and $\proj_+,\act_+$ (respectively, $\proj_-,\act_-$) are the restrictions of $\proj,\act$ to $\R^{1|1}_{\ge 0}\subset \R^{1|1}$ (respectively, $\R^{1|1}_{\le 0}\subset\R^{1|1}$).
%Restrictions of the above provide source and target data for the families of bordisms~\eqref{eq:IM} and~\eqref{eq:LM}. We start with the families~$\sI^\pm_M$. 
%where the 2-commutativity data comes from intersecting the choice of collar in the definition of the family $\sI^+_M$ which the choice of collar for the families of objects $\spt^+_M$. 
The other structure maps in $1|1\EBord(M)$ lead to similar 2-commuting diagrams.

\begin{prop}\label{prop:stI}
Let $\u$ be the unit functor in $1|1\EBord(M)$. There is a 2-commutative diagrams of stacks
\beq
&&\begin{tikzpicture}[baseline=(basepoint)];
\node (A) at (0,0) {$\Map(\R^{1|1},M)$};
\node (B) at (0,-1.5) {$\R^{1|1}_{\ge 0}\times \Map(\R^{1|1},M)$};
\node (C) at (7,0) {$\Ob(1|1\EBord(M))$};
\node (D) at (7,-1.5) {$\Mor(1|1\EBord(M))$};
\node (E) at (3.25,-.75) {$\twocommute$};
\draw[->] (A) to node [above] {$\spt^+_M$} (C);
\draw[->,right hook-latex] (A) to node [left] {$i_0\times\id_{\Map(\R^{1|1},M)}$} (B);
\draw[->] (C) to node [right] {$\u$} (D);
\draw[->] (B) to node [above] {$\sI^+_M$} (D);
\path (0,-.75) coordinate (basepoint);
\end{tikzpicture} \label{eq:unitcyl}
\eeq
where the left vertical arrow is inclusion along $0\in \R^{1|1}_{\ge 0}$. Let $\c$ be the composition functor. There is a 2-commutative diagrams of stacks
\beq
\begin{tikzpicture}[baseline=(basepoint)];
\node (A) at (0,0) {$\R^{1|1}_{\ge 0}\times \R^{1|1}_{\ge 0} \times\Map(\R^{1|1},M)$};
\node (B) at (8,0) {$\Mor(1|1\EBord(M))^{[2]}$};
\node (C) at (0,-1.5) {$\R^{1|1}_{\ge 0} \times\Map(\R^{1|1},M)$};
\node (D) at (8,-1.5) {$\Mor(2|1\EBord(M))$};
\node (E) at (4,-.75) {$\twocommute$};
\draw[->] (A) to node [above] {$\sI^+_M \circ p_2 \times  \sI^+_M \circ p_1$} (B);
\draw[->] (A) to node [left] {${\sf m}$} (C);
\draw[->] (C) to node [below] {$\sI^+_M$} (D);
\draw[->] (B) to node [right] {$\c$} (D);
\path (0,-.75) coordinate (basepoint);
\end{tikzpicture} \label{eq:comp2iso}
\eeq
where the top horizontal arrow is determined by the family $\sI^+_M$ and the maps
\beq
&&p_i\colon \R^{1|1}_{\ge 0}\times \R^{1|1}_{\ge 0} \times\Map(\R^{1|1},M) \to \R^{1|1}_{\ge 0} \times\Map(\R^{1|1},M),\nonumber\\ 
&&p_1\big((s,\eta),(t,\theta),\Phi\big)=(t,\theta,\Phi),\qquad p_2\big((s,\eta),(t,\theta),\Phi\big)=(s,\eta,\Phi\circ T^{-1}_{t,\theta})\label{eq:composablesuperpath}
\eeq
and ${\sf m}$ is inherited from multiplication in $\R^{1|1}$,
\beq\label{eq:sfm}
{\sf m}\big((s,\eta),(t,\theta),\Phi\big)=(s+t+\eta\theta,\eta+\theta,\Phi).
\eeq
There are analogous diagrams to~\eqref{eq:unitcyl} and~\eqref{eq:comp2iso} for the family of bordisms classified by~$\sI^-_M$. 
%There are analogous diagrams for $\sI^-_M$ and $\spt^-_M$. Equivalently, we have isomorphisms over~$S$ 
%\beq
%\u(\spt^\pm,\phi)=(I_{0,0}^\pm,\phi), \quad \c((I^\pm_{t',\theta'},\phi\circ T_{t,\theta}^{-1}),(I^\pm_{t,\theta},\phi))&=&(I^\pm_{t'+t+\theta'\theta,\theta'+\theta},\phi)
%%\quad \s(I_{t,\theta}^\pm,\Phi)&=&(\spt^\pm,\Phi), \quad \t(I_{t,\theta}^\pm,\Phi)= (\spt^\pm,\Phi\circ T_{t,\theta}^{-1}),\\ 
%\eeq
%that are natural in~$S$.
%\beq
%&&\begin{tikzpicture}[baseline=(basepoint)];
%\node (A) at (0,0) {$\R^{1|1}_{\ge 0} \times \Map(\R^{1|1},M)$};
%\node (B) at (7,0) {$\Mor(1|1\EBord(M))$};
%\node (C) at (0,-1.5) {$\Map(\R^{1|1},M)\times \Map(\R^{1|1},M)$};
%\node (D) at (7,-1.5) {$\Ob(1|1\EBord(M))\times \Ob(1|1\EBord(M))$};
%\node (E) at (3.25,-.75) {$\twocommute$};
%\draw[->] (A) to node [above] {$\sI^+_M$} (B);
%\draw[->] (A) to node [left] {$\proj^+\coprod \act^+$} (C);
%\draw[->] (C) to node [below] {$\spt^+_M\times \spt^+_M$} (D);
%\draw[->] (B) to node [right] {$\s\times \t$} (D);
%\path (0,-.75) coordinate (basepoint);
%\end{tikzpicture} \label{eq:cyl2comm}
%\eeq

%Let $\sC_0$ denote the restriction of $\sC$ to $\mathcal{L}_0(M)=(\R_{>0}\coprod \R_{>0})\times \Map(\R^{0|1},M)\xhookrightarrow{0} \mathcal{C}_0(M)$, included along $\{0\}\coprod \{0\} \in \HH^{2|1}_+\coprod \HH^{2|1}_-$. Then there is a 2-commutative diagram 
\end{prop}

\bp
The commutativity of the diagram~\eqref{eq:unitcyl} is almost immediate; the only subtlety comes from the choice of collars in the construction of the family $\sI^+_M$. Indeed, rather than a strictly commutative diagram this choice leads to 2-commutativity data gotten from the isomorphism between objects of $\Mor(1|1\EBord(M))$ from shrinking the collar in~\eqref{eq:collardataI}. 

Turning to the diagram~\eqref{eq:comp2iso}, the definition of the $S$-family $(I_{t,\theta}^+,\Phi)$ of bordisms implies that the pair~\eqref{eq:composablesuperpath} is indeed in the fibered product. The construction of 2-commutativity data comes from unpacking the definition of the functor $\c$, i.e., composition of bordisms. Consider the diagram 
\beq\label{eq:triplepath}
S\times \R^{0|1}\coprod S\times \R^{0|1}\coprod S\times \R^{0|1} \xhookrightarrow{i_0 \coprod i_{t,\theta}\coprod i_{(s,\eta)\cdot (t,\theta)}} S\times \R^{1|1} \xrightarrow{\Phi} M.
\eeq
By forgetting any of the three inclusions of $S\times \R^{0|1}$, one obtains three possible $S$-points of the source of $\sI^+_M$ via the input datum~\eqref{eq:morphism11}: forgetting $i_{(s,\eta)\cdot (t,\theta)}$ gives the $S$-family of bordisms $(I^+_{t,\theta},\Phi)$, forgetting $i_{t,\theta}$ gives the $S$-family of bordisms $(I^+_{(s,\eta)\cdot (t,\theta)},\Phi)$, and forgetting $i_0$ gives a family that is isomorphic to $(I^+_{s,\eta},\Phi\circ T^{-1}_{t,\theta})$ (where the isomorphism is specified by $T_{t,\theta}$). These three $S$-families of bordisms gotten from~\eqref{eq:triplepath} coincide with the two projections out of the fibered product and the value of the functor~$\c$. The only subtle aspect in commutativity of the diagram~\eqref{eq:comp2iso} is the choice of collars in the definitions of these families of bordisms, i.e., the choices in the construction of the functor $\sI^+_M$. However, since any pair of choices lead to isomorphic bordisms (by intersecting the collars), we obtain canonical 2-commutativity data in~\eqref{eq:comp2iso}.
% also uses shrinking of collars together with the isomorphism of super Euclidean pairs at the end of Example~\ref{ex:Eucpair} and compatibility between composition in $1|1\EBord(M)$ and concatenation of superpaths. 

The arguments involving $\spt^-_M$ and $\sI^-_M$ are identical. 
\ep
\begin{rmk} In terms of $S$-families of bordisms, the top horizontal map in~\eqref{eq:comp2iso} is 
$$
\big((s,\eta),(t,\theta),\Phi\big)\mapsto \big((I^+_{s,\eta},\Phi\circ T^{-1}_{t,\theta}),(I^+_{t,\theta},\Phi)\big). 
$$
\end{rmk}
\begin{notation} Let $\sI^\pm_0$ denote the composition $\u\circ \spt^\pm_M$, 
$$
\sI^\pm_0\colon \Map(\R^{1|1},M)\to \Mor(1|1\EBord(M)).
$$
For an $S$-point $\Phi$ of the source above, let $(I_0^\pm,\Phi)$ denote the associated $S$-family of (identity) $1|1$-Euclidean bordisms. 
\end{notation}

The following generalizes the adjunction between super semigroups coming from the last relation in~\eqref{eq:11generatoreasy}.
%The following shows that $\sL^+_M$ gives an adjunction between the superpath categories $\sI^+_M$ and $\sI^-_M$. 

\begin{prop}\label{prop:adjunctionbordisms}
There is a 2-commutative diagrams of stacks
\beq
&&\begin{tikzpicture}[baseline=(basepoint)];
\node (A) at (0,0) {$\R^{1|1}_{\ge 0}\times \R^{1|1}_{\ge 0}\times \R^{1|1}_{\le 0}\times \Map(\R^{1|1},M)$};
\node (B) at (9,0) {$\Mor(1|1\EBord(M))^{[2]}$};
\node (C) at (0,-1.5) {$\R^{1|1}_{\ge 0}\times \Map(\R^{1|1},M)$};
\node (D) at (9,-1.5) {$\Mor(1|1\EBord(M))$};
\node (E) at (3.25,-.75) {$\twocommute$};
\draw[->] (A) to node [above] {$\sL^+_M\circ p_2 \times (\sI^+_M\circ p_1 \coprod \sI^-_M\circ p_3)$} (B);
\draw[->] (A) to (C);
\draw[->] (C) to node [below] {$\sL^+_M$} (D);
\draw[->] (B) to node [right] {$\c$} (D);
\path (0,-.75) coordinate (basepoint);
\end{tikzpicture} \label{eq:adjunctioncyl}
\eeq
where the top vertical arrow is determined by the families $\sL^+_M$, $\sI^+_M$ and $\sI^-_M$ and the maps
\beq
&&p_i\colon \R^{1|1}_{\ge 0}\times \R^{1|1}_{\ge 0}\times \R^{1|1}_{\le 0}\times \Map(\R^{1|1},M)\to \R^{1|1}_{\ge 0} \times \Map(\R^{1|1},M), \quad i=1,2\nonumber\\
&&p_1((t,\theta),(s,\eta),(r,\rho),\Phi)= (t,\theta,\Phi), \quad p_2((t,\theta),(s,\eta),(r,\rho),\Phi)= (s,\eta,\Phi\circ T^{-1}_{t,\theta}) \label{eq:targetofadj}\\
%\Big((L^+_{s,\eta},\Phi\circ T^{-1}_{t,\theta}),\big((I^+_{t,\theta},\Phi)\coprod (I^-_{r,\rho},\Phi\circ T^{-1}_{t,\theta}\circ T^{-1}_{s,\eta})\big)\Big)\\
&&p_3\colon \R^{1|1}_{\ge 0}\times \R^{1|1}_{\ge 0}\times \R^{1|1}_{\le 0}\times \Map(\R^{1|1},M)\to \R^{1|1}_{\le 0} \times \Map(\R^{1|1},M)\nonumber\\
&&p_3((t,\theta),(s,\eta),(r,\rho),\Phi)= (r,\rho,\Phi\circ T^{-1}_{t,\theta}\circ T^{-1}_{s,\eta})
%\Big((L^+_{s,\eta},\Phi\circ T^{-1}_{t,\theta}),\big((I^+_{t,\theta},\Phi)\coprod (I^-_{r,\rho},\Phi\circ T^{-1}_{t,\theta}\circ T^{-1}_{s,\eta})\big)\Big)
\eeq
and the left vertical arrow is determined by the multiplication in $\R^{1|1}$,
$$
((t,\theta),(s,\eta),(r,\rho),\Phi)\mapsto ((r,\rho)^{-1}\cdot (s,\eta)\cdot (t,\theta),\Phi). 
$$
There is a similar 2-commutative diagram for $\sL^-_M$. 
\end{prop}
\bp
From the definitions of source and target data, the top map in~\eqref{eq:adjunctioncyl} does indeed yield an object in the fibered product. The 2-commutativity data comes from unpacking the definition of the functor $\c$. Similar to before, consider the diagram
\beq
\resizebox{.9\textwidth}{!}{$
S\times \R^{0|1}\coprod S\times \R^{0|1}\coprod S\times \R^{0|1}\coprod S\times \R^{0|1} \xhookrightarrow{i_0 \coprod i_{t,\theta}\coprod i_{(s,\eta)\cdot (t,\theta)}\coprod i_{(r,\rho)^{-1}\cdot (s,\eta)\cdot (t,\theta)}} S\times \R^{1|1} \xrightarrow{\Phi} M.$}\nonumber
\eeq
In this case, forgetting any pair of factors in the coproduct gives the input data to an $S$-family of bordisms. For example, forgetting the middle two factors corresponds to the $S$-family of bordisms $(L^+_{(r,\rho)^{-1}\cdot (s,\eta)\cdot (t,\theta)},\Phi)$, which is the claimed composition.  The bordisms in~\eqref{eq:targetofadj} can similarly be extracted by forgetting the appropriate pair of inclusions and applying an appropriate isometry determined by a translation of $S\times \R^{1|1}$. 
%where it is important here to keep track of the partitions of the collars when considering source and targets. 
\ep

\begin{rmk} In terms of $S$-families of bordisms, the top horizontal map in~\eqref{eq:adjunctioncyl} is
$$
((t,\theta),(s,\eta),(r,\rho),\Phi)\mapsto \Big((L^+_{s,\eta},\Phi\circ T^{-1}_{t,\theta}),\big((I^+_{t,\theta},\Phi)\coprod (I^-_{r,\rho},\Phi\circ T^{-1}_{t,\theta}\circ T^{-1}_{s,\eta})\big)\Big).
$$
\end{rmk}

The following shows that the family of bordisms $\sI^\pm_M$ factors as a composition involving $\sL^\pm_M$ and $\sR^\pm_M$. The proof uses similar techniques as above. 
%This generalizes the dualizability relation in the topological bordism category. 

\begin{prop}\label{prop:duality}
There is a 2-commutative diagrams of stacks
\beq
&&\begin{tikzpicture}[baseline=(basepoint)];
\node (A) at (0,0) {$\R^{1|1}_{\ge 0}\times \R^{1|1}_{>0}\times \Map(\R^{1|1},M)$};
\node (B) at (10,0) {$\Mor(1|1\EBord(M))^{[2]}$};
\node (C) at (0,-1.5) {$\R^{1|1}_{\ge 0}\times \Map(\R^{1|1},M)$};
\node (D) at (10,-1.5) {$\Mor(1|1\EBord(M))$};
\node (E) at (3.25,-.75) {$\twocommute$};
\draw[->] (A) to node [above] {$(\sL_M^+\circ p_1\coprod \sI_0^+\circ q_2)\times (\sI_0^+\circ q_1\coprod \sR^+_M\circ p_2)$} (B);
\draw[->] (A) to node [left] {${\sf m}$} (C);
\draw[->] (C) to node [below] {$\sI^+_M$} (D);
\draw[->] (B) to node [right] {$\c$} (D);
\path (0,-.75) coordinate (basepoint);
\end{tikzpicture} \label{eq:thickdiagram}
\eeq
where the top vertical arrow is determined by the families $\sL^+_M$, $\sR^+_M$, $\sI^\pm_0$, and the maps 
\beq
&&p_i\colon \R^{1|1}_{\ge 0}\times \R^{1|1}_{>0}\times \Map(\R^{1|1},M)\to \R^{1|1}_{>0}\times \Map(\R^{1|1},M)\nonumber\\
&&p_1((t,\theta),(s,\eta),\Phi)=(t,\theta,\Phi),\quad p_2((t,\theta),(s,\eta),\Phi)=(s,\eta,\Phi\circ T^{-1}_{t,\theta}), \label{eq:targetofdual1}\\
&&q_i\colon \colon \R^{1|1}_{\ge 0}\times \R^{1|1}_{>0}\times \Map(\R^{1|1},M)\to  \Map(\R^{1|1},M)\nonumber\\
&&\quad q_1((t,\theta),(s,\eta),\Phi)=\Phi, \quad q_2((t,\theta),(s,\eta),\Phi)=\Phi\circ T_{t,\theta}^{-1}\circ T_{s,\eta}^{-1}\nonumber
%&&\mapsto \Big((L^+_{s,\eta},\Phi)\coprod (I^+_0,\Phi\circ T^{-1}_{t,\theta} \circ T^{-1}_{s,\eta})\Big) \circ \Big((I_0^+,\Phi)\coprod (R_{s,\eta}^+,\Phi\circ T^{-1}_{t,\theta}) \Big) 
\eeq
and ${\sf m}$ is~\eqref{eq:sfm}. There is a similar 2-commutative diagram factoring $\sI^-_M$ as a composition.
\end{prop}

\begin{rmk} In terms of $S$-families of bordisms, the top horizontal arrow in~\eqref{eq:thickdiagram} is 
$$
((t,\theta),(s,\eta),\Phi)\mapsto \Big((L^+_{t,\theta},\Phi)\coprod (I^+_0,\Phi\circ T^{-1}_{t,\theta} \circ T^{-1}_{s,\eta})\Big) , \Big((I_0^+,\Phi)\coprod (R_{s,\eta}^+,\Phi\circ T^{-1}_{t,\theta}) \Big) . 
$$
\end{rmk}

\begin{prop}\label{prop:symmetry}
Let $\sL^\pm_M\circ\sigma$ denote the functor $\sL^\pm_M$ post-composed with the isomorphism in $\Mor(1|1\EBord(M))$ that exchanges the order of components of the target, i.e.,
$$
\t(\sL^\pm_{t,\theta}\circ\sigma,\Phi)=(\spt^\mp,\Phi\circ T_{t,\theta}^{-1})\coprod (\spt^\pm,\Phi).
$$
There is a 2-commutative triangle
\beq
&&\begin{tikzpicture}[baseline=(basepoint)];
\node (A) at (0,0) {$\R^{1|1}_{\ge 0} \times \Map(\R^{1|1},M)$};
\node (B) at (7,-.75) {$\Mor(1|1\EBord(M))$};
\node (C) at (0,-1.5) {$\R^{1|1}_{\le 0} \times \Map(\R^{1|1},M)$};
\node (E) at (3.25,-.75) {$\twocommute$};
\draw[->] (A) to node [above] {$\sL^+_M\circ \sigma$} (B);
\draw[->] (A) to (C);
\draw[->] (C) to node [below] {$\sL^-_M$} (B);
\path (0,-.75) coordinate (basepoint);
\end{tikzpicture} \label{eq:symmetrybordism}
\eeq
for the left vertical arrow
\beq
\R^{1|1}_{\ge 0} \times \Map(\R^{1|1},M)&\to& \R^{1|1}_{\le 0} \times \Map(\R^{1|1},M),\nonumber\\
(t,\theta,\phi)&\mapsto& (-t,-\theta,\phi\circ T_{t,\theta}^{-1}). \label{eq:2morphsigma}
\eeq
\end{prop} 

\bp
Starting with the data that determines $(L_{t,\theta}^+,\Phi)$, 
\beq\nonumber
S\times \R^{0|1}\coprod S\times \R^{0|1} \xhookrightarrow{i_0 \coprod i_{t,\theta}} S\times \R^{1|1} \xrightarrow{\Phi} M.
\eeq
exchanging the factors of the coproduct and acting on $S\times \R^{1|1}$ by $T^{-1}_{t,\theta}=T_{(-t,-\theta)}$ provides the data determining $(L_{-t,-\theta}^-,\Phi\circ T_{t,\theta}^{-1})$.
\ep

%\begin{lem}\label{lem:sourcetarget1}
%The source and target functors ${\sf s,t}\colon \Mor(1|1\EBord(M))\to \Ob(1|1\EBord(M))$ applied to~\eqref{eq:morphisms11EB} yield objects isomorphic to coproducts from Definition~\ref{defn:11EBspt}:
%\beq\label{eq:sourcetargetIt}
%\begin{array}{rcl}
%\s(I_{t,\theta}^+,\Phi)=(\spt^+,\Phi), && \t(I_{t,\theta}^+,\Phi)\simeq (\spt^+,\Phi\circ T_{t,\theta}^{-1}),\\ 
%\t(I_{t,\theta}^-,\Phi)= (\spt^-,\Phi), && \s(I_{t,\theta}^-,\Phi)\simeq(\spt^-,\Phi\circ T_{t,\theta}^{-1}), \\
%\t(L_{t,\theta}^+,\Phi)=\emptyset, && \s(L_{t,\theta}^+,\Phi)\simeq (\spt^+,\Phi)\coprod (\spt^-,\Phi\circ T_{t,\theta}^{-1})\\
%\t(L_{t,\theta}^-,\Phi)=\emptyset, && \s(L_{t,\theta}^-,\Phi)\simeq (\spt^-,\Phi\circ T_{t,\theta}^{-1})\coprod (\spt^+,\Phi),\\
%\s(R_{t,\theta}^+,\Phi)= \emptyset, && \t(R_{t,\theta}^+,\Phi)\simeq (\spt^+,\Phi\circ T_{t,\theta}^{-1})\coprod (\spt^-,\Phi),\\
%\s(R_{t,\theta}^-,\Phi)= \emptyset, && \t(R_{t,\theta}^+,\Phi)\simeq (\spt^-,\Phi)\coprod (\spt^+,\Phi\circ T_{t,\theta}^{-1}).
%\end{array}
%\eeq
%Furthermore, $(I_{0}^\pm,\Phi)$ is the identity bordism on $(\spt^\pm,\Phi)$, where $(t,\theta)=0\colon S\to \R^{1|1}$ is the constant map to~$0\in \R^{1|1}$. \end{lem}

We recall the super Euclidean loop space $\mathcal{L}^{1|1}(M)$ from Definition~\ref{defn:sEucloop}. By \cite[Lemma 2.12]{DBEChern}, there is a functor
$$
\mathcal{L}^{1|1}(M)\to \Mor(1|1\EBord(M)),
$$
sending a super Euclidean loop in $M$ to a bordism over $M$ whose source and target are the empty family of $1|1$-Euclidean manifolds. 

\begin{prop}\label{prop:tracerelations}
There is 2-commutative diagrams of stacks
\beq
&&\begin{tikzpicture}[baseline=(basepoint)];
\node (A) at (0,0) {$(\R^{1|1}_{\le 0}\times \R^{1|1}_{> 0}\times \Map(\R^{1|1},M))^\Z$};
\node (B) at (8,0) {$\Mor(1|1\EBord(M))^{[2]}$};
\node (C) at (0,-1.5) {$\mathcal{L}^{1|1}(M)$};
\node (D) at (8,-1.5) {$\Mor(1|1\EBord(M))$};
\node (E) at (3.25,-.75) {$\twocommute$};
\draw[->] (A) to node [above] {$\sL^-_M\circ p_1\times \sR^+_M\circ p_2$} (B);
\draw[->] (A) to (C);
\draw[->] (C) to (D);
\draw[->] (B) to node [right] {$\c$} (D);
\path (0,-.75) coordinate (basepoint);
\end{tikzpicture} \label{eq:tracerelations}
\eeq
where the $\Z$-fixed subspace has $S$-points $(t,\theta)\in \R^{1|1}_{\le 0}(S)$, $(s,\eta)\in \R^{1|1}_{>0}(S)$ and $\Phi\colon S\times \R^{1|1}\to M$ satisfying
\beq\label{eq:periodic1}
\Phi\circ T_{t,\theta}^{-1}=\Phi\circ T_{s,\eta}^{-1}\iff \Phi=\Phi\circ (T_{(s,\eta)\cdot (t,\theta)^{-1}})^{n},\quad n\in \Z
\eeq
the top vertical arrow is determined by the families $\sL^-_M$ and $\sR^+_M$ and the maps
\beq\label{eq:targetofdual}
&&p_1((t,\theta),(s,\eta),\Phi)=(t,\theta,\Phi), \qquad p_2((t,\theta),(s,\eta),\Phi)=(s,\eta,\Phi)
\eeq
and the left vertical arrow is determined by
$$
((t,\theta),(s,\eta),\Phi)\mapsto ((s,\eta)\cdot (t,\theta)^{-1},\Phi)\in \big(\R^{1|1}_{>0}\times \Map(\R^{1|1},M)\big)(S)\supset \mathcal{L}^{1|1}(M)(S)
$$
which determines an $S$-point of $\mathcal{L}^{1|1}(M)$ by the condition~\eqref{eq:periodic1}. There is a diagram analogous to~\eqref{eq:tracerelations} involving $\sL^-_M$ and $\sR^+_M$. 
\end{prop}
\bp
The compatibility $\Phi\circ T_{t,\theta}^{-1}=\Phi\circ T_{s,\eta}^{-1}$ in~\eqref{eq:periodic1} guarantees that the target of $(R^-_{s,\eta},\Phi)$ agrees with the source of~$(L^+_{t,\theta},\Phi)$. Hence, the bordisms are indeed composable and their composition is necessarily a bordism from the empty set to itself. Unpacking the definition of the functor $\c$ and examining the periodicity condition~\eqref{eq:periodic1}, one finds that this composition is indeed the claimed super Euclidean loop in~$M$ with circumference $(\ell,\lambda)=(s,\eta)\cdot (t,\theta)^{-1}$. 
\ep
\begin{rmk} In terms of $S$-families of bordisms, the top horizontal arrow in~\eqref{eq:tracerelations} is 
$$
((t,\theta),(s,\eta),\Phi)\mapsto \Big((L^-_{t,\theta},\Phi),(R^+_{s,\eta},\Phi) \Big).
$$
\end{rmk}

\subsection{Restriction of additional structures}

Consider the restriction of extra structures from Lemma~\ref{lem:extrastrEB} to the families of bordisms above.

\begin{prop}\label{prop:restrictRR}
The real structure on $1|1\EBord(M)$ sends $S$-families of bordisms to $\overline{S}$-families of bordisms specified by
\beq
\RR(\spt^\pm,\phi)&\simeq&(\spt^\pm,\overline{\phi}\circ \rr)\nonumber\\
\RR(I^\pm_{t,\theta},\Phi)&\simeq&(I^\pm_{\overline{t},\overline{\theta}},\overline{\Phi}\circ \rr) \nonumber\\
\RR(L^\pm_{t,\theta},\Phi)&\simeq&(L^\pm_{\overline{t},\overline{\theta}},\overline{\Phi}\circ \rr) \nonumber\\
\RR(R^\pm_{t,\theta},\Phi)&\simeq&(R^\pm_{\overline{t},\overline{\theta}},\overline{\Phi}\circ \rr ) \nonumber
\eeq
where (in an abuse of notation) $\rr\colon \overline{S}\times \R^{1|1}\to \overline{S}\times \overline{\R}{}^{1|1}$ and $\rr\colon \overline{S}\times \R^{0|1}\to \overline{S}\times \overline{\R}{}^{0|1}$ are the maps determined by the real structure~\eqref{eq:real11} on $\R^{1|1}$, $(\overline{t},\overline{\theta})\colon \overline{S}\to \overline{\R}{}^{1|1}$ is value of the conjugation functor on the map $(t,\theta)\colon S\to \R^{1|1}$, and $\overline{\Phi}\colon \overline{S}\times \overline{\R}{}^{1|1}\to M$, $\overline{\phi}\colon \overline{S}\times\overline{\R}{}^{0|1}\to M$ are the maps conjugate to $\Phi$ and $\phi$ that use the canonical real structure $\overline{M}\simeq M$ on the ordinary manifold $M$.
\end{prop}
\bp The real structure on $1|1\EBord(M)$ evaluated on any of the above families is determined by the diagram
\beq
&&\begin{tikzpicture}[baseline=(basepoint)];
\node (AA) at (5,1.5) {$\overline{M}$};
\node (AAA) at (5.7,1.5) {$\simeq M$};
\node (A) at (0,0) {$\overline{S}\times \overline{\R}{}^{0|1}$};
\node (B) at (5,0) {$\overline{S}\times \overline{\R}{}^{1|1}$};
\node (C) at (10,0) {$\overline{S}\times \overline{\R}{}^{0|1}$};
\node (D) at (0,-1.5) {$\overline{S}\times \R^{0|1}$};
\node (E) at (5,-1.5) {$\overline{S}\times \R^{1|1}$};
\node (F) at (10,-1.5) {$\overline{S}\times \R^{0|1}$};
\draw[->,right hook-latex] (A) to node [above] {$i_{\overline{t},\overline{\theta}}$} (B);
\draw[->] (D) to node [left] {$\rr$} (A);
\draw[->] (E) to node [left] {$\rr$} (B);
\draw[->] (F) to node [left] {$\rr$} (C);
\draw[->,right hook-latex] (D) to (E);
\draw[->,left hook-latex] (F) to (E);
\draw[->] (B) to node [left] {$\overline{\Phi}$} (AA);
\draw[->] (C) to node [above] {$i_{\overline{0}}$} (B);
\path (0,-.75) coordinate (basepoint);
\end{tikzpicture} \nonumber
\eeq
where the middle row is the value of the conjugation functor on a bordism over $M$, and the isomorphisms with the lower row use the real structure on the $1|1$-Euclidean geometry. Together with the fact that the real structure preserves the subspaces $\R^{1|1}_{\ge 0},\R^{1|1}_{\le 0}\subset \R^{1|1}$, one can choose isomorphisms as indicated in the proposition. These isomorphisms are unique up to the choice of collar data for the source and target. 
\ep

\begin{prop}\label{prop:orientation11restrict}
The orientation reversal automorphism of $1|1\EBord(M)$ takes values
\beq
\Or(\spt^\pm,\phi)&\simeq&(\spt^\mp,\phi\circ \orr^{-1})\nonumber\\
\Or(I^\pm_{t,\theta},\Phi)&\simeq&(I^\mp_{-t,-i\theta},\Phi\circ \orr^{-1}) \nonumber\\
\Or(L^\pm_{t,\theta},\Phi)&\simeq&(L^\mp_{-t,-i\theta},\Phi\circ \orr^{-1}) \nonumber\\
\Or(R^\pm_{t,\theta},\Phi)&\simeq&(R^\mp_{-t,-i\theta},\Phi\circ \orr^{-1}) \nonumber
\eeq
where (in an abuse of notation) $\orr\colon S\times \R^{1|1}\to S\times \R^{1|1}$ is determined by the map~\eqref{eq:supertimereverse}.
\end{prop}
\bp The orientation reversal automorphism of $1|1\EBord(M)$ evaluated on any of the above families is determined by the diagram
\beq
&&\begin{tikzpicture}[baseline=(basepoint)];
\node (AA) at (5,1.5) {$M$};
\node (A) at (0,0) {$S\times \R^{0|1}$};
\node (B) at (5,0) {$S\times \R^{1|1}$};
\node (C) at (10,0) {$S\times \R^{0|1}$};
\node (D) at (0,-1.5) {$S\times \R^{0|1}$};
\node (E) at (5,-1.5) {$S\times \R^{1|1}$};
\node (F) at (10,-1.5) {$S\times \R^{0|1}$};
\draw[->,right hook-latex] (A) to node [above] {$i_{t,\theta}$} (B);
\draw[->] (A) to node [left] {$\orr$} (D);
\draw[->] (B) to node [left] {$\orr$} (E);
\draw[->] (C) to node [left] {$\orr$} (F);
\draw[->,right hook-latex] (D) to node [above] {$i_{-t,-i\theta}$} (E);
\draw[->,left hook-latex] (F) to node [above] {$i_0$} (E);
\draw[->] (B) to node [left] {$\Phi$} (AA);
\draw[->] (C) to node [above] {$i_0$} (B);
\path (0,-.75) coordinate (basepoint);
\end{tikzpicture} \nonumber
\eeq
where the inclusions in the lower row are the unique ones making the diagram commute. Together with the fact that orientation reversal exchanges the subspaces $\R^{1|1}_{\ge 0},\R^{1|1}_{\le 0}\subset \R^{1|1}$, one can choose isomorphisms as in the statement. These isomorphisms are unique up to the choice of collar data for the source and target. 
\ep

\begin{prop}\label{prop:RG11restrict}
For $\mu\in \R_{>0}$, the renormalization group automorphism of $1|1\EBord(M)$ takes values
\beq
\RG_\mu(\spt^\pm,\phi)&\simeq&(\spt^\pm,\phi\circ \rg_\mu^{-1})\nonumber\\
\RG_\mu(I^\pm_{t,\theta},\Phi)&\simeq&(I^\pm_{\mu^2t,\mu\theta},\Phi\circ \rg_\mu^{-1}) \nonumber\\
\RG_\mu(L^\pm_{t,\theta},\Phi)&\simeq&(L^\pm_{\mu^2 t,\mu\theta},\Phi\circ \rg_\mu^{-1}) \nonumber\\
\RG_\mu(R^\pm_{t,\theta},\Phi)&\simeq&(R^\pm_{\mu^2t,\mu\theta},\Phi\circ \rg_\mu^{-1}) \nonumber
\eeq
where (in an abuse of notation) $\rg_\mu\colon S\times \R^{1|1}\to S\times \R^{1|1}$ is determined by the map~\eqref{eq:superRG}.
\end{prop}

\bp
The argument is analogous to the proof of Proposition~\ref{prop:orientation11restrict}, using this time that $\rg_\mu$ preserves the subspaces $\R^{1|1}_{\ge 0},\R^{1|1}_{\le 0}\subset \R^{1|1}$.
\ep

\subsection{Small bordisms and nearly constant superpaths}

\begin{defn}\label{defn:smallbordisms}
Consider the families of bordisms associated to the maps of stacks
\beq\label{eq:LM0}
\begin{array}{c} 
\spt^\pm_M\colon \Map(\R^{0|1},M)\to \Ob(1|1\EBord(M))\\
{\sf i}^+_M\colon \R^{1|1}_{\ge 0}\times \Map(\R^{0|1},M)\to \Mor(1|1\EBord(M))\\ 
{\sf i}^-_M\colon \R^{1|1}_{\le 0}\times \Map(\R^{0|1},M)\to \Mor(1|1\EBord(M))\\
{\rm l}^+_M\colon \R^{1|1}_{\ge 0}\times \Map(\R^{0|1},M)\to \Mor(1|1\EBord(M))\\ 
{\rm l}^-_M\colon \R^{1|1}_{\le 0}\times \Map(\R^{0|1},M)\to \Mor(1|1\EBord(M))\\
{\sf r}^+_M\colon \R^{1|1}_{\ge 0}\times \Map(\R^{0|1},M)\to \Mor(1|1\EBord(M))\\ 
{\sf r}^-_M\colon \R^{1|1}_{\le 0}\times \Map(\R^{0|1},M)\to \Mor(1|1\EBord(M))\end{array}
\eeq
by restriction of~\eqref{eq:sptM},~\eqref{eq:IM} and~\eqref{eq:LM} to the (finite-dimensional) subspace 
$$
\Map(\R^{0|1},M)\subset \Map(\R^{1|1},M),\qquad \phi_0\mapsto (S\times \R^{1|1}\xrightarrow{p}S\times \R^{0|1}\xrightarrow{\phi_0} M).
$$
We refer to the families~\eqref{eq:LM0} as \emph{small bordisms} in $M$.
\end{defn}

In particular, the families~\eqref{eq:LM0} determine maps of stacks
\beq\label{eq:superpathtoEB}
&&\Ob(\sP_0(M))\to \Ob(1|1\EBord(M)),\quad \Mor(\sP_0(M))\to \Mor(1|1\EBord(M)).
\eeq

\begin{prop}\label{prop:FTrestrict1}
The maps of stacks~\eqref{eq:superpathtoEB} are part of the data of an internal functor~\eqref{eq:maininclude} from the super Lie category of nearly constant superpaths in~$M$ (viewed as a category internal to stacks) to the $1|1$-Euclidean bordism category over $M$. Under this functor, the restriction of the structures $\Or, \RG_\mu$ and $\RR$ in Lemma~\ref{lem:extrastrEB} are compatible with the structures on $\sP_0(M)$ in Lemma~\ref{lem:structureonsP}. 
\end{prop}

\bp
The data of an internal functor is reviewed in Definition~\ref{defn:internalfunctor}. In addition to the maps of stacks~\eqref{eq:superpathtoEB}, we require compatibility data for composition and units. This is given by Proposition~\ref{prop:stI}. Since the 2-isomorphism data in this proposition is canonical (coming from intersecting collars of bordisms), the coherence property demanded by these 2-isomorphisms is automatic. This completes the construction of the internal functor~\eqref{eq:maininclude}.

The compatibility of this functor with the additional structures follows from Propositions~\ref{prop:restrictRR}, \ref{prop:orientation11restrict} and~\ref{prop:RG11restrict}: the isomorphisms in these propositions restricted to the appropriate families in~\eqref{eq:LM0} supply 2-commutativity data, 
\beq
\label{eq:RGOrRR}
&&\begin{tikzpicture}[baseline=(basepoint)];
\node (A) at (0,0) {$\sP_0(M)$};
\node (B) at (3,0) {$1|1\EBord(M)$};
\node (C) at (0,-1.5) {$\sP_0(M)$};
\node (D) at (3,-1.5) {$1|1\EBord(M)$};
\draw[->] (A) to (B);
\draw[->] (A) to node [left] {$\RG_\mu,\Or$} (C);
\draw[->] (C) to  (D);
\draw[->] (B) to node [right] {$\RG_\mu,\Or$} (D);
\node (E) at (1.5,-.75) {$\twocommute$};
\path (0,-.75) coordinate (basepoint);
\end{tikzpicture}\quad 
\begin{tikzpicture}[baseline=(basepoint)];
\node (A) at (0,0) {$\sP_0(M)$};
\node (B) at (3,0) {$1|1\EBord(M)$};
\node (C) at (0,-1.5) {$\overline{\sP_0(M)}$};
\node (D) at (3,-1.5) {$1|1\EBord(M)$};
\draw[->] (A) to (B);
\draw[->] (A) to node [left] {$\RR$} (C);
\draw[->] (C) to  (D);
\draw[->] (B) to node [right] {$\RR$} (D);
\node (E) at (1.5,-.75) {$\twocommute$};
\path (0,-.75) coordinate (basepoint);
\end{tikzpicture}
\eeq
witnessing the claimed compatibility. 
\ep

\begin{rmk}
Composition and disjoint union of small bordisms generates a subcategory $1|1\ebord(M)\subset 1|1\EBord(M)$ of small bordisms in $M$. Morally, the categories $1|1\eft^n(M)$ in Theorems~\ref{thm:cocycle} and~\ref{thm:index} considers degree~$n$ field theories defined on these subcategories. However, it is difficult to formulate a precise notion of generators for a (sub)category in the framework from~\cite{ST11}, and so instead we formulate our results in terms of the further subcategory $\sP_0(M)\subset 1|1\ebord(M)$. 
\end{rmk}

\section{$1|1$-Euclidean field theories of degree~$n$}\label{sec:degreendefns}
In this section we review Stolz and Teichner's degree~$n$ twist and degree~$n$ twisted field theories over~$M$ from~\cite{ST11}. In brief, these are internal natural transformations~\eqref{eq:twistedEFT}
%\beq
%&&\begin{tikzpicture}[baseline=(basepoint)];
%\node (A) at (0,0) {$1|1\EBord(M)$};
%\node (B) at (5,0) {$\TA,$};
%\node (C) at (2.5,0) {$E \Downarrow$};
%\draw[->,bend left=15] (A) to node [above] {$\one$} (B);
%\draw[->,bend right=15] (A) to node [below] {$\twist^{\otimes n}$} (B);
%\path (0,-.1) coordinate (basepoint);
%\end{tikzpicture},\qquad n\in \Z\label{eq:twistedEFT}
%\eeq
for a $\otimes$-invertible functor $\twist$ called the \emph{degree twist}.
% that factors through the functor $1|1\EBord(M)\to 1|1\EBord(\pt)$ determined by $M\to \pt$. 
A complete construction of $\twist$ has yet to appear, but expected properties of~$\twist$ are described in~\cite{ST04,HST,ST11} for three variations on the definition of the bordism category~$1|1\EBord(M)$ and Morita bicategory~$\TA$. With all of these slightly differing past formalisms and others under current development, our goal below is provide a conceptual description of the degree~$n$ twist that can be parsed in any framework for field theories. Our description relies on the cobordism hypothesis for 2-dimensional extended topological field theories. Assuming this hypothesis, one can calculate the values of the degree~$n$ twist and twisted field theories for the families of bordisms constructed in the previous section. 
%In \S\ref{Sec:invertible}-\ref{Sec:invertible2} we sketch a construction of the degree~$n$ twist that relies on the cobordism hypothesis for topological field theories. This hypothesis should hold for any reasonable definition of (extended) field theory, and hence the salient features of the degree~$n$ twist should not depend on the nuances of a $1|1$-Euclidean bordism category; see Hypothesis~\ref{hyp:twist}. 
The punchline is that degree~$n$ field theories are (non-topological) boundary conditions for the 2-dimensional topological field theory determined by the Morita invertible (and hence fully-dualizable) algebra~$\cCl_n$, the $n$th Clifford algebra. The Morita equivalences between $\cCl_n\otimes \cCl_m$ and $\cCl_{n+m}$ then give the desired (graded) multiplication on degree~$n$ field theories via the tensor product. 

\begin{rmk}\label{rmk:genrelate} 
%We briefly review previous work on the degree twist. 
%All candidate constructions of nontrivial functors out of Stolz and Teichner's super Euclidean bordism categories depend on proposed generators and relations presentations. In low-dimensions (i.e., $1|1$ and $2|1$) the basic geometry of this statement is more-or-less clear, e.g., see~\eqref{diag:supergenerators} and~\eqref{eq:11generatoreasy}, as well as~\cite[6.7]{HST}, \cite[3.2]{ST11}, and the proof of~\cite[Proposition 4.16]{ST11}. 
%However, a precise statement has yet to appear, and it depends sensitively on the chosen framework. 
%%Presentations have proven difficult to formulate and prove in the framework for field theories in~\cite{ST11}. 
%The desired examples of functors out of the bordism category constructed via generators and relations include the
%% field theory from super parallel transport along a superconnection (discussed) and the 
% degree~$n$ twist~\eqref{eq:maintwist}. 
%%In these cases the values on super intervals are expected to determine the theory, e.g.,~\cite[Remark 3.8]{TwistAugusto}. 
%
The values of the degree~$n$ twist for a $d$-dimensional Euclidean field theories without supersymmetry are listed in~\cite[\S5.3]{ST11}. 
%Roughly, this non-supersymmetric twist comes from quantizing a theory of free fermions on Euclidean bordisms.
%\footnote{There remain some details to pin down in this sketch: the data analogous to the map~\eqref{Eq:twisttrace} when $d=2$ is fixed in the last paragraph of the proof of the periodicity theorem~\cite[Theorem 1.16]{ST11}. Different choices lead to periodicities of order 24, 48, or $\infty$. In particular, there are non-isomorphic degree twists when $d=2$, and it is unknown which one is the ``right" one to generalize for the $2|1$-dimensional supersymmetric twist.} 
%It is not clear how to extend this construction to super Euclidean bordisms. 
% though there remains both additional data to specify and condition to check to verify this sketch. 
One possible way to extend this construction to (supersymmetric) $1|1$-Euclidean bordisms involves a generators and relations presentation for the bordism category. Such presentations have been described in~\cite[6.7]{HST}, \cite[3.2]{ST11}, and the proof of~\cite[Proposition 4.16]{ST11}. However, complete statements and proofs have yet to appear. 
In terms of these expected generators and relations, $1|1$-Euclidean field theories of degree~$n$ are equivalent to a quantum mechanical theory with the $n$th Clifford algebra $\Cl_n$ acting by symmetries in the sense of~\S\ref{eq:CliffordlinearFT}, see~\cite[\S3.2]{ST04} and~\cite[Main Theorem]{HST}. The construction of the degree~$n$ twist outlined below is compatible with these prior expectations. 

\end{rmk}
%previously described structures. 
%when~$d=1$ that relies on a framework for fully-extended geometric field theories that includes topological cobordism hypothesis for topological field theories. 
%Ultimately, the details depend sensitively on the chosen framework for geometric field theories. 

%The definition~\eqref{eq:twistedEFT} of a $1|1$-Euclidean field theory over $M$ of degree~$n$ is complete modulo the construction of the degree~$n$ twist, discussed in the following remark. 

\subsection{The values of twisted field theories}
Twisted geometric field theories are defined as internal natural transformations~\eqref{eq:twistedGFT} for an internal functor $\twist$ from a geometric bordism category to the internal Morita category; see~\S\ref{sec:twistedFT} for a review. We unpack the values of a twisted $1|1$-Euclidean field theory on the families of bordisms from the previous section; this expands upon the discussion~\cite[page 52]{ST11}. 
%The maps $\proj_\pm$ and $\act_\pm$ used in the description below are defined immediately after~\eqref{eq:actproj}. 

\begin{prop}\label{prop:dataoftwistedFT} Given an internal functor $\twist\colon 1|1\EBord(M)\to \TA$, a $\twist$-twisted field theory~$E$ determines data
\begin{enumerate}
\item bundles of left modules
\beq\nonumber
\begin{array}{c}
\nonumber _{\twist(\spt^+_M)}E(\spt^+_M)\to \Map(\R^{1|1},M),\quad \nonumber _{\twist(\spt^-_M)}E(\spt^-_M)\to \Map(\R^{1|1},M)
% \in \Mor(\TA),
\end{array}
\eeq
where $\twist(\spt^\pm_M)$ are super algebra bundles over $\Map(\R^{1|1},M)$;
\item maps of bimodule bundles over $\R^{1|1}_{\ge 0}\times \Map(\R^{1|1},M)$ 
\beq
&&E(\sI^+_M)\colon \twist(\sI^+_M)\otimes_{\proj_+^*\twist(\sI^+_M)} \proj_+^*E(\spt^+_M)\to \act_+^*E(\spt^+_M),\nonumber\\
&& E(\sL^+_M)\colon \twist(\sL^+_M)\otimes_{\proj_+^*\twist(\spt^+_M)\otimes \act_+^*\twist(\spt^-_M)} \proj_+^*E(\spt^+_M)\otimes \act_+^*E( \spt^-_M)\to \underline{\C} ,\nonumber
\eeq
maps of bimodule bundles over $\R^{1|1}_{> 0}\times \Map(\R^{1|1},M)$ 
\beq
&&E(\sR^+_M)\colon \twist(\sR^+_M)\to \proj_+^*E(\spt^-_M) \otimes \act_+^*E(\spt^+_M), \nonumber
 \eeq
 maps of bimodule bundles over  $\R^{1|1}_{\le 0}\times \Map(\R^{1|1},M)$,
 \beq
&&E(\sI^-_M)\colon \twist(\sI^-_M)\otimes_{\proj_-^*\twist(\sI^-_M)} \proj_-^*E(\spt^-_M)\to \t^*E(\spt^-_M),\nonumber\\
&& E(\sL^-_M)\colon \twist(\sL^-_M)\otimes_{\proj_-^*\twist(\spt^-_M)\otimes \act_-^*\twist(\spt^+_M)} \proj_-^*E(\spt^-_M)\otimes \act_-^*E( \spt^+_M)\to \underline{\C} ,\nonumber
\eeq
and maps of bimodule bundles over  $\R^{1|1}_{< 0}\times \Map(\R^{1|1},M)$,
\beq
&&E(\sR^-_M)\colon \twist(\sR^-_M)\to \proj_-^*E(\spt^+_M)\otimes \act_-^*E(\spt^-_M), \nonumber
 \eeq
for the projection and action maps $\proj_\pm$ and $\act_\pm$ defined immediately after~\eqref{eq:actproj};
\item a section of the dual of the equivariant vector bundle $\twist(\mathcal{L}^{1|1}(M))$ over the stack~$\mathcal{L}^{1|1}(M)$
$$
E(\mathcal{L}^{1|1}(M))\colon \twist(\mathcal{L}^{1|1}(M)) \to \underline{\C}\nonumber
$$
that is invariant under the action of super Euclidean isometries on $\mathcal{L}^{1|1}(M)$. 
\end{enumerate}
% involving the families~\eqref{eq:somefamilies}, and the 
%including those in Lemma~\ref{lem:relations}. 
\end{prop}
\bp
By Definition~\ref{defn:internalfunctor} of an internal functor, the twist determines maps of stacks
$$
\Ob(1|1\EBord(M))\to \Ob(\TA),\qquad \Mor(1|1\EBord(M))\to \Mor(\TA)
$$
which can be evaluated on the objects of the source stacks defined in~\S\ref{sec:superpathbordism}. Hence, the twist determines bundles:
\begin{itemize}
\item $\twist(\spt^\pm_M)\to \Map(\R^{1|1},M)$ of super algebras;
\item $\twist(\sI^+_M)\to \R^{1|1}_{\ge 0}\times \Map(\R^{1|1},M)$ of $\proj_+^*\twist(\spt^+_M)$-$\act_+^*\twist(\spt^+_M)$-bimodules;
\item $\twist(\sI^-_M)\to \R^{1|1}_{\le 0}\times \Map(\R^{1|1},M)$ of $\proj_-^*\twist(\spt^-_M)$-$\act_-^*\twist(\spt^-_M)$-bimodules;
\item $\twist(\sL^+_M)\to \R^{1|1}_{\ge 0}\times \Map(\R^{1|1},M)$ of right $\proj_+^*\twist(\spt^+_M)\otimes\act_+^*\twist(\spt^-_M)$-modules;
\item $\twist(\sL^-_M)\to \R^{1|1}_{\le 0}\times \Map(\R^{1|1},M)$ of right $\proj_-^*\twist(\spt^-_M)\otimes\act_-^*\twist(\spt^+_M)$-modules;
\item $\twist(\sR^+_M)\to \R^{1|1}_{> 0}\times \Map(\R^{1|1},M)$ of left $\proj_+^*\twist(\spt^-_M)\otimes\act_+^*\twist(\spt^+_M)$-modules;
\item $\twist(\sR^-_M)\to \R^{1|1}_{< 0}\times \Map(\R^{1|1},M)$ of left $\proj_-^*\twist(\spt^+_M)\otimes\act_-^*\twist(\spt^-_M)$-modules; 
\item and $\twist(\mathcal{L}^{1|1}(M))\to\mathcal{L}^{1|1}(M)$, a vector bundle.
\end{itemize}
The additional data $\mu$ in~\eqref{eq:internalfunctordiag} of an internal functor provides maps of bimodule bundles associated with compositions of bordisms; we spell out the data inherited from Propositions~\ref{prop:stI}, \ref{prop:adjunctionbordisms},~\ref{prop:duality},~\ref{prop:symmetry}, and~\ref{prop:tracerelations}. We obtain a map of $\twist(\spt^+_M)$-bimodules,
% (something like the following, from Proposition~\ref{prop:stI})
\beq\nonumber
&&p_2^* \twist(\sI^+_M)\otimes_{p_2^*\twist(\spt^+_M)} p_1^* \twist(\sI^+_M)\to  {\sf m}^*\twist(\sI^+_M),\nonumber\\ 
&&{\rm over}\quad  \R^{1|1}_{\ge0}\times \R^{1|1}_{\ge 0}\times \Map(\R^{1|1},M) \nonumber
\eeq
where $p_1,p_2$ are the maps in~\eqref{eq:composablesuperpath}. There is a similar isomorphism of $\twist(\spt^-_M)$ bimodules over $\R^{1|1}_{\le 0}\times \R^{1|1}_{\le 0}\times \Map(\R^{1|1},M)$. We also obtain a map of right $\twist(\spt^+_M)\otimes \twist(\spt^-_M)$-module bundles
\beq\label{eq:twistadjunction}
&&p_2^* \twist(\sL^+_M)\otimes_{p_2^*\twist(\spt^+_M)\otimes p_3^*\twist(\spt^-_M)} (p_1^* \twist(\sI^+_M)\otimes p_3^*\twist(\sI^-_M))\to {\sf m}^*\twist(\sL^+_M),\nonumber \\ 
&&{\rm over} \quad \R^{1|1}_{\ge0}\times \R^{1|1}_{\ge 0}\times\R^{1|1}_{\le 0}\times \Map(\R^{1|1},M) \nonumber
\eeq
where $p_1,p_2,p_3$ are the maps in~\eqref{eq:targetofadj}. There is a similar isomorphism involving $\twist(\sL^-_M)$. There is a map of vector bundles
\beq\label{Eq:twisttrace}
&&p_1^*\twist(\sL^+_M)\otimes_{p_1^*\twist(\spt^+_M)\otimes p_2^*\twist(\spt^-_M)} p_2^*\twist(\sR^-_M)\to \c^*\twist(\mathcal{L}^{1|1}(M)),\quad {\rm over } \quad \mathcal{L}^{1|1}(M)
\eeq
for $p_1,p_2$ the maps in~\eqref{eq:targetofdual}. There is an isomorphism of right $\proj_+^*\twist(\spt^+_M)\otimes \act_+^*\twist(\spt^-_M)$-module bundles
\beq\label{eq:twistsymemtry}
\inv^*\twist(\sL^-_M)\simeq \twist(\sL^+_M)\quad {\rm over} \quad \R^{1|1}_{\ge 0}\times \Map(\R^{1|1},M) 
\eeq
where $\inv$ is the map~\eqref{eq:2morphsigma}. The isomorphisms above associated to composition are subject to further coherence conditions, as dictated by the definition of an internal functor.

A twisted field theory is defined as an internal natural transformation (see Definition~\ref{rmk:twistedcategory}), and hence is the data of a map of symmetric monoidal stacks with involution $\eta\colon \Ob(1|1\EBord(M))\to \Mor(\TA)$ and compatibility data~\eqref{diag:nu}. The map of stacks~$\eta$ hence determines morphisms over $\Map(\R^{1|1},M)$
\beq\nonumber
\begin{array}{c}
\nonumber E(\spt^\pm_M) \in \Mor(\TA), \quad {\rm with} \quad \t(E(\spt^\pm_M))= \twist(\spt^\pm_M), \ \s(E(\spt^\pm_M))= \underline{\C}.
\end{array}
\eeq
By the definition of $\TA$, this provides the bundles of left modules in~(1): the right action by $\underline{\C}$ is no additional data.  The data $\rho$ of a twisted field theory yields the data in (2) and (3) by the argument from~\cite[page 52]{ST11}. For example, $\rho$ evaluated on the family $\sI^\pm_M$ is a map of bimodules
$$
E(\sI^+_M)\colon \twist(\sI^+_M)\otimes_{\proj_+^*\twist(\sI^+_M)} \proj^*E(\spt^+_M)\to \act_+^*E(\spt^+_M)
%\twist(\sI^\pm_M)\otimes_{\s^*\twist(\sI^\pm_M)} \s^*E(\spt^\pm_M)\to \t^* E(\spt^\pm_M)\otimes_{\t^*\underline{\C}} \underline{\C}\simeq \t^* E(\spt^\pm_M),
$$
which is the first map of bimodules listed above. 
\ep

\begin{prop} \label{prop:maincompositionstatement}
The data in Proposition~\ref{prop:dataoftwistedFT} satisfy compatibility properties with respect to composition of bordisms. 
In particular, we have an equality of bimodule maps 
% (something like the following, from Proposition~\ref{prop:stI})
\beq\label{eq:comppathbimodule}
p_1^* E(\sI^+_M)\circ p_2^* E(\sI^+_M)= {\sf m}^*E(\sI^+_M),\quad {\rm over}\quad  \R^{1|1}_{\ge0}\times \R^{1|1}_{\ge 0}\times \Map(\R^{1|1},M) 
\eeq
where $p_1,p_2$ are the maps~\eqref{eq:composablesuperpath}. There is further an equality of bimodule maps
\beq\label{eq:bimoduleadjunction}
&&
\resizebox{.9\textwidth}{!}{$p_2^* E(\sL^+_M)\circ (p_1^* E(\sI^+_M)\otimes p_3^*E(\sI^-_M))= {\sf m}^*E(\sL^+_M), \ {\rm over} \ \R^{1|1}_{\ge0}\times \R^{1|1}_{\ge 0}\times\R^{1|1}_{\le 0} \times \Map(\R^{1|1},M) $}
\eeq
where $p_1,p_2,p_3$ are the maps in~\eqref{eq:targetofadj}. We have the equality of dual sections of $\twist(\mathcal{L}^{1|1}(M))$ 
\beq\label{Eq:bimoduletrace}
p_1^*E(\sL^+_M)\circ p_2^*E(\sR^-_M)={\sf m}^*E(\mathcal{L}^{1|1}(M)),\quad {\rm over } \quad \mathcal{L}^{1|1}(M)
\eeq
for $p_1,p_2$ the projections to the factors in~\eqref{eq:targetofdual}. Proposition~\ref{prop:symmetry} yields the equality of maps
\beq\label{eq:bimodulesymemtry}
\inv^*E(\sL^-_M)=E(\sL^+_M)\quad {\rm over} \quad \R^{1|1}_{\ge 0}\times \Map(\R^{1|1},M) 
\eeq
where $\inv$ is the map~\eqref{eq:2morphsigma}. Analogous statements hold after exchanging the $+$ and $-$ decorations for the families of bordisms.

\end{prop}
\bp
The claim follows from Definition~\ref{rmk:twistedcategory} of an internal natural transformation, where the data in Proposition~\ref{prop:dataoftwistedFT} must satisfy compatibility properties for composition in the internal category~$1|1\EBord(M)$ and invariance properties for isomorphisms in the stack $\Mor(1|1\EBord(M))$. We emphasize that the relations demanded above are now \emph{conditions} that certain maps of vector bundles be equal rather than additional data.

In more detail,~\eqref{eq:comppathbimodule} follows from compatibility with composition for the families $\sI^+_M$ as computed in Proposition~\ref{prop:stI}. The equality~\eqref{eq:bimoduleadjunction} follows from compatibility for composition of the families $\sL^+_M$ and $\sI^+_M\coprod \sI^-_M$ computed in Proposition~\ref{prop:adjunctionbordisms}. The equality~\eqref{Eq:bimoduletrace} follows from compatibility with the composition of $\sL^+_M$ and $\sR^-_M$ computed in Proposition~\ref{prop:tracerelations}. Finally,~\eqref{eq:bimodulesymemtry} follows from the isomorphism of families of bordisms computed in Proposition~\ref{prop:symmetry}. 
%For Propositions~\ref{prop:stI} and~\ref{prop:duality}, have a commuting diagram
%\beq
%\begin{tikzpicture}[baseline=(basepoint)];
%\node (A) at (0,0) {$\R^{1|1}_{\ge 0}\times \R^{1|1}_{\ge 0} \times\Map(\R^{1|1},M)$};
%\node (B) at (8,0) {$\Mor(1|1\EBord(M))^{[2]}$};
%\node (C) at (0,-1.5) {$\R^{1|1}_{\ge 0} \times\Map(\R^{1|1},M)$};
%\node (D) at (8,-1.5) {$\Mor(2|1\EBord(M))$};
%\node (E) at (4,-.75) {$\twocommute$};
%\draw[->] (A) to  (B);
%\draw[->] (A) to node [left] {${\sf m},p_1,p_2$} (C);
%\draw[->] (C) to (D);
%\draw[->] (B) to node [right] {$\c,p_1,p_2$} (D);
%\path (0,-.75) coordinate (basepoint);
%\end{tikzpicture} \nonumber
%\eeq
%for projection to the first and second morphisms in the composition versus their composite. So can pullback along these and compare over $\R^{1|1}_{\ge 0}\times \R^{1|1}_{\ge 0} \times\Map(\R^{1|1},M)$. Similar story for Proposition~\ref{prop:adjunctionbordisms}, only the relevant family is $\R^{1|1}_{\ge 0}\times \R^{1|1}_{\ge 0}\times \R^{1|1}_{\le 0} \times\Map(\R^{1|1},M)$. Finally for Proposition~\ref{prop:tracerelations} the point is that the dual section of the line bundle needs to be equal to the dual section coming from the tensor product of module maps. 
%
%For Proposition~\ref{prop:symmetry}, have compatibility for the pullback of the bimodule map along~\eqref{eq:2morphsigma}. 
%
\ep

\subsection{Reflection and real structures}\label{sec:reflrealEFT}
Next we describe reflection and real structures for twists and twisted field theories. We require a little notation: extend the conjugation functor on supermanifolds to a conjugation functor on generalized supermanifolds that sends a presheaf $\mathcal{F}$ to the presheaf
$$
\overline{\mathcal{F}}(S):=\mathcal{F}(\overline{S}).
$$
We will phrase reflection positivity and reality data using the map of sheaves
\beq
 \Map(\R^{1|1},M)&\to& \overline{\Map(\R^{1|1},M)},\label{eq:firstr} \\
(\phi\colon S\times \R^{1|1}\to M)&\mapsto &(\overline{S}\times \R^{1|1}\xrightarrow{\id_{\overline{S}}\times \rr} \overline{S}\times\overline{\R}{}^{1|1}\xrightarrow{\overline{\phi}} \overline{M}\simeq M),\nonumber
\eeq
which also determines maps
\beq
\R^{1|1}_{\ge 0} \times \Map(\R^{1|1},M)&\to& \overline{\R}^{1|1}_{\ge 0}\times \overline{\Map(\R^{1|1},M)},\label{eq:firstrr} \\
\R^{1|1}_{\le 0} \times \Map(\R^{1|1},M)&\to& \overline{\R}^{1|1}_{\le 0}\times \overline{\Map(\R^{1|1},M)}.\label{eq:firstrrr}
\eeq
In an abuse of notation, the maps~\eqref{eq:firstr}, \eqref{eq:firstrr} and~\eqref{eq:firstrrr} will all be denoted by $\rr$. Proposition~\ref{prop:restrictRR} gives 2-commuting data for the diagrams in stacks
\beq
%\label{eq:moreRG}\nonumber
&&\begin{tikzpicture}[baseline=(basepoint)];
\node (A) at (0,0) {$\Map(\R^{1|1},M)$};
\node (B) at (4,0) {$\Ob(1|1\EBord(M))$};
\node (C) at (0,-1.5) {$\overline{\Map(\R^{1|1},M)}$};
\node (D) at (4,-1.5) {$\Ob(1|1\EBord(M))$};
\draw[->] (A) to node [above] {$\spt^\pm_M$} (B);
\draw[->] (A) to node [left] {$\rr$} (C);
\draw[->] (C) to node [below] {$\overline{\spt}^\pm_M$}  (D);
\draw[->] (B) to node [right] {$\RR_0$} (D);
\node (E) at (2,-.75) {$\twocommute$};
\path (0,-.75) coordinate (basepoint);
\end{tikzpicture}\quad 
\begin{tikzpicture}[baseline=(basepoint)];
\node (A) at (0,0) {$\R^{1|1}_{\ge 0}\times \Map(\R^{1|1},M)$};
\node (B) at (4,0) {$\Mor(1|1\EBord(M))$};
\node (C) at (0,-1.5) {$\overline{\R}^{1|1}_{\ge 0}\times  \overline{\Map(\R^{1|1},M)}$};
\node (D) at (4,-1.5) {$\Mor(1|1\EBord(M))$};
\draw[->] (A) to node [above] {$\sI^+_M$} (B);
\draw[->] (A) to node [left] {$\rr$} (C);
\draw[->] (C) to node [below] {$\overline{\sI}^+_M$}  (D);
\draw[->] (B) to node [right] {$\RR_1$} (D);
\node (E) at (2,-.75) {$\twocommute$};
\path (0,-.75) coordinate (basepoint);
\end{tikzpicture}\nonumber
\eeq
with similar diagrams for $\sI^-_M$, $\sL^\pm_M$ and $\sR^\pm_M$. Next define the map
\beq
\Map(\R^{1|1},M)&\to& \Map(\R^{1|1},M),\label{eq:firstor} \\
(\phi\colon S\times \R^{1|1}\to M)&\mapsto &(S\times \R^{1|1}\xrightarrow{\id_{\overline{S}}\times \orr^{-1}} S\times\overline{\R}{}^{1|1}\to M),\nonumber
\eeq
which also determines maps
\beq
\R^{1|1}_{\ge 0} \times \Map(\R^{1|1},M)&\to& \R^{1|1}_{\le 0} \times \Map(\R^{1|1},M)\label{eq:firstorr}\\
\R^{1|1}_{\le 0} \times \Map(\R^{1|1},M)&\to& \R^{1|1}_{\ge 0} \times \Map(\R^{1|1},M)\label{eq:firstorrr}
\eeq
The maps~\eqref{eq:firstor}, \eqref{eq:firstorr} and~\eqref{eq:firstorrr} will all be denoted by $\orr$. 
Proposition~\ref{prop:orientation11restrict} gives 2-commuting data 
\beq
%\label{eq:moreRG}\nonumber
&&\begin{tikzpicture}[baseline=(basepoint)];
\node (A) at (0,0) {$\Map(\R^{1|1},M)$};
\node (B) at (4,0) {$\Ob(1|1\EBord(M))$};
\node (C) at (0,-1.5) {$\Map(\R^{1|1},M)$};
\node (D) at (4,-1.5) {$\Ob(1|1\EBord(M))$};
\draw[->] (A) to node [above] {$\spt^\pm_M$} (B);
\draw[->] (A) to node [left] {$\orr$} (C);
\draw[->] (C) to node [below] {$\spt^\mp_M$}  (D);
\draw[->] (B) to node [right] {$\Or_0$} (D);
\node (E) at (2,-.75) {$\twocommute$};
\path (0,-.75) coordinate (basepoint);
\end{tikzpicture}\quad 
\begin{tikzpicture}[baseline=(basepoint)];
\node (A) at (0,0) {$\R^{1|1}_{\ge 0}\times \Map(\R^{1|1},M)$};
\node (B) at (4,0) {$\Mor(1|1\EBord(M))$};
\node (C) at (0,-1.5) {$\R^{1|1}_{\le 0}\times \Map(\R^{1|1},M)$};
\node (D) at (4,-1.5) {$\Mor(1|1\EBord(M))$};
\draw[->] (A) to node [above] {$\sI^+_M$} (B);
\draw[->] (A) to node [left] {$\orr$} (C);
\draw[->] (C) to node [below] {$\sI^-_M$}  (D);
\draw[->] (B) to node [right] {$\Or_1$} (D);
\node (E) at (2,-.75) {$\twocommute$};
\path (0,-.75) coordinate (basepoint);
\end{tikzpicture}\nonumber
\eeq
with similar diagrams for $\sI^-_M$, $\sL^\pm_M$ and $\sR^\pm_M$. 

We recall from Definition~\ref{ex:realGM} that reflection structures for twists and twisted field theories are defined in terms of $\Z/2$-equivariance data in the sense of Definition~\ref{defn:equivariancedata} for the $\Z/2$-action on $1|1\EBord(M)$ determined by $\RR\circ \Or$. 

%and observe that the maps of stacks $\spt^\pm_M$ and $\spt^\mp_M\circ \Or$ are isomorphic (and hence the associated families of bordisms are also isomorphic). 
%The notion of a reflection structures is given in Definition~\ref{defn:RP0}. 

\begin{lem}\label{lem:RPtwist}
A reflection structure for a twist $\twist\colon 1|1\EBord(M)\to \TA$ determines an invertible $\twist(\spt^\pm_M)$-$\orr^*\rr^*\overline{\twist(\spt^\mp_M))}$-bimodule $\mathscr{R}_\twist^\pm\to \Map(\R^{1|1},M)$
and isomorphisms of bimodules (where the bimodule on the right is the identity bimodule)
\beq\label{eq:involutiondata}
\mathscr{R}_\twist^\pm\otimes \orr^*\rr^*\overline{\mathscr{R}}_\twist^\mp\simeq \twist(\sS^\pm_M)
\eeq
and isomorphisms of bimodules
\beq
&&\twist(\sI^\pm_M)\otimes_{\act_\pm^*\twist(\spt^\pm_M)} \act_\pm^*\mathscr{R}_\twist^\pm \xrightarrow{\sim} \proj_\pm^*\mathscr{R}_\twist^\pm \otimes_{\orr^*\rr^*\overline{\proj}_\mp^*\overline{\twist(\spt^\mp_M)}} \orr^*\rr^*\overline{\twist(\sI^\mp_M)},\nonumber
\eeq
\beq
&&\twist(\sL^\pm_M)\otimes_{\proj_\pm^*\twist(\spt^\pm_M)\otimes \act_\pm^*\twist(\spt^\mp_M)}(\proj_\pm^*\mathscr{R}_\twist^\pm\otimes \act_\pm^*\mathscr{R}_\twist^\mp) \xrightarrow{\sim} \orr^*\rr^*\overline{\twist(\sL^-_M)},\nonumber
\eeq
over $\R^{1|1}_{\ge 0}\times \Map(\R^{1|1},M)$  or $\R^{1|1}_{\le 0}\times \Map(\R^{1|1},M)$, and 
\beq
&&\twist(\sR^+_M)\xrightarrow{\sim}(\proj_\pm^*\mathscr{R}_\twist^\mp \otimes \act_\pm^*\mathscr{R}_\twist^\pm)\otimes_{\orr^*\rr^*(\overline{\proj}_\mp^*\overline{\twist(\spt^\pm_M)}\otimes \overline{\act}_\mp^*\overline{\twist(\spt^\mp_M)})}  \orr^*\rr^*\overline{\twist(\sR^-_M)},
\eeq
over $\R^{1|1}_{> 0}\times \Map(\R^{1|1},M)$ or $\R^{1|1}_{<0}\times \Map(\R^{1|1},M)$ satisfying an involutive property, where above we have used the notation for values of the twist listed in the proof of Proposition~\ref{prop:dataoftwistedFT}. 
\end{lem} 
\bp
%By Definition~\ref{defn:RP0}, a reflection structure is $\Z/2$-equivariance data in the sense of Definition~\ref{defn:equivariancedata}. 
Equivariance data demands an internal natural isomorphism, which (from Definition~\ref{rmk:twistedcategory}) is a map of stacks 
$$
\phi\colon \Ob(1|1\EBord(\pt))\to \Mor(\TA)^\times
$$
valued in invertible bimodules, together with coherence data $\rho$ satisfying properties. Applying $\phi$ to the family $\spt^\pm_M$, we obtain the claimed bundle of invertible bimodules. Applying $\rho$ to the families $\sI^\pm_M$, $\sL^\pm_M$ and $\sR^\pm_M$ from~\S\ref{sec:superpathbordism}, a reflection structure specifies data of isomorphisms of bimodules involving the values of the twist listed in the proof of Proposition~\ref{prop:dataoftwistedFT}. The other $\Z/2$-equivariance datum is an isomorphism $\beta$ between internal natural isomorphisms, see~\eqref{eq:equivariantfunctor2}; this determines the isomorphism of bimodules~\eqref{eq:involutiondata}. 
\ep

%\begin{rmk} In practice, a reflection structure for a twist comes from a $*$-structure on the algebras $\twist(\spt^\pm)$, so that the invertible bimodule is $\mathscr{R}_\twist^\pm=\twist(\spt^\pm)$ for the left action of the algebra on itself and right action determined by the $*$-structure. 
%\end{rmk}

\begin{lem} \label{lem:RPFT}
Fix a reflection structure for a twist $\twist\colon 1|1\EBord(M)\to \TA$. Then a reflection structure for a $\twist$-twisted field theory is an isomorphism of bundles of left $\twist(\spt^\pm_M)$-modules over $\Map(\R^{1|1},M)$
%
%
%is the data of isomorphisms of $\cCl_n$- and $\cCl_{-n}$-modules, respectively 
%\beq
%\V_+\stackrel{\sim}{\to} \overline{(r_0\circ \Or_0)^*\V}_- \qquad \V_-\stackrel{\sim}{\to} \overline{(r_0\circ \Or_0)^*\V}_+\label{eq:RPdataonspt}
%\eeq
\beq
E(\spt^\pm_M) \stackrel{\sim}{\to} \mathscr{R}_\twist^\pm \otimes_{\orr^*\rr^*\overline{\twist(\spt^\mp_M)}}  \orr^*\rr^*\overline{E(\spt^\mp_M))}.\label{eq:RPdataonspt}
%\overline{E(r_0 \circ \Or_0(\spt^\pm_M))}\simeq \rr^*\orr^*\overline{E(\spt^\mp_M)} 
\eeq
These isomorphisms are required to satisfy compatibility conditions for the families $\sL^\pm(M)$, $\sR^\pm(M)$, and~$\sI^\pm(M)$, namely commutative squares, 
\beq
\resizebox{.95\textwidth}{!}{$
\begin{tikzpicture}[baseline=(basepoint)];
\node (A) at (0,0) {$\twist(\sI^\pm_M)\otimes_{\proj_\pm^*\twist(\spt^\pm_M)} \proj_\pm^*E(\spt^\pm_M)$};
\node (B) at (8,0) {$\act_\pm^*E(\spt^\pm_M)$};
\node (C) at (0,-1.5) {$\orr^*\rr^*\Big(\overline{\twist(\sI^\mp_M)} \otimes_{\overline{\proj}_\mp^*\overline{\twist(\spt^\mp_M)}} \overline{\proj}_\mp^*\overline{E(\spt^\mp_M)}\Big)$};
\node (D) at (8,-1.5) {$\orr^*\rr^*\overline{\act}_\mp^* \overline{E(\spt^\mp_M)}$};
\draw[->] (A) to node [above] {$E(\sI^\pm_M)$} (B);
\draw[->] (A) to node [left] {$\simeq$} (C);
\draw[->] (C) to node [below] {$\orr^*\rr^*\overline{E(\sI^\mp_M)}$} (D);
\draw[->] (B) to node [right] {$\simeq$} (D);
\path (0,-.75) coordinate (basepoint);
\end{tikzpicture}$}\label{eq:selfadjointdiagram}\\
\resizebox{.95\textwidth}{!}{$
\begin{tikzpicture}[baseline=(basepoint)];
\node (A) at (0,0) {$\twist(\sL^\pm_M)\otimes_{\proj_\pm^*\twist(\spt^\pm_M)\otimes \act_\pm^*\twist(\spt^\mp_M)} \proj_\pm^*E(\spt^\pm_M)\otimes\act_\pm^*E(\spt^\mp_M)$};
\node (B) at (9,0) {$\underline{\C}$};
\node (C) at (0,-1.5) {$\orr^*\rr^*(\overline{\twist(\sL^\mp_M)} \otimes_{\overline{\proj}_\mp^*\overline{\twist(\spt^\mp_M)}\otimes \overline{\act}_\mp^*\overline{\twist(\spt^\pm_M)}} \overline{\proj}_\mp^*\overline{E(\spt^\mp_M)}\otimes \overline{\act}_\mp^*\overline{E(\spt^\pm_M)}$};
\node (D) at (9,-1.5) {$\orr^*\rr^* \underline{\overline{\C}}$};
\draw[->] (A) to node [above] {$E(\sL^+_M)$} (B);
\draw[->] (A) to node [left] {$\simeq$} (C);
\draw[->] (C) to node [below] {$\orr^*\rr^*\overline{E(\sL^\mp_M)}$} (D);
\draw[->] (B) to node [right] {$\simeq$} (D);
\path (0,-.75) coordinate (basepoint);
\end{tikzpicture}$}\label{eq:RPdataonL0m}
\eeq
over $\R^{1|1}_{\ge 0}\times \Map(\R^{1|1},M)$ and $\R^{1|1}_{\le 0}\times \Map(\R^{1|1},M)$, and commutative squares over $\R^{1|1}_{> 0}\times \Map(\R^{1|1},M)$ and $\R^{1|1}_{<0}\times \Map(\R^{1|1},M)$
\beq\resizebox{.95\textwidth}{!}{$
\begin{tikzpicture}[baseline=(basepoint)];
\node (A) at (0,0) {$\underline{\C}$};
\node (B) at (9,0) {$\twist(\sR^\pm_M)\otimes_{\proj_\pm^*\twist(\spt^\mp_M)\otimes \act_\pm^*\twist(\spt^\pm_M)} \proj_\pm^*E(\spt^\mp_M)\otimes\act_\pm^*E(\spt^\pm_M)$};
\node (C) at (0,-1.5) {$\orr^*\rr^* \underline{\overline{\C}}$};
\node (D) at (9,-1.5) {$\orr^*\rr^*(\overline{\twist(\sR^\mp_M)} \otimes_{\overline{\proj}_\mp^*\overline{\twist(\spt^\pm_M)}\otimes \overline{\act}_\mp^*\overline{\twist(\spt^\mp_M)}} \overline{\proj}_\mp^*\overline{E(\spt^\pm_M)}\otimes \overline{\act}_\mp^*\overline{E(\spt^\mp_M)}$};
\draw[->] (A) to node [above] {$E(\sR^\pm_M)$} (B);
\draw[->] (A) to node [left] {$\simeq$} (C);
\draw[->] (C) to node [below] {$\orr^*\rr^*\overline{E(\sR^\mp_M)}$} (D);
\draw[->] (B) to node [right] {$\simeq$} (D);
\path (0,-.75) coordinate (basepoint);
\end{tikzpicture}$}\nonumber
\eeq
where the vertical arrows use data from Lemma~\ref{lem:RPtwist} and~\eqref{eq:RPdataonspt}, though we have suppressed the tensor products with the bimodules $\mathscr{R}_\twist^\pm$ that mediate between $\twist(\spt^\pm)$-modules and $\orr^*\rr^*\overline{\twist(\spt^\mp)}$-modules to avoid (further) cluttering the notation. 
\end{lem}
\bp
Reflection positivity data for a twisted field theory is defined as an isomorphism between internal natural transformations, see Definition~\ref{defn:equivariancedata2}. From Definition~\ref{rmk:twistedcategory}, this is an isomorphism between the maps of stacks $\Ob(1|1\EBord(M))\to \Mor(\TA)$ that are part of the data of the internal natural transformation. Evaluating this isomorphism between maps of stacks on the families $\spt^\pm_M$ yields the data~\eqref{eq:RPdataonspt}. 

For the data~\eqref{eq:RPdataonspt} to determine an isomorphism between internal natural transformations, there are compatibility requirements for each object of $\Mor(1|1\EBord(M))$. Evaluating on the families $\sI^\pm_M$, $\sL^\pm_M$ and $\sR^\pm_M$ recovers the claimed commutative squares. 
\ep

We may consider real structures for twists and twisted field theories in the sense of Definition~\ref{defn:Real0}, using the real structure on $1|1\EBord(M)$ from Lemma~\ref{lem:extrastrEB}, inherited from the real structure on the $1|1$-Euclidean geometry determined by~\eqref{eq:real11}. 

\begin{lem}\label{lem:Realtwist}
A real structure for a twist $\twist\colon 1|1\EBord(M)\to \TA$ determines an invertible $\twist(\spt^\pm_M)$-$\rr^*\overline{\twist(\spt^\pm_M))}$-bimodule $\mathscr{R}_\twist^\pm\to \Map(\R^{1|1},M)$ with involutive data analogous to~\eqref{eq:involutiondata} and isomorphisms of bimodules analogous to those in Lemma~\ref{lem:RPtwist}. 
\end{lem}
\bp
The proof is analogous to that of Lemma~\ref{lem:RPtwist}. 
\ep

\begin{lem} \label{lem:RrealTw}
Fix a real structure for a twist $\twist\colon 1|1\EBord(M)\to \TA$. Then a real structure for a $\twist$-twisted field theory is an isomorphism of bundles of left $\twist(\spt^\pm_M)$-modules 
%
%
%is the data of isomorphisms of $\cCl_n$- and $\cCl_{-n}$-modules, respectively 
%\beq
%\V_+\stackrel{\sim}{\to} \overline{(r_0\circ \Or_0)^*\V}_- \qquad \V_-\stackrel{\sim}{\to} \overline{(r_0\circ \Or_0)^*\V}_+\label{eq:RPdataonspt}
%\eeq
\beq
E(\spt^\pm_M) \stackrel{\sim}{\to} \mathscr{R}_\twist^\pm \otimes_{\rr^*\overline{\twist(\spt^\pm_M)}}  \rr^*\overline{E(\spt^\pm_M))}\label{eq:RPdataonspt2}
%\overline{E(r_0 \circ \Or_0(\spt^\pm_M))}\simeq \rr^*\orr^*\overline{E(\spt^\mp_M)} 
\eeq
over $\Map(\R^{1|1},M)$. These isomorphisms are required to satisfy compatibility properties analogous to those in Lemma~\ref{lem:RPFT}. 
\end{lem}

\bp
The proof is analogous to that of Lemma~\ref{lem:RPFT}. 
\ep

\subsection{Invertible twists}\label{Sec:invertible}

The twists of interest for the conjectures~\eqref{eq:conjecture} are the degree~$n$ twists for $n\in\Z$. These are determined by an invertible twist
\beq\label{eq:thedegreetwist}
d|1\EBord(\pt)\xrightarrow{\twist} \TA^\times,\qquad d=1,2
\eeq
which then gives the degree~$n$ twists $\twist^{\otimes n}\colon d|1\EBord(M)\to \TA^\times$ by taking tensor powers and precomposing with the canonical functor $d|1\EBord(M)\to d|1\EBord(\pt)$ determined by the map $p\colon M\to \pt$. The twists $\twist^{\otimes n}$ give a $\Z$-grading on the associated collection of degree~$n$ twisted field theories for $n\in \Z$. This grading is compatible with the tensor product of field theories: the tensor product of a degree~$n$ and degree~$m$ field theory is a field theory of degree~$n+m$. Under the conjectured maps~\eqref{eq:conjecture}, the tensor product is expected to be compatible with the product in cohomology. 

When $d=1$,~\eqref{eq:thedegreetwist} determines an invertible object of the (usual) Morita bicategory, 
\beq\label{eq:valueoftwist}
A:=\twist(\spt^+)
\eeq
gotten from evaluatation on the family of bordisms $\spt^+$ over $\pt=\Map(\R^{1|1},\pt)$ from Definition~\ref{defn:11EBspt}, using that an algebra bundle over $\pt$ is just an algebra. By the cobordism hypothesis,~$A$ determines a 2-dimensional topological field theory.  Such a topological theory also determines an invertible twist via the composition
\beq\label{eq:degreetwistsketch}
1|1\EBord(\pt)\to 1\Bord(\pt)\hookrightarrow 2\Bord(\pt) \xrightarrow{{\sf TFT}(A)} \TA
\eeq
where the first arrow is a forgetful functor from $1|1$-Euclidean bordisms to 1-dimensional spin bordisms, and the second arrow includes 1-dimensional bordisms into the (fully extended, spin) 2-dimensional bordism category. The arrow labeled ${\sf TFT}(A)$ takes the fully-extended spin topological field theory determined by the invertible (and hence dualizable) algebra~$A$. 

%~\eqref{eq:thedegreetwist} is also determined by~$A$. 
%%For twists that arise as restriction of topological theories, the algebra $A$ should completely determine the twist. 
%%A more conceptual description of the degree~$n$ twist comes from viewing the assignments~\eqref{eq:degnvalues} as the restriction of a 2-dimensional topological field theory. Morally, 
%Indeed, consider the restriction of the invertible 2-dimensional topological field theory determined by~$A$ along the composition

\begin{hyp}\label{hyp1}
When $d=1$, an invertible twist $\twist$ is equivalent to the invertible twist~\eqref{eq:degreetwistsketch} constructed from the value~\eqref{eq:valueoftwist} of $\twist$ on $\spt^+$. 
\end{hyp}

Assuming this hypothesis, we can compute the values of the twist on various families of bordisms over~$M=\pt$ (compare Proposition~\ref{prop:dataoftwistedFT})
\beq
A=\twist(\spt^+),\quad A^\op=\twist(\spt^-),\quad _AA_A=\twist(\sI^+),\quad _{A^\op}{A^\op}_{A^\op}=\twist(\sI^-), \nonumber\\ 
A/[A,A]=\twist(\mathcal{L}^{1|1}(\pt))),\quad _\C A_{A\otimes A^\op}=\twist(\sL^+),\quad _\C A_{A^\op\otimes A}=\twist(\sL^-), \\\
_{A^\op\otimes A}A_\C=\twist(\sR^+),\quad  _{A\otimes A^\op}A_\C=\twist(\sR^-), \nonumber
\eeq
where the above are interpreted as bimodule bundles over the corresponding family parameter for the bordisms involved. For example, $\twist(\sI^+)$ is the trivial bundle over $\R^{1|1}_{\ge 0}$ with fiber $A$ as a bimodule over itself. The twist also requires the data of isomorphisms between various compositions of the above bordisms; this data is canonical in the case above, e.g., concatenation of superpaths corresponds to the isomorphisms
$$
A\otimes_A A\simeq A, \qquad A^\op\otimes_{A^\op} A^\op\simeq A^\op,
$$
and the relation amongst bordisms in Proposition~\ref{prop:tracerelations} corresponds to the canonical isomorphism
$$
A\otimes_{A\otimes A^\op} A\simeq A/[A,A]. 
$$
%Over~$\C$, any Morita invertible super algebra is equivalent to~$\C$ or~$\cCl_1$, and hence the twist~\eqref{eq:valueoftwist} is equivalent to one determined by $\C$ or $\cCl_1$. 

%The calculations in~\cite{ST04} of $1|1$-Euclidean field theories of degree~$n$ are consistent with the above values wit
%Alas, a characterization of 2-dimensional topological field theories in this context requires a version of the cobordism hypothesis within Stolz and Teichner's framework for geometric theories. However, assuming the usual story, the values~\eqref{eq:degnvalues} are all algebraically determined by the invertible (and hence fully dualizable) 

We refer to Definitions~\ref{defn:RP0} and~\ref{defn:Real0} for reflection and real data for a twist. Again by the cobordism hypothesis, reflection and real data for the twist determined by~\eqref{eq:valueoftwist} amounts to additional structure on~$A$. By Lemma~\ref{lem:RPtwist}, a reflection structure for the degree twist can be specified by an anti-linear anti-involution
$$
A=\twist(\spt^+)\xrightarrow{\sim} \overline{\twist(\spt^-)}=\overline{A}{}^\op,
$$
which allows us to view $A$ itself as a $A$-$\overline{A}^\op$ bimodule. 
This is the data of a $*$-superalgebra structure, which is equivalent to the data of an ordinary $*$-structure on~$A$; see~\S\ref{sec:superstar}. Similarly, Lemma~\ref{lem:Realtwist} shows that a real structure for the twist is determined by an anti-linear involution
$$
A=\twist(\spt^+)\xrightarrow{\sim} \overline{\twist(\spt^+)}=\overline{A}.
$$
i.e., a real structure on~$A$. 
%The see~\cite[\S2.C]{Sam}. 

%, respectively. As~\eqref{eq:degreetwistsketch} is determined by an invertible topological theory, this additional structure comes from an additional structure on the value of this topological field theory on the point. As we explain below in \S? below, a reflection structure is equivalent to a $*$-structure on the Clifford algebra, whereas real data is the same as a real structure on the Cliffford algebra. 

%Indeed, $\cCl_n$ is an invertible (and hence fully dualizable) object in the Morita category and so determines a fully extended 2-dimensional topological field theory by the cobordism hypothesis. 
%Indeed, the coherence data above comes from dualizability data. 

%\begin{rmk}
%An obvious generalization of Definition~\ref{defn:degreen} gives twists $\twist_V$ indexed by real inner product spaces~$V$ and determined by $T_V(\spt^+)=\cCl(V)$ with $T_V\otimes T_W\simeq T_{V\oplus W}$. 
%\end{rmk}

\subsection{The degree twist and degree~$n$ field theories}\label{Sec:invertible2}

For $\Z/2$-graded algebras over $\C$, any $\otimes$-invertible algebra is Morita equivalent to $\C$ or $\cCl_1$. Hence, assuming Hypothesis~\ref{hyp1}, the degree twist is necessarily equivalent to the twist determined by $\cCl_1$.

\begin{hyp}\label{hyp:twist}
The degree twist $\twist=\twist^{\otimes 1}$ is determined under~\eqref{eq:valueoftwist} by the Clifford algebra,~$\cCl_1$. A reflection structure for the degree twist is fixed by the standard $*$-structure~\eqref{eq:starCl} on~$\cCl_1$. A real structure for the twist is fixed by the standard real structure on $\cCl_1$. 
\end{hyp}

Assuming Hypothesis~\ref{hyp:twist}, the degree~$n$ twist is determined by the $n$th Clifford algebra~$\cCl_n$, since $\cCl_{-1}$ is the $\otimes$-inverse to $\cCl_1$ and $\cCl_{\pm 1}^{\otimes n}\simeq \cCl_{\pm n}$ for $n\ge 0$. The real and reflection structure on the degree twists uniquely determines a real and reflection structure on the degree~$n$ twist.

\begin{rmk}
The 2-dimensional spin topological field theory associated to a Clifford algebra has been studied in detail by Gunningham~\cite{Sam}. 
\end{rmk}

Propositions~\ref{prop:dataoftwistedFT} and~\ref{prop:maincompositionstatement} extract the following data from a degree~$n$ theory. 

\begin{cor} \label{cor:degreenvalues}
Assuming Hypothesis~\ref{hyp:twist}, a degree~$n$ twisted field theory $E$ determines
\begin{itemize}
\item bundles of left $\cCl_{\pm n}$-modules
$$
\V_\pm :=E(\spt^\pm_M)\to \Map(\R^{1|1},M);
$$
\item maps of bundles of $\cCl_{\pm n}$-modules over $\R^{1|1}_{\ge 0}\times \Map(\R^{1|1},M)$ and $\R^{1|1}_{\le 0}\times \Map(\R^{1|1},M)$ 
\beq
&&E(\sI^\pm_M)\colon \proj_\pm^*E(\spt^\pm_M)\to \act_\pm^*E(\spt^\pm_M),\nonumber
\eeq
compatible with composition~\eqref{eq:comppathbimodule};
\item maps of vector bundles over $\R^{1|1}_{\ge 0}\times \Map(\R^{1|1},M)$ and $\R^{1|1}_{\le 0}\times \Map(\R^{1|1},M)$ 
\beq\label{eq:Clnliearpairing}
&& E(\sL^\pm_M)\colon \cCl_{\pm n} \otimes_{\cCl_{\pm n} \otimes \cCl_{\mp n}} \Big(\proj_\pm^*E(\spt^\pm_M)\otimes \act_\pm^*E( \spt^\mp_M)\Big)\to \underline{\C} ,
\eeq
satisfying the adjunction relation~\eqref{eq:bimoduleadjunction} and the symmetry relation~\eqref{eq:bimodulesymemtry}; 

\item maps of left $\cCl_{\mp n}\otimes \cCl_{\pm n}$-modules over $\R^{1|1}_{> 0}\times \Map(\R^{1|1},M)$ and $\R^{1|1}_{< 0}\times \Map(\R^{1|1},M)$
\beq
&&E(\sR^\pm_M)\colon \cCl_{\mp n} \to \proj_\pm^*E(\spt^\mp_M) \otimes \act_\pm^*E(\spt^\pm_M);\nonumber
 \eeq
\item a section of the dual of the trivial line bundle with fiber $\cCl_n/[\cCl_n,\cCl_n]$ over the stack $\mathcal{L}^{1|1}(M)$
$$
E(\mathcal{L}^{1|1}(M))\colon \cCl_n/[\cCl_n,\cCl_n] \to \underline{\C}\nonumber
$$
that is invariant under the action of super Euclidean isometries on $\mathcal{L}^{1|1}(M)$ and agrees with the value gotten from~\eqref{Eq:bimoduletrace}. 
\end{itemize}
\end{cor}

Lemma~\ref{lem:RPFT} extracts data from a reflection structure on a degree~$n$ theory. 

\begin{cor}\label{cor:RPFT}
Assuming Hypothesis~\ref{hyp:twist}, a reflection structure on a degree~$n$ twisted field theory $E$ determines
%
%
%is the data of isomorphisms of $\cCl_n$- and $\cCl_{-n}$-modules, respectively 
%\beq
%\V_+\stackrel{\sim}{\to} \overline{(r_0\circ \Or_0)^*\V}_- \qquad \V_-\stackrel{\sim}{\to} \overline{(r_0\circ \Or_0)^*\V}_+\label{eq:RPdataonspt}
%\eeq
\begin{itemize}
\item An isomorphism of bundles of $\cCl_{\pm n}$-modules over $\Map(\R^{1|1},M)$ 
\beq
\mathcal{V}_\pm \stackrel{\sim}{\to} \orr^*\rr^*\overline{\mathcal{V}_{\mp}},\label{eq:RPdataonspt3}
%\overline{E(r_0 \circ \Or_0(\spt^\pm_M))}\simeq \rr^*\orr^*\overline{E(\spt^\mp_M)} 
\eeq
relative to the $*$-structure on $\cCl_{\pm n}$;
\item a commutative square of bundles of $\cCl_{\pm n}$-modules over $\R^{1|1}_{\ge 0}\times \Map(\R^{1|1},M)$ and $\R^{1|1}_{\le 0}\times \Map(\R^{1|1},M)$,
\beq
&&\begin{tikzpicture}[baseline=(basepoint)];
\node (A) at (0,0) {$\proj_\pm^*\mathcal{V}_\pm$};
\node (B) at (6,0) {$\act_\pm^*\mathcal{V}_\pm$};
\node (C) at (0,-1.5) {$\orr^*\rr^* \overline{\proj}_\mp^*\overline{\mathcal{V}_\mp}$};
\node (D) at (6,-1.5) {$\orr^*\rr^* \overline{\act}_\mp^*\overline{\mathcal{V}_\mp}$};
\draw[->] (A) to node [above] {$E(\sI^\pm_M)$} (B);
\draw[->] (A) to node [left] {$\simeq$} (C);
\draw[->] (C) to node [below] {$\orr^*\rr^*\overline{E(\sI^\mp_M)}$} (D);
\draw[->] (B) to node [right] {$\simeq$} (D);
\path (0,-.75) coordinate (basepoint);
\end{tikzpicture}\label{eq:selfadjointdiagram2}
\eeq
\item a commutative square of vector bundles over $\R^{1|1}_{\ge 0}\times \Map(\R^{1|1},M)$ and $\R^{1|1}_{\le 0}\times \Map(\R^{1|1},M)$,
\beq
&&\begin{tikzpicture}[baseline=(basepoint)];
\node (A) at (0,0) {$\cCl_{\pm n} \otimes_{\cCl_{\pm n} \otimes \cCl_{\mp n}} \big( \proj_\pm^*\mathcal{V}_\pm\otimes\act_\pm^*\mathcal{V}_\mp\big)$};
\node (B) at (8,0) {$\underline{\C}$};
\node (C) at (0,-1.5) {$\orr^*\rr^*\Big(\overline{\cCl}_{\mp n} \otimes_{\overline{\cCl}_{\mp n} \otimes \overline{\cCl}_{\pm n}} \big(\overline{\proj}_\mp^*\overline{\mathcal{V}_\mp}\otimes \overline{\act}_\mp^*\overline{\mathcal{V}_\pm}\big)\Big)$};
\node (D) at (8,-1.5) {$\orr^*\rr^* \underline{\overline{\C}}$};
\draw[->] (A) to node [above] {$E(\sL^\pm_M)$} (B);
\draw[->] (A) to node [left] {$\simeq$} (C);
\draw[->] (C) to node [below] {$\orr^*\rr^*\overline{E(\sL^\mp_M)}$} (D);
\draw[->] (B) to node [right] {$\simeq$} (D);
\path (0,-.75) coordinate (basepoint);
\end{tikzpicture}\label{eq:RPdataonL0m2}
\eeq
%\item commutative squares over $\R^{1|1}_{> 0}\times \Map(\R^{1|1},M)$ and $\R^{1|1}_{<0}\times \Map(\R^{1|1},M)$
%\beq
%\begin{tikzpicture}[baseline=(basepoint)];
%\node (A) at (0,0) {$\underline{\C}$};
%\node (B) at (8,0) {$\twist(\sR^\pm_M)\otimes_{\proj^*\twist(\spt^\mp_M)\otimes \act^*\twist(\spt^\pm_M)} \proj^*E(\spt^\mp_M)\otimes\act^*E(\spt^\pm_M)$};
%\node (C) at (0,-1.5) {$\orr^*\rr^* \underline{\overline{\C}}$};
%\node (D) at (8,-1.5) {$\orr^*\rr^*(\overline{\twist(\sR^\mp_M)} \otimes_{\overline{\proj}^*\overline{\twist(\spt^\pm_M)}\otimes \overline{\act}^*\overline{\twist(\spt^\mp_M)}} \overline{\proj}^*\overline{E(\spt^\pm_M)}\otimes \overline{\act}^*\overline{E(\spt^\mp_M)}$};
%\draw[->] (A) to node [above] {$E(\sR^\pm_M)$} (B);
%\draw[->] (A) to node [left] {$\simeq$} (C);
%\draw[->] (C) to node [below] {$\orr^*\rr^*\overline{E(\sR^\mp_M)}$} (D);
%\draw[->] (B) to node [right] {$\simeq$} (D);
%\path (0,-.75) coordinate (basepoint);
%\end{tikzpicture}\nonumber
%\eeq
\end{itemize}
\end{cor}
Lemma~\ref{lem:RrealTw} extracts data from a real structure on a degree~$n$ theory. 

\begin{cor}\label{cor:RFT}
Assuming Hypothesis~\ref{hyp:twist}, a real structure on a degree~$n$ twisted field theory $E$ determines
%
%
%is the data of isomorphisms of $\cCl_n$- and $\cCl_{-n}$-modules, respectively 
%\beq
%\V_+\stackrel{\sim}{\to} \overline{(r_0\circ \Or_0)^*\V}_- \qquad \V_-\stackrel{\sim}{\to} \overline{(r_0\circ \Or_0)^*\V}_+\label{eq:RPdataonspt}
%\eeq
\begin{itemize}
\item an isomorphism of bundles of $\cCl_{\pm n}$-modules over $\Map(\R^{1|1},M)$ 
\beq
\mathcal{V}_\pm \stackrel{\sim}{\to} \rr^*\overline{\mathcal{V}_{\mp}},\label{eq:Realdataonspt}
%\overline{E(r_0 \circ \Or_0(\spt^\pm_M))}\simeq \rr^*\orr^*\overline{E(\spt^\mp_M)} 
\eeq
relative to the real structure on $\cCl_{\pm n}$;
\item a commutative square of bundles of $\cCl_{\pm n}$-modules over $\R^{1|1}_{\ge 0}\times \Map(\R^{1|1},M)$ and $\R^{1|1}_{\le 0}\times \Map(\R^{1|1},M)$,
\beq
&&\begin{tikzpicture}[baseline=(basepoint)];
\node (A) at (0,0) {$\proj_\pm^*\mathcal{V}_\pm$};
\node (B) at (6,0) {$\act_\pm^*\mathcal{V}_\pm$};
\node (C) at (0,-1.5) {$\rr^* \overline{\proj}_\mp^*\overline{\mathcal{V}_\pm}$};
\node (D) at (6,-1.5) {$\rr^* \overline{\act}_\mp^*\overline{\mathcal{V}_\pm}$};
\draw[->] (A) to node [above] {$E(\sI^\pm_M)$} (B);
\draw[->] (A) to node [left] {$\simeq$} (C);
\draw[->] (C) to node [below] {$\rr^*\overline{E(\sI^\pm_M)}$} (D);
\draw[->] (B) to node [right] {$\simeq$} (D);
\path (0,-.75) coordinate (basepoint);
\end{tikzpicture}\nonumber
\eeq
\item a commutative square of vector bundles over $\R^{1|1}_{\ge 0}\times \Map(\R^{1|1},M)$ and $\R^{1|1}_{\le 0}\times \Map(\R^{1|1},M)$,
\beq
&&\begin{tikzpicture}[baseline=(basepoint)];
\node (A) at (0,0) {$\cCl_{\pm n} \otimes_{\cCl_{\pm n} \otimes \cCl_{\mp n}} \big( \proj_\pm^*\mathcal{V}_\pm\otimes\act_\pm^*\mathcal{V}_\mp\big)$};
\node (B) at (8,0) {$\underline{\C}$};
\node (C) at (0,-1.5) {$\rr^*\Big(\overline{\cCl}_{\pm n} \otimes_{\overline{\cCl}_{\pm n} \otimes \overline{\cCl}_{\mp n}} \big(\overline{\proj}_\mp^*\overline{\mathcal{V}_\pm}\otimes \overline{\act}_\mp^*\overline{\mathcal{V}_\mp}\big)\Big)$};
\node (D) at (8,-1.5) {$\rr^* \underline{\overline{\C}}$};
\draw[->] (A) to node [above] {$E(\sL^\pm_M)$} (B);
\draw[->] (A) to node [left] {$\simeq$} (C);
\draw[->] (C) to node [below] {$\rr^*\overline{E(\sL^\pm_M)}$} (D);
\draw[->] (B) to node [right] {$\simeq$} (D);
\path (0,-.75) coordinate (basepoint);
\end{tikzpicture}\nonumber
\eeq
over $\R^{1|1}_{\ge 0}\times \Map(\R^{1|1},M)$ and $\R^{1|1}_{\le 0}\times \Map(\R^{1|1},M)$.
%\item commutative squares over $\R^{1|1}_{> 0}\times \Map(\R^{1|1},M)$ and $\R^{1|1}_{<0}\times \Map(\R^{1|1},M)$
%\beq
%\begin{tikzpicture}[baseline=(basepoint)];
%\node (A) at (0,0) {$\underline{\C}$};
%\node (B) at (8,0) {$\twist(\sR^\pm_M)\otimes_{\proj^*\twist(\spt^\mp_M)\otimes \act^*\twist(\spt^\pm_M)} \proj^*E(\spt^\mp_M)\otimes\act^*E(\spt^\pm_M)$};
%\node (C) at (0,-1.5) {$\orr^*\rr^* \underline{\overline{\C}}$};
%\node (D) at (8,-1.5) {$\orr^*\rr^*(\overline{\twist(\sR^\mp_M)} \otimes_{\overline{\proj}^*\overline{\twist(\spt^\pm_M)}\otimes \overline{\act}^*\overline{\twist(\spt^\mp_M)}} \overline{\proj}^*\overline{E(\spt^\pm_M)}\otimes \overline{\act}^*\overline{E(\spt^\mp_M)}$};
%\draw[->] (A) to node [above] {$E(\sR^\pm_M)$} (B);
%\draw[->] (A) to node [left] {$\simeq$} (C);
%\draw[->] (C) to node [below] {$\orr^*\rr^*\overline{E(\sR^\mp_M)}$} (D);
%\draw[->] (B) to node [right] {$\simeq$} (D);
%\path (0,-.75) coordinate (basepoint);
%\end{tikzpicture}\nonumber
%\eeq
\end{itemize}
\end{cor}

We observe that the data and structures in Corollaries~\ref{cor:degreenvalues},~\ref{cor:RPFT}, and~\ref{cor:RFT} can be restricted along the inclusions 
\beq\label{eq:restrictalongthis}
M\subset \Pi TM\simeq \Map(\R^{0|1},M)\subset \Map(\R^{1|1},M),
\eeq
yielding bundles of $\cCl_{\pm n}$-bundles 
$$
V_\pm:=\mathcal{V}_\pm|_M,\qquad \Omega^\bullet(M;V_\pm)\simeq \mathcal{V}_\pm|_{\Pi TM}.
$$
where the isomorphism follows from Lemma~\ref{lem:sVect}. In particular, the restriction of~\eqref{eq:Clnliearpairing} and~\eqref{eq:RPdataonspt3} combine to give a map
$$
\cCl_{\pm n}\otimes_{\overline{\cCl}_{\pm n}\otimes \cCl_{\pm n}} \Big(\overline{V}_\pm \otimes V_\pm\Big)\simeq \cCl_{\mp n}\otimes_{\cCl_{\mp n}\otimes \cCl_{\pm n}} \Big(V_\mp \otimes V_\pm\Big) \to \underline{\C}
$$
which determines a pairing for which the Clifford actions are adjoint,
\beq\label{eq:Lzeropairing3}\label{eq:Cliffadjoint} 
&&\langle-,-\rangle \colon \overline{V}_+ \otimes_{\cCl_n} V_+ \to \underline{\C},\qquad \langle \omega^*\cdot x,y\rangle =(-1)^{|\omega||x|}\langle x,\omega\cdot y\rangle, \ \omega\in \cCl_n.
\eeq

\begin{lem}\label{lem:hermitian}
The pairing~\eqref{eq:Lzeropairing3} is hermitian in the $\Z/2$-graded sense
\beq
\langle x,y\rangle =(-1)^{|x||y|}\overline{\langle y,x\rangle}.\label{eq:sesquilinear}
\eeq
\end{lem}
\bp
The argument uses similar ideas to the first part of \cite[Corollary 6.25]{HST}, though the signs here are different as we are working in the oriented setting (versus the unoriented framework of \cite{HST}). 
Equation~\eqref{eq:sesquilinear} follows from commutativity of the diagram 
\beq\label{eq:hermitianargument}
\begin{tikzpicture}[baseline=(basepoint)];
\node (AA) at (0,1.5) {$\overline{\cCl}_{\pm n}\otimes_{\overline{\cCl}_{\pm n}\otimes\cCl_{\pm n}} (\orr^*\rr^*\overline{\proj}_\pm^*\overline{\V}_\pm\otimes \act_\mp^*\V_\pm)$};
\node (A) at (0,0) {$\cCl_{\mp n}\otimes_{\cCl_{\mp n}\otimes \cCl_{\pm n}}(\proj_\mp^*\V_{\mp}\otimes \act_\mp^*\V_\pm)$};
\node (B) at (8,0) {$\underline{\C}$};
\node (C) at (0,-1.5) {$\orr^*\rr^*\Big(\overline{\cCl}_{\pm n}\otimes_{\overline{\cCl}_{\pm n}\otimes \overline{\cCl}_{\mp n}} (\overline{\proj}_\pm^*\overline{\V}_\pm\otimes \overline{\act}_\pm^*\overline{\V}_\mp)\big)$};
\node (D) at (8,-1.5) {$\orr^*\rr^*\overline{\underline{\C}}$};
\node (DD) at (4,-3.5) {$\orr^*\rr^*\big(\overline{\cCl}_{\mp n}\otimes_{\overline{\cCl}_{\mp n}\otimes \overline{\cCl}_{\pm n}} (\overline{\act}_\pm^*\overline{\V}_\mp\otimes \overline{\proj}_\pm^*\overline{\V}_\pm)$};
\draw[->] (AA) to (A);
\draw[->] (A) to node [above] {$E(\sL^\mp_M)$} (B);
\draw[->] (A) to (C);
\draw[->] (C) to node [below] {$\overline{E(\sL^\pm_M)}$} (D);
\draw[->] (B) to node [right] {$\simeq$} (D);
\draw[->] (C) to node [below] {$\sigma$} (DD);
\draw[->] (DD) to node [right=10pt] {$\overline{\inv}^*\overline{E(\sL^\mp_M)}$} (D);
\draw[->, bend left=10] (AA) to node [above=8pt] {restricts to $\langle-,-\rangle$} (B);
\path (0,0) coordinate (basepoint);
\end{tikzpicture}
\eeq
where the $*$-structure on $\cCl_{\pm n}$ is used implicitly throughout. As indicated, the upper right arrow restricts to the pairing $\langle-,-\rangle$ by definition. 
Commutativity of the lower triangle follows from~\eqref{eq:bimodulesymemtry}.  Commutativity of the inner square is the condition ~\eqref{eq:RPdataonL0m2} on reflection positivity data.
%together with Lemma~\ref{lem:orientationrelations} part (3) restricted to $M\subset L_0^\pm(M)$. 
The composition of arrows starting from the top left and using the lower triangle is the right hand side of~\eqref{eq:sesquilinear}, and hence the restriction of this diagram along~\eqref{eq:restrictalongthis} implies the claimed formula. 
%Finally, consider the relation~\eqref{eq:relation4} and~\eqref{eq:relation5} for $(L^\pm_0,pr)$ (so $(t,\theta)$, respectively $(s,\eta)$, is zero). Since $E(I_0^\pm,pr)$ is the identity module map (and hence invertible), it is invertible in a neighborhood of~$0$. This shows that $E(L_0^\pm,pr)$ is nondegenerate. 
\ep

\begin{rmk} We comment in more detail about the differing signs in the above argument and \cite[Corollary 6.25]{HST}.
First, in the unoriented framework of \cite{HST} there is a $\C$-linear isomorphism $\V_\pm\simeq \orr^*\V_\mp$. This results in a $\C$-bilinear pairing through a definition analogous to~\eqref{eq:Cliffadjoint}. Second, $\orr^2=(-1)^{\sf F}$ whereas $\overline{(\rr\circ \orr)}\circ (\rr\circ \orr)=\id$. This extra factor of $(-1)^{\sf F}$ results in a pairing that is naively symmetric in \cite[Corollary 6.25]{HST}, as compared with the graded symmetry of~\eqref{eq:sesquilinear}. 
\end{rmk}

We refer to~\eqref{eq:positivitypairing} for a discussion of the signs involved when formulating the positivity property for a $\Z/2$-graded hermitian pairing. 
\begin{defn}\label{defn:RPFT} A $1|1$-Euclidean field theory of degree~$n$ with reflection positivity data is \emph{reflection positive} if the pairing~\eqref{eq:Lzeropairing3} is a (graded) positive-definite hermitian pairing on the super vector bundle~$V_+$ over $M$. 
% and the $\cCl_n$-action on $\V_+$ restricts to a $*$-representation with respect to this hermitian pairing. 
\end{defn}

\begin{rmk}
For a reflection positive field theory,~\eqref{eq:Cliffadjoint} implies that the $\cCl_n$-action on $V_+$ is automatically a $*$-representation relative to the hermitian inner product.
% meaning the $\cCl_n\simeq \overline{\cCl}_n^\op$-action on $\overline{V}_+$ is adjoint (via~\eqref{eq:Lzeropairing2}) to the $\cCl_n$-action on $V_+$.
\end{rmk}

\begin{rmk} \label{rmk:RPFT}
We contrast Definition~\ref{defn:RPFT} with the definition of positivity in~\cite{HST}, where the $1|1$-Euclidean bordism category includes the orientation reversing map as an isometry. By functoriality, this orientation-reversing isometry must be sent to a $\C$-linear map. In particular, the positivity condition in~\cite[Definition 6.26 and 6.49]{HST} is not equivalent to reflection positivity above, where orientation-reversal corresponds to a $\C$-antilinear map. However, after one imposes real structures on field theories, positivity in the sense of~\cite{HST} implies the real reflection-positive structure in the sense above. 
%Indeed, equivariance for $r\circ T$ (reflection positivity) and $r$ (reality) end up being the same as equivariance for $T$ (unoriented) and $r$ (reality), and the positivity property from reflection positivity then matches the positivity property in~\cite{HST} (modulo some additional ingrdients there involving a positive bordism category~\cite[Remark~74]{HST}). 
%Our choice of definition for positivity comes from both its relationship with established ideas in physics, and also because it is the only option that seems to generalize to the $2|1$-dimensional setting (where the $2|1$-dimensional super Euclidean geometry no longer has a real structure). 
% where the functor $r$ and $\Or$ do not exist independently, but their composition $r\circ \Or$ does. 
%This is all to say that the ubiquity of involutions on $1|1\EBord(M)$ allows for several different (and interrelated) connections with $\KO$-cocycles. 
\end{rmk}

\section{Extracting geometric data from a representations of superpaths}\label{sec:repsP}

%Given a $1|1$-dimensional, degree~$n$, real reflection positive field theory over $M$ we extract two piece of geometric data: (1) a Clifford linear, real, self-adjoint superconnections gotten from the restriction to constant superpaths in $M$ (Proposition~\ref{prop1}) and (2) a differential form gotten from restriction to constant super circles in $M$, following~\cite[\S2]{DBEChern}. These data are compatible, owing to the fact that the trace of the values on superpaths recovers the values on super circles. 
In this section we consider the restrictions~\eqref{eq:twistedrestriction} of degree~$n$ field theories along the functor~$\iota$ constructed in~\S\ref{sec:iota}. From this we extract Clifford linear superconnections with additional structure and property.

\subsection{Twisted representations of super Lie categories}
We begin with some general definitions of representations of super Lie categories. The strict 2-category of super Lie categories, smooth functors, and smooth natural transformations is recalled in Definition~\ref{defn:superLiecat}. Throughout, a vector bundle is taken to be a locally free sheaf of nuclear Fr\'echet spaces; we make no finite rank assumption. The following generalizes the standard definition of a representation of a Lie groupoid, e.g., see~\cite[Definition~1.7.1]{Mackenzie}. 

\begin{defn} \label{defn:repofSL}
A \emph{representation} of a super Lie category ${\sf C}$ is the data of a vector bundle $\V\to \Ob({\sf C})$ and a map of vector bundles $\rho\colon \s^*\V\to \t^*\V$ over $\Mor({\sf C})$. This data is required to satisfy the condition $\u^*\rho=\id_\V$ on $\Ob({\sf C})$ and $p_1^*\rho\circ p_2^*\rho=\m^*\rho$ on $\Mor({\sf C})\times_{\Ob({\sf C})}\Mor({\sf C})$. An \emph{isomorphism} $(\V,\rho)\simeq (\V',\rho')$ between representations is an isomorphism of vector bundles $\V\to \V'$ over $\Ob({\sf C})$ for which $\rho$ and $\rho'$ satisfy a compatibility property over $\Mor({\sf C})$.
\end{defn}

\begin{ex} A supermanifold $S$ determines a super Lie category ${\sf C}$ with $\Ob({\sf C})=\Mor({\sf C})=S$, i.e., a discrete super Lie category. In this case, a representation of ${\sf C}$ is the same data as a vector bundle over $S$.
\end{ex}

\begin{ex}\label{ex:supersemigroup}
Let ${\sf C}$ be the super Lie category with a single object and $\Mor({\sf C)}=\R^{1|1}_{\ge 0}$, where composition is inherited from the super semigroup structure on $\R^{1|1}_{\ge 0}$. Then a representation of ${\sf C}$ as a super Lie category is equivalent to a super semigroup representation~\eqref{eq:N1super}. 
\end{ex}

%\begin{ex}
%Generalizing the previous example, for a superconnection $\A$ on a vector bundle $V\to M$, the formula~\eqref{eq:basicrepofsuperpath} from the introduction gives a representation for the super Lie category $\{\proj,\act\colon \R^{1|1}_{\ge 0}\times \Map(\R^{0|1},M)\rightrightarrows \Map(\R^{0|1},M)\}$. Together with the conjugate bundle
%$$
%\Pi TV\coprod \Pi T\overline{V}\to \Pi TM\coprod \Pi TM\simeq \Map(\R^{0|1},M)\coprod \Map(\R^{0|1},M)=\Ob(\sP_0(M))
%$$
%and the conjugate of~\eqref{eq:basicrepofsuperpath}, one obtains a representation of $\sP_0(M)$ from a superconnection. 
%\end{ex}

\begin{defn}\label{defn:SLtwist}
A \emph{twist} $\twist$ for a super Lie category ${\sf C}$ is the data of:
\begin{enumerate}
\item $\twist_0\to \Ob({\sf C})$ a bundle of algebras;
\item $\twist_1 \to \Mor({\sf C})$ a bundle of invertible $\t^*\twist_0$-$\s^*\twist_0$-bimodules;
\item an isomorphism of bimodule bundles $\epsilon\colon \u^*\twist_1\stackrel{\sim}{\to} \twist_0$ over $\Ob({\sf C})$ where $\twist_0$ is regarded as a bimodule over itself; 
\item a map of bimodule bundles $\mu\colon p_1^*\twist_1\otimes_{p_2^*\t^*\twist_0} p_2^*\twist_1\stackrel{\sim}{\to} \m^*B\twist_1$ over the fibered product of pairs of composable morphisms, $\Mor({\sf C})\times_{\Ob({\sf C})}\Mor({\sf C})$.
%\begin{enumerate}
%\item $A_+\to \Map(\R^{0|1},M)$ and $A_-\to \Map(\R^{0|1},M)$ algebra bundles;
%\item $B_+\to \R^{1|1}_{\ge 0}\times \Map(\R^{0|1},M)$, $B_-\to \R^{1|1}_{\ge 0} \times \Map(\R^{0|1},M)$ bundles of $A_\pm$-bimodules;
%\item maps of bimodules 
%$$
%p_1^*B_\pm\otimes_{A_\pm} p_2^*B_\pm\to m^*B_\pm
%$$
%over $\R^{1|1}_{\ge 0}\times \R^{1|1}_{\ge 0}\times \Map(\R^{0|1},M)$ and $\R^{1|1}_{\le 0}\times \R^{1|1}_{\le 0}\times \Map(\R^{0|1},M)$, where $m$ is composition of superpaths, and $p_1,p_2$ are the projections 
%$$
%\R^{1|1}_{\ge 0}\times \R^{1|1}_{\ge 0}\times \Map(\R^{0|1},M)\to \R^{1|1}_{\ge 0}\times \Map(\R^{0|1},M)\quad \R^{1|1}_{\le 0}\times \R^{1|1}_{\le 0}\times \Map(\R^{0|1},M)\to \R^{1|1}_{\le 0}\times \Map(\R^{0|1},M)
%$$
\end{enumerate}
The maps of $\twist_0$-bimodules are required to satisfy further coherence conditions for composition with the identity and associativity of composition (see Lemma~\ref{lem:fancyrep}). An \emph{isomorphism} $\twist\simeq \twist'$ between twists is an invertible $\twist_0$-$\twist_0'$ bimodule $\mathscr{B}$ over $\Ob({\sf C})$ and an isomorphism $\twist_1\otimes_{\s^*\twist_0} \s^*\mathscr{B}\simeq \t^*\mathscr{B}\otimes_{\t^*\twist_0'}\twist_1'$ of bimodules over $\Mor({\sf C})$ satisfying a property over $\Mor({\sf C})\times_{\Ob({\sf C})}\Mor({\sf C})$.
\end{defn}

\begin{rmk}
When ${\sf C}$ is a Lie groupoid and the fibers of $\twist_0\to \Ob({\sf C})$ are central simple algebras, the above definition reduces to Freed's notion of an invertible algebra bundle over~${\sf C}$ \cite[Definition 1.59]{Freedalg}. In particular, gerbes provide examples of twists for the (\v{C}ech) groupoid associated to an open cover~\cite[Example 1.74]{Freedalg}. Invertible algebra bundles over quotient Lie groupoids determine twists for (equivariant) K-theory~\cite[Definition 1.78]{Freedalg}, see also Karoubi~\cite{KaroubiFrench} and Donovan--Karoubi \cite{DonovanKaroubi}. All the examples of twistings below will have algebra bundles whose fibers are central simple. 
\end{rmk} 

\begin{rmk}
Definition~\ref{defn:SLtwist} is set up so that Lemma~\ref{lem:fancyrep} holds using \cite[Definition 5.1]{ST11} for the internal Morita category $\TA$. As a consequence, Definition~\ref{defn:SLtwist} does not play well with equivalences of super Lie categories, see Remark~\ref{rmk:badcat}. There are various well-known ways to fix this issue. 
%In short, one needs an appropriate stackification of Definition~\ref{defn:Morita} of the internal Morita category. 
For example, define a new version of twisting consisting of an equivalence of super Lie categories ${\sf C}'\xrightarrow{\sim} {\sf C}$ and a twist for ${\sf C}'$ in the above sense; this is essentially the definition of twistings for the K-theory of topological groupoids in~\cite{FHTI}. In all the examples below, twistings fit into the more naive framework of Definition~\ref{defn:SLtwist}.
\end{rmk}

\begin{defn}\label{def:twrep}
Given a twist for a super Lie category ${\sf C}$, a \emph{twisted representation is}
\begin{enumerate}
\item $\eta\colon \V\to \Ob({\sf C})$ a bundle of left $\twist_0$-modules;
\item $\rho\colon \twist_1\otimes_{\s^*\twist_0} \s^*\V\to \t^*\V$ a map of bundles of left $\t^*\twist_0$-modules over $\Mor({\sf C})$
%a map of $A$-bimodules $p_1^*B\otimes_Ap_2^*B\to B$ over $\Mor(\sP(M))\times_{\Ob(\sP(M))}\Mor(\sP(M))$, where $
%p_1,p_2,m\colon \Mor(\sP(M))\times_{\Ob(\sP(M))}\Mor(\sP(M))\to \Mor(\sP(M))$ are the projections and composition. 
%\begin{enumerate}
%\item $A_+\to \Map(\R^{0|1},M)$ and $A_-\to \Map(\R^{0|1},M)$ algebra bundles;
%\item $B_+\to \R^{1|1}_{\ge 0}\times \Map(\R^{0|1},M)$, $B_-\to \R^{1|1}_{\ge 0} \times \Map(\R^{0|1},M)$ bundles of $A_\pm$-bimodules;
%\item maps of bimodules 
%$$
%p_1^*B_\pm\otimes_{A_\pm} p_2^*B_\pm\to m^*B_\pm
%$$
%over $\R^{1|1}_{\ge 0}\times \R^{1|1}_{\ge 0}\times \Map(\R^{0|1},M)$ and $\R^{1|1}_{\le 0}\times \R^{1|1}_{\le 0}\times \Map(\R^{0|1},M)$, where $m$ is composition of superpaths, and $p_1,p_2$ are the projections 
%$$
%\R^{1|1}_{\ge 0}\times \R^{1|1}_{\ge 0}\times \Map(\R^{0|1},M)\to \R^{1|1}_{\ge 0}\times \Map(\R^{0|1},M)\quad \R^{1|1}_{\le 0}\times \R^{1|1}_{\le 0}\times \Map(\R^{0|1},M)\to \R^{1|1}_{\le 0}\times \Map(\R^{0|1},M)
%$$
\end{enumerate}
with the condition that the module maps (2) are compatible with composition in ${\sf C}$. Given an isomorphism between twists, and isomorphism between twisted representations $(\V,\rho)\simeq (\V',\rho')$ is an isomorphism of left modules $\mathscr{B}\otimes_{\twist_0} \V\simeq \V'$ over $\Ob({\sf C})$ with a compatibility condition over $\Mor({\sf C})$. 
\end{defn}

\begin{defn} \label{defn:algebratwist}
For any super Lie category ${\sf C}$ and super algebra~$A$, define the twist $\twist_A$ given by the trivial bundle of algebras over ${\rm Ob}({\sf C})$ with fiber $A$, and trivial bundle of bimodules over $\Mor({\sf C})$ with fiber $A$ regarded as a bimodule over itself. The isomorphism $\epsilon$ is the identity, and $\mu$ is determined by the isomorphism $A\otimes_A A\simeq A$. We observe that Morita equivalent algebras give rise to isomorphic twists. For $A=\C$, the twist $\twist_\C$ is the \emph{trivial twist}, which we will also denote by $\one\colon {\sf C}\to \TA$. 
\end{defn}

The following is immediate from the definitions, compare~\cite[Lemma~5.7]{ST11}. 

\begin{lem}
A twisted representation for the trivial twist determines a representation of ${\sf C}$ in the sense of Definition~\ref{defn:repofSL}. 
\end{lem}

%In brief, twisted representations are representations valued in modules over an algebra. Hence, this specializes to ordinary representations when the algebra is~$\C$; compare~\cite[Lemma 5.7]{ST11}. 
\begin{lem} \label{lem:fancyrep}
A twist for a super Lie category ${\sf C}$ is the same as an internal functor $\twist\colon {\sf C}\to \TA$, and a $\twist$-twisted representation is equivalent to an internal natural transformation 
\beq
&&\begin{tikzpicture}[baseline=(basepoint)];
\node (A) at (0,0) {${\sf C}$};
\node (B) at (4,0) {$\TA,$};
\node (C) at (2.1,0) {$E \Downarrow$};
\draw[->,bend left=12] (A) to node [above] {$\one$} (B);
\draw[->,bend right=12] (A) to node [below] {$\twist$} (B);
\path (0,0) coordinate (basepoint);
\end{tikzpicture}\nonumber
\eeq
where $\TA$ is the internal Morita category (see Definition~\ref{defn:Morita}). Isomorphisms between twists in the sense of Definition~\ref{defn:SLtwist} are equivalent to internal natural isomorphisms between $\twist,\twist'\colon {\sf C}\to \TA$. Isomorphisms between twisted representations in the sense of Definition~\ref{def:twrep} are the same as isomorphisms between internal natural transformations $E\simeq E'$. 
\end{lem}

\bp
The claim follows from unpacking the definitions. Super Lie categories can be regarded as categories internal to stacks; see \S\ref{sec:STBord} below. Using that $\Ob({\sf C})$ and $\Mor({\sf C})$ are representable stacks, the data of the internal functor~$\twist$ (see Definition~\ref{defn:internalfunctor}) is the same as a twist in the sense of Definition~\ref{defn:SLtwist}. The data of an internal natural transformation~$E$ (see Definition~\ref{rmk:twistedcategory}) is equivalent to a twisted representation $(\V,\rho)$ in the sense of Definition~\ref{def:twrep}. The statements for isomorphisms between twists and isomorphisms between twisted representations follow similarly. 
%The twisted representation theory of super Lie categories defined above corresponds to an internal natural transformation
\ep

\begin{ex} 
Let ${\sf C}$ be the super Lie category from Example~\ref{ex:supersemigroup}. Consider the twist from Definition~\ref{defn:algebratwist} relative to the algebra $\cCl_n$. 
%\to \pt$ and the trivial bimodule bundle $_{\cCl_n}(\cCl_n)_{\cCl_n}\times \R^{1|1}_{\ge 0}\to \R^{1|1}_{\ge 0}$ determines a twist for~${\sf C}$. 
Then a twisted representation of ${\sf C}$ is a Clifford linear super semigroup representation $\R^{1|1}_{\ge 0}\to \End_{\cCl_n}(\V)$ for a $\cCl_n$-module~$\V$.
\end{ex}

\subsection{Real representations}
Lemma~\ref{lem:fancyrep} allows us to import Stolz and Teichner's notions of Real representations (e.g.,~\cite[\S6.8]{HST}) to the setting of super Lie categories. This turns out to provide a generalization of Atiyah's Real vector bundles, see Example~\ref{Ex:Atiyah} below. To start, we observe that the conjugation functor on supermanifolds (see Definition~\ref{defn:realsmfld}) induces a strict involution on the strict 2-category of super Lie categories. We use the following notation and terminology for the values of this involution. 

\begin{defn} For a super Lie category ${\sf C}$, the \emph{conjugate} super Lie category, denoted $\overline{\sf C}$, is defined by applying the conjugation functor on supermanifolds to the objects, morphisms, and structure maps of ${\sf C}$. Similarly, for a smooth functor $F\colon {\sf C}\to {\sf D}$ between super Lie categories, let $\overline{F}\colon \overline{\sf C}\to \overline{\sf D}$ be the \emph{conjugate} smooth functor between conjugate categories, and for a smooth natural transformation $\eta\colon F\Rightarrow G$, let $\overline{\eta}\colon \overline{F}\Rightarrow \overline{G}$ denote the \emph{conjugate} smooth natural transformation between the conjugate smooth functors. 
\end{defn}

The following is a special case of a real structure on a category internal to stacks, see Definition~\ref{defn:interinvol}.

\begin{defn}\label{defn:RealLie} A \emph{Real} structure on a super Lie category ${\sf C}$ is a smooth functor $\RRR\colon {\sf C}\to \overline{\sf C}$ and a smooth natural isomorphism $\eta\colon \id_{\sf C}\Rightarrow \overline{\RRR}\circ \RRR$. A super Lie category with a real structure is a \emph{Real super Lie category}, denoted $({\sf C},\RRR,\eta)$. 
\end{defn}

The data of the natural isomorphism $\eta\colon \id_{\sf C}\Rightarrow \overline{\RRR}\circ \RRR$ will often be trivial in examples below, and so we often omit it from the notation. Real (twisted) representations of a super Lie category are then defined in terms of $\Z/2$-equivariance data for the $\Z/2$-action on ${\sf C}$ given by a Real structure, and the $\Z/2$-action on $\TA$ given by complex conjugation of vector bundles. 
%We recall that for categories internal to stacks, there are isomorphisms between internal natural transformations coming from isomorphisms between maps of stacks; see Definition~\ref{rmk:twistedcategory}. 
The following is a special case of Definitions~\ref{defn:equivariancedata} and~\ref{defn:equivariancedata2}.

\begin{defn}\label{defn:Realrep}
For ${\sf C}$ be a super Lie category with Real structure $\RRR$, a \emph{Real structure} for a twist $\twist$ is an internal natural isomorphism~$\phi_\twist$
\beq
&&\begin{tikzpicture}[baseline=(basepoint)];
\node (A) at (0,0) {${\sf C}$};
\node (C) at (0,-1) {$\overline{\sf C}$};
\node (B) at (4,0) {$\TA$};
\node (D) at (4,-1) {$\TA,$};
\node (E) at (2,-.5) {$\phi_\twist \ \twocommute$};
\draw[->] (A) to node [above] {$\twist$} (B);
\draw[->] (A) to node [left] {$\RRR$} (C);
\draw[->] (C) to node [below] {$\overline{\twist}$} (D);
\draw[->] (B) to node [right] {$(\overline{\phantom{A}})$} (D);
\path (0,-.5) coordinate (basepoint);
\end{tikzpicture}\nonumber
\eeq
and an isomorphism $\beta_\twist$ of internal natural transformations $\phi_\twist^2\circ \eta_{\sf C}\simeq \eta_{\TA}=\id$, see~\eqref{eq:equivariantfunctor2}. The data $(\twist,\phi_\twist,\beta_\twist)$ is a \emph{Real twist}. A \emph{Real structure} for a $\twist$-twisted representation $E$ is a Real structure for $\twist$ and an isomorphism between internal natural transformations $\phi_\twist \circ E\simeq E\circ \phi_\one$ compatible with $\beta_\twist$ using the canonical real structure for the internal functor~$\one$. 

% $\sigma\colon \eta\circ( \id_{(\overline{\phantom{A}})}\circ\overline{E}\circ \id_{\RRR})) \circ c\simeq E\colon \one \to \twist$ in the diagram
%\beq
%&&\begin{tikzpicture}[baseline=(basepoint)];
%\node (A) at (0,0) {${\sf C}$};
%\node (C) at (0,-1) {$\overline{\sf C}$};
%\node (B) at (4,0) {$\TA$};
%\node (D) at (4,-1) {$\TA,$};
%\node (E) at (2,-.5) {$ \Uparrow \ \eta$};
%\draw[->] (A) to node [above] {$\Downarrow E$} (B);
%\draw[->,bend left=25] (A) to node [above] {$\one$} (B);
%\draw[->] (A) to node [left] {$\RRR$} (C);
%\draw[->] (C) to node [below] {$\Uparrow \overline{E}$} (D);
%\draw[->,bend right=25] (C) to node [below] {$\overline{\one}$} (D);
%%\draw[->,bend right=55] (C) to node [below] {$(\rr^{-1})^*\one$} (D);
%\draw[->] (D) to node [right] {$(\overline{\phantom{A}})$} (B);
%\path (0,-.5) coordinate (basepoint);
%\end{tikzpicture}\label{eq:nattransest}
%\eeq
%using the identity $\one= (\overline{\phantom{A}})\circ \overline{\one}\circ \RRR$, where $\sigma$ is further required to satisfy the involutive property~\eqref{eq:involproperty}. 
% since both functors classify the trivial algebra bundle over $\Ob({\sf C})$ and trivial bimodule bundle over $\Mor({\sf C})$
\end{defn}

\begin{lem}\label{lem:RealData} For a Real super Lie category $({\sf C},\RRR)$ and twist $\twist$, a Real structure on $\twist$ determines data
\begin{enumerate}
\item[(i)]
a bundle $N \to \Ob({\sf C})$ of invertible $\twist_0$-${\RRR_0}^*\overline{\twist_0}$  bimodules with an isomorphism
\beq\label{eq:dataofbeta}
N\otimes_{\twist_0} {\overline{\RRR}_0}^*\overline{N}\simeq \twist_0
\eeq
of bimodule bundles over $\Ob({\sf C})$;
\item[(ii)] an isomorphism of $\t^*\twist_0$-$\s^*\RRR_0^*\overline{\twist_0}$ bimodule bundles over $\Mor({\sf C})$
$$
\twist_1\otimes_{\s^*\twist_0} \s^*N\xrightarrow{\sim} \t^*N\otimes_{\t^*\RRR_0^*\overline{\twist_0}} \RRR_1^*\overline{\twist_1}
$$
\end{enumerate}
satisfying a compatibility condition on $\Mor({\sf C})\times_{\Ob({\sf C})}\Mor({\sf C})$. 

Given a Real twist, a Real twisted representation determines additional data 
\begin{enumerate}
\item[(iii)] an isomorphism of left $\twist_0$-modules $N\otimes_{\RRR_0^*\overline{\twist_0}} \RRR_0^*\overline{\V}\xrightarrow{\sim} \V$ over $\Ob({\sf C})$
\end{enumerate}
satisfying the property that the diagram of $\t^*\RRR_0\overline{\twist_0}$-$\s^*\twist_0$ bimodule bundles over $\Mor({\sf C})$ commutes
\beq
&&\begin{tikzpicture}[baseline=(basepoint)];
\node (A) at (-4,0) {$\twist_1\otimes_{\s^*\twist_0} \s^*(N\otimes_{\RRR_0^*\overline{\twist_0}} \RRR_0^*\overline{\V})$};
\node (B) at (-4,-1.5) {$\t^*N\otimes_{\t^*\RRR_0^* \overline{\twist_0}} \RRR_1^*(\overline{\twist_1} \otimes_{\s^*\overline{\twist_0}} \s^*\overline{\V})$};
\node (C) at (2,0) {$\twist_1\otimes_{\s^*\twist_0} \s^*\V$};
\node (D) at (2,-1.5) {$\t^*\V$};
\node (E) at (-1,-3) {$\t^*(N\otimes_{\RRR_0^*\overline{\twist_0}} \RRR_0^*\overline{\V})$};
\draw[->] (C) to node [right] {$\rho$} (D);
\draw[->] (A) to node [left] {${\rm {\rm (ii)}}\otimes \id$} (B);
\draw[->] (A) to node [above] {$\id\otimes {\rm (iii)}$} (C);
\draw[->] (B) to node [left=5pt] {$\id\otimes \RRR_1^*\overline{\rho}$} (E);
\draw[->] (E) to node [right=5pt] {$\t^*{\rm (iii)}$} (D);
\path (0,-.75) coordinate (basepoint);
\end{tikzpicture}\label{diag:reflectionpositivitycondition}
\eeq
with maps defined by the data above as indicated.
\end{lem}

\bp
By Definition~\ref{rmk:twistedcategory}, the first piece of data of the internal natural transformation $\phi_\twist$ is a map of stacks
\beq\label{eq:bimoduleN}
N\colon \Ob({\sf C})\to \Mor(\TA), \nonumber
\eeq
which by the Yoneda lemma is equivalent to a bundle of bimodules over the supermanifold $\Ob({\sf C})$. The source of this bimodule is the algebra bundle $\RRR_0^*\overline{\twist}_0$, and the target is $\twist_0$. As the internal natural transformation is invertible, this bundle of bimodules is also required to be invertible. The isomorphism~\eqref{eq:dataofbeta} is the datum $\beta_\twist$. This provides the data in part (i) of the lemma. Part (ii) supplies the datum denoted $\rho$ in~\eqref{diag:nu} for an internal natural transformation. The condition on the data (i) and (ii) similarly follows from the definition of an internal natural transformation over $\Mor({\sf C})$ and $\Mor({\sf C})\times_{\Ob({\sf C})}\Mor({\sf C})$. 

The data of an isomorphism between internal natural transformations is an isomorphism between their corresponding maps of stacks $\Ob({\sf C})\to \Mor(\TA)$; see Definition~\ref{rmk:twistedcategory}. In this case, this is precisely the datum~(iii). The diagram~\eqref{diag:reflectionpositivitycondition} is the coherence condition so that this isomorphism of maps of stacks determines an isomorphism between internal natural transformations. 
\ep

\begin{ex}\label{Ex:Atiyah}
For an ordinary manifold $X$ regarded as a super Lie category, an involution $\alpha\colon X\to X$ determines a Real structure on~$X$ by setting $\RRR$ to be the composition
$$
X\xrightarrow{\alpha} X\simeq \overline{X}
$$ 
for the canonical real structure on $X$ coming from complex conjugation of smooth functions. Taking the trivial twist with its canonical Real structure, a Real representation of $(X,\alpha)$ is then the same data as a Real vector bundle in the sense of Atiyah~\cite{AtiyahReal}. 
\end{ex}

\subsection{Degree~$n$ representations of superpaths} 

Now we specialize to the super Lie category $\sP_0(M)$ of nearly constant superpaths in~$M$. 

\begin{defn}\label{defn:Cliffordtwist} The \emph{degree~$n$ twist} for the super Lie category $\sP_0(M)$ is the coproduct of trivial algebra bundles on objects,
$$
\underline{\cCl}_n\coprod \underline{\cCl}_{-n}\to \Map(\R^{0|1},M)\coprod \Map(\R^{0|1},M)=\Ob(\sP_0(M))
$$
and the coproduct of trivial bimodule bundles over morphisms
$$
\underline{\cCl}_n\coprod \underline{\cCl}_{-n} \to \R^{1|1}_{\ge 0} \times \Map(\R^{0|1},M)\coprod \R^{1|1}_{\le 0}\times \Map(\R^{0|1},M)=\Mor(\sP_0(M)).
$$
The remaining data for this twist is specified by the identity isomorphism $\epsilon=\id\colon \cCl_{\pm n}\to \cCl_{\pm n}$ and the canonical isomorphisms of Clifford bimodules, $
\mu\colon \cCl_{\pm n}\otimes_{\cCl_{\pm n}}\cCl_{\pm n}\stackrel{\sim}{\to} \cCl_{\pm n}$ (for $\epsilon$ and $\mu$ in the notation of Definition~\ref{defn:SLtwist}). The required properties in this case are trivially satisfied. 
\end{defn}

%\begin{rmk} We observe that an $n$-dimensional metrized vector bundle $V\to M$ determines a twist for the super Lie category $\sP_0(M)$ given by the bundles of algebras $\cCl(V)\coprod \cCl(V)^\op\to \Ob(\sP_0(M))$ with the identity bimodules over morphisms. A choice of spin structure affords an isomorphism from this twist to the degree~$n$ twist, using Stolz and Teichner's definition of spin structure as a bundle of invertible $\cCl(V)$-$\cCl_n$-bimodules \cite[Definition 2.3.1]{ST04}. When $V$ does not have a spin structure, modules over the bundle $\cCl(V)$ determine classes in twisted K-theory, e.g., see Karoubi~\cite{KaroubiFrench}, Donovan--Karoubi \cite{DonovanKaroubi} and Freed--Hopkins--Teleman~\cite[\S3.6]{FHTI}. Hence, twisted representations of $\sP_0(M)$ offer one inroad to a variant of the first conjecture in~\eqref{eq:conjecture} involving twisted K-theory; see also~\cite{TwistAugusto}. \end{rmk}

\begin{defn}\label{defn:degreensP}
A \emph{degree $n$ representation} of $\sP_0(M)$ is a twisted representation for the degree~$n$ twist. 
\end{defn}

\begin{prop} \label{prop:EFT1}
Assuming Hypothesis~\ref{hyp:twist}, the restriction of a degree~$n$ field theory along~\eqref{eq:twistedrestriction} determines a degree~$n$ representation of $\sP_0(M)$. 
\end{prop}
\bp
The value of the restriction on objects and morphisms of $\sP_0(M)$ is immediate from Corollary~\ref{cor:degreenvalues}, after further restricting the data there to the subspace
\beq\label{eq:constantpathsubspace}
\Map(\R^{0|1},M)\subset \Map(\R^{1|1},M),
\eeq
which includes along the $S$-families of maps that factor through the quotient
$$
S\times \R^{1|1}\to (S\times \R^{1|1})/\R\simeq S\times \R^{0|1}\to M. 
$$
The compatibility with source, target, and units is clear; compatibility with composition comes from restricting~\eqref{eq:comppathbimodule} along~\eqref{eq:constantpathsubspace}. 
\ep

\begin{prop}\label{prop:superconn}
A degree~$n$ representation of $\sP_0(M)$ 
%on a  vector bundle 
%$$
%\V_+\coprod \V_-\to \Map(\R^{0|1},M)\coprod \Map(\R^{0|1},M)=\Ob(\sP_0(M))
%$$ 
determines $\cCl_{\pm n}$-linear superconnections $\A_\pm$ on super vector bundles $V_\pm \to M$.
\end{prop}

%\begin{enumerate}
%\item a bundle $V_+\to M$ of $\cCl_n$-modules and a bundle $V_-\to M$ of $\cCl_n^\op$-modules;
%\item Clifford linear super semigroup representations satisfying a Leibniz-type rule
%\beq
%\rho_+\colon \R^{1|1}_{\ge 0}\to \End(\Omega^\bullet(M;V_+)), && \rho_-\colon \R^{1|1}_{\le 0} \to \End(\Omega^\bullet(M;V_-)) \label{eq:semigroupmaps}\\
%\rho_\pm(t,\theta) (\alpha v)=(\alpha-\theta d\alpha)\rho_\pm (t,\theta)(v),&& \alpha\in \Omega^\bullet(M), \ v\in \Omega^\bullet(M;V_\pm).\label{eq:Leibniz} 
%\eeq
%
%%where we have identified $C^\infty(\Map(\R^{0|1},M))\simeq \Omega^\bullet(M)$ and $\Gamma(\Map(\R^{0|1},M);V_\pm)\simeq \Omega^\bullet(M;V_\pm)$. 
%\end{enumerate}
%\end{prop}
\bp
Following Definition~\ref{def:twrep}, let $\V=\V_+\coprod \V_-$ denote the value of the representation on objects, where $\V_\pm\to \Map(\R^{0|1},M)$ is a bundle of $\cCl_{\pm n}$-modules and $\Ob(\sP_0(M))=\Map(\R^{0|1},M)\coprod \Map(\R^{0|1},M)$. Using Lemma~\ref{lem:sVect}, we make the identification
\beq\label{eq:vbidentification}
\Gamma(\Map(\R^{0|1},M);\V_\pm)\simeq \Omega^\bullet(M;V_\pm),\qquad V_\pm=\V_\pm|_M
\eeq
for vector bundles $V_\pm\to M$ gotten by restriction along the canonical inclusion $M\hookrightarrow \Map(\R^{0|1},M)$. 

The map of left $A$-modules in Definition~\ref{def:twrep} translates in this case to maps $\rho_\pm$ of $\cCl_{\pm n}$-module bundles 
\beq
&&\rho_+\colon \proj_+^*\V_+\to \act_+^*\V_+\ {\rm over} \ \R_{\ge 0}^{1|1}\times \Map(\R^{0|1},M) \nonumber\\
&&\rho_-\colon \proj_-^*\V_-\to \act_-^*\V_-\ {\rm over} \ \R_{\le 0}^{1|1}\times \Map(\R^{0|1},M)\nonumber
\eeq
for the projection and action maps gotten by restriction of~\eqref{eq:actproj}. 
%Identifying $C^\infty(\Map(\R^{0|1},M))\simeq \Omega^\bullet(M)$ and applying Lemma~\ref{lem:action}, the above vector bundle maps are maps of $\Omega^\bullet(M)$-modules for the action on $\s^*V_\pm$ by pulling back along the projection, and the action on $\t^*V_\pm$ twisted by the restriction of~\eqref{eq:actdeRham}. 
%We similarly use the same notation $\rho_\pm$ for vector bundle maps involving the $\Omega^\bullet(M)$-modules $\Omega^\bullet(M;V_\pm)$. A
Using~\eqref{eq:vbidentification} and applying Lemma~\ref{lem:action}, $\rho_\pm$ are maps of $\Omega^\bullet(M)$-modules for the action on the source and target gotten by pulling back along the projection or the restriction of~\eqref{eq:actdeRham}. This gives the description 
\beq
%\rho_+\in C^\infty(\R^{1|1}_{\ge 0},\End_{\cCl_{\pm n}}(\Omega^\bullet(M;V_+))), && \rho_-\in C^\infty(\R^{1|1}_{\le 0},\End_{\cCl_{\pm n}}(\Omega^\bullet(M;V_-)))\label{eq:semigroupmaps}\\
%\rho_+ \colon \R^{1|1}_{\ge 0}\to \End_{\cCl_{n}}(\Omega^\bullet(M;V_+)),&& \rho_- \colon \R^{1|1}_{\le 0}\to \End_{\cCl_{-n}}(\Omega^\bullet(M;V_-))  \label{eq:semigroupmaps}\\
\rho_\pm(t,\theta,\alpha v)=(\alpha-\theta d\alpha)\rho_\pm (t,\theta,v),&& \alpha\in \Omega^\bullet(M), \ v\in \Omega^\bullet(M;V_\pm)\label{eq:Leibniz} 
\eeq
for the standard coordinates $t$ and $\theta$ on $\R^{1|1}_{\ge 0}$ and $\R^{1|1}_{\le 0}$.
%where $\End_{\cCl_{\pm n}}(\Omega^\bullet(M;V_+)))$ denotes the presheaf on supermanifolds determined by vector space of endomorphisms of $\Omega^\bullet(M;V_\pm)$ commuting with the $\cCl_{\pm n}$-actions; see Example~\ref{ex:tvspresheaf}. 
Compatibility with composition promotes $\rho_+$ and $\rho_-$ to super semigroup representations, since multiplication in $\R^{1|1}_{\ge 0}$ and $\R^{1|1}_{\le 0}$ corresponds to composition in~$\sP_0(M)$. 

%The maps~\eqref{eq:semigroupmaps} are equivalently functions 

Taylor expanding~\eqref{eq:Leibniz}  in the odd variable $-\theta\in C^\infty(\R^{1|1}_{\ge 0})$, define
\beq\label{eq:AB}
\rho_\pm(t,\theta,-)=A_\pm(t)- \theta B_\pm(t)
\eeq
%\beq
%&& A_\pm (t)\in \End_{\cCl_{\pm n}}(\Omega^\bullet(M;V_\pm))^\ev,\ B_\pm (t)\in \End_{\cCl_{\pm n}}(\Omega^\bullet(M;V_\pm))^\odd 
%\eeq
where $A_\pm(t)\colon \Omega^\bullet(M;V_\pm)\to \Omega^\bullet(M;V_\pm)$ is an even map and $B_\pm (t)\colon \Omega^\bullet(M;V_\pm)\to \Omega^\bullet(M;V_\pm)$ is odd, and both commute with the Clifford actions. Define
\beq
&&\A_\pm:=B_\pm(0)\colon \Omega^\bullet(M;V_\pm)\to \Omega^\bullet(M;V_\pm).\label{eq:superconnconstrr}
\eeq
Since $B_\pm (t)$ is odd, so is $\A_\pm$. Furthermore, because $B_\pm(t)$ is $\cCl_{\pm n}$-linear, so is $\A_\pm$. Finally, combining~\eqref{eq:Leibniz} and~\eqref{eq:AB} yields the pair of equations
$$
A_\pm(t)(\alpha v)=\alpha A_\pm(t)(v),\qquad B_\pm(t)(\alpha v)=d\alpha A_\pm(t)(v)+(-1)^{|\alpha|}\alpha B_\pm(t)(v).
$$
Using that $A_\pm(0)=\id_{V_\pm}$, the value of the second equation at $t=0$ is the Leibniz rule for the superconnections~$\A_\pm$. 
\ep

Under the above translation, the renormalization group flow corresponds to a standard structure in index theory. 

\begin{defn}[{\cite[page~267]{BGV}}]
Define the \emph{(Getzler) rescaling} of a Clifford linear superconnection by the $\R_{>0}$-action~\eqref{eq:RGsconn}.
\end{defn}

The action of the renormalization group on $1|1$-Euclidean bordisms is computed in Lemma~\ref{prop:RG11restrict} and restricts to the category of nearly constant superpaths via~\eqref{eq:RGOrRR}.

\begin{prop} \label{prop:Getzler}
The renormalization group flow on a degree~$n$ field theory restricts along~\eqref{eq:twistedrestriction} to determine the rescaling of the superconnections extracted in Proposition~\ref{prop:superconn}. 
\end{prop}

\bp
Using~\eqref{eq:RGOrRR}, we can compute using the action of the renormalization group on $\sP_0(M)$. This action was computed in~\eqref{eq:RGonsP}. In particular, we have an isomorphism
$$
\V_\pm\xrightarrow{\sim}\RG_\mu^*\V_\pm,\qquad s\mapsto \mu^{\deg/2}s,\qquad s\in \Gamma(\V_\pm)
$$
viewing the source and target as $\Omega^\bullet(M)$-modules as in~\eqref{eq:vbidentification}, where the target module structure is twisted by pulling back along the RG-action. Hence, under the above identification a representation $\rho_\pm(t,\theta)$ of $\sP_0(M)$ is sent under $\RG_\mu$ to $\mu^{-\deg/2}\rho(\mu^2t,\mu\theta)\mu^{\deg/2}$. From the  definition~\eqref{eq:superconnconstrr} of the superconnections $\A_+$ and $\A_-$ in terms of the representations $\rho_\pm(t,\theta)$ this yields the formula~\eqref{eq:RGsconn}. 
\ep

\subsection{Degree~$n$ representations with a pairing}
The restriction of the map~\eqref{eq:2morphsigma} determines
$$
\inv\colon \R^{1|1}_{\ge 0}\times \Map(\R^{0|1},M)\to \R^{1|1}_{\le 0}\times \Map(\R^{0|1},M),\quad (t,\theta,\phi)\mapsto (-t,-\theta,\phi\circ T^{-1}_{t,\theta}).
$$

\begin{lem} The pullback of $\rho_-$ along $\inv$ gives a map of vector bundles 
$$
\inv^*\rho_-\colon \act_+^*\V_-\to \proj_+^*\V_-
$$
over $\R^{1|1}_{\ge 0}\times \Map(\R^{0|1},M)$ commuting with the $\cCl_{-n}$-action. 
\end{lem}
\bp We observe that the triangles commute,
\beq
\begin{tikzpicture}[baseline=(basepoint)];
\node (A) at (0,0) {$\R^{1|1}_{\ge 0}\times \Map(\R^{0|1},M)$};
\node (C) at (2,-2) {$\Map(\R^{0|1},M)$};
\node (B) at (4,0) {$\R^{1|1}_{\le 0}\times \Map(\R^{0|1},M)$};
\draw[->] (A) to node [above] {$\inv$} (B);
\draw[->] (A) to node [left=5pt] {$\act_+$} (C);
\draw[->] (B) to node [right=5pt] {$\proj_-$} (C);
\path (0,-.5) coordinate (basepoint);
\end{tikzpicture}\quad
\begin{tikzpicture}[baseline=(basepoint)];
\node (A) at (0,0) {$\R^{1|1}_{\ge 0}\times \Map(\R^{0|1},M)$};
\node (C) at (2,-2) {$\Map(\R^{0|1},M)$};
\node (B) at (4,0) {$\R^{1|1}_{\le 0}\times \Map(\R^{0|1},M)$};
\draw[->] (A) to node [above] {$\inv$} (B);
\draw[->] (A) to node [left=5pt] {$\proj_+$} (C);
\draw[->] (B) to node [right=5pt] {$\act_-$} (C);
\path (0,-.5) coordinate (basepoint);
\end{tikzpicture}\nonumber
\eeq
using that 
$$
\act_\pm(t,\theta,\phi)=\phi\circ T^{-1}_{t,\theta},\qquad \proj_\pm(t,\theta,\phi)=\phi.
$$
Hence, we have
$$
\act_+^*\V_-\simeq \inv^*\proj_-^*\V_-\xrightarrow{\inv^*\rho_-} \inv^*\act_-^*\V_-\simeq \proj_+\V_-
$$
so $\inv^*\rho_-$ uniquely determines a map with the claimed source and target. 
\ep

\begin{defn} \label{defn:pairing}
A \emph{pairing} on a degree~$n$ representation of $\sP_0(M)$ is a vector bundle map over $\Map(\R^{0|1},M)$
\beq\label{eq:sPpairing}
L\colon \cCl_n\otimes_{\cCl_{-n}\otimes \cCl_n} (\V_-\otimes \V_+)\to \underline{\C}
\eeq
with the property that the diagram of vector bundles over $\R^{1|1}_{\ge 0}\times \Map(\R^{0|1},M)$ commutes 
\beq
&&\begin{tikzpicture}[baseline=(basepoint)];
\node (A) at (0,0) {$\cCl_n\otimes_{\cCl_{-n}\otimes \cCl_n}(\act_+^*\V_-\otimes \proj_+^*\V_+)$};
\node (C) at (0,-1.5) {$\proj_+^*(\cCl_n\otimes_{\cCl_{-n}\otimes \cCl_n}(\V_-\otimes \V_+))$};
\node (B) at (7,0) {$\act^*_+(\cCl_n\otimes_{\cCl_{-n}\otimes \cCl_n}(\V_-\otimes \V_+))$};
\node (D) at (7,-1.5) {$\underline{\C}.$};
\draw[->] (A) to node [above] {$\id\otimes \rho_+$} (B);
\draw[->] (A) to node [left] {$\inv^*\rho_-\otimes \id$} (C);
\draw[->] (C) to node [below] {$\proj_+^*L$} (D);
\draw[->] (B) to node [right] {$\act_+^*L$} (D);
\path (0,-.5) coordinate (basepoint);
\end{tikzpicture}\label{eq:adjunction}
\eeq
\end{defn}

%\beq
%&&\begin{tikzpicture}[baseline=(basepoint)];
%\node (A) at (-4,0) {$\inv^*\proj^*\V_-\otimes \proj^*\V_+$};
%\node (B) at (2,0) {$\inv^*\act^*\V_-\otimes \proj^*\V_+$};
%\node (C) at (-7,-1.5) {$\act^*\V_-\otimes \proj^*\V_+$};
%\node (D) at (2,-1.5) {$\proj^*(\V_-\otimes \V_+)$};
%\node (E) at (-4,-3) {$\act^*(\V_-\otimes \V_+)$};
%\node (F) at (2,-3) {$\underline{\C}$};
%\draw[->] (A) to node [above] {$\inv^*\rho_-\otimes \id$} (B);
%\draw[->] (C) to node [left] {$\simeq$} (A);
%\draw[->] (C) to node [left] {$\id\otimes\rho_+$} (E);
%\draw[->] (E) to node [below] {$\act^*\langle-,-\rangle$} (F);
%\draw[->] (D) to node [right] {$\proj^*\langle-,-\rangle$} (F);
%\draw[->] (B) to node [right] {$\simeq$} (D);
%\path (0,-.75) coordinate (basepoint);
%\end{tikzpicture}\nonumber
%\eeq

\begin{lem}\label{prop:EFT2} Assuming Hypothesis~\ref{hyp:twist}, the restriction of a degree~$n$ field theory along~\eqref{eq:twistedrestriction} equips the resulting degree~$n$ representation of $\sP_0(M)$ with a pairing. 
\end{lem}

\bp The datum~\eqref{eq:sPpairing} follows from restricting $E(\sL^-_M)$ in Corollary~\ref{cor:degreenvalues} along~\eqref{eq:constantpathsubspace}. The property~\eqref{eq:adjunction} follows from~\eqref{eq:bimoduleadjunction} in Proposition~\ref{prop:maincompositionstatement}. 
\ep

\begin{prop} \label{eq:proppairing}
Given a pairing on a degree~$n$ representation of $\sP_0(M)$, the Clifford linear superconnections $\A_+,\A_-$ from Proposition~\ref{prop:superconn} satisfy the adjunction formula,
\beq\label{eq:initialadjunction}
L(\A_-v,w)+(-1)^{|v|}L( v,\A_+w)=dL( v,w). 
\eeq
\end{prop}
\bp
Using the descriptions of $\proj^*_+\V_\pm$ and $\act^*_+\V_\pm$ in terms of modules over $\Omega^\bullet(M;C^\infty(\R^{1|1}_{\ge 0}))$ and $\Omega^\bullet(M;C^\infty(\R^{1|1}_{\le 0}))$ (see Lemma~\ref{lem:stmap}) and the notation from Proposition~\ref{prop:superconn} we compute
\beq
(\proj_+^*L)((\inv^*\rho_-)(t,\theta,v),w)&=&L(\rho_-(-t,-\theta,v),w)=L( (A_-(-t)+\theta B_-(-t))v,w)\nonumber\\
&=&(\act^*_+L)( v,\rho_+(t,\theta,w))\nonumber\\
&=&L( v,\rho_+(t,\theta,w)\rangle+\theta d\langle v,\rho_+(t,\theta,w))\nonumber\\
&=&L( v,(A_+(t)-\theta B_+(t))w\rangle+\theta d\langle v,(A_+(t)-\theta B_+(t))w), \nonumber
\eeq
where the equalities above are between elements of $\Omega^\bullet(M;C^\infty(\R^{1|1}_{\ge 0}))$. 
Collecting the terms involving $\theta$, setting $t=0$, and using the definition of the superconnections $\A_+$ and $\A_-$, the result follows. 
%The sign on $\theta d$ comes from the fact that the pairing takes values in $\act^*\underline{\C}$, and then we are taking the isomorphism back to $\underline{\C}$ coming from the action by the inverse, hence the flip in sign.
\ep

\subsection{Reflection structures and real structures}
Using the functors $\RR$ and $\Or$ from Lemma~\ref{lem:structureonsP}, $(\sP_0(M),\RR\circ \Or)$ is a Real super Lie category in the sense of Definition~\ref{defn:RealLie}, where $\overline{(\RR\circ \Or)}\circ (\RR\circ \Or)=\id$ (and so $\eta$ is the identity natural isomorphism). 

\begin{defn}\label{defn:reflectionstructure}
A \emph{reflection structure} for the degree~$n$ twist is a Real structure for the twist relative to $(\sP_0(M),\RR\circ \Or)$.
\end{defn}

We refer to \S\ref{sec:superstar} for the definition of a $*$-superalgebra.

\begin{lem} A $*$-superalgebra structure on $\cCl_{\pm n}$ determines a reflection structure for the degree~$n$ twist. 
\end{lem}
\bp
By Definition~\ref{defn:Realrep} and Lemma~\ref{lem:RealData}, a Real structure is an invertible $\twist_0-\Or_0^*\RR_0^*\overline{\twist_0}$-bimodule
over $\Ob(\sP_0(M))\simeq \Map(\R^{0|1},M)\coprod \Map(\R^{0|1},M)$ with additional involutive data. For the degree~$n$ twist it suffices to consider the case $M=\pt$. Such an invertible bimodule can be constructed from isomorphisms of algebras 
\beq\label{eq:whatisstar}
\cCl_{\pm n}\simeq \overline{\cCl}_{\mp n}=\overline{\cCl}_{\pm n}^\op
\eeq
satisfying an involutive property. Indeed, the above allows us to view $\cCl_n$ as a $\overline{\cCl}_{\pm n}^\op$-$\cCl_n$ bimodule. An involutive algebra isomorphism~\eqref{eq:whatisstar} is precisely the data of a $*$-superalgebra structure, proving the lemma.
\ep

Hereafter, we fix a reflection structure for the degree~$n$ twist from the following $*$-superalgebra structure on $\cCl_{\pm n}$.

\begin{defn} Define a $*$-superalgebra structure on $\cCl_n$ given by the $\C$-antilinear extension of the map on generators,
\beq\label{eq:superstarCl}
f_j^*=-if_j,\qquad e_j^*=ie_j.
\eeq
using the notation of generators for the Clifford algebra from~\eqref{eq:Clifford}. 
\end{defn}

\begin{defn}\label{defn:rssP}
Suppose we are given a degree~$n$ representation $E$ of $\sP_0(M)$ with a pairing~$L$. For the reflection structure for the degree~$n$ twist fixed by~\eqref{eq:superstarCl}, define a \emph{reflection structure} on $(E,L)$ as a Real twisted representation (see Definition~\ref{defn:Realrep}) with the additional property that the the pairing,
$$
\langle-,-\rangle\colon \Or_0^*\RR_0^*\overline{\V}_+\otimes_{\cCl_n} \V_+\simeq \cCl_n\otimes_{\overline{\cCl}_n \otimes \cCl_n}(\Or_0^*\RR_0^*\overline{\V}_+\otimes \V_+)\simeq 
\cCl_n\otimes_{\cCl_{-n} \otimes \cCl_n}(\V_-\otimes \V_+)\xrightarrow{L}\underline{\C}
$$
is hermitian in the $\Z/2$-graded sense, where the second isomorphism above makes use of the Real structure on the representation. A degree~$n$ representation with a reflection structure is \emph{reflection positive} if the above pairing restricts along $M\subset \Map(\R^{0|1},M)$ to a positive pairing (again in the $\Z/2$-graded sense).
%diagram commutes
%\beq
%&&\begin{tikzpicture}[baseline=(basepoint)];
%\node (A) at (0,0) {$\cCl_n\otimes_{\overline{\cCl}_n \otimes \cCl_n}(\orr^*\rr^*\overline{\V}_+\otimes \V_+)$};
%\node (C) at (0,-1.5) {$\cCl_n\otimes_{\cCl_{-n} \otimes \cCl_n}(\V_-\otimes \V_+)$};
%\node (B) at (7,0) {$\overline{\cCl}_n\otimes_{\overline{\cCl}_n \otimes \overline{\cCl}_{-n}}(\orr^*\rr^*(\overline{\V}_+\otimes \overline{\V}_-)$};
%\node (D) at (7,-1.5) {$\overline{\cCl}_n\otimes_{\overline{\cCl}_n \otimes \overline{\cCl}_{-n}}(\orr^*\rr^*(\overline{\V}_-\otimes \overline{\V}_+)$};
%\node (E) at (0,-3) {$\underline{\C}$};
%\node (F) at (7,-3) {$\orr^*\rr^*\overline{\underline{\C}}$};
%\draw[->] (A) to node [above] {$\id\otimes \rp_+$} (B);
%\draw[->] (A) to node [left] {$\rp_-^{-1}\otimes \id$} (C);
%\draw[->] (C) to node [left] {$L$} (E);
%\draw[->] (D) to node [right] {$\orr^*\rr^*\overline{L}$} (F);
%\draw[->] (E) to node [below] {$\simeq$} (F);
%\draw[->] (B) to node [right] {$\sigma$} (D);
%\path (0,-1.5) coordinate (basepoint);
%\end{tikzpicture}\label{eq:symmetryofpairing}
%\eeq
%I think we need to spell out $\orr^*\rr^*\overline{\V}_+$ as a module over $\Omega^\bullet(M)$: same vector space as section of $\V_+$ but module structure is twisted by algebra map. 
\end{defn}

\begin{prop} \label{prop:EFT3}
Assuming Hypothesis~\ref{hyp:twist}, the restriction of a reflection structure on a degree~$n$, $1|1$-Euclidean field theory along~\eqref{eq:twistedrestriction} determines a reflection structure on a degree~$n$ representation of $\sP_0(M)$. The restriction of a degree~$n$ reflection positive $1|1$-Euclidean field theory along~\eqref{eq:twistedrestriction} determines a degree~$n$ reflection positive representation of $\sP_0(M)$. 
\end{prop}
\bp
Hypothesis~\ref{hyp:twist} states that the reflection structure for the degree~$n$ twist $\twist^{\otimes n}\colon 1|1\EBord(M)\to \TA$ comes from the $*$-superalgebra structure on $\cCl_n$. If we then further examine the data in Proposition~\ref{lem:RPtwist} and Corollary~\ref{cor:degreenvalues}, the restriction of a reflection structure along~\eqref{eq:constantpathsubspace} evidentially determines a Real degree~$n$ representation for the real structure on $\sP_0(M)$ given by $\RR\circ \Or$. The symmetry of the pairing in Definition~\ref{defn:rssP} follows from Lemma~\ref{lem:hermitian}. Finally, we observe that the positivity condition in Definition~\ref{defn:RPFT} is equivalent to the positivity condition in Definition~\ref{defn:rssP}. 
\ep

\begin{prop} \label{prop:selfadjoint} 
The $\cCl_{n}$-linear superconnection $\A_+$ associated with a degree~$n$ reflection positive representation is self-adjoint with respect to the pairing from Definition~\ref{defn:rssP}.
\end{prop}

\bp
A reflection structure gives isomorphisms of vector bundles over $\Map(\R^{0|1},M)$, 
\beq\label{eq:rsisos}
\V_\mp\stackrel{\sim}{\to} \Or_0^*\RR_0^*\overline{\V}_\pm.
\eeq 
Relative to these isomorphisms, we have the identities 
$$
\rho_\pm(t,\theta)=(\orr^*\rr^*\overline{\rho}_\mp)(t,\theta)=i^\deg\overline{\rho}_\mp(-t,i\theta).
$$
Taylor expanding the above expression (as in~\eqref{eq:AB}), we find that the superconnection $\A_-$ is sent to $-i^{\deg+1}\overline{\A_+}$ under the isomorphism~\eqref{eq:rsisos}. Then by the adjunction formula~\eqref{eq:initialadjunction} and the definition of the hermitian pairing~\eqref{eq:sesquilinear}, we obtain the desired self-adjointness property for $\A_+$. 
\ep

We may also consider Real degree~$n$ representations of $\sP_0(M)$ for the Real structure given by $\RR$. This leads to a discussion analogous to the one for reflection structures. 

\begin{lem}
For the Real super Lie category $(\sP_0(M),\RR)$, a Real structure for the degree~$n$ twist is determined by a real structure on $\cCl_n$. 
\end{lem}

\begin{defn} \label{def:realCl}
Define an antilinear involution $\cCl_n\to \overline{\cCl}_n$ by declaring the generators of the Clifford algebra in~\eqref{eq:Clifford} to be real, and take the extension to a $\C$-antilinear map. 
\end{defn}

 \begin{defn}\label{defn:rrsP}
Suppose we are given a degree~$n$ representation $E$ of $\sP_0(M)$ with a pairing~$L$. For the Real structure for the degree~$n$ twist fixed by Definition~\ref{def:realCl}, define a \emph{real, reflection positive degree~$n$ representation of $\sP_0(M)$} as a Real structure on $E$ relative to $\RR$ (in the sense of Definition~\ref{defn:Realrep}), a reflection structure on $E$ satisfying the reflection positivity condition, and the further condition that the pairing satisfies the reality condition
\beq\label{eq:realcondition}
\langle x,y\rangle=(-1)^{|\overline{x}|} \RR_0^*\overline{\langle \overline{x},\overline{y}\rangle}
\eeq
for $x,y\in \Gamma(V_+)$, and $\overline{x},\overline{y}$ their image under the isomorphism $V_+\xrightarrow{\sim}\RR_0^*\overline{V}_+$. 

%diagram commutes
%\beq
%&&\begin{tikzpicture}[baseline=(basepoint)];
%\node (A) at (0,0) {$\cCl_n\otimes_{\overline{\cCl}_n \otimes \cCl_n}(\orr^*\rr^*\overline{\V}_+\otimes \V_+)$};
%\node (C) at (0,-1.5) {$\cCl_n\otimes_{\cCl_{-n} \otimes \cCl_n}(\V_-\otimes \V_+)$};
%\node (B) at (7,0) {$\overline{\cCl}_n\otimes_{\overline{\cCl}_n \otimes \overline{\cCl}_{-n}}(\orr^*\rr^*(\overline{\V}_+\otimes \overline{\V}_-)$};
%\node (D) at (7,-1.5) {$\overline{\cCl}_n\otimes_{\overline{\cCl}_n \otimes \overline{\cCl}_{-n}}(\orr^*\rr^*(\overline{\V}_-\otimes \overline{\V}_+)$};
%\node (E) at (0,-3) {$\underline{\C}$};
%\node (F) at (7,-3) {$\orr^*\rr^*\overline{\underline{\C}}$};
%\draw[->] (A) to node [above] {$\id\otimes \rp_+$} (B);
%\draw[->] (A) to node [left] {$\rp_-^{-1}\otimes \id$} (C);
%\draw[->] (C) to node [left] {$L$} (E);
%\draw[->] (D) to node [right] {$\orr^*\rr^*\overline{L}$} (F);
%\draw[->] (E) to node [below] {$\simeq$} (F);
%\draw[->] (B) to node [right] {$\sigma$} (D);
%\path (0,-1.5) coordinate (basepoint);
%\end{tikzpicture}\label{eq:symmetryofpairing}
%\eeq
%I think we need to spell out $\orr^*\rr^*\overline{\V}_+$ as a module over $\Omega^\bullet(M)$: same vector space as section of $\V_+$ but module structure is twisted by algebra map. 
\end{defn}

\begin{prop} \label{prop:reality}
Assuming Hypothesis~\ref{hyp:twist}, the restriction of a real, reflection positive, degree~$n$, $1|1$-Euclidean field theory along~\eqref{eq:twistedrestriction} determines a real, reflection positive degree~$n$ representation of $\sP_0(M)$. 
\end{prop}
\bp
The reality data on the degree~$n$ representation $E$ of $\sP_0(M)$ is clear. It remains to check the condition on the pairing. This follows from the argument in the proof of \cite[Theorem 6.48]{HST}, where the relation between the reflection and real structures in~\eqref{Eq:Orrelation} introduces the sign $(-1)^{|\overline{x}|}$. 
%We recall that Hypothesis~\ref{hyp:twist} states that the reflection structure for the degree~$n$ twist $\twist^{\otimes n}\colon 1|1\EBord(M)\to \TA$ comes from the $*$-superalgebra structure on $\cCl_n$. If we then further examine the data in Proposition~\ref{lem:RPtwist} and Corollary~\ref{cor:degreenvalues}, the restriction of a reflection structure along~\eqref{eq:constantpathsubspace} evidentially determines a Real degree~$n$ representation for the real structure on $\sP_0(M)$ given by $\RR\circ \Or$. The symmetry of the pairing in Definition~\ref{defn:rssP} follows from Lemma~\ref{lem:hermitian}. Finally, we observe that the positivity condition in Definition~\ref{defn:RPFT} is equivalent to the positivity condition in Definition~\ref{defn:rssP}. 
\ep

\begin{prop} \label{prop:realA} 
A real, reflection positive degree~$n$ representation of $\sP_0(M)$ equips the vector bundle $V=V_+$ with a real structure such that the pairing~\eqref{eq:realcondition} comes from an ordinary (i.e., not $\Z/2$-graded) real pairing on the underlying real vector bundle and the self-adjoint superconnection $\A=\A_+$ is a real superconnection.
\end{prop}

\bp
A Real structure for a degree~$n$ representation with respect to the involution $\RR$ gives isomorphisms of bundles of $\cCl_{\pm n}$-modules over $\Map(\R^{0|1},M)$, 
\beq\label{eq:risos}
\V_\pm\stackrel{\sim}{\to} \RR^*_0\overline{\V}_\pm,
\eeq 
relative to the chosen real structure on $\cCl_n$. This gives an ordinary real structure on $V_+,V_-\to M$. The super semigroup representations satisfy the identity
\beq\label{eq:realdataonC}
\rho_\pm(t,\theta)=(\RR^*_1\overline{\rho}_\pm)(t,\theta).
\eeq
Taylor expanding the above expression (as in~\eqref{eq:AB}), the associated superconnections are real. The condition~\eqref{eq:realcondition} implies that the pairing is real for pairings between even sections, zero for pairings between even and odd sections, and purely imaginary for pairings between odd sections. The translation between a super Hilbert space and a $\Z/2$-graded Hilbert space is (compare~\cite[page~91]{strings1})
\beq\label{Eq:metric}
( v,w):=\left\{\begin{array}{ll} \langle v,w\rangle & v,w \in \Gamma(V_+)^\ev \\ i^{-1}\langle v,w\rangle & v,w\in \Gamma(V_+)^\odd \\ 0 & {\rm else} \end{array}\right.
\eeq
and so we see that the condition~\eqref{eq:realcondition} leads to a real pairing~\eqref{Eq:metric} in the usual sense.
\ep

\subsection{Trace class representations} 
%reference Lemma~\ref{lem:spanofinclude}?

We recall the super Euclidean loop space $\mathcal{L}^{1|1}(M)$ from Definition~\ref{defn:sEucloop} with its action by $\R^{1|1}\rtimes \Z/2$. We observe that the loop space and its equivariant structure are natural in the manifold $M$, and hence we have a map of stacks
\beq\label{eq:sloopover}
\mathcal{L}^{1|1}(M)\sq (\R^{1|1}\rtimes \Z/2) \to \mathcal{L}^{1|1}(\pt)\sq (\R^{1|1}\rtimes \Z/2) = \R^{1|1}_{>0}\sq (\R^{1|1}\rtimes \Z/2 ).
\eeq

\begin{defn} \label{defn:Pfaffian}
The \emph{Pfaffian line bundle} $\Pf\to \mathcal{L}^{1|1}(M)\sq \R^{1|1}\rtimes \Z/2$ is the pullback along \eqref{eq:sloopover} of the odd line bundle classified by
$$
\R^{1|1}_{>0}\sq (\R^{1|1}\rtimes \Z/2) \xrightarrow{p} \pt\sq \Z/2=\pt\sq \{\pm 1\}\to \pt\sq \C^\times
$$
where $p$ is determined by the projections $\R^{1|1}_{>0}\to \pt$ and $\R^{1|1}\rtimes \Z/2\to \Z/2$, and we take the standard $\C^\times$-action on $\Pi \C$, the trivial odd line bundle over $\pt$. The Pfaffian line bundle is natural for maps of stacks $\mathcal{L}^{1|1}(M)\sq (\R^{1|1}\rtimes \Z/2) \to\mathcal{L}^{1|1}(N)\sq (\R^{1|1}\rtimes \Z/2)$ induced by maps of smooth manifolds $M\to N$. 
\end{defn}

\begin{rmk} 
Identify the quotient stack $\R_{>0}\sq \R\rtimes \Z/2$ with the stack of metrized circles with periodic spin structure, where $\Z/2$ acts by the spin flip on the spinor bundle, and $\R$ acts by rotation of the circle. The Dirac operator on the circle for the chosen spin structure gives an $\R\rtimes \Z/2$-equivariant family of operators over $\R_{>0}$. The terminology in Definition~\ref{defn:Pfaffian} comes from the fact that the restriction of $\Pf$ along the canonical inclusion $\R_{>0}\sq \R\rtimes \Z/2\subset \R^{1|1}_{>0}\sq \R^{1|1}\rtimes \Z/2$ is the Pfaffian line bundle of the family of Dirac operators. This Pfaffian line bundle appears in Stolz and Teichner's sketch of the degree~$n$ twist for the Euclidean bordism category without supersymmetry; see~\cite[\S5.3]{ST11}. 
\end{rmk}

\begin{prop} For $n\in \Z$, let $\Pf^{\otimes n}$ denote the $n$th tensor power of the Pfaffian line bundle over $\mathcal{L}^{1|1}(M)\sq (\R^{1|1}\rtimes \Z/2)$. Sections are in bijection with
$$
\Gamma(\mathcal{L}_0^{1|1}(M)\sq (\R^{1|1}\rtimes \Z/2);\Pf^{\otimes n}) \simeq \left\{Z,Z_\ell\in \Omega^\bullet(M;C^\infty(\R_{>0}))\mid \begin{array}{c} dZ=0, \ \partial_\ell Z=dZ_\ell \\ \deg(Z)=n \mod 2, \\ \deg(Z_\ell)=n-1\mod 2\end{array}\right\}.
$$
\end{prop} 

\begin{proof} 
Under the identification 
$$
C^\infty(\R^{1|1}_{>0}\times \Map(\R^{0|1},M))\simeq \Omega^\bullet(M;C^\infty(\R^{1|1}_{>0}))
$$
sections of $\Pf^{\otimes n}$ can be identified with a differential form valued in $C^\infty(\R^{1|1}_{>0})$ with equivariance properties for the action by $\R^{1|1}\rtimes \Z/2$. Choosing coordinates $(\ell,\lambda)$ on $\R^{1|1}_{>0}$ and Taylor expanding a section $s$ in the odd variable $\lambda$, consider the parameterization 
$$
s=\ell^{\deg/2}Z+2i\lambda \ell^{(\deg+1)/2}Z_\ell,\qquad Z,Z_\ell\in\Omega^\bullet(M;C^\infty(\R_{>0})).
$$
By \cite[Proposition~1.4]{DBEChern}, $\R^{1|1}$-invariance is equivalent to the conditions $dZ=0$ and $\partial_\ell Z=dZ_\ell$. The degree condition on $(Z,Z_\ell)$ comes from equivariance for the $\Z/2$-action which acts by $(-1)^\deg$ on differential forms and $(\ell,\lambda)\mapsto (\ell,\pm \lambda)$. The result follows.
\ep

Define the element $\Gamma_n$ of the $n$th Clifford algebra
\beq\label{eq:Gamma}
\Gamma_n:=\left\{\begin{array}{ll} 2^{-n/2}f_1f_2\dots f_n & n \ge 0 \\  2^{n/2}e_1e_2\dots e_n & n<0\end{array}\right.
\eeq
in the notation from~\eqref{eq:Clifford}. 
The \emph{Clifford super trace} of a $\cCl_n$-linear endomorphism $T\colon V\to V$ is 
$$
\sTr_{\cCl_n}(T):=\sTr(\Gamma_n\circ T),
$$
see~\S\ref{sec:supertrace} below. 

\begin{defn}\label{defn:tracestru}
A \emph{trace class representation of $\sP_0(M)$ of degree~$n$} is the data of a degree~$n$ representation $E$ of $\sP_0(M)$ and a section $(Z(E),Z(E)_\ell)\in \Gamma(\mathcal{L}^{1|1}_0(M);\Pf^{\otimes -n})$ satisfying the property 
\beq\label{eq:traceclass}\label{eq:constpartitionfun2}
&&
\resizebox{.9\textwidth}{!}{$Z(E)=\ell^{-\deg/2}\Tr_{\cCl_n}(\rho_+|_{\R_{>0}\times \Map(\R^{0|1},M)})\in C^\infty(\R_{>0}\times \Map(\R^{0|1},M))\simeq \Omega^\bullet(M;\C^\infty(\R_{>0}))$}
\eeq
where $\Tr_{\cCl_n}$ denotes the Clifford super trace. 
\end{defn}

\begin{lem} \label{lem:4periodic}
For a trace class degree~$n$ representation, $Z(E)\in \Omega_\cl^\bullet(M;C^\infty(\R_{>0}))$ is a closed differential form in degrees equal to the parity of~$n$ (i.e., either even or odd degree). When $E$ is a real field theory, $Z(E)$ is a real differential form. When $E$ is real and reflection positive, $Z(E)$ is real and concentrated in degrees equal to $n$ mod 4.  \end{lem}

\bp
For the first, note that the element $\Gamma_n\in \Cl_n\subset \cCl_n$ has parity equal to the parity of $n$, and so the endomorphism $\Gamma_n\circ \rho_+$ also has parity~$n$ mod~2. Hence, it follows for any degree~$n$ representation for which the trace~\eqref{eq:traceclass} is defined, the resulting function is even or odd according to this parity. The real structure preserves $\Gamma_n\in \cCl_{n}$, and the reality condition~\eqref{eq:realdataonC} implies that $Z(E)$ be equal to its conjugate. Since $\Gamma_n\in \Cl_n$ determines a real trivialization (see~\eqref{eq:Gamma}), this amounts to demanding that $Z(E)\in C^\infty(\R_{>0}\times \Map(\R^{0|1},M))\simeq \Omega^\bullet(M;C^\infty(\R_{>0}))$ is real as a differential form, for the standard conjugation action on $\C$-valued differential forms. When $E$ is also reflection positive, we observe that $\Gamma_n\mapsto i^n\Gamma_n$ under the $*$-superalgebra structure~\eqref{eq:superstarCl}. This imposes an additional invariance property on $Z(E)$, implying that it takes values in 4-periodic differential forms.
% refining the 2-periodic version above implied by the parity condition on $Z(E)$.
\ep

\begin{lem}
Assuming Hypothesis~\ref{hyp:twist}, a choice of orientation on $\R^n$ determines an isomorphism $\Pf^{\otimes n} \simeq \twist^{\otimes n}(\mathcal{L}^{1|1}(M)\sq (\R^{1|1}\rtimes \Z/2))$ between the Pfaffian line bundle and the valued of the degree~$n$ twist on $\mathcal{L}^{1|1}(M)\sq (\R^{1|1}\rtimes \Z/2)$. In particular, a degree~$n$ field theory determines a section of the dual of $\Pf^{\otimes n}$.
\end{lem}

\bp
It suffices to take $M=\pt$. Then the value of the twist on $\mathcal{L}^{1|1}(\pt)$ is the line bundle with fiber 
%By Corollary~\ref{cor:degreenvalues}, a degree~$n$ field theory~$E$ determines a map of line bundles over $\R_{>0}\times \Map(\R^{0|1},M)$ 
$$
\cCl_{n}\otimes_{\cCl_{n}\otimes \cCl_{- n}}\cCl_n\simeq \cCl_n/[\cCl_n,\cCl_n]. 
$$
A choice of orientation $\Gamma_n\in \Lambda^{\rm top}(\C^n)\simeq \cCl_n/[\cCl_n,\cCl_n]$ trivializes this line, and so after this choice it remains to compute the equivariance structure. Under Hypothesis~\ref{hyp:twist}, the $\R^{1|1}$-action is trivial, and the $\Z/2$-action is by parity. This completes the proof. 
%
%\beq\label{eq:constpartitionfun}
%&&E(C_\ell,\eval)\colon \cCl_n/[\cCl_n,\cCl_n] \to \underline{\C},  
%\eeq
%where the source denotes the trivial bundle with fiber $\cCl_n/[\cCl_n,\cCl_n]$, and we have used
%\beq\label{eq:constpartitionfun2}
%E(C_\ell,\eval)\in C^\infty(\R_{>0}\times \Map(\R^{0|1},M))\simeq \Omega^\bullet(M;C^\infty(\R_{>0})),
%\eeq
%and hence a differential form on $M$ with values in $C^\infty(\R_{>0})$. 
\ep

\begin{prop} \label{prop:trace}
Assuming Hypothesis~\ref{hyp:twist}, the restriction of a degree~$n$ field theory along~\eqref{eq:twistedrestriction} determines a trace class degree~$n$ representation of $\sP_0(M)$.
\end{prop} 

\begin{proof}
The argument is similar to \cite[Lemma~3.20]{ST11}. In the terminology of \cite{STTraces}, the morphisms in the image of the composition
\beq\label{eq:thickfamily}
\R_{>0}\times \Map(\R^{0|1},M)\hookrightarrow \Mor(\sP_0(M))\to \Mor(1|1\EBord(M))
\eeq
are \emph{thick}: this follows from Proposition~\ref{prop:duality}. The categorical trace of a thick morphism is an endomorphism of the monoidal unit. Proposition~\ref{prop:tracerelations} shows that this categorical trace is the image of the composition
$$
\R_{>0}\times \Map(\R^{0|1},M)\to \mathcal{L}_0^{1|1}(M)\to \Mor(1|1\EBord(M)).
$$
Applying a degree~$n$ field theory to these families, Corollary~\ref{cor:degreenvalues} shows that the categorical trace of the family~\eqref{eq:thickfamily} is sent to the composition 
$$
\cCl_n\otimes_{\cCl_n\otimes \cCl_n} \cCl_n\xrightarrow{\id_{\cCl_n}\otimes R_{t,0}^+}\cCl_n\otimes_{\cCl_n\otimes \cCl_n} (\V_+\otimes \V_-)\xrightarrow{\sigma}\cCl_n\otimes_{\cCl_n\otimes \cCl_n} (\V_-\otimes \V_+) \xrightarrow{L_0} \underline{\C}.
$$
The element $\Gamma_n$ provides an identification $\cCl_n\otimes_{\cCl_n\otimes \cCl_n} \cCl_n\simeq \underline{\C}$ with the trivial line bundle, and the above composition becomes the Clifford super trace. Finally, using that nuclear Fr\'echet spaces have the approximation property~\cite[page~109]{Schaefer}, we find that the categorical trace above is independent of the choice of thickener~\cite[Theorem 1.7]{STTraces}, completing the proof.
\ep

\section{Energy cutoffs and the families index}\label{sec:index}

In this section we prove Theorem~\ref{thm:cocycle} via an index bundle construction applied to the self-adjoint superconnections extracted from field theories in Proposition~\ref{prop:geodata}. The general index bundle construction for a superconnection is reviewed in~\S\ref{sec:KOindex}. When comparing with objects from index theory, we refer to~\S\ref{eq:superlinear} for the translation between the super sign conventions in the previous sections with the standard sign conventions from differential geometry for $\Z/2$-graded inner products, self-adjoint superconnections, and Clifford actions. 

\subsection{Energy cutoffs}\label{sec:cutoffs}

%The failure of existence and uniqueness of solutions to differential equations (see Remark~\ref{rmk:SParinfinite}) necessitates the following definition to ensure that the Chern character of the superconnection extracted from a field theory is compatible with its partition function. 

%Our present goal is to define a subcategory $1|1\EFT^n_\eff(M)\subset 1|1\EFT^n(M)$ on which the the index construction applied to the superconnection from Proposition~\ref{prop:geodata} and differential form $Z(E)(1)\in \Omega^\bullet_\cl(M)$ (from evaluating the partition function~\eqref{eq:constpartitionfun2} at $1\in \R_{>0}$) yields a class in $\dKO^n(M)$. 

\begin{defn}\label{def:hypothesis}
Let $E$ denote a degree~$n$ representation of $\sP_0(M)$, and let $\A_\pm$ denote the superconnections afforded by Lemma~\ref{prop:FTrestrict1}. Then $E$ is \emph{infinitesimally generated} if we have the equality of differential forms 
\beq\label{eq:partitionisChern}
\Tr_{\cCl_n}(\rho_+|_{\R_{>0}\times \Map(\R^{0|1},M)})=\sTr_{\Cl_n}(\exp(-t\A^2))\in \Omega^\bullet(M;C^\infty(\R_{>0})),
\eeq
compare~\eqref{eq:constpartitionfun2}. 
%\beq\label{eq:hypothesis}
%\rho_\pm(t,\theta)=\exp(-t\A_\pm^2\pm \theta \A_\pm). 
%\eeq
An infinitesimally generated representation \emph{admits smooth energy cutoffs} if furthermore the superconnections $\A_\pm$ admits smooth index bundles in the sense of Definition~\ref{defn:superconncutoff}. 
\end{defn}

%{could add a CS form condition to the infinitesimal generation, then this becomes and equivalence, maybe a remark citing Bertram that the CS condition is implied by the trace condition in finite rank.}

\begin{rmk}
We comment on the conditions imposed in Definition~\ref{def:hypothesis}. Although one expects a field theory to be determined by infinitesimal data (implying the condition~\eqref{eq:partitionisChern}), this is not automatic in the formalism~\cite{ST11} due to the failure of existence and uniqueness of solutions to differential equations in topological vector spaces; see Remark~\ref{rmk:SParinfinite}. 
% asserts that the partition function is determined by infinitesimal data. 
The existence of smooth index bundles is more subtle. For a family of Dirac operators, the existence of an index bundle follows from ellipticity, using that eigensections of the Dirac Laplacian are smooth, e.g., see~\cite[Proposition~9.10]{BGV} or~\cite[III, Theorem~5.8]{LM}. As such, the existence of smooth cutoffs can be thought of as an ellipticity property for field theories. In higher dimensions, one can formulate similar smooth cutoff conditions for subspaces~$\mathcal{H}^{\vec\lambda}\subset \mathcal{H}$ indexed by both energies and momenta, with the number of cutoff indices equal to the dimension of the field theory. 
%Condition (4) restricts the growth rates of the energies, where the integral~\eqref{eq:improperintegral} measures the change in the partition function from ``integrating out" the higher energy modes. These types of integrals are common in Wilson's effective approach to quantum field theory with evident generalizations in higher dimensions, e.g., see~\cite[\S1.3]{costbook} for a mathematical introduction. For the Bismut superconnection, the convergence of this integral follows from the growth rate of eigenvalues of the Dirac Laplacian. 
\end{rmk}

Definition~\ref{def:hypothesis} is the final ingredient in the groupoid $1|1\eft^n(M)$, see Definition~\ref{defn:effective}. For simplicity, we will often refer to an object as $E\in 1|1\eft^n(M)$, emphasizing the degree~$n$ representation of $\sP_0(M)$ and suppressing the remaining data. The simplest examples of objects in $1|1\eft^n(M)$ come from superconnections on finite rank bundles.

%With the above preamble in place, we return to the discussion of an index map for field theories. 
%We emphasize that we assume the geometric data previously extracted from a field theory has been translated to follow the standard conventions from index theory rather than the super sign conventions; see~\S\ref{eq:superlinear}. 

% The following uses Definition~\ref{defn:STfieldtheories} for the category $1|1\EFT^n(M)$.
 
%\begin{prop}\label{prop:EFTtoeft}
%Assuming Hypothesis~\ref{hyp:twist} there is a partially-defined functor induced by restriction to $\sP_0(M)$
%$$
%1|1\EFT^n(M)\dashrightarrow 1|1\eft^n(M).
%$$
%\end{prop}
%\bp
%This follows from Propositions~\ref{prop:EFT1},~\ref{prop:EFT2},~\ref{prop:EFT3}, \ref{prop:reality}, \ref{prop:trace}, and the fact that the equality~\eqref{eq:partitionisChern} and existence of cutoffs are properties of a field theory (rather than additional data). 
%\ep

%The subgroupoids $1|1\EFT^n_\eff(M)\subset 1|1\EFT^n(M)$ and $1|1\eft^n_\eff(M)\subset 1|1\eft^n(M)$ are full: for isomorphic objects $E,E'$, we see $E$ admits cutoffs if and only if $E'$ admits cutoffs.

%

The following is the evident degree~$n$ generalization of Proposition~\ref{prop:dataofeft}. 
\begin{prop}\label{prop:dataofeft2}
There is an equivalence of categories
\beq\nonumber
{\rm sPar}\colon \Vect^{\A}_{\Cl_n,{\rm ind}}(M) \xrightarrow{\sim} 1|1\eft^n(M)
\eeq
from the category of real, $\Cl_n$-linear self-adjoint superconnections admitting smooth index bundles to the category~$1|1\eft^n(M)$. 
\end{prop}

\bp
For a $\Cl_n$-linear superconnection $\A$ on a real super vector bundle $V\to M$, consider the representation of $\sP_0(M)$ defined on the vector bundles with sections
$$
\Gamma(\mathcal{V}_+):=\Omega^\bullet(M;V_\C),\qquad \Gamma(\mathcal{V}_-):=\Omega^\bullet(M;\overline{V}_\C)
$$
and super semigroup actions given by 
\beq\label{eq:expof}
e^{-t\A^2+\theta\A}\colon \R^{1|1}_{\ge 0} \to \End(\Omega^\bullet(M;V_\C)),\qquad e^{-t\overline{\A}^2+\theta\overline{\A}}\colon \R^{1|1}_{\le 0} \to \End(\Omega^\bullet(M;\overline{V}_\C)).
\eeq
A pairing for this representation is inherited from the hermitian pairing $\overline{V}_\C\otimes V_\C\to \underline{\C}$ built in the usual way from the underlying metric on~$V$. The real and reflection structures are evident from the construction, as is the reflection positivity condition. The trace structure uses that $\exp(-t\A^2)\in \Omega^\bullet(M;\End(V_\C))$ is trace class (from the existence of smooth index bundles), and $Z_\ell$ is determined by the Chern--Simons form for the 1-parameter family of superconnections gotten from Getzler rescaling of~$\A$. Existence of cutoffs follows by assumption. This constructs the functor ${\rm sPar}$ in~\eqref{eq:dataofeft} on objects. The construction on morphisms is straightforward: using~\ref{eq:expof}, isomorphisms between superconnections are sent to isomorphisms between representations of $\sP_0(M)$. The functor is full and faithful since isomorphisms of representations of $\sP_0(M)$ are determined by isomorphisms of vector bundles. Finally, Proposition~\ref{prop:geodata} verifies that ${\rm sPar}$ is essentially surjective, and hence an equivalence. 
\ep
\begin{cor}\label{prop:florin1}
The groupoid of self-adjoint $\Cl_n$-linear superconnections on finite-rank real vector bundles admits a functor ${\rm sPar}$ 
\beq\label{eq:littlesPar}
&&\begin{tikzpicture}[baseline=(basepoint)];
\node (A) at (0,0) {$\Vect^\A_{\cCl_n,{\rm finite}}(M)$};
\node (B) at (4,0) {$1|1\eft^n(M)$};
\draw[->] (A) to node [above] {\small sPar}  (B);
%\draw[->,dashed] (A) to (C);
%\draw[->,right hook-latex] (C) to (B);
\path (0,-.05) coordinate (basepoint);
\end{tikzpicture} \qquad (V,\A)\mapsto E_{V,\A}
\eeq
determined by super parallel transport along nearly constant superpaths.
\end{cor}
\bp
This follows immediately from Proposition~\ref{prop:dataofeft} by the existence of smooth index bundles in finite dimensions, e.g., see~\cite[\S9.1]{BGV}. 
\ep

%\begin{rmk} 
%To compare with objects from index theory, we switch from the super sign conventions in the previous section to the standard sign conventions from differential geometry for $\Z/2$-graded inner products, self-adjoint superconnections, and Clifford actions; see~\S\ref{eq:superlinear} for this translation. We view the Clifford algebras $\Cl_n$ as real $\Z/2$-graded $*$-algebras. 
%\end{rmk}

For an object $E\in 1|1\eft^n(M)$, let $\mathcal{H}\to M$ denote the fiberwise Hilbert completion of the real vector bundle underlying $E$ relative to the real inner product (see Proposition~\ref{prop:geodata}); $\mathcal{H}$ is typically only a continuous bundle on~$M$. 

\begin{defn}
Given $E\in 1|1\eft^n(M)$ and $\lambda\in \R_{>0}$, the \emph{$\lambda$-cutoff theory} $E_{<\lambda}\in 1|1\eft^n(U_\lambda)$ is defined as the functor~\eqref{eq:littlesPar} applied to the super vector bundle $\mathcal{H}^{<\lambda}\to U_\lambda$ with Clifford linear self-adjoint connection $\nabla^{\mathcal{H}^{<\lambda}}$, i.e., 
$$
E_{<\lambda}:=E_{\mathcal{H}^{<\lambda},\nabla^{\mathcal{H}^{<\lambda}}} \in 1|1\eft^n(U_\lambda)
$$ 
in the notation of Corollary~\ref{prop:florin1}. 
\end{defn}

\begin{lem} 
The partition function of $E$ and the partition function of the $\lambda$-cutoff theory are related by
$$
Z(E)|_{U_\lambda}=Z_{E^{<\lambda}}+d\eta_\lambda\in \Omega^\bullet(U_\lambda;C^\infty(\R_{>0})).
$$
\end{lem}
\bp
This follows from the construction of the $\eta$-forms~\eqref{eq:improperintegral}.
\ep
\subsection{Cutoffs and the RG flow}

The renormalization group flow interpolates between different cutoff theories, as the following lemma makes precise. 

\begin{prop}\label{prop:RGcutoff}\label{lem:RGcutoff}  The existence of cutoffs is preserved under the renormalization group flow: for $E\in 1|1\eft^n(M)$, the $\lambda$-cutoff of~$\RG^\mu(E)$ is canonically isomorphic to the $\mu\lambda$-cutoff of $E$
\beq\label{eq:RGcutoff}
\RG^\mu(E)_{<\lambda}\simeq E_{<\mu\lambda}\in 1|1\eft^n(U_{\mu\lambda}).
\eeq
\end{prop}

\bp
%We first prove a version of the statement replacing the groupoids $1|1\EFT^n(M)$ by reflection positive functors that may not admit cutoffs. By Proposition~\ref{prop:FTrestrict1}, it suffices to calculate the affect of the renormalization group flow on $\sP_0(M)$. 
%%Starting with a degree~$n$ field theories in the upper left corner of the diagram, 2-commutativity of the diagram in Lemma~\ref{lem:equalRG} shows that there is a canonical isomorphism between the corresponding vector bundles in the lower right hand corner. To understand the affect on superconnections we again use  Lemma~\ref{lem:equalRG}. 
%The map $\Mor(\sP_0(M))\to \Mor(\sP_0(M))$ is determined by 
%$$
%(t,\theta,x,\psi)\mapsto (\mu^2t,\mu\theta,x,\mu^{-1}\psi),\qquad (t,\theta)\in (\R^{1|1}_{\ge 0}),\ (x,\psi)\in \Pi TM_\C(S)\simeq \Map(\R^{0|1},M)(S).
%$$
%% in the definition~\eqref{eq:superconnconstrr} of $\A_\pm$. This yields the formula~\eqref{eq:RGsconn}. 
%Taylor expanding the semigroup representation associated to a field theory as in~\eqref{eq:AB} and pulling back along the above action, we find
%\beq
%\RG^\mu_1(\rho_+(t,\theta))&=&\RG^\mu_1(A_+(t)-\theta B_+(t))=\RG^\mu_1(A_+(t))-\mu \theta \cdot \RG^\mu_1(B_+(t))\nonumber\\
%&=&\mu^{-\deg}A_+(t)-\mu \theta (\mu^{-\deg}  B_+(t)).\nonumber
%\eeq
%Since $\A_+=B_+(0)$, this yields $\A_k\mapsto \mu^{1-k}\A_k$ for $\A=\sum_k \A_k$. This is exactly~\eqref{eq:RGsconn}. 
Using Proposition~\ref{prop:RGcutoff}, the description of the $\RG$-action in terms of Getzler rescaling allows us to analyze cutoffs. Unraveling Definition~\ref{defn:effective}, the $\mu\lambda$-cutoff for $E$ is the same data as the $\lambda$-cutoff for $\RG^\mu(E)$: the cutoff bundles are the \emph{same} subspace of $\mathcal{H}\to U_\lambda$ with the \emph{same} connection (note $\A_1$ is invariant under the $\RG$-flow). Hence, there is indeed a canonical isomorphism~\eqref{eq:RGcutoff}, and furthermore $\RG^\mu(E)$ satisfies properties (1) and (2) of Definition~\ref{defn:effective}. Property (3) follows from the Getzler rescaling description of the $\RG$-action, where the partition function of $\RG^\mu(E)$ is computed by the Getzler-rescaled superconnection. Finally, by a change of variables one identifies the integral~\eqref{eq:improperintegral} for~$E$ and $\lambda\mu$ with the integral for $\RG^\mu (E)$ and $\lambda$. This completes the proof. 
%By~\eqref{eq:RGeqn}, this acts on superconnections through pre- and post-composition with the degree endomorphism $\mu^{-\deg}$, with an extra factor of $\mu$ from the action on $\theta$
\ep

\begin{defn}
An object $E\in 1|1\eft^n(M)$ \emph{flows to zero} under the renormalization group flow if for all $\lambda\in \R_{>0}$ and compact submanifolds $K\subset M$, there exists $\mu\in \R_{>0}$ such that $\RG^\mu(E)|_K$ has the zero theory as its $\lambda$-cutoff, i.e., $K\subset U_\lambda$ and $\mathcal{H}^{<\lambda}|_K=\{0\}$. 
\end{defn}

%\begin{rmk} Dropping condition (3) from Definition~\ref{defn:effective} gives a version of effective field theories with an index map valued in $\Cl_n$-bundles (without superconnection). This leads to a cocycle map valued in (non-differential) $\KO^n(M)$. 
%\end{rmk}

%\begin{rmk}
%The terminology \emph{effective} comes from Wilsonian effective field theory, e.g., see~\cite[\S1.3]{costbook} for a mathematical exposition of the key ideas. 
%\end{rmk}

%We recall the incarnations of the renormalization group flow: Definition~\ref{defn:sRG} as an action on $\sP_0(M)$ and 
%We recall Definition~\ref{defn:RG} of the renormalization group flow as an action on field theories induced from an $\R_{>0}$-action on $1|1\EBord(M)$ or $1|1\ebord(M)$. 
%These actions were shown to be compatibility in Lemma~\ref{lem:equalRG}. 

%We also let $\RG^\mu$ denote the precomposition action of this functor on the category of degree~$n$ representations of $\sP(M)$, and call this action the \emph{renormalization group flow}. 
%

%\begin{rmk}
%A closely related option for defining when a field theory flows to zero under the RG flow is to demand that the supersemigroup have zero limit,
%\beq
%\lim_{\mu\to \infty} \rho_+(\mu^2 t,\mu\theta)=0. \label{eq:limitiszero}
%\eeq
%This definition can be quite technical to check in practice. In the geometric examples, the vanishing of the limit~\eqref{eq:limitiszero} is equivalent to invertibility of the degree zero part of the associated superconnection~\cite[Theorem~9.26]{BGV}. The following lemma makes contact with this. 
%\end{rmk}

\begin{lem}\label{lem:RGlimit}  
An object $E\in 1|1\eft^n(M)$ flows to zero under the renormalization group flow if and only if the degree zero component of its associated superconnection $\A_{0}$ is invertible. 
\end{lem}
\bp
It suffices to prove the lemma when $M$ is compact. Let $\nu\in \R_{>0}$ be the minimum eigenvalue of $\A_{0}^2$. If $\nu>0$ (i.e., $\A_{0}$ is invertible), for any $\lambda$ we may choose some $\mu\gg 0$ such that $\mu\nu>\lambda$. By Lemma~\ref{lem:RGcutoff} the $\lambda$-cutoff theory of $\RG^\mu(E)$ is the zero theory. Conversely, if there exists some $\mu$ such that $\RG^\mu(E)$ has the zero bundle as its cutoff bundle, then $\mu\nu$ must be larger than $\lambda$, and hence $\nu\ne 0$ and $\A_0$ is invertible.
\ep 
\subsection{Proof of Theorems~\ref{thm:cocycle} and~\ref{thm:index}}

\begin{defn} For $E\in 1|1\eft^n(M)$, choose a discrete subset $\Lambda\subset \R_{>0}$ such that $\{U_\lambda\}_{\lambda\in \Lambda}$ is a locally finite cover. Then define the differential cocycle 
\beq\label{eq:IndexLambda}
{\rm Index}(E,\Lambda)=\{\mathcal{H}^{<\lambda},\nabla^{\mathcal{H}^{<\lambda}},g_{\lambda\mu},e_{\lambda\mu},\eta_\lambda\}_{\lambda\in \Lambda} \in \dKO^n(M),
\eeq
in the description of differential $\KO$-theory from \cite[\S7]{DBEBis}.
\end{defn}

In brief, the formalism of \cite[\S7]{DBEBis} expands on the ideas of Atiyah--Bott--Shapiro to identify $(\mathcal{H}^{<\lambda},\nabla^{\mathcal{H}^{<\lambda}})$ with the local data of a differential $\KO$-class with Pontraygin form $\sTr_{\Cl_n}(e^{-(\nabla^{\mathcal{H}^{<\lambda}})^2})$. The inclusions~\eqref{eq:compatibilitydata} and odd endomorphisms~\eqref{eq:compatibilitydata2} provide the additional data to refine these local Pontryagin forms into a global \v{C}ech--de~Rham cocycle representing the Pontryagin character of ${\rm Index}(E,\Lambda)$. The forms $\eta_\lambda$ are a \v{C}ech--de~Rham coboundary mediating between this cocycle and the global differential form $Z(E)$.

%The proofs of our main results now follow quickly from the definitions of cutoffs. 

\begin{proof}[Proof of Theorem~\ref{thm:cocycle}]
The index map~\eqref{eq:IndexLambda} applied to the self-adjoint superconnection underlying $E\in 1|1\eft^n(M)$ gives a map of the desired type, but a priori depends on a choice of subset $\Lambda\subset \R_{>0}$ determining a choice of locally finite subcover~$\{U_\lambda\}_{\lambda\in \Lambda}$. To show that the differential KO-class is independent of this choice, one applies the same argument that proves the differential families index is well-defined~\cite[\S8.3]{DBEBis}. To review this argument, for a different choice $\Lambda'\subset \R_{>0}$, we can compare the resulting $\Cl_{-n}$-bundles on the mutual refinement $\{U_\lambda\}_{\lambda\in (\Lambda\bigcup \Lambda')}$. For $\lambda<\lambda'$, we obtain inclusions~\eqref{eq:compatibilitydata} of $\Cl_{-n}$-bundles with commuting $\Cl_{-1}$-actions on the orthogonal complement to the inclusion. This gives the data of a stable equivalence of $\Cl_{-n}$-bundles. By~\cite[Lemma~5.9]{DBEBis}, the $\Cl_{-n}$-bundles are therefore concordant and the differential KO-classes agree. This shows that~\eqref{eq:IndexLambda} is independent of the choice of $\Lambda$, and the index map is well-defined. 

%For a choice of locally finite subcover of $\{U_\lambda\}_{\lambda\in \R_{>0}}$ of $M$ associated with a discrete subset $\Lambda\subset \R_{>0}$, define the cocycle map
%\beq
%1|1\eft^n_\eff(M)\to \dKO^n(M),\label{eq:cocycle map}
%\eeq
%via Definition~\ref{defn:effective}. This determines a differential cocycle by definition. Different choices of cover associated with $\Lambda,\Lambda'\subset \R_{>0}$ give $\Cl_n$-bundles that can be compared on the mutual refinement~$\{U_\lambda\}_{\Lambda\cup \Lambda'}$. By the same argument that shows the differential index of a family of spin manifolds is well-defined, one sees that the $\Cl_n$-bundles with superconnection are concordant 

By Lemma~\ref{lem:RGlimit}, objects in $1|1\eft^n(M)$ that are sent to zero under the RG flow are also sent to the zero class under the composition
$$
1|1\eft^n(M)\to \dKO^n(M)\to \KO^n(M),
$$
since the degree zero part of their associated superconnection is invertible, and hence the underlying KO-class is trivial. This completes the proof. \ep

The following verifies that Theorem~\ref{thm:cocycle} is compatible with Fei Han's construction of the Chern form of a $1|1$-Euclidean field theory built from Dumitrescu's super parallel transport. 
\begin{prop}\label{prop:fei}
The composition of the first two arrows
\beq\nonumber
\Vect^\A_{\cCl_n,{\rm finite}}(M)\stackrel{{\rm sPar}}{\to}1|1\eft^n(M)\stackrel{{\rm Index}}{\longrightarrow} \dKO^n(M)\stackrel{{\rm curv}}{\to}\Omega^\bullet_\cl(M)
\eeq
sends a Clifford linear superconnection to the expected differential $\KO$-class, and hence the above composition sends a $\Cl_n$-linear superconnection $\A$ to the differential form $\sTr_{\Cl_n}(e^{-
\A^2})$ representing its Chern character. 
\end{prop}
\bp
%After possibly choosing an exhaustion function, 
Without loss of generality, we may assume $M$ is compact. Choosing a cutoff valued~$\lambda$ larger than the maximum eigenvalue of $\A_0^2$, the index construction sends a vector bundle with Clifford linear superconnection to the associated differential class in $\dKO^n(M)$, and the curvature of this class is as stated above. This completes the proof. 
\ep

\begin{proof}[Proof of Theorem~\ref{thm:index}]
Using Proposition~\ref{prop:dataofeft}, for a family $\pi\colon X\to M$ of Riemannian spin manifolds and $E\to X$ a real vector bundle the $\Cl_{-n}$-linear superconnection~\eqref{eq:Bismutsuper} acting on the Fr\'echet bundle $V=\pi_*(\bS\times E)\otimes \C$ determines an object of $1|1\eft^{-n}(B)$ via the representation of positively oriented superpaths,
$$
\rho_+(t,\theta)=\exp(-t\A^2+\theta \A).
$$
The pairing and reflection positivity data are inherited from the hermitian pairing on $\pi_*(\bS\times E)\otimes \C$. The real structure comes from complex conjugation on $\pi_*(\bS\times E)\otimes \C$. The standard properties of the Bismut superconnection reviewed in~\S\ref{sec:dKO} imply the existence of a differential index bundle, so we do indeed obtain an object in $1|1\eft^{-n}(M)$. Finally, applying the index map~\eqref{eq:IndexLambda} to this object of $1|1\eft^{-n}(B)$ is precisely the (real) family Clifford index~\cite[Theorem~1.3]{DBEBis}. This completes the proof. \ep

%\begin{defn} Define the subcategory category $1|1\EFT^n_\eff(M)\subset 1|1\EFT^n(M)$ as the (homotopy) pullback in groupoids
%\beq
%\begin{tikzpicture}[baseline=(basepoint)];
%\node (A) at (0,0) {$1|1\EFT_\eff^n(M)$};
%\node (B) at (3,0) {$1|1\eft^n(M)$};
%\node (AA) at (0,-1) {$1|1\EFT^n(M)$};
%\node (BB) at (3,-1) {$1|1\eft^n(M),$};
%\draw[->] (A) to (B);
%\draw[->] (AA) to  (BB);
%\draw[->] (A) to (AA);
%\draw[->] (B) to (BB);
%\path (0,-.75) coordinate (basepoint);
%\end{tikzpicture}\nonumber
%\eeq
%using that pullbacks preserve monomorphisms and therefore $1|1\EFT^n_\eff(M)$ is equivalent to a subcategory of $1|1\EFT^n(M)$. 
%\end{defn}
%

\appendix 

\section{Super hermitian structures, {$*$}-superalgebras, and supermanifolds}\label{eq:superlinear}

%Hermitian pairings, real structures, and $*$-structures in super algebra end up being equivalent to the usual notions from ordinary (ungraded) algebra, with translations coming from judicious factors of $i=\sqrt{-1}$~\cite[\S4]{DM}. We review these translations below, with a focus on the examples of interest for this paper. 

\subsection{Real structures and hermitian pairings}\label{sec:superherm}
The category of \emph{super vector spaces} is the category of $\Z/2$-graded vector spaces over $\C$ equipped with the graded tensor product. The adjective ``super" refers to the sign in the braiding isomorphism for this tensor product
$$
\sigma\colon V\otimes W\stackrel{\sim}{\to} W\otimes V \qquad v\otimes w\mapsto (-1)^{|v||w|} w\otimes v \quad v\in V, \ w\in W
$$
where $v$ and $w$ are homogeneous elements of degree $|v|,|w|\in \Z/2$. 
Let $V^\ev\oplus V^\odd$ denote the direct sum decomposition of $V$ into its even and odd subspaces. The \emph{grading involution} is the linear map $(-1)^{\sf F}\colon V \to V$ that acts by $+1$ on $V^\ev$ and~$-1$ on $V^\odd$. 

%The \emph{grading involution} $(-1)^{\sf F}\colon V\to V$ acts by $+1$ on the even subspace $V^\ev\subset V$ and $-1$ on the odd subspace~$V^\odd\subset V$. 

There is an involution on the category of super vector spaces that on objects reverses complex structures, $V\mapsto \overline{V}$, i.e., $\overline{V}$ is gotten from $V$ by precomposing the action of $\C$ with complex conjugation. 

\begin{defn}
A \emph{real structure} on a super vector space $V$ is an isomorphism $r\colon V\stackrel{\sim}{\to} \overline{V}$ of super vector spaces such that $\overline{r}\circ r=\id_V$. The \emph{real subspace} associated with a real structure is the fixed point set $V_\R\subset V$ for $r$. 
\end{defn}

%For example, complex conjugation gives a real structure $(\overline{\phantom{A}})\colon \C\to \overline{\C}$ on $\C$ with real subspace the standard inclusion $\R\subset \C$. 

\begin{defn}\label{defn:sesquilinear} A \emph{hermitian form} on a super vector space $V$ is the data of a map $\langle-,-\rangle\colon \overline{V}\otimes V\to \C$ making the diagram commute, 
\beq
&&\begin{tikzpicture}[baseline=(basepoint)];
\node (A) at (0,0) {$\overline{V}\otimes V$};
\node (B) at (4,0) {$\C$};
\node (C) at (0,-1.25) {$V\otimes \overline{V}$};
\node (D) at (4,-1.25) {$\overline{\C}$};
\draw[->] (A) to node [above] {$\langle-,-\rangle$} (B);
\draw[->] (A) to node [left] {$\sigma$} (C);
\draw[->] (C) to node [below] {$\overline{\langle-,-\rangle}$} (D);
\draw[->] (B) to node [right] {$(\overline{\phantom{A}})$} (D);
\path (0,-.75) coordinate (basepoint);
\end{tikzpicture}
\eeq
or equivalently, satisfying the formula $\langle x,y\rangle =(-1)^{|x||y|}\overline{\langle y,x\rangle}$ for homogeneous elements. A hermitian form is \emph{positive} if 
\beq
\langle x,x\rangle>0, & x\ne 0 \ {\rm even} \quad {\rm and} \quad i^{-1}\langle x,x\rangle>0, & x\ne 0 \ {\rm odd.}\label{eq:positivitypairing}
\eeq
A \emph{real} hermitian form on a super vector space with real structure is a hermitian form that when restricted to $(V_\R)^\ev\subset V$ takes values in $\R\subset \C$ and when restricted to $(V_\R)^\odd\subset V$ takes values in $i^{-1}\R\subset \C$. 
\end{defn}
%Hermitian form on a real super vector space is \emph{real} if the restriction of the pairing to the real subspace is real on the even part and pure imaginary on the odd part. 

A positive hermitian form in the above sense determines a positive pairing $(-,-)$ in the usual sense on the underlying (ungraded) vector space $V$ by the formula
\beq
(x,y)=\left\{\begin{array}{ll} 0 & |x|\ne |y| \\ \langle x,y\rangle & x,y \ {\rm even} \\  i\langle x,y\rangle & x,y \ {\rm odd} \end{array}\right.\label{eq:ordinarypairing}
\eeq
with the property that the even and odd subspaces are orthogonal. 
When $V$ is endowed with a real structure and a real hermitian form,~\eqref{eq:ordinarypairing} determines a real inner product on~$V_\R$.

\begin{defn}
Let $V$ be a super vector space equipped with a positive hermitian form~$\langle-,-\rangle$. The \emph{super adjoint} of a linear map $T\colon V\to V$ is characterized by
\beq
\langle x,T y\rangle =(-1)^{|x||T|}\langle T^*x,y\rangle. \label{eq:superadjoint}
\eeq
\end{defn}

In the translation to the ordinary pairing~\eqref{eq:ordinarypairing}, we have
\beq
&&T^*=\left\{\begin{array}{ll} T^\dagger & T \ {\rm even} \\ iT^\dagger & T\ {\rm odd}\end{array}\right.\label{eq:superadjointtrans}
\eeq
where $(-)^\dagger$ denotes the usual adjoint with respect to the positive pairing $(-,-)$.

\subsection{$*$-superalgebras}\label{sec:superstar}
A \emph{super algebra} is an algebra object in super vector spaces. Super algebras are the objects of a monoidal category, where the tensor product of super algebras~$A$ and $B$ is the super vector space $A\otimes B$ with multiplication
$$
(a\otimes b)\cdot (a'\otimes b')=(-1)^{|b||a'|} aa'\otimes bb'. 
$$
The \emph{opposite} of a super algebra $A$ is a super algebra $A^\op$ with the same underlying super vector space and multiplication 
\beq
a\cdot_{\op} b:=(-1)^{|a||b|}b\cdot a\label{eq:opposite}
\eeq
where $b\cdot a$ is the multiplication in $A$. A \emph{$*$-superalgebra} is a super algebra together with a homomorphism $(-)^*\colon A\to \overline{A}{}^\op$ that squares to the identity on $A$; we emphasize that an anti-homomorphism in the graded world carries a sign, and hence
$$
(ab)^*=(-1)^{|a||b|}b^*a^*.
$$
A \emph{real structure} on a super algebra is an involutive isomorphism $r\colon A\stackrel{\sim}{\to} \overline{A}$, and $A_\R\subset A$ is the fixed subalgebra. Given a super algebra with real structure, a $*$-structure is \emph{real} if $(-)^*$ preserves the subalgebra $(A_\R)^\ev\subset A$ and sends the subspace $(A_\R)^\odd\subset A$ to $i\cdot (A_\R)^\odd$.

\begin{ex}
Let $V$ be a super vector space with hermitian inner product. Then the super algebra $\End(W)$ has the structure of a $*$-superalgebra determined by the super adjoint~\eqref{eq:superadjoint}. A real structure on $W$ determines a real structure on $\End(W)$ for which this $*$-structure is real. 
\end{ex}

\begin{ex}
The Clifford algebras $\cCl_n$ defined in~\eqref{eq:Clifford} have the structure of $*$-superalgebras using~\eqref{eq:superstarCl}. Complex conjugation gives $\cCl_n$ the structure of a real super algebra for which this $*$-structure is real. 
\end{ex}

A second notion of $*$-structure on a super algebra comes from endowing the underlying (ungraded) algebra with a $*$-structure; we refer to this as a \emph{$\Z/2$-graded $*$-algebra}. This is equivalent to an antilinear involution $(-)^\dagger$ of $A$ satisfying
$$
(xy)^\dagger=y^\dagger x^\dagger.
$$
One can pass back and forth between the two versions of $*$-structure by inserting some appropriate factors of $i$~\cite[page 90]{strings1}. We review this translation in the two examples of interest. 

\begin{ex} \label{ex:superstar}
For a super vector space $V$ with positive hermitian form $\langle-,-\rangle$, the super algebra $\End(V)$ has the structure of a $*$-superalgebra via the super adjoint. It has the structure of a $\Z/2$-graded $*$-algebra via the adjoint with respect to the (ungraded) inner product determined by~\eqref{eq:ordinarypairing}. These two notion of $*$-structure are related by~\eqref{eq:superadjointtrans}. 
\end{ex}
%
%\begin{rmk} Given a super algebra $A$, there are two notions of $*$-structure: (1) a $*$-structure on the underlying (ungraded) algebra of $A$, and (2) a $*$-superalgebra structure on~$A$. In both cases the relevant structure is an involution, $\C$-antilinear, antihomomoprhism, but there are two choices
%$$
%(1): \ (xy)^\dagger=y^\dagger x^\dagger \qquad (2): \ (xy)^*=(-1)^{|x||y|}y^*x^*
%$$ 
%corresponding to whether one takes the opposite of $A$ in the ungraded or graded sense, respectively. 

\begin{ex}\label{ex:Cliffordstar}
In the case of the Clifford algebras $\cCl_n$, we obtain a $\Z/2$-graded $*$-algebra via the $\C$-antilinear extension of the map on generators,
\beq
f_j^\dagger=-f_j,\qquad e_j^\dagger=e_j.\label{eq:starCl}
\eeq
This differs from the $*$-superalgebra structure~\eqref{eq:superstarCl} by a factor of $i$. 
\end{ex}

\subsection{Modules over super $*$-algebras}\label{sec:supermod}
Given a super algebra $A$, let ${}_A\Mod$ and $\Mod_A$ denote the categories of left and right $A$-modules, respectively, whose objects are super vector spaces with a (graded) $A$-action. These categories have symmetric monoidal structures from direct sum of $A$-modules. Similarly, let ${}_A\Mod_B$ denote the category of $A-B$-bimodules. There are canonical equivalence of categories 
\beq
\Mod_{A^\op\otimes B^\op}\simeq {}_{A}\Mod_{B^{\op}} \simeq {}_{A\otimes B}\Mod\label{eq:leftright}
\eeq 
where a left $A\otimes B$-bimodule $W$ is sent to the $A-B^\op$-bimodule with the same underlying super vector space and action
$$
a\cdot w\cdot b:=(-1)^{|w||b|} (a\otimes b) \cdot w
$$
with similar formulas defining a right $A^\op\otimes B^\op$-module. 
%where a $B-A$-bimodule $W$ is sent to the $A\otimes B^\op$-module with the same underlying super vector space and action
%$$
%(a\otimes b) \cdot w:=(-1)^{|m||a|} b\cdot w\cdot a.
%$$
We often use notation like ${}_AW$ for an object of ${}_A\Mod$, and ${}_AW_B$ for an object of ${}_A\Mod_B$. 

\begin{ex}\label{ex:leftright}
A \emph{Clifford module} is a graded module over a Clifford algebra $\cCl_n$. The canonical isomorphism $\cCl_n^\op\simeq \cCl_{-n}$ leads to equivalences of categories 
\beq
{}_{\cCl_n\otimes \cCl_{-m}}\Mod\simeq {}_{\cCl_n}\Mod_{\cCl_m}\simeq \Mod_{\cCl_{-n}\otimes \cCl_m}\label{eq:Cliffmods}
\eeq
which will allow us to translate any Clifford (bi)module into a left Clifford module.
\end{ex}

\begin{defn} A module $V$ over a  $*$-superalgebra $A$ is \emph{self-adjoint} if the structure map $A\to \End(V)$ is a $*$-homomorphism for a specified hermitian form on $V$. \end{defn}

%Super algebras $A$ and $B$ are \emph{Morita equivalent}  if there exist bimodules ${}_AW_B$ and ${}_BV_A$ such that ${}_AW\otimes_B V_A\simeq {}_AA_A$ and ${}_BV\otimes_A W_B\simeq {}_BB_B$, where $A$ and $B$ are regarded as bimodules over themselves. For Morita equivalent algebras, there are equivalences of categories
%$$
%{}_A\Mod\to {}_B\Mod,\qquad W\mapsto {}_BV\otimes_A W
%$$
%induced by the indicated functor, with inverse functor ${}_AW\otimes_B$. 
\begin{rmk} When $A=\cCl_n$ with the $*$-superalgebra structure~\eqref{eq:superstarCl}, the above definition is equivalent to~\cite[Definition~3.3]{BGV} via the translation to $\Z/2$-graded $*$-algebras in Examples~\ref{ex:superstar} and~\ref{ex:Cliffordstar}. 
\end{rmk}

\begin{ex}\label{ex:Morita1}
Define a left $\cCl_{1}\otimes \cCl_{-1}$-module whose underlying super vector space is $\C^{1|1}$ with action determined by
\beq
f\mapsto\left[\begin{array}{cc} 0 & 1 \\ -1 & 0\end{array}\right],\qquad e \mapsto\left[\begin{array}{cc} 0 & 1 \\ 1 & 0\end{array}\right].\label{eq:Cl11bimod}
\eeq
The assignment~\eqref{eq:Cl11bimod} induces an isomorphism of super algebras, $\cCl_{1}\otimes \cCl_{-1}\simeq \End(\C^{1|1})$. In other words, $\cCl_{1}\otimes \cCl_{-1}$ is Morita equivalent to $\C$. 
This $\cCl_{1}\otimes \cCl_{-1}$-action is self-adjoint for the graded hermitian pairing on $\C^{1|1}$ coming from the standard hermitian pairing on $\C^2$ under~\eqref{eq:ordinarypairing}. 
%Using~\eqref{eq:isosofCliff}, we find that $\Cl_n\otimes \Cl_m$ is Morita equivalent to $\Cl_{n+m}$ for any $n,m\in \Z$. From~\eqref{eq:Cliffmods}, we also see that~\eqref{eq:Cl11bimod} determines an invertible $\Cl_1-\Cl_1$-bimodule. The isomorphism $\Cl_{1,1}\simeq \End(\R^{1|1})$ shows that there are two isomorphism classes of invertible $\Cl_1-\Cl_1$-bimodules, corresponding under Morita equivalence to the $\R$-module $\R$ with its two possible gradings. The left $\Cl_{1,1}$-module from~\eqref{eq:Cl11bimod} corresponds to $\Cl_1$ as a bimodule over itself. Finally, the automorphisms of $\Cl_1$ as a bimodule over itself is the group $\Z/2$, acting by $\{\pm 1\}$ on $\R^{1|1}\simeq \Cl_1$. 
\end{ex}

%Opposite algebras are taken in the graded sense, so this is a super $*$-algebra structure, e.g., see~\cite[page 89]{strings1}.  

%Let $\rho\colon A\to \End(W)$ endow $W$ with the structure of a left $A$-module for a super $*$-algebra $A$. Endow the real vector space underlying $W$ with a positive definite inner product such that $W^\ev$ and $W^\odd$ are orthogonal subspaces of $W$. Then $W$ is a \emph{self-adjoint} $A$-module if

%In the following we will restrict attention to the class of \emph{reflection positive} field theories; this feature is present in~\cite{ST04,HST}, but not~\cite{ST11}, so sketch the definition presently. 
\subsection{Super traces}\label{sec:supertrace}
The \emph{super trace} of an endomorphism $T\colon V\to V$ is
\beq
\sTr(T):= {\sf Tr}((-1)^{\sf F}\circ T)\label{eq:supetrace}
% D=\left[\begin{array}{cc} a & b \\ c & d \end{array}\right],\label{eq:supetrace}
\eeq
where ${\sf Tr}$ is the ordinary trace. 
%Throughout, $\R^{n|m}$ will denote the super vector space with~$(\R^{n|m})^\ev=\R^n$ and~$(\R^{n|m})^\odd=\R^m$. 
The \emph{super commutator} of elements in a super algebra is
\beq
[a,b]=ab-(-1)^{|a||b|}ba\qquad a,b\in A.\label{eq:supercommutator}
\eeq
%and it satisfies a graded version of the Jacobi identity. 
A \emph{super trace} on a super algebra is a map of super vector spaces 
\beq
\sTr_A\colon A\to \C,\qquad \sTr_A([a,b])=0\label{eq:superalgsupertrace}
\eeq
that vanishes on super commutators, and hence is determined by a linear map~$A/[A,A]\to \C$. 

A super trace on a super algebra $A$ induces a super trace on the category of (finite-dimensional) $A$-modules. For right $A$-modules, this super trace comes from the composition
\beq
\End_A(W)\stackrel{\sim}{\leftarrow} W\otimes_A W^\vee \to A/[A,A]\stackrel{\sTr_A}{\to}  \C \label{eq:HattoriStallings}
\eeq
where $W^\vee:={\rm Hom}_A(W,A)$ is the dual (left) $A$-module, and the map to $A/[A,A]$ is the Hattori--Stallings trace induced by the evaluation pairing,
$$
{\rm eval}\colon W\otimes_\C W^\vee \to A. 
$$
One obtains traces for bimodules (and left modules) from traces on $A\otimes B^\op$ using~\eqref{eq:leftright}. 

Define the \emph{Clifford super trace} 
\beq
\sTr_{\Cl_n}\colon \End(W)\to \C,\qquad \sTr_{\Cl_n}(T):=\sTr(\Gamma_n \circ T\colon W\to W)\label{eq:Cliffsuper}
\eeq
where $\Gamma_n$ is defined in~\eqref{eq:Gamma}. 
%\beq
%\Gamma_n=2^{-n/2}e_1\cdot e_2\cdots e_n\in \cCl_n.\label{eq:Gamma}
%\eeq
%e.g., see~\cite[Definition 3.2.16]{ST04}. 

%\begin{rmk}
%Without fixing the choice~\eqref{eq:Gamma} for $\Gamma$, the Clifford super trace takes values in the determinant line $\det(\C^n)=\Lambda^{\rm top}\C^n\simeq \cCl_n/[\cCl_n,\cCl_n]$~\cite[Proposition 3.2.1]{BGV}. This line is even or odd line, depending on the parity of $n$. The choice of trivialization $\Gamma$ trivializes this line.
%\end{rmk}

\subsection{(Real) supermanifolds}

Let $\SMfld$ denote the category of supermanifolds with structure sheaves defined over $\C$; we refer to~\cite[\S4.1]{ST11} for a brief introduction and~\cite{DM} for a more comprehensive one. We require a slight generalization of the usual definition that allows for supermanifolds with boundary. Analogously to the usual case, these are defined as Hausdorff, second countable locally ringed spaces that are locally isomorphic to $\HH^{n|m}$, where $\HH^{n|m}$ is the restriction of the structure sheaf of $\R^{n|m}$ to the half-space $\HH^n:=\{(x_1,\dots, x_n)\in \R^n\mid x_1\ge 0\}.$

Throughout, $C^\infty$ will denote the structure sheaf of a supermanifold $S$, which is a sheaf of $\Z/2$-graded complex algebras. Global sections are denoted $C^\infty(S)$. We refer to sections as \emph{functions} on the supermanifold. Let $S_{\red}$ denote the reduced manifold of~$S$, i.e., the locally ringed space with the same underlying topological space and structure sheaf $C^\infty/{\rm nil}$ where ${\rm nil}$ is the sheaf of nilpotent functions, or equivalently, the sheaf of ideals generated by odd sections of $C^\infty$. 

\begin{ex}\label{ex:Roneonege}
%The two main examples of supermanifolds with (nonempty) boundary considered in this paper are $\R^{1|1}_{\ge 0}$ and $\R^{1|1}_{\le 0}$. 
Let $\R^{1|1}$ denote the supermanifold whose underlying topological space is $\R$, and structure sheaf is the sheaf of sections of $\Lambda^\bullet \underline{\C}^\vee$, where~$\underline{\C}$ is the trivial complex line bundle on $\R$. Then define $\R^{1|1}_{\ge 0}$ and $\R^{1|1}_{\le 0}$ by restricting the structure sheaf of $\R^{1|1}$ to $\R_{\ge 0}\subset \R$ and $\R_{\le 0}\subset\R$. 
\end{ex}

\begin{ex}
For a complex vector bundle $E\to N$ over an ordinary manifold $N$, let $\Pi E$ denote the supermanifold whose structure sheaf of $\Z/2$-graded algebras is given by 
$$
U\mapsto C^\infty(U):=\Gamma(U;\Lambda^\bullet E^\vee),\qquad U\subset N. 
$$
Every supermanifold $S$ is isomorphic to $\Pi E$ for some vector bundle $E\to N=S_{\rm red}$~\cite{Batchelor}.
\end{ex}

In computations, we frequently identify a supermanifold with its functor of points, i.e., its image under the Yoneda embedding into presheaves (of sets) on $\SMfld$. We use the notation $N(S)$ to denote the $S$-points of a supermanifold~$N$. 

\begin{ex}
The supermanifold $\R^{n|m}$ is the locally ringed space with underlying topological space $(\R^{n|m})_{\rm red}=\R^n$ and structure sheaf the sheaf of sections of $\Lambda^\bullet (\underline{\C}^m)^\vee$, the trivial rank $m$ complex vector bundle on $\R^n$. The $S$-points of $\R^{n|m}$ are
$$
\R^{n|m}(S)\simeq \{t_1,\dots, t_n\in C^\infty(S)^\ev, \ \theta_1,\dots, \theta_m \in C^\infty(S)^\odd \mid (f_i)_{\rm red} =\overline{(f_i)}_{\rm red}\in C^\infty(S_{\rm red})\}
$$
where the reality condition on the $f_i$ is imposed on the reduced manifold of $S$. 
\end{ex}

Many definitions in supermanifolds can be phrased in terms of the functor of points. For example, a \emph{super Lie group} is a group object in supermanifolds. 
Sometimes supermanifolds are first defined in terms of their functor of points and then proven to be representable. 

\begin{ex}
Define $\Map(\R^{0|1},M)$ as the presheaf on supermanifolds whose value on $S$ is the set of maps $S\times \R^{0|1}\to M$; this supermanifold is representable~$\Map(\R^{0|1},M)\simeq \Pi TM_\C$, e.g., see~\cite[Proposition~3.1]{HKST}. Consider the precomposition action of $\R^{0|1}$ on itself by translation,
$$
\R^{0|1}\times \Map(\R^{0|1},M)\to \Map(\R^{0|1},M). 
$$
Differentiation at $0\in \R^{0|1}$ gives an odd vector field on $\Map(\R^{0|1},M)$. Under the isomorphism $C^\infty(\Map(\R^{0|1},M))\simeq C^\infty(\Pi TM_\C)\simeq \Omega^\bullet(M)$, this vector field is the de~Rham differential, acting on forms as an odd derivation. 
\end{ex}

%The functor of points description can reveal super geometry that is perhaps less apparent in the standard point of view. 
Following the terminology from \cite{HST,ST11}, a \emph{generalized supermanifold} is a presheaf on the site of supermanifolds (that typically is not representable). A map between generalized supermanifolds is a map between presheaves. 

\begin{ex} \label{ex:mappingsheaf}
Given supermanifolds $\Sigma$ and $X$, define the mapping object 
$$
\Map(\Sigma,X)\colon \SMfld^\op\to {\sf Set},\qquad \Map(\Sigma,X)(S):=\SMfld(S\times \Sigma,X)
$$
as a generalized supermanifold with the indicated $S$-points. 
\end{ex}

\begin{ex}[{\cite[Example~5.5]{HST}}]\label{ex:tvspresheaf}
Given a $\Z/2$-graded topological vector space $\V$, define a generalized supermanifold by the $S$-points 
$$
\V\colon \SMfld^\op\to {\sf Set},\qquad  \V(S):=(C^\infty(S)\otimes \V)^\ev
$$ 
where the tensor product is the $\Z/2$-graded projective tensor product. 
\end{ex}

\begin{defn}\label{defn:realsmfld}
There is an involution on the category of supermanifolds, 
\beq\label{eq:conjugationfunctor}
(\overline{\phantom{-}})\colon \SMfld\to \SMfld,\qquad S\mapsto \overline{S}
\eeq
sending a supermanifold to the same locally ringed space but whose sheaf has complex numbers acting through complex conjugation. A \emph{real structure} on a supermanifold $S$ is an isomorphism $r_S\colon S\stackrel{\sim}{\to} \bar S$ such that $\overline{r}\circ r=\id_S$. A map $f\colon S'\to S$ between supermanifolds with real structure is \emph{real} if $r_S\circ f=f\circ r_{S'}$. 
\end{defn}

\begin{ex} \label{ex:real}
For a complex vector bundle $E\to N$, a real structure on $E$ (as a vector bundle) determines a real structure on $\Pi E$ (as a supermanifold) 
$$
C^\infty(\Pi E)\simeq \Gamma(N;\Lambda^\bullet E^\vee)\simeq \Gamma(N;(\Lambda E_\R^\vee)\otimes \C) \qquad (E_\R)\otimes \C\simeq E. 
$$
Indeed, complex conjugation determines a $\C$-antilinear map $C^\infty(\Pi E)\to C^\infty(\Pi E)$ which provides an isomorphism of supermanifolds $\Pi E\stackrel{\sim}{\to} \overline{\Pi E}$. 
%In the case of $\Pi TM_\C \simeq \Map(\R^{0|1},M)$, this real structure is complex conjugation of $\C$-valued differential forms on $M$, viewed as a $\C$-antilinear map and hence an isomorphism $\Map(\R^{0|1},M)\stackrel{\sim}{\to} \overline{\Map(\R^{0|1},M)}$. 
\end{ex}

Given a vector bundle $V\to M$, we obtain a vector bundle $p^*V\to \Map(\R^{0|1},M)$ by pulling back along the projection $p\colon \Map(\R^{0|1},M)\simeq \Pi TM\to M$. 
%A superconnection~\eqref{eq:ordsuperconn} on $V$ can therefore be reinterpreted as a (differential) operator on sections of $p^*V\to \Map(\R^{0|1},M)$. 
%Furthermore, any super vector bundle over $\Map(\R^{0|1},M)$ can be given such a description. 

%Given a vector bundle $W\to \Map(\R^{0|1},M)$, we can restrict along the inclusion of the reduced manifold $i\colon M\hookrightarrow \Map(\R^{0|1},M)$ to obtain a vector bundle $i^*W\to M$. Conversely, given a vector bundle $V\to M$ we can pullback along the projection $p\colon \Map(\R^{0|1},M)\to M$ to obtain a vector bundle $p^*V\to \Map(\R^{0|1},M)$. 

\begin{lem} \label{lem:sVect}For any vector bundle $\V\to \Map(\R^{0|1},M)$, there is an isomorphism $\V \simeq p^*V$ where $V\to M$ is a vector bundle on $M$.
\end{lem}

\bp
Let $i\colon M\hookrightarrow \Map(\R^{0|1},M)$ denote the inclusion of the reduced manifold. Given a vector bundle $\V\to \Map(\R^{0|1},M)$, our candidate for $V\to M$ in the statement is $V:=i^*\V$. Indeed, the composition $i\circ p$ is homotopic to the identity on $\Map(\R^{0|1},M)$. For example, consider the homotopy $\R\times \Pi TM_\C\to \Pi TM_\C$ that rescales the odd fibers of $\Pi TM_\C\to M$ by a scalar $t\in \R$; in the functor of points this map is $(t,x,\psi)\mapsto (x,t\psi)$ for $t\in \R(S)$ and $(x,\psi)\in \Pi TM_\C(S)\simeq \Map(\R^{0|1},M)(S)$. Now we use that pullbacks along homotopic maps yield isomorphic vector bundles,
%(the standard partition of unity argument in the category of manifolds proceeds identically in the category of supermanifolds). 
concluding that $p^*V=p^*i^*\V\simeq (i\circ p)^*\V$ is isomorphic to~$\V$. 
\ep

\section{Review of Stolz and Teichner's Euclidean field theories}\label{sec:ebord}

In this section we review Stolz and Teichner's definition of twisted geometric field theory following~\cite{ST11}. Their central construction is a \emph{geometric bordism category}. This takes as input a \emph{rigid geometry}~$(G,\M)$ and outputs a functor from smooth manifolds to categories internal to symmetric monoidal stacks on the site of supermanifolds, 
$$
(G,\M)\Bord(-)\colon \Mfld\to \left\{\begin{array}{c} {\rm Categories\ internal\ to} \\ {\rm symmetric\ monoidal\ stacks}\end{array}\right\},\qquad  M\mapsto (G,\M)\Bord(M),
$$
where $(G,\M)\Bord(M)$ is the \emph{category of $(G,\M)$-bordisms over $M$}. Geometric field theories are then defined as functors out of $(G,\M)\Bord(M)$, and twisted field theories are certain natural transformations between functors out of $(G,\M)\Bord(M)$.
% as in~\eqref{eq:twistedEFT}. 

We also introduce \emph{reflection positivity data} within this formalism (see~\S\ref{sec:additionalstructuresGM}), following the motivational discussion from \S\ref{sec:SUSYQMmot}. Although this structure is not part of Stolz and Teichner's definition from~\cite{ST11}, it is implicit in their previous work~\cite[Definition~2.1.1, Remark 2.1.2]{ST04} and seems to be 
% and is related to the definition of \emph{positivity} that is crucial to the connection with the KO-spectrum in~\cite[\S6]{HST}, see also~\cite[Definition~1.1.32 part (iii)]{Cheung}.
%It is entirely possible that a version of the $d=1$ conjecture~\eqref{eq:conjecture} holds without reflection positivity, and so in some sense this condition may not be necessary. However, in dimension~$1|1$ reflection positive theories connect better with classical objects from index theory, and all the examples from physics have this structure. Furthermore, this structure seems necessary for the $d=2$ conjecture in~\eqref{eq:conjecture}; 
crucial for the $d=2$ case of conjecture~\eqref{eq:conjecture}; see Remark~\ref{rmk:KTate}.

%So for the time being, it seems reasonable to include this additional structure when investigating the conjectures~\eqref{eq:conjecture}. 

\begin{rmk} \label{rmk:extenddown}
The $d=2$ case of conjecture~\eqref{eq:conjecture} requires a generalization of the definitions below to fully-extended $2|1$-dimensional Euclidean field theories of degree~$n$. In particular, this requires that one enhance the $2|1$-Euclidean bordism category to (some flavor of) bicategory that encodes cutting and gluing of $2|1$-Euclidean manifolds with corners. The framework for such a definition is a subject of current investigation, and likely to depart in various ways from prior approaches. Indeed, extending  the definitions in~\cite{ST11} to bicategories internal to symmetric monoidal stacks seems rather unwieldy. At present, smooth bordism categories inspired by iterated Segal spaces seem more tractable, e.g., see~\cite{ClaudiaTheo,GradyPavlov,LudewigStoffel,BEPTFT}. 
\end{rmk}

\subsection{Rigid geometries} \label{sec:rigidgeometry}
Inspired by Thurston~\cite[Ch.~3.3]{Thurston}, a \emph{rigid geometry} is a triple $(G,\M,\M^c)$ where $\M$ is a supermanifold, $G$ is a super Lie group acting on $\M$, and $\M^c\subset \M$ is a codimension~1 submanifold. The supermanifold $\M$ is the \emph{local model} for the geometry, the group $G$ is the \emph{isometry group}, and $\M^c\subset \M$ is a local model for (collared) boundaries. The data of a rigid geometry is often abbreviated to the pair $(G,\M)$.

\begin{defn}[Sketch of {\cite[Definitions~2.33 and~4.4]{ST11}}]\label{defn:rigid} A \emph{family of $(G,\M)$-supermanifolds} is data
    \begin{enumerate}
        \item a smooth submersion $p\colon Y\to S$;
        \item a maximal atlas $\{U_i\}$ of $Y$ equipped with open embeddings $\varphi_i \colon U_i \hookrightarrow S \times \M$ over~$S$;
        \item transition data $g_{ij}\colon p(U_i\bigcap U_j)\to G$. 
    \end{enumerate}
The embeddings $\varphi_i$ are required to be compatible with transition data $g_{ij}$ (using the $G$-action on $\M$) and the transition data must satisfy a further cocycle condition. 

An \emph{isometry} between $(G,\M)$-families is a pullback $f\colon Y\to Y'$ over a base change $S\to S'$ with the property that the restriction of $f$ to the atlases determining the $(G,\M)$-structures is given by the action of~$G$ on open subspaces of $\M$. The collection of $(G,\M)$-supermanifolds and isometries defines a stack $(G,\M)$-$\SMfld$ on the site of supermanifolds; see~\cite[\S2.8]{ST11} and~\cite[Remark A.23]{Powerops}.

A \emph{family of $(G,\M)$-pairs} is a $(G,\M)$-family $Y\to S$ with a codimension~1 submanifold $Y^c\subset Y$ over $S$ with the property that in a chart for the $(G,\M)$-structure, the image of $U_i\bigcap Y^c\subset U_i$ under  $\varphi_i$ lands in the subspace $S\times \M^c\subset S\times \M$. Isometries between $(G,\M)$-pairs are isometries of $(G,\M)$-supermanifolds $Y\to Y'$ such that $Y^c \subset Y$ is sent to $Y'^c\subset Y'$. This determines a stack on the site of supermanifolds whose objects are $(G,\M)$-pairs. 
\end{defn}

We also refer to~\cite[\S6.3]{HST} for additional discussion of rigid geometries. 

%\begin{rmk}
%The definition of a $(G,\M)$-pair in~\cite[Definition 2.33]{ST11} disallows certain bordisms asserted to exist in~\cite[\S3.2]{ST11}. The above definition differs from \cite{ST11} to resolve this issue. In brief, in \cite[\S3.2]{ST11} the source of the ``left cylinder" is a super Euclidean pair where the core $Y^c\subset Y$ is not a codimension~1 sub-supermanifold, but rather a coproduct of maps $Y^c\coprod Y^c\to Y$ for the same sub-supermanifold $Y^c\subset Y$. This bordism leads to a pairing on the vector space of the associated field theory, and hence it is important to allow for these more general $(G,\M)$-pairs as defined above (compare~\S\ref{sec:superpathbordism}  above). 
%%It is possible that one should require a further condition on the images $Y^c_\alpha\subset Y$, but in the examples we have considered no such conditions are needed. 
%\end{rmk}

%The objects and morphisms above determine stacks $(G,\M)$-$\SMfld$ on the site of supermanifolds whose objects are families of $(G,\M)$-supermanifolds. Similarly, there is a stack of \emph{$(G,\M)$-pairs} whose objects are $S$-families of supermanifolds $Y\to S$ with a $(G,\M)$-structure together with a codimension~1 $S$-family $Y^c\to Y$ that is locally of the form $\M^c\subset \M$. 

\begin{ex}[Euclidean rigid geometry]\label{ex:modelgeom}
Let $\A^d$ be affine $d$-space, regarded as a Riemannian spin manifold. Its group of spin structure preserving isometries is $\R^d\rtimes \Spin(d)$, where $\Spin(d)$ acts on $\R^d$ through the double cover $\Spin(d)\to \SO(d)$. The \emph{$d$-dimensional Euclidean geometry} is defined by the triple $(G,\M,\M^c)=(\R^d\rtimes \Spin(d),\A^d,\A^{d-1})$, where $\A^{d-1}\subset \A^d$ comes from the standard inclusion of a totally geodesic hyperplane; after a choice of basepoint this is the standard inclusion $\R^{n-1}\subset \R^n$. The stack $(\R^d\rtimes \Spin(d),\A^d)$-$\SMfld$ of $d$-dimensional Euclidean manifolds has as objects fiber bundles $Y\to S$ of $d$-manifolds endowed with flat Riemannian metrics and spin structures. The stack of $d$-dimensional Euclidean pairs has as objects fiber bundles $Y\to S$ of $d$-manifolds endowed with flat Riemannian metrics and spin structures together with a family of totally geodesic codimension~1 submanifolds $Y^c\subset Y$ over~$S$. 
\end{ex}
\begin{lem}\label{lem:GMaddstructure1}
Fix a rigid geometry $(G,\M)$. Given a homomorphism $g\colon G\to G$ and a $G$-equivariant map $f\colon \M\to \M$ relative to $g$ that preserves the image of the inclusion $\M^c\subset \M$, there is a uniquely determined map of stacks
$$
F\colon (G,\M)\hbox{-}\SMfld \to (G,\M)\hbox{-}\SMfld.
$$
\end{lem}
\begin{proof}[Sketch of Proof.]
The proof is the same argument as~\cite[Lemma 6.19]{HST}. Suppose we are given an $S$-family of $(G,\M)$-manifolds $Y\to S$. We obtain a new $(G,\M)$-family $f(Y\to S)$ given by the same fibration $Y\to S$, but where local charts $U_i\subset Y$ are modified by postcomposition with $f$,
$$
Y\supset U\stackrel{\varphi_i}{\hookrightarrow} \M\times S\stackrel{f\times \id_S}{\to} \M\times S. 
$$
Because $f$ is $G$-equivariant relative to $g$, the gluing data $g^{-1}g_{ij}g$ determines a $(G,\M)$-structure on $Y$. Similarly, because $f$ preserves the image of $\M^c\subset \M$, it determines an endofunctor on the category of $(G,\M)$-pairs. Since the underlying map on the family $Y\to S$ is the identity map, these endofunctors of $(G,\M)$-manifolds determine maps of stacks.
\ep

\begin{defn}[{\cite[Definition 6.18]{HST}}]\label{defn:realFT}
%The image of a rigid geometry $(G,\M)$ under the conjugation functor (see ) determines a (possibly distinct) rigid geometry $(\overline{\M},\overline{G})$. 
A \emph{real structure} on a rigid geometry $(G,\M)$ is a real structure on $\M^c$, $\M$ and $G$ such that all the maps defining the rigid geometry are real in the sense of Definition~\ref{defn:realsmfld}. 
\end{defn}

%\begin{rmk}
%The Euclidean geometry defined above leads to the definition of Euclidean field theory in Stolz and Teichner's terminology; this is in keeping with Thurston's notion of Euclidean geometry. However, we caution that the terminology conflicts with a different common usage. In physics, \emph{Euclidean field theories} are defined on spacetimes with Riemannian (as opposed to Lorentzian) metrics, but---in contrast to Stolz and Teichner's definition---there is typically no condition that the metric be flat. Of course, for 1-dimensional field theories (the focus of this paper) 
%\end{rmk}
%\begin{ex}\label{ex:modelgeo}
%In this paper, the main example of interest are the 1-dimensional and $1|1$-Euclidean geometries, which are determined by the triples 
%$$
%(G,\M,\M^c)=(\R,\A,\pt)\quad {\rm respectively} \quad (\R^{1|1}\rtimes \Z/2,\A^{1|1},\A^{0|1}).
%$$ 
%The stack of $1$-dimensional Euclidean manifolds has as objects fiber bundles $Y\to S$ whose fibers are endowed with flat Riemannian metrics. The stack of $1|1$-Euclidean manifolds has as object fiber bundles $Y\to S$ whose fibers are endowed with a super Euclidean structure; see~\cite[\S4.2]{ST11}. 
%\end{ex}

\subsection{Categories internal to stacks} \label{sec:internalcategories}

\begin{defn}[Sketch of~{\cite[Definition~2.13]{ST11}}]\label{defn:internalcat}
A \emph{category ${\sf C}$ internal to symmetric monoidal stacks} is the data of symmetric monoidal stacks $\Ob({\sf C})$ and $\Mor({\sf C})$ of objects and morphisms with source, target, unit, and composition functors
$$
\s,\t\colon \Mor({\sf C})\to \Ob({\sf C})\quad \u\colon \Ob({\sf C})\to \Mor({\sf C}),\quad {\sf m}\colon \Mor({\sf C})\times_{\Ob({\sf C})} \Mor({\sf C})\to \Mor({\sf C})
$$
along with additional coherence data (e.g., associator and unitor) satisfying compatibility properties.
\end{defn} 

\begin{defn}[Sketch of {\cite[Definition~2.18]{ST11}}]\label{defn:internalfunctor}
For ${\sf C}$, ${\sf D}$ categories internal to symmetric monoidal stacks, an \emph{internal functor} $F=(F_0,F_1,\mu,\varepsilon)\colon {\sf C}\to {\sf D}$ is defined as morphisms of stacks $F_0\colon \Ob({\sf C})\to \Ob({\sf D})$, $F_1\colon \Mor({\sf C})\to \Mor({\sf D})$, and 2-commuting data, 
\beq
&&\begin{tikzpicture}[baseline=(basepoint)];
\node (A) at (0,0) {$\Mor({\sf C})\times_{\Ob({\sf C})} \Mor({\sf C})$};
\node (B) at (4,0) {$\Mor({\sf C})$};
\node (C) at (0,-1.5) {$\Mor({\sf D})\times_{\Ob({\sf D})}\Mor({\sf D})$};
\node (D) at (4,-1.5) {$\Mor({\sf D})$};
\node (E) at (2,-.75) {$\mu \ \twocommute$};
\draw[->] (A) to node [above] {${\sf m}$} (B);
\draw[->] (A) to node [left] {$F_1\times F_1 $} (C);
\draw[->] (C) to node [below] {${\sf m}$} (D);
\draw[->] (B) to node [right] {$F_1$} (D);
\path (0,-.75) coordinate (basepoint);
\end{tikzpicture}
\qquad \begin{tikzpicture}[baseline=(basepoint)];
\node (A) at (0,0) {$\Ob({\sf C})$};
\node (B) at (4,0) {$\Ob({\sf D})$};
\node (C) at (0,-1.5) {$\Mor({\sf C})$};
\node (D) at (4,-1.5) {$\Mor({\sf D}).$};
\node (E) at (2,-.75) {$\varepsilon \ \twocommute$};
\draw[->] (A) to node [above] {$F_0$} (B);
\draw[->] (A) to node [left] {$\u$} (C);
\draw[->] (C) to node [below] {$F_1$} (D);
\draw[->] (B) to node [right] {$\u$} (D);
\path (0,-.75) coordinate (basepoint);
\end{tikzpicture}\label{eq:internalfunctordiag}
\eeq
These data are required to satisfy further properties involving the associators and unitors. 
\end{defn}

\begin{defn}[Sketch of {\cite[Definition~2.19]{ST11}}]\label{rmk:twistedcategory}
An \emph{internal natural transformation} $(\eta,\rho)\colon F\Rightarrow G$ between internal functors $F,G\colon {\sf C}\to {\sf D}$ is a morphism of stacks $\eta\colon \Ob({\sf C})\to \Mor({\sf D})$ and 2-commuting data $\rho$
\beq
&&\begin{tikzpicture}[baseline=(basepoint)];
\node (A) at (0,0) {$\Mor({\sf C})$};
\node (B) at (6,0) {$\Mor({\sf D})\times_{\Ob({\sf D})}\Mor({\sf D})$};
\node (C) at (0,-1.5) {$\Mor({\sf D})\times_{\Ob({\sf D})}\Mor({\sf D})$};
\node (D) at (6,-1.5) {$\Mor({\sf D}).$};
\node (E) at (3,-.75) {$\rho \ \twocommute$};
\draw[->] (A) to node [above] {$\eta\circ \t \times F_1$} (B);
\draw[->] (A) to node [left] {$G_1\times \eta\circ \s $} (C);
\draw[->] (C) to node [below] {${\sf m}$} (D);
\draw[->] (B) to node [right] {${\sf m}$} (D);
\path (0,-.75) coordinate (basepoint);
\end{tikzpicture}\label{diag:nu}
\eeq
These data are required to satisfy additional coherence properties. Define an \emph{isomorphism between internal natural transformations} $(\eta,\rho)\Rightarrow (\eta',\rho')$ as the the data of 2-morphism $\eta\Rightarrow \eta'$ between maps of stacks with respect to which~$\rho$ and~$\rho'$ are compatible. 
%Compatibility uses the diagram~\eqref{diag:nu} and the isomorphism $\eta\Rightarrow \eta'$, and require a pair of 2-morphisms between maps of stacks to be equal.
\end{defn}

\begin{defn}\label{defn:superLiecat}
A \emph{super Lie category} is a category internal to supermanifolds whose source and target maps are submersions. For a super Lie category ${\sf C}$, let $\Ob({\sf C})$ and $\Mor({\sf C})$ denote the supermanifolds of objects and morphisms, and let
$$
\u\colon \Ob({\sf C})\to \Mor({\sf C}), \ \s,\t,\colon\Mor({\sf C})\to \Ob({\sf C}),\ \m\colon\Mor({\sf C})\times_{\Ob({\sf C})}\Mor({\sf C})\to \Mor({\sf C})
$$
denote the unit, source, target, and composition maps, respectively. A \emph{smooth functor} $F\colon {\sf C}\to {\sf D}$ is a pair of maps of supermanifolds $F_0\colon \Ob({\sf C})\to \Ob({\sf D})$, $F_1\colon \Mor({\sf C})\to \Mor({\sf D})$ satisfying the axioms for a functor between categories. A \emph{smooth natural transformation} $\eta\colon F\Rightarrow G$ is a map of supermanifolds $\eta\colon \Ob({\sf C})\to \Mor({\sf D})$ satisfying the axioms for a natural transformation between functors. Super Lie categories, smooth functors, and smooth natural transformations form a strict 2-category. 
\end{defn}

By regarding a supermanifold as a (representable) stack, any super Lie category determines a category internal to stacks (see Definition~\ref{defn:internalcat}). By the Yoneda lemma for stacks, functors and natural transformations between super Lie categories are equivalent to functors and natural transformations between their associated categories internal to stacks. This also allows us to consider functors from a super Lie category to a category internal to stacks. 

\begin{ex} Any supermanifold determines a discrete super Lie category, i.e., a super Lie category with only identity morphisms. \end{ex}

\begin{rmk}\label{rmk:badcat}
Categories internal to stacks, internal functors, and internal natural transformations are not as well-behaved as one might hope. For example, an internal functor ${\sf C}\to {\sf D}$ that is (appropriately) essentially surjective and fully faithful need not have an inverse functor. For example, let ${\sf C}$ be the \v{C}ech groupoid of an open cover of $X$ and ${\sf D}$ be the (discrete) internal category determined by~$X$; the evident internal functor ${\sf C}\to {\sf D}$ has an inverse internal functor if and only if the cover admits a section. The ``correct" framework for categories internal to stacks therefore requires a localization at an appropriate class of weak equivalences, analogous to the localization relating the 2-category of Lie groupoids with geometric stacks~\cite{Pronk}. For the manipulations in this paper, it suffices to work with categories internal to stacks before such a localization. 
\end{rmk}

%\begin{ex}\label{ex:supersemigroup} The super semigroup $\R^{1|1}_{\ge 0}$ from~\S\ref{sec:SUSYQM} determines a super Lie category ${\sf C}$ with $\Ob({\sf C})=\pt$ and $\Mor({\sf C})=\R^{1|1}_{\ge 0}$. The unit is the inclusion of $0\in \R^{1|1}_{\ge 0}$ and composition is determined by multiplication in $\R^{1|1}_{\ge 0}$. 
%\end{ex}

The following notions of internal involution and equivariance data expand on~\cite[Definition 6.47]{HST} and \cite[2.6]{ST11}. The ambient involution $\SMfld\to \SMfld$ is taken in the sense of~\cite[Definition~B.1]{FreedHopkins}, and in our applications will either be the identity involution or the conjugation functor~\eqref{eq:conjugationfunctor} on supermanifolds.

\begin{defn}\label{defn:interinvol}
Fix an involution $\SMfld\to \SMfld$ on the category of supermanifolds. For a category ${\sf C}$ internal to (symmetric monoidal) stacks on the site of supermanifolds, an \emph{internal involution} is the data of an internal functor $\tau\colon {\sf C} \to {\sf C}$ covering the fixed involution of supermanifolds, an internal natural isomorphism $\eta \colon \id_{\sf C}\Rightarrow \tau\circ\tau$, and an isomorphism $\delta$ between internal natural isomorphisms $\tau\Rightarrow \tau$ in the diagram
\beq
\begin{tikzpicture}[baseline=(basepoint)];
\node (A) at (0,0) {${\sf C}$};
\node (B) at (3,0) {${\sf C}$};
\node (C) at (6,0) {${\sf C}$};
\node (D) at (9,0) {${\sf C}$};
\node (E) at (3,.5) {$\eta \ \Downarrow$};
\node (E) at (6,-.5) {$\eta \ \Uparrow$};
\draw[->] (A) to node [above] {$\tau$} (B);
\draw[->,bend left] (A) to node [above] {$\id_{\sf C}$} (C);
\draw[->,bend right] (B) to node [below] {$\id_{\sf C}$} (D);
\draw[->] (B) to node [above] {$\tau$} (C);
\draw[->] (C) to node [above] {$\tau$} (D);
\path (0,-.75) coordinate (basepoint);
\end{tikzpicture}\nonumber
\eeq
The isomorphism $\delta$ between internal natural transformations satisfies a compatibility condition with $\eta$ for the quadruple composition $\tau\circ\tau\circ\tau\circ\tau$. 
\end{defn}

The following generalizes \cite[Definition~B.6]{FreedHopkins} to Stolz and Teichner's internal categories.

\begin{defn} \label{defn:equivariancedata}
Let $({\sf C},\tau_{\sf C},\eta_{\sf C},\delta_{\sf C})$ and $({\sf D},\tau_{\sf D},\eta_{\sf D},\delta_{\sf D})$ be internal categories with internal involution, and let $F\colon {\sf C}\to {\sf D}$ be a functor. \emph{Equivariance data} for $F$ is an internal natural transformation~$\phi\colon F\circ \tau_{\sf C}\Rightarrow \tau_{\sf D}\circ F$, i.e., 2-commuting data
\beq
\begin{tikzpicture}[baseline=(basepoint)];
\node (A) at (0,0) {${\sf C}$};
\node (B) at (4,0) {${\sf D}$};
\node (C) at (0,-1.5) {${\sf C}$};
\node (D) at (4,-1.5) {${\sf D}$};
\node (E) at (2,-.75) {$\phi \ \twocommute$};
\draw[->] (A) to node [above] {$F$} (B);
\draw[->] (A) to node [left] {$\tau_{\sf C}$} (C);
\draw[->] (C) to node [below] {$F$} (D);
\draw[->] (B) to node [right] {$\tau_{\sf D}$} (D);
\path (0,-.75) coordinate (basepoint);
\end{tikzpicture}\label{eq:equivariantfunctor}
\eeq
%$$
% F\xrightarrow{\eta_{\sf C}} F\circ \tau_{\sf C}\circ \tau_{\sf C} \xrightarrow{\phi^2} \tau_{\sf D}\circ \tau_{\sf D} 
%$$
%
%\beq
%\begin{tikzpicture}[baseline=(basepoint)];
%\node (A) at (0,0) {$F$};
%\node (B) at (4,0) {$F\circ \tau_{\sf C}\circ \tau_{\sf C}$};
%\node (C) at (0,-1.5) {$\tau_{\sf D}\circ \tau_{\sf D}\circ F$};
%\draw[->] (B) to node [above] {$\eta_{\sf C}$} (A);
%\draw[->] (C) to node [left] {$\eta_{\sf C}$} (A);
%\draw[->] (C) to node [left] {$\phi^2$} (B);
%\path (0,-.75) coordinate (basepoint);
%\end{tikzpicture}\label{eq:equivariantfunctor}
%\eeq
and an isomorphism $\beta\colon \phi^2\circ \eta_{\sf C}\rightarrow \eta_{\sf D}$ between internal natural transformations in the diagram
\beq
\begin{tikzpicture}[baseline=(basepoint)];
\node (A) at (0,0) {${\sf C}$};
\node (B) at (4,0) {${\sf D}$};
\node (C) at (0,-1.5) {${\sf C}$};
\node (D) at (4,-1.5) {${\sf D}$};
\node (F) at (0,-3) {${\sf C}$};
\node (G) at (4,-3) {${\sf D}$};
\node (E) at (2,-.75) {$\phi \ \twocommute$};
\draw[->,bend right=50] (A) to node [left] {$\id_{\sf C}$} (F);
\draw[->,bend left=50] (B) to node [right] {$\id_{\sf D}$} (G);
\draw[->] (A) to node [above] {$F$} (B);
\draw[->] (A) to node [left] {$\tau_{\sf C}$} (C);
\draw[->] (C) to node [below] {$F$} (D);
\draw[->] (B) to node [right] {$\tau_{\sf D}$} (D);
\draw[->] (C) to node [left] {$\tau_{\sf C}$} (F);
\draw[->] (F) to node [below] {$F$} (G);
\draw[->] (D) to node [right] {$\tau_{\sf D}$} (G);
\node (H) at (2,-2.25) {$\phi \ \twocommute$};
\node (I) at (-.45,-1.5) {$\stackrel{\eta_{\sf C}}{\Rightarrow}$};
\node (I) at (4.45,-1.5) {$\Leftarrow$};
\path (0,-1.5) coordinate (basepoint);
\end{tikzpicture}\label{eq:equivariantfunctor2}
\eeq
%where (in an abuse of notation) we implicitly use $\tau_{\sf C}$ and $\tau_{\sf D}$ to regard $\phi\circ \phi$ as a natural transformation between the functor $F\simeq \id_{\sf D}\circ F \circ \id_{\sf C}$ and itself. 
where $\beta$ is compatible with $\delta_{\sf C}$ and $\delta_{\sf D}$. 
\end{defn}

Categories internal to stacks also have a notion of equivariance data for internal natural transformations, defined as follows. 

\begin{defn}  \label{defn:equivariancedata2}
Let  $({\sf C},\tau_{\sf C},\eta_{\sf C},\delta_{\sf C})$ and $({\sf D},\tau_{\sf D},\eta_{\sf D},\delta_{\sf D})$ be internal categories with internal involution, and let $F,G\colon {\sf C}\to {\sf D}$ be internal functors with equivariance data $(\phi_F,\beta_F)$ and $(\phi_G,\beta_G)$. \emph{Equivariance data} for an internal natural transformation $E\colon F\Rightarrow G$ is an isomorphism $\sigma\colon \phi_G\circ E\simeq E\circ \phi_F$  of natural transformations compatible with~$\beta_F$ and~$\beta_G$.
%\beq
%&&\begin{tikzpicture}[baseline=(basepoint)];
%\node (A) at (0,0) {${\sf C}$};
%\node (B) at (5,0) {${\sf D}$};
%\node (C) at (0,-1.5) {${\sf C}$};
%\node (D) at (5,-1.5) {${\sf D}$};
%\node (E) at (2,-.75) {$\phi_G \ \twocommute$};
%\node (AA) at (2.5,.3) {$\Downarrow E$};
%\node (BB) at (2.5,-1.8) {$\Uparrow E$};
%\draw[->,bend left=25] (A) to node [above] {$F$} (B);
%\draw[->] (A) to node [below] {$G$} (B);
%\draw[->] (A) to node [left] {$\tau_{\sf C}$} (C);
%\draw[->] (C) to node [above] {$G$} (D);
%\draw[->,bend right=25] (C) to node [below] {$F$} (D);
%\draw[->] (B) to node [right] {$\tau_{\sf D}$} (D);
%\path (0,-.75) coordinate (basepoint);
%\end{tikzpicture}\label{diag:RPtwisted}
%\eeq
%with the property that the composition of internal natural transformations 
%\beq\label{eq:involproperty}
%&&\phi_G\circ E\xrightarrow{\phi_G\star \sigma} \phi_G\circ\phi_G\circ E\xrightarrow{\beta_G} E
%\eeq
%equals $\sigma$. 
\end{defn}

\subsection{Stolz and Teichner's geometric bordism categories}\label{sec:STBord}
Given a smooth manifold~$M$ and rigid geometry $(G,\M)$, we now recall the definition of $(G,\M)\Bord(M)$ as a category internal to symmetric monoidal stacks.
% this definition is closely related to (but distinct from) \cite[Definitions~51, 52, 68, 69]{HST}.

\begin{defn}[{\cite[Definitions~2.21, 2.46, end of 2.48, and 4.4]{ST11}}]\label{defn:GBord}
Define a category $(G,\M)\Bord(M)$ internal to symmetric monoidal stacks as follows. 

\vspace{.1in}

\noindent \underline{An object of the stack $\Ob((G,\M)\Bord(M))$ is an equivalence class of the data:}
\begin{enumerate}
\item a $(G,\M)$-pair $Y^c\subset Y\to S$ where $Y^c\to S$ is proper;
\item a partition into two components $Y^+\coprod Y^-\simeq Y\setminus Y^c\subset Y$; and
\item a map of supermanifolds $\phi\colon Y\to M$.
\end{enumerate}
The supermanifold $Y^c$ is called the \emph{core}. The supermanifold $Y\supset Y^c$ in which the core sits is the \emph{collar}. The collar is partitioned into an \emph{outgoing} component $Y^+$ and an \emph{incoming} component~$Y^-$.
% see~\cite[Figure~1]{ST11} for a picture. 
%Isomorphisms between  objects in $\Ob((G,\M)\Bord(M))$ are roughly given by spans $Y \hookleftarrow W \hookrightarrow Y'$ of open inclusions over $S$ that are compatible with the $(G,\M)$-pair structure and the map to $M$. 
We apply the following equivalence relation to these data: 
\beq\label{eq:equivalencerelation}
&&(Y^c\subset Y, Y^\pm,\phi)\sim (Y^c\subset Y, Y^\pm,\phi') \iff (\phi)|_{W^+\bigcup Y^c}=(\phi')|_{W^+\bigcup Y^c}
\eeq
where $W\subset Y$ is an open submanifold containing $Y^c\subset Y$, and $W^+=W\bigcap Y^+$. In short, we only remember the data of the map $\phi$ on a half-open germ of a neighborhood of $Y^c\subset Y$ of the outgoing collar. We often abbreviate the data of an object of $\Ob((G,\M)\Bord(M))$ to $\phi\colon Y\to M$. 

\vspace{.1in}

\noindent\underline{A morphism in $\Ob((G,\M)\Bord(M))$}, denoted $(Y^c_0\subset Y_0,Y^\pm_0,\phi_0)\to (Y^c_1\subset Y_1,Y^\pm_1,\phi_1)$ is an equivalence class of open submanifolds $V_0\subset Y_0$ and $V_1\subset Y_1$ containing the cores and isometries $f\colon V_0\to V_1$ in the commutative diagram
\beq
&&\begin{tikzpicture}[baseline=(basepoint)];
\node (AA) at (-2,0) {$Y_0^c$};
\node (BB) at (-2,-1) {$Y_1^c$};
\node (A) at (0,0) {$V_0$};
\node (B) at (2,0) {$Y_0$};
\node (C) at (0,-1) {$V_1$};
\node (D) at (2,-1) {$Y_1$};
\node (E) at (4,-.5) {$M$.};
\draw[->,right hook-latex] (AA) to (A);
\draw[->,right hook-latex] (BB) to (C);
\draw[->,right hook-latex] (A) to (B);
\draw[->] (A) to node [left] {$f$} (C);
\draw[->,right hook-latex] (C) to (D);
\draw[->] (AA) to node [left] {$f|_{Y_0^c}$} (BB);
\draw[->] (B) to node [above] {$\phi_0$} (E);
\draw[->] (D) to node [below] {$\phi_1$} (E);
\path (0,-.75) coordinate (basepoint);
\end{tikzpicture}\nonumber
\eeq
These data are subject to the equivalence relation that $f\colon V_0\to V_1$ and $f'\colon V_0'\to V_1'$ represent the same isomorphism if $f|_{V_0\bigcap V_0'}=f'|_{V_0\bigcap V_0'}$. 

%roughly given by morphisms of $(G,\M)$-pairs compatible with the map to $M$, where morphisms are regarded as equal subject to an equivalence relation analogous to~\eqref{eq:equivalencerelation}. 
Disjoint union endows $\Ob((G,\M)\Bord(M))$ with a symmetric monoidal structure. 

\vspace{.1in}

\noindent \underline{An object of the stack $\Mor((G,\M)\Bord(M))$ is given by the data:}
\begin{enumerate}
\item a $(G,\M)$-family $\Sigma\to S$,
\item a map of supermanifolds $\Phi\colon \Sigma\to M$,
\item a pair of objects of $\Ob((G,\M)\Bord(M))$, denoted $\phi_{\rm in}\colon Y_{\rm in}\to M$, $\phi_{\rm out}\colon Y_{\rm out}\to M$, 
\item spans of smooth maps $Y_{\out} \hookleftarrow W_{\out} \stackrel{i_{\out}}{\to} \Sigma\stackrel{i_{\inn}}{\leftarrow}  W_{\inn} \hookrightarrow Y_{\inn}$, where $W_{\out/\inn}\hookrightarrow Y_{\out/\inn}$ are open embeddings. 
\end{enumerate}
Define the notation
$$
W^{\pm}_{\rm in/out}:= Y^{\pm}_{\rm in/out}\bigcap W_{\rm in/out},\quad i^-_{\rm in}=i_{\rm in}|_{W^-_{\rm in}},\quad  i^-_{\rm out}=i_{\rm out}|_{W^-_{\rm out}}.
$$
Then $\Sigma^c:=\Sigma\setminus (i_{\rm in}(W_{\rm in}^-)\cup i_{\rm out}(W_{\rm out}^+))$ is the \emph{core}. The above data are require to satisfy
\begin{enumerate}
\item the image $W_{\rm in/out}\hookrightarrow Y_{\rm in/out}$ contains a neighborhood of the core $Y^c_{\rm in/out}$;
\item $\Sigma^c\to S$ is proper;
\item $i^-_{\rm in}$ and $i^-_{\rm out}$ are isometries onto their image in $\Sigma\setminus i_{\rm out}(W_{\rm out}^+\cup Y_{\rm out}^c)$; and
\item the maps $\phi_{\rm in/out}\colon Y_{\rm in/out} \to M$ are compatible with $\Phi\colon \Sigma\to M$. 
\end{enumerate}
Similar to~\eqref{eq:equivalencerelation}, we impose an equivalence relation on objects
\beq\label{eq:equivalencerelation2}
&&(\Sigma,Y_\inn,Y_\out,\Phi_1)\sim (\Sigma,Y_\inn,Y_\out,\Phi_2)\iff (\Phi_1)|_{W^+_\out\bigcup \Sigma^c}=(\Phi_2)|_{W^+_\out\bigcup \Sigma^c}
\eeq
so again the equivalence class of $\Phi$ is only the data of a half-open neighborhood of the core $\Sigma^c\to M$ inside the outgoing collar. 

\vspace{.1in}

%A useful picture is~\cite[Figure~2]{ST11}, though beware a typo: the underbrace denoting~$\Sigma^c$ should extend to $i_0(Y_0^c)$. 
\noindent\underline{A morphism in $\Mor((G,\M)\Bord(M))$} is an equivalence classes of open sub supermanifolds $\Omega_i\subset \Sigma_i$, $V_{i,\inn}\subset Y_{i,\inn}$, $V_{i,\out}\subset Y_{i,\out}$ satisfying
$$
\Sigma_i^c\subset \Omega_i\subset \Sigma_i,\quad Y_{i,\inn}^c\subset V_{i,\inn}\subset W_{i,\inn}\bigcap i^{-1}_\inn(\Omega_i)\subset Y_{i,\inn},\quad Y_{i,\out}^c\subset V_{i,\out}\subset W_{i,\out}\bigcap i^{-1}_\out(\Omega_i)\subset Y_{i,\out}
$$ 
and isometries $f_\inn \colon V_{0,\inn}\to V_{1,\inn}$, $f_\out \colon V_{0,\out}\to V_{1,\out}$ and $F\colon \Omega_0\to \Omega_1$ in the commutative diagram of supermanifolds over $M$
\beq
&&\begin{tikzpicture}[baseline=(basepoint)];
\node (A) at (-2,0) {$Y_{0,\out}^c$};
\node (B) at (-2,-1) {$Y_{1,\out}^c$};
\node (C) at (0,0) {$V_{0,\out}$};
\node (D) at (2,0) {$\Omega_0$};
\node (E) at (0,-1) {$V_{1,\out}$};
\node (F) at (2,-1) {$\Omega_1$};
\node (G) at (6,0) {$Y_{0,\inn}^c$};
\node (H) at (6,-1) {$Y_{1,\inn}^c$};
\node (I) at (4,0) {$V_{0,\inn}$};
\node (J) at (4,-1) {$V_{1,\inn}$};
\draw[->] (G) to node [right] {$f_\inn|_{Y_{0,\out}^c}$} (H);
\draw[->] (A) to node [left] {$f_\out|_{Y_{0,\inn}^c}$} (B);
\draw[->,left hook-latex] (C) to (A);
\draw[->,left hook-latex] (E) to (B);
\draw[->] (E) to (F);
\draw[->] (C) to (D);
\draw[->] (I) to (D);
\draw[->] (J) to (F);
\draw[->] (C) to node [left] {$f_\out$} (E);
\draw[->] (I) to node [left] {$f_\inn$} (J);
\draw[->] (D) to node [left] {$F$} (F);
\draw[->,right hook-latex] (I) to (G);
\draw[->,right hook-latex] (J) to (H);
\path (0,-.75) coordinate (basepoint);
\end{tikzpicture}\label{eq:dataofmorphismofmorphism}
\eeq
where we also require $F$ preserves the core of the bordism, $F|_{\Sigma^c_0}\colon \Sigma^c_0\to \Sigma_1^c\subset \Omega_1$. These data are equivalent~$(f,g,F)\sim (f',g',F')$ when
$$
f_\inn |_{V_{0,\inn}\bigcap V_{0,\inn}'}=f'_\inn |_{V_{0,\inn}\bigcap V_{0,\inn}'},\quad f_\out |_{V_{0,\out}\bigcap V_{0,\out}'}=f'_\inn |_{V_{0,\out}\bigcap V_{0,\out}'},\quad F|_{\Omega_0\bigcap \Omega_0'}=F'|_{\Omega_0\bigcap \Omega_0'}.
$$
Disjoint union endows $\Mor((G,\M)\Bord(M))$ with a symmetric monoidal structure. We often abbreviate the data of an object of $\Mor((G,\M)\Bord(M))$ to $Y_{\rm in}\coprod Y_{\rm out} \to \Sigma\xrightarrow{\Phi} M$, or even more succinctly to $\Phi\colon \Sigma\to M$. 

\vspace{.1in}

\noindent \underline{Unit, source, and target functors} are given by the assignments on objects
$$
\u(Y\to M)=(Y\coprod Y\stackrel{{\rm fold}}{\to} Y\to M), \qquad \s(Y_{\rm in}\coprod Y_{\rm out} \to \Sigma\to M)=(Y_{\rm in}\to M)
$$
$$
\t(Y_{\rm in}\coprod Y_{\rm out} \to \Sigma\to M)=(Y_{\rm out}\to M),
$$
with similar assignments on morphisms. For any $S$, let $\emptyset\in \Ob((G,\M)\Bord(M))(S)$ and $\emptyset\in \Mor((G,\M)\Bord(M))(S)$ denote the empty object over~$S$, i.e., the symmetric monoidal unit over~$S$ for these symmetric monoidal stacks.

\vspace{.1in}

\noindent \underline{Composition} is given by gluing bordisms along collars; see~\cite[page~20]{ST11} and \cite[Definition~6.5]{HST}. 
\end{defn}

We now address a small technical issue. The categorical foundations of \cite{ST11} (following~\cite{Pseudocategories}) use strict fibered products of stacks in the definition of internal categories, see Definition~\ref{defn:internalcat} above. This framework leads to problems whenever one looks to replace an internal category with one whose objects and morphisms are equivalent (but not isomorphic) stacks. For example, generators and relations arguments involve such equivalences which are not isomorphisms. 
%The main issue is that an internal functor that induces an equivalence on object and morphism stacks fails to give an equivalence of internal categories. 
A simple remedy to this technical issue is to modify Definition~\ref{defn:internalcat}, requiring the source and target functors to be fibrations of stacks. This implies that the strict fibered product agrees with the homotopy fibrered product. 
%We also remark that the construction of the composition functor for $(G,\M)\Bord(M)$ in~\cite[page~20]{ST11} implicitly uses the homotopy fibered product: the gluing construction involves modifying collars of bordisms, which leads to isomorphic (but not equal) objects. 
Fortuitously, the source and target functors in $(G,\M)\Bord(M)$ do indeed satisfy this property.
% resolving all of these technicalities. 

\begin{lem} \label{lem:isofibration}
The source and target functors 
$$
\s,\t\colon \Mor((G,\M)\Bord(M))\to \Ob((G,\M)\Bord(M))
$$ 
are fibrations of stacks.\end{lem}
\bp
We recall that a map of stacks is a fibration if it induces a fibration of groupoids over each supermanifold $S$, e.g., see~\cite{Hollander}. We prove the statement of the lemma for the functor $\s$; the case for $\t$ is similar. 
It suffices to show that for any isomorphism~$f\colon \s(\Sigma_0,\Phi_0)\to (Y_1,\phi_1)$ in $\Ob((G,\M)\Bord(M'))(S)$ where
$$
(Y_1,\phi_1)\in \Ob((G,\M)\Bord(M'))(S), \quad (\Sigma_0,\Phi_0)\in \Mor((G,\M)\Bord(M'))(S),
$$ 
there exists an isomorphism $(f,g,F)\colon (\Sigma_0,\Phi_0)\to (\Sigma_1,\Phi_1)$ in $\Mor((G,\M)\Bord(M'))(S)$ with $\s(f,g,F)=f$. We recall that the notation $\Sigma_0$ for a cobordism is an abbreviation for the data 
\beq\label{eq:olddata}
(Y_{0,\out} \hookleftarrow W_{0,\out} \stackrel{i_{0,\out}}{\hookrightarrow} \Sigma_0\stackrel{i_{0,\inn}}{\hookleftarrow}  W_{0,\inn} \hookrightarrow Y_{0,\inn},\Phi_0\colon \Sigma_0\to M).
\eeq
Similarly, $f\colon \s(\Sigma_0,\Phi_0)\to (Y_1,\phi_1)$ is an abbreviation for the data of open sub supermanifolds $V_0\subset Y_0$, $V_1\subset Y_1$ containing the core and an isometry $f\colon V_0\to V_1$ preserving the cores. 

Using the notation above, we define $(\Sigma_1,\Phi_1)$ by the data 
\beq\label{eq:newdata}
(Y_{0,\out} \hookleftarrow W_{0,\out} \stackrel{i_{0,\out}}{\hookrightarrow} \Sigma_0\stackrel{i_{0,\inn}\circ f^{-1}}{\hookleftarrow}  f(W_{0,\inn}\bigcap V_0) \hookrightarrow Y_{1},\Phi_0\colon \Sigma_0\to M).
\eeq
where the relevant properties are inherited from~\eqref{eq:olddata}. 
Next define an isomorphism $(\Sigma_0,\Phi_0)\to (\Sigma_1,\Phi_1)$ by taking $V_{0,\out}=V_{1,\out}=W_{\out}$, $\Omega_0=\Omega_1=\Sigma_0$, $V_{0,\inn}=W_{0,\inn}\bigcap V_0$, $V_{1,\inn}=f(W_{0,\inn}\bigcap V_0)$, $f_\out=\id_{W_{\out}}$, $F=\id_{\Sigma_0}$ and $f_\inn=f|_{W_{0,\inn}\bigcap V_0}$. By construction, the diagram~\eqref{eq:dataofmorphismofmorphism} commutes with this data and the requisite properties are satisfied to define an isomorphism in $\Mor(1|1\EBord(M))$. The source functor applied to this isomorphism agrees with $(f,V_0,V_1)$ on restriction to $V_0\bigcap W_{0,\inn}$ and $f(V_0\bigcap W_{0,\inn})$, so this morphism is equal to~$f$ because morphisms are defined as equivalence classes. This completes the proof. 
%Possibly after shrinking neighborhoods defining $f$ (which by definition gives the same isomorphism $f$), we may define the desired lift $(\Sigma_1,\Phi_1)$ in terms of the same data as $(\Sigma_0,\Phi_1)$, but where the map $i_{1,\inn}\colon W_{1,\inn}\to \Sigma$ is modified by $f^{-1}$. Inspecting the diagram~\eqref{eq:dataofmorphismofmorphism}, taking $g=\id$ and $F=\id$ gives an isomorphism with $\s(f,g,F)=f$. 
%To promote this to a fibration of symmetric monoidal stacks, we require these lifts to be compatible with the symmetric monoidal structure, i.e., disjoint union. This is clear from the construction. 
\ep

The following is evident from Definition~\ref{defn:GBord}. 
\begin{lem}\label{lem:naturalityofGBord}
A map $f\colon M\to M'$ of smooth manifolds induces an internal functor 
$$
f_*\colon (G,\M)\Bord(M)\to (G,\M)\Bord(M'),\qquad \begin{array}{ccc} (\phi\colon Y\to M) & \mapsto & (f\circ \phi\colon Y\to M') \\ (\Phi\colon \Sigma\to M)&\mapsto & (f\circ \Phi\colon \Sigma\to M')\end{array}
$$ 
and these internal functors are compatible with composition of maps between manifolds. 
\end{lem}

%\begin{rmk}
%The conditions~\eqref{eq:equivalencerelation} and~\eqref{eq:equivalencerelation2} 
%\end{rmk}

\subsection{Twisted geometric field theories}\label{sec:twistedFT}

\begin{defn}[Sketch of {\cite[Definitions~2.24, 2.47, \S4.3]{ST11}}]\label{defn:TV}
Define the category~$\TV$ internal to symmetric monoidal stacks as follows. The stack $\Ob(\TV)$ classifies bundles of $\Z/2$-graded, nuclear Fr\'echet spaces over~$\C$, i.e., locally free sheaves of nuclear Fr\'echet super vector spaces. The stack $\Mor(\TV)$ has objects maps between vector bundles, and as morphisms commuting squares 
\beq
&&\begin{tikzpicture}[baseline=(basepoint)];
\node (A) at (0,0) {$V$};
\node (B) at (3,0) {$W$};
\node (C) at (0,-1) {$V'$};
\node (D) at (3,-1) {$W'$};
\draw[->] (A) to (B);
\draw[->] (A) to node [left] {$\simeq $} (C);
\draw[->] (C) to (D);
\draw[->] (B) to node [right] {$\simeq$} (D);
\path (0,-.75) coordinate (basepoint);
\end{tikzpicture}\label{diag:unitary}
\eeq
whose horizontal arrows are maps of vector bundles and whose vertical arrows are fiberwise isomorphisms of vector bundles, i.e., pullbacks along base changes. There are functors ${\sf s,t}\colon \Mor(\TV)\to \Ob(\TV)$ from restricting to the left and right halves of the square~\eqref{diag:unitary}. Finally, endow $\Ob(\TV)$ and $\Mor(\TV)$ with symmetric monoidal structures from the projective tensor product of topological vector bundles.  
\end{defn}

\begin{rmk} The above deviates from the definition in~\cite{ST11} in two ways. First, we only work with vector bundles rather than general sheaves of topological vector spaces (i.e., we impose a locally free condition). Second, we specialize from general locally convex topological vector spaces to nuclear Fr\'echet spaces. We make this choice in order to guarantee that certain formal categorical structures in the functorial definition of quantum field theory agree with standard structures in physics: Fr\'echet spaces have the \emph{approximation property}~\cite[page~109]{Schaefer}, leading to a well-behaved theory for traces of nuclear operators, see~\cite[\S4.5]{STTraces}. It also happens that all the desired constructions in index theory result in vector bundles of nuclear Fr\'echet spaces, so this restriction does not any exclude any of the desired examples. We refer to~\cite[Appendix~2]{costbook} for a succinct overview of nuclear Fr\'echet spaces. Stolz and Teichner's motivation to work with more general sheaves of topological vector spaces is explained in Remark~\ref{rmk:concordance}. 
\end{rmk}

%Definition~\ref{defn:TV} naturally leads to an internal Morita category~$\TA$.
\begin{defn}[Sketch of {\cite[Definition~5.1]{ST11}}] \label{defn:Morita}
Define the \emph{internal Morita category}~$\TA$ as a category internal to symmetric monoidal stacks as follows. The stack $\Ob(\TA)$ classifies locally trivial bundles of nuclear Fr\'echet super algebras and isomorphisms of algebra bundles. The stack $\Mor(\TA)$ classifies locally trivial bundles of nuclear Fr\'echet bimodules over topological algebras. Composition comes from the projective tensor product of bimodule bundles,
$$
\m \colon \Mor(\TA)\times_{\Ob(\TA)}\Mor(\TA)\to \Mor(\TA),\qquad \m(B,B')=B\otimes_A B'.
$$ 
Viewing an algebra bundle as a bimodule over itself constructs the unit map $\u\colon \Ob(\TA)\to \Mor(\TA)$. Endow $\Ob(\TA)$ and $\Mor(\TA)$ with monoidal structures from the projective tensor product of algebra and bimodule bundles over the structure sheaf of the base. 
\end{defn}

%\begin{rmk} Following the same type of argument as Lemma~\ref{lem:isofibration}, the source maps in $\TV$ and $\TA$ can be seen to be fibrations of stacks. Hence, we may identify the strict fibered product for composition in these internal categories with the weak fibered product.
%\end{rmk}

\begin{ex}\label{ex:one} Given any category ${\sf C}$ internal to symmetric monoidal stacks, there is an internal functor
$$
\one \colon {\sf C}\to \TA
$$
defined as follows. For $Y\in \Ob({\sf C})(S)$, assign the trivial algebra bundle on~$S$, i.e., the sheaf of algebras whose global sections are $C^\infty(S)$. For $\Sigma\in \Mor({\sf C})(S)$, assign the trivial bimodule bundle over $S$, i.e., the sheaf of bimodules whose global sections are $C^\infty(S)$ as a bimodule over $C^\infty(S)$. The 2-commuting data $\mu$ is determined by the canonical isomorphism $C^\infty(S)\otimes_{C^\infty(S)} C^\infty(S)\simeq C^\infty(S)$, while the 2-commuting data~$\epsilon$ comes from the canonical identification of the algebra $C^\infty(S)$ with the identity bimodule over itself. 
\end{ex}

\begin{defn}[{\cite[Definitions~2.25,~2.48, 4.12]{ST11}}]\label{defn:FTmain}
A \emph{geometric field theory over $M$} is an internal functor
\beq
E\colon (G,\M)\Bord(M)\to \TV.\label{eq:GFT}
\eeq
\end{defn}

\begin{defn}[{\cite[Definition~5.2]{ST11}}]\label{defn:FTmain2}
A \emph{twisted geometric field theory over $M$} is an internal natural transformation~$E$
\beq
&&\begin{tikzpicture}[baseline=(basepoint)];
\node (A) at (0,0) {$(G,\M)\Bord(M)$};
\node (B) at (5,0) {$\TA$};
\node (C) at (2.5,0) {$E \Downarrow$};
\draw[->,bend left=15] (A) to node [above] {$\one$} (B);
\draw[->,bend right=15] (A) to node [below] {$\twist$} (B);
\path (0,0) coordinate (basepoint);
\end{tikzpicture}\label{eq:twistedGFT}
\eeq
where $\twist$ is a \emph{twist functor} and $\one$ is defined in Example~\ref{ex:one}. Twisted field theories are the objects of a groupoid whose morphisms are isomorphisms between internal natural transformations. 
\end{defn}

\begin{rmk} An (untwisted) geometric field theory in the sense of Definition~\ref{defn:FTmain} determines a twisted field theory in the sense of Definition~\ref{defn:FTmain2} for the trivial twist functor $\twist=\one$~\cite[Lemma~5.7]{ST11}. 
\end{rmk}

%\begin{rmk} We caution that the above definition of twisted field theory is unrelated to \emph{topological twists} in supersymmetric quantum field theory, e.g., the $A$- and $B$-twists in mirror symmetry. Instead, Stolz and Teichner's terminology reflects the fact that twists in the sense of Definition~\ref{defn:FTmain2} are expected to lead to twisted cohomology groups under the conjectures~\eqref{eq:conjecture}, with the $\Z$-grading on cohomology being a special kind of twist. 
%\end{rmk}

\subsection{Additional structures}\label{sec:additionalstructuresGM}

Additional structures in the rigid geometry $(G,\M)$ lead to additional structures on the internal bordism category. 

\begin{lem}\label{lem:GMaddstructure}
Fix a rigid geometry $(G,\M)$. Given a homomorphism $g\colon G\to G$ and a $G$-equivariant map $f\colon \M\to \M$ relative to $g$ that preserves the image of the inclusion $\M^c\subset \M$, there is a uniquely determined internal functor
$$
F\colon (G,\M)\Bord(M)\to (G,\M)\Bord(M). 
$$
\end{lem}
\bp
The endofunctors of $(G,\M)$-manifolds from Lemma~\ref{lem:GMaddstructure1} determine maps of stacks,
$$
F_0\colon \Ob((G,\M)\Bord(M))\to \Ob((G,\M)\Bord(M)), \quad (Y\to S,\phi)\mapsto (f(Y\to S),\phi) $$$$ F_1\colon \Mor((G,\M)\Bord(M))\to \Mor((G,\M)\Bord(M)),\quad (\Sigma\to M,\Phi)\mapsto (f(\Sigma\to M),\Phi),
$$
by taking the same collar data on bordisms, but different $(G,\M)$-structures. 
By construction, these maps of stacks are (strictly) compatible with source, target, unit, and composition functors defining $(G,\M)\Bord(M)$. 
\ep

%\begin{rmk}\label{rmk:autofromcenter} There is an isomorphic description of the functor $F$ above in the special case that $Y$ or $\Sigma\subset S\times \M$ is an open submanifold. In this case, $f$ determines a (non-identity) isomorphism $f|_Y\colon Y\to Y'\subset \M\times S$ or $f|_{\Sigma}\colon \Sigma\to \Sigma'\subset \M\times S$ of submanifolds, and so we obtain the assignments
%$$
%(Y\to S,\phi)\mapsto (Y'\to S,\phi\circ (f^{-1}\times \id_S)|_{Y'}),\quad (\Sigma \to S,\Phi)\mapsto (\Sigma'\to S,\Phi\circ (f^{-1}\times \id_S)|_{Y'}).
%$$
%% an isomorphic description of the functor~$F$ is $F(\Sigma,\phi)=(f(\Sigma),\phi|\circ (f^{-1}\times \id_S)|_{f(\Sigma)}$. In the examples below, this leads to an isomorphic description of the values of the functor~$F$ when restricted to open submanifolds $Y,\Sigma\subset \M\times S$. 
%%In practice, there is often a isometry $\tilde{f}\colon f(\Sigma\to M)\to (\Sigma\to M)$, in which case it is more convenient to view the functor $F$ defined in the previous lemma as leading 
%
%\end{rmk}
%Geometric bordism categories often have additional structures that allow one to refine the notion of field theory. We describe a few of these structures that are important below.
%These structures typically are inherited from additional symmetries in the rigid geometry. 

\begin{ex}[Bordism categories with flips]\label{ex:flips}
In the context of Lemma~\ref{lem:GMaddstructure}, suppose that $\fl \colon \M\to \M$ is an involution commuting with the action of $G$ and preserving $\M^c\subset \M$. This endows every family of $(G,\M)$-manifolds with an involution. In practice, an involution $\fl\colon \M\to \M$ comes from a \emph{spin flip} automorphism, e.g., see~\cite[\S6.7, Lemma 6.39]{HST},~\cite[Example~2.39 and \S4.1]{ST11} and~\eqref{eq:spinflipfor11}. Following~\cite[\S2.6]{ST11}, such an involution~$\fl$ allows one to regard the bordism category as a category internal to symmetric monoidal stacks with involution. 
%For the rigid geometries relevant to the Stolz--Teichner program, such spin flip automorphisms exist canonically, and field theories are always taken to be flip preserving relative to this structure. 
\end{ex}

%The parity involution $(-1)^F$ on algebras, bimodules, and bimodule maps allows one to similarly regard $\TA$ as a category internal to symmetric monoidal stacks with involution. 
Next we observe that every super vector bundle or algebra bundle has a canonical automorphism coming from the parity involution on super vector spaces,
\beq
(-1)^{\sf F}\colon V\to V,\qquad (-1)^{\sf F}\colon A\to A\label{eq:flip1}
\eeq
acting by the identity on even sections and $-1$ on odd sections. This allows one to regard $\TV$ and $\TA$ as categories internal to symmetric monoidal stacks with involution. Stolz and Teichner incorporate these involutions to their field theories as follows.

%\begin{defn}[{\cite[\S2.6 and Definition 2.48]{ST11}}]
%A twisted geometric field theory is \emph{flip preserving} if it comes with a preferred extension to a natural transformation for categories internal to symmetric monoidal stacks with involution, where the involution on $(G,\M)\Bord(M)$ comes from $\fl$ and the involution on $\TA$ comes from the parity involution. 
%\end{defn}

\begin{defn}[{\cite[Definition 6.44]{HST}, \cite[\S2.6]{ST11}}]\label{defn:flips}
%The spin flip~\eqref{eq:spinflip} and parity involution~\eqref{eq:flip1} allow us to regard $\TV$, $\TA$, and $1|1\EBord(M)$ as categories internal to stacks with involution. 
Fix a rigid geometry $(G,\M)$ with an involution $\fl\colon \M\to \M$ as in Example~\ref{ex:flips}.  A \emph{flip-preserving geometric field theory} is a functor~\eqref{eq:GFT} internal to the category of symmetric monoidal stacks with involution. Similarly, one defines \emph{flip-preserving twists} and \emph{flip-preserving twisted geometric field theories} as functors and natural transformations internal to symmetric monoidal stacks with involution. 
\end{defn} 

Whenever a bordism category is endowed with a flip, we will always take (twisted) geometric field theories to mean \emph{flip-preserving} (twisted) geometric field theories. 

\begin{ex}[Orientation reversal] \label{ex:orientationGM}
Fix an orientation on the manifold $\M_{\rm red}$, and suppose that $\orr \colon \M\to \M$ is a $G$-equivariant map with the property that $\orr_{\rm red}\colon \M_{\rm red}\to \M_{\rm red}$ is an orientation reversing map and $\orr^2=\fl$ is the spin flip. Then Lemma~\ref{lem:GMaddstructure} determines an \emph{orientation reversal functor}
\beq
\Or\colon (G,\M)\Bord(M)\to (G,\M)\Bord(M)\label{eq:orientationreversalGM}
\eeq
with a natural isomorphism $\Or\circ \Or\Rightarrow \id$ (via the spin flip). Indeed, $\Or$ determines an internal involution in the sense of Definition~\ref{defn:interinvol}. 
\end{ex}

In 1-dimensional field theories (i.e., quantum mechanics) orientation-reversal is the same as time-reversal, which is typically implemented on the space of states as a $\C$-antilinear map. This structure is part of \emph{reflection positivity}, which can be adapted to functorial field theories by requiring equivariance data relative to actions generated by orientation reversal~\eqref{eq:orientationreversalGM} and the functor on vector spaces coming from complex conjugation, e.g., see~\cite[\S3]{FreedHopkins}. However complex conjugation is more subtle in the context of supermanifolds, as we explain presently. There are conjugation functors,
\beq
&&(\overline{\phantom{-}}) \colon \TV\to \TV, \ \ (V\to S)\mapsto (\overline{V}\to \overline{S})\qquad (\overline{\phantom{-}})\colon \TA\to \TA, \ \ (A\to S)\mapsto (\overline{A}\to \overline{S})\label{eq:realTVS}
\eeq
sending a vector bundle, algebra bundle, etc., to the same bundle but where the complex numbers act through conjugation; when the input is a bundle over a supermanifold~$S$, the output is a bundle over the conjugate supermanifold, $\bar{S}$  (see Definition~\ref{defn:realsmfld}). In other words, the natural internal involution of $\TA$ (see Definition~\ref{defn:interinvol}) covers the conjugation functor on supermanifolds. For this reason it does not make sense to ask for equivariance data mixing~\eqref{eq:orientationreversalGM} and~\eqref{eq:realTVS}, as they are internal functors covering different functors on supermanifolds ($\Or$ covers the identify functor). This is repaired by demanding a real structure on the relevant bordism category, which can be gotten from a real structure on the rigid geometry in the sense of Definition~\ref{defn:realFT}. 

\begin{ex}[Bordism categories with real structures, {\cite[\S8]{HST}}]\label{ex:realGM}
Given an $S$-family of $(G,\M)$-supermanifolds $(Y\to S)$, a real structure on the rigid geometry $(G,\M)$ gives a unique $(G,\M)$-stucture for the $\overline{S}$-family $(\overline{Y}\to \overline{S})$. This extends to an internal functor 
\beq
&&\RR\colon (G,\M)\Bord(M)\to (G,\M)\Bord(M),\qquad (S\leftarrow Y\to M)\mapsto (\overline{S}\leftarrow \overline{Y}\to \overline{M}\simeq M).\label{eq:realstru}
\eeq
The identification $\overline{M}\simeq M$ uses that $M$ is an ordinary manifold, and hence has a canonical real structure when regarded as a supermanifold. Via the forgetful functor $(G,\M)\Bord(M)\to \SMfld$ that sends an $S$-family of bordisms to $S$, the internal functor~\eqref{eq:realstru} covers the conjugation functor on supermanifolds. We observe that $\RR$ determines an internal involution in the sense of Definition~\ref{defn:interinvol}. 
\end{ex}

The functor $\one\colon (G,\M)\Bord\to \TA$ has canonical equivariance data for any internal involution on the source.

\begin{defn}\label{defn:RP0}
Let $(G,\M)$ be a rigid geometry with a chosen real structure and orientation reversing automorphism. 
A \emph{reflection structure} for a twist $\twist$ is equivariance data for the involutions $\RR\circ \Or$ and $(\overline{\phantom{-}})$.
A \emph{reflection structure} for a twisted field theory $E\colon \one\to \twist$ is equivariance data for the internal natural transformation $E$. 
% natural a choice of isomorphism between this pair of internal natural transformations $\one\Rightarrow\twist$.
% (see Definition~\ref{rmk:twistedcategory}). 
\end{defn}

\begin{rmk}
%The above definition closely follows Freed--Hopkins~\cite[\S3]{FreedHopkins}, though in the context of supermanifolds one needs a real structure on the rigid geometry as the functor $(\overline{\phantom{-}})$ covers complex conjugation on the category of supermanifolds. In short, 
Roughly, a reflection structure stipulates that orientation reversal of bordisms correspond to complex conjugation of vector spaces (which is data). In tandem with this structure is the \emph{reflection positivity property}, which demands that a certain hermitian pairing is positive. Suggestively we write this pairing as
\beq\label{eq:ahermpairing}
\overline{E(Y)}\otimes E(Y) \stackrel{h\otimes \id}{\simeq} E(\overline{Y})\otimes E(Y) \stackrel{E(\Sigma)}{\to} \C
\eeq
where the isomorphism $h$ comes from the reflection structure, and $\Sigma$ is a bordism with boundary $\overline{Y}\coprod Y$ where $Y$ is the orientation reversal of $Y$. Then the positivity property is the positivity of the pairing~\eqref{eq:ahermpairing}. 
%The pairing itself comes from a bordism $\Sigma\colon \overline{Y}\coprod Y\to \emptyset$, where $\overline{Y}$ is $Y$ with the opposite orientation. 
We provide the precise property for $1|1$-Euclidean field theories in Definition~\ref{defn:RPFT}.

%There are examples (e.g., in dimension~2) where neither $r$ nor $\Or$ can be defined as indicated above, but instead only the composition $\Or\circ r$ is well-defined as an internal endofunctor on the category of bordisms. This makes reflection positivity somewhat more delicate in these examples. 
\end{rmk}

\begin{defn}[Compare~{\cite[Definition 6.47]{HST}}]\label{defn:Real0}
Let $(G,\M)$ be a rigid geometry with a chosen real structure, and $\twist\colon (G,\M)\Bord(M)\to \TA$. A \emph{real structure} for $\twist$ is equivariance data for the involutions $\RR$ and $(\overline{\phantom{-}})$. Given a twisted field theory $E\colon \one\to \twist$, a \emph{real structure} on $E$ is equivariance data for the internal natural transformation. 
%
%
% the data $(E,\eta,\rho)$ where $E$ is a geometric field theory and $(\eta,\rho)$ is an internal natural isomorphism 
%\beq
%&&\begin{tikzpicture}[baseline=(basepoint)];
%\node (A) at (0,0) {$(G,\M)\Bord(M)$};
%\node (B) at (4,0) {$\TV$};
%\node (C) at (0,-1.5) {$(G,\M)\Bord(M)$};
%\node (D) at (4,-1.5) {$\TV.$};
%\node (E) at (2,-.75) {$(\eta,\rho) \ \Uparrow$};
%\draw[->] (A) to node [above] {$E$} (B);
%\draw[->] (A) to node [left] {$r$} (C);
%\draw[->] (C) to node [below] {$E$} (D);
%\draw[->] (D) to node [right] {$(\overline{\phantom{-}})$} (B);
%\path (0,-.75) coordinate (basepoint);
%\end{tikzpicture}\label{diag:real}
%\eeq
%where the functor~$r$ is given by~\eqref{eq:realstru} and $(\overline{\phantom{-}})$ is from~\eqref{eq:realTVS}. Real twists and real twisted field theories are defined completely analogously to~\eqref{diag:RPtwisted}. 
\end{defn}

%The following is inspired by~\cite[\S3]{FreedHopkins}. 

\begin{ex}[Renormalization group flow]\label{ex:RGgeneral}
In many cases, there is a nontrivial $\R_{>0}$-action on $\M$ that commutes with the $G$-action. For example, the Euclidean rigid geometries from Example~\ref{ex:modelgeom} admit such a structure from the standard dilation action on $\R^d$. Applying Lemma~\ref{lem:GMaddstructure}, this leads to the action of the \emph{renormalization group},
\beq
\RG^\mu\colon (G,\M)\Bord(M)\to (G,\M)\Bord(M)\qquad \mu\in \R_{>0},\label{eq:generalRG}
\eeq
which then acts on the category of (twisted) field theories by precomposition with the internal functor~\eqref{eq:generalRG}. 
\end{ex}

\section{The families index in KO-theory}\label{sec:KOindex}

\subsection{Superconnections}\label{sec:cutoffconstrution}

A super vector bundle is taken to be a locally free sheaf of $\Z/2$-graded nuclear Fr\'echet spaces; we make no finite rank assumption. 

\begin{defn}[{\cite{Quillensuperconn}}]\label{defn:superconn}
Let $V\to M$ be a complex vector bundle. A \emph{superconnection}~$\A$ on $V$ is an odd, $\C$-linear map of Fr\'echet spaces satisfying a Leibniz rule
\beq
&&\A\colon \Omega^\bullet(M;V)\to \Omega^\bullet(M;V),\quad \A(\alpha\cdot v)=d\alpha \cdot v+(-1)^{|\alpha|} \alpha \cdot \A v\label{eq:ordsuperconn}
\eeq
where $v\in \Omega^\bullet(M;V)$ and $\alpha\in \Omega^\bullet(M)$. 
\end{defn}

Using the $\Z$-grading on forms, a superconnection can be written as 
\beq
\A=\sum \A^{[k]},\qquad \A^{[k]}\colon \Omega^\bullet(M;V)\to \Omega^{\bullet+k}(M;V)\label{eq:superconnectionZ}
\eeq 
where $\A_1=\nabla$ is an ordinary connection on $V$, and $\A_{2k}\in \Omega^{2k}(M;\End(V)^{\odd})$, $\A_{2k+1}\in \Omega^{2k+1}(M;\End(V)^{\ev})$ are endomorphism valued forms. 

\begin{defn}
Suppose that $\langle-,-\rangle\colon\overline{V}\otimes V\to \underline{\C}$ is a (super) hermitian pairing, where~$\otimes$ is the projective tensor product. Then a superconnection $\A$ on $V$ is \emph{self-adjoint} if 
\beq
\langle \overline{\A}_1x,y\rangle+(-1)^{|x||y|}\langle x,\A_1y\rangle=d\langle x,y\rangle,\qquad \A_{k}^*=i^{k+1}\overline{\A}_k,\label{eq:supersa}
\eeq
where $\A_k^*$ is the super adjoint of $\A_k$. 
%where $\langle-,-\rangle\colon \overline{V}\otimes V\to \underline{\C}$ is the super hermitian pairing associated to $(-,-)\colon \overline{V}\otimes_0 V\to \underline{\C}$ using~\eqref{eq:ordinarypairing}, and the super adjoints relative to $\langle-,-\rangle$ are compared to adjoints relative to $(-,-)$ using~\eqref{eq:superadjointtrans}.
\end{defn}

\begin{rmk}
Translating from super adjoints to ordinary adjoints as in~\eqref{eq:superadjointtrans}, the above is equivalent to the standard definition of self-adjoint superconnection~\cite[page 117]{BGV}. 
\end{rmk}

For a superconnection $\A$, define the 1-parameter family,
\beq
&&\A(u):=u^{-1/2}\A^{[0]}+\A^{[1]}+u^{1/2}\A^{[2]}+\dots+ u^{j/2}\A^{[j]}+\dots\quad\nonumber
%&&\A(u)\colon \Omega^\bullet(M;V)[u^{1/2},u^{-1/2}]\to \Omega^\bullet(M;V)[u^{1/2},u^{-1/2}] \nonumber
\eeq

\begin{defn}\label{defn:chernform}
The \emph{Chern form} of a $\Cl_n$-linear superconnnection is the Clifford super trace (which need not exist for infinite rank bundles)
\beq\label{eq:theChernchar}
\Ch(\A):=\sTr_{\cCl_n}(e^{-u\A(u)^2})\in \Omega^\bullet(M;\C[u^{\pm 2}])).
\eeq
\end{defn}
The absence of fractional powers of~$u$ in $\Ch(\A)$ follows from the fact that the super trace of an odd endomorphism is zero. The self-adjoint property further implies that the trace only has even powers of~$u$.

\begin{lem}
When it exists, the Chern form of a $\Cl_n$-linear superconnection $\A$ is closed and has total degree $n$.
\end{lem}

\subsection{Clifford modules and KO-theory}
We overview a model for the KO-spectrum due to~\cite[\S1.3-1.4]{Cheung} and~\cite[Theorem~7.1]{HST}, see also~\cite[\S4]{DBEBis}. 

\begin{defn}[{\cite[\S2]{Segalclassifying}}]
A \emph{topological category} is a category in which the objects and morphisms are topological spaces for which the structure maps (source, target, unit, composition) are continuous. For a topological category $\mathcal{C}$, let~$\mathcal{C}_0$ denote the space of objects and $\mathcal{C}_1$ the space of morphisms. 
\end{defn}

Let $\mathcal{H}_n$ denote an infinite-dimensional $\Z/2$-graded separable Hilbert space with a self-adjoint action by $\Cl_n$ with the property that every irreducible representation of $\Cl_n$ appears with infinite multiplicity in $\mathcal{H}_n$. Let $\mathcal{B}(\mathcal{H}_n)$ denote the algebra of $\Cl_n$-linear bounded operators in the norm topology. 

\begin{defn} 
Let $\mathcal{M}_n$ be the category whose objects are finite-dimensional $\Cl_{n}$-submodules $W\subset \mathcal{H}_n$ and whose morphisms are
\beq
&&{\rm Mor}(W_0,W_1):=\left\{ \begin{array}{cl} \begin{array}{c} \Cl_{n}\otimes \Cl_{-1}-{\rm action\ on\ } W_0^\perp \\ {\rm extending\ the }\ \Cl_{n}{\rm-action}\end{array}  & {\rm if} \ W_0\subseteq W_1\subset \mathcal{H}_n \\
\null\\
\emptyset & {\rm else} \end{array}\right.
\eeq
where $W_0^\perp\subseteq W_1$ denotes the orthogonal complement of the subspace $W_0\subseteq W_1$. Composition for a nested inclusion $W_0\subseteq W_1\subseteq W_2$ uses the direct sum of Clifford modules to obtain a $\Cl_{n}\otimes \Cl_{-1}$-module structure on the orthogonal complement of $W_0\subseteq W_2$. Endow the objects and morphisms of $\mathcal{M}_n$ with the subspace topology for the embeddings
\beq
{\rm Obj}(\mathcal{M}_n)\hookrightarrow \mathcal{B}(\mathcal{H}_n),\qquad {\rm Mor}(\mathcal{M}_n)\hookrightarrow \mathcal{B}(\mathcal{H}_n)^{\times 3} \label{eq:bddops}
\eeq
where the first inclusion identifies a submodule $W_0\subset \mathcal{H}_n$ with the (finite-rank) projection~$\mathcal{H}_n\to W_0$. For the second inclusion, a triple $(W_0,W_1,e)$ determines a pair of finite rank projection operators and a bounded operator~$e$. 
% gotten by extending $e$ to an operator on $\mathcal{H}_n$ via the zero operator on the orthogonal complement of $W_0^\perp\subset \mathcal{H}_n$. 
\end{defn} 

\begin{defn} The \emph{nerve} of a topological category $\mathcal{C}$ is the simplicial space whose $k$-simplices are $N_k\mathcal{C}=\mathcal{C}_k$ for $k=0,1$, and for $k>1$ are length $k$ chains of composable morphisms
\beq
N_k\mathcal{C}&=&\{x_0\stackrel{f_1}{\to}x_1\stackrel{f_2}{\to}x_2\stackrel{f_3}{\to}\cdots \stackrel{f_k}{\to}x_k\mid x_i\in \mathcal{C}_0, f_i\in \mathcal{C}_1\}\nonumber\\&=&\mathcal{C}_1\times_{\mathcal{C}_0}\cdots \times_{\mathcal{C}_0} \mathcal{C}_1,\nonumber 
\eeq
topologized as the fibered product.
% where the maps determining the fibered product are source and target maps,~$s,t\colon \mathcal{C}_1\rightrightarrows \mathcal{C}_0$. 
The face maps in  $N_\bullet \mathcal{C}$ are determined by composing morphisms, and degeneracies are determined by the unit map in $\mathcal{C}$. 
\end{defn}

%Face maps compose a pair of morphisms, giving a chain of length~$k-1$. Degeneracy map inserts an identity morphism $x_i\to x_i$, giving a chain of length~$k+1$. 

The \emph{geometric realization} of a simplicial space $Z_\bullet$ is 
\beq
| Z_\bullet | :=(\coprod_n Z_n\times \Delta^n)/{\sim}\qquad (f^*x,t)\sim (x,f_*t) \quad \forall f\colon [l]\to [k]\label{eq:fatrealizae}
\eeq
where $\Delta^n$ is the standard $n$-simplex, and $f\colon [l]\to [k]$ is an order-preserving map. Geometric realization is natural: a map of simplicial spaces $F\colon Z_\bullet\to Y_\bullet$ determines a map between realizations, $|F|\colon |Z_\bullet|\to |Y_\bullet|$. 

\begin{prop}[{\cite[Theorem~7.1]{HST}, \cite[\S1.3-1.4]{Cheung}}]\label{prop:ABScat}
The classifying space ${\rm B} \mathcal{M}_n$ represents the functor $\KO^{n}(-)$. 
\end{prop}
Given a topological space $Z$ and an open cover $\{U_i\}_{i\in I}$ of $Z$, define $U_\sigma=\bigcap_{i\in \sigma} U_i$ as the intersection for $\sigma\subset I$ a nonempty finite subset.

\begin{defn} \label{ex:opencover0}
For $Z$ a topological space and $\{U_i\}_{i\in I}$ an open cover, the \emph{\v{C}ech category}, denoted $\check{{\rm C}}(U_i)$, is the topological category whose objects are $\check{{\rm C}}(U_i)_0=\coprod_{\sigma} U_\sigma$ and whose morphisms are the space $\check{{\rm C}}(U_i)_1=\coprod_{\sigma\subseteq \tau} U_\tau$. The source and target maps are determined by the projection and inclusion $U_\tau\subseteq U_\sigma$, the unit is induced by the identity map $U_\sigma \to U_\sigma$, and composition comes from nested inclusions. 
\end{defn}

\begin{lem}[{\cite[Proposition~4.12]{Segalclassifying}, \cite[Theorem~2.1]{DuggerIsaksen}}]\label{lem:Cechrealize}
%The nerve of the topological category $\check{C}(U_i)$ from Example~\ref{ex:opencover0} has as $k$-simplices the nested inclusions 
%$$
%N_k\check{C}(U_i)=\coprod_{\sigma_0\subseteq \sigma_1\subseteq \cdots \subseteq \sigma_k} U_{\sigma_k}. 
%$$
%Segal shows that 
For $Z$ a topological space and $\{U_i\}_{i\in I}$ an open cover, the canonical map 
\beq
{\rm B} \check{C}(U_i)\stackrel{\sim}{\to} Z,\label{eq:Segalwe}
\eeq 
is a homotopy equivalence. 
% when the cover is numerable~\cite[Proposition~4.12]{Segalclassifying}, and Dugger--Isaksen prove the map is a weak equivalence for an arbitrary open cover~\cite[Theorem~2.1]{DuggerIsaksen}. 
%When the cover $\{U_i\}$ is countable, this is a result of Segal~\cite[Proposition~4.12]{Segalclassifying}.
% The general case is proved by Dugger and Isaksen~\cite[Theorem~2.1]{DuggerIsaksen}. 
%possibly also~\cite[\S{A}]{segal_gamma} for the equivalence between fat realizations and ordinary ones, using that the nerve of a cover is good. 
 \end{lem}

\begin{defn}[{\cite[\S4]{DBEBis}}]
A \emph{$\Cl_n$-bundle} on a manifold $M$ is data 
\begin{enumerate}
\item an open cover $\{U_i\}_{i\in I}$ of $M$,
\item for each finite subset $\sigma \subset I$, a metrized super vector bundle $W_\sigma\to U_\sigma$ with $\Cl_n$-action such that the fibers of $W_\sigma$ are self-adjoint $\Cl_n$-modules;
\item for each $\sigma\subset \tau$, a $\Cl_n$-equivariant, injective map of super vector bundles $g_{\sigma\tau}\colon W_\tau\hookrightarrow W_\sigma$, together with a $\Cl_{n}\otimes \Cl_{-1}$-action on $W_\tau^\perp\subset W_\sigma$ on the orthogonal complement of the image extending the $\Cl_n$-action.
\end{enumerate}
For $\sigma\subset \tau\subset \rho $, we require an equality of maps $g_{\sigma\tau}\circ g_{\tau\rho}=g_{\sigma\rho}\colon W_\rho\hookrightarrow W_\sigma$, and  that the $\Cl_{-1}$-action on $W_\rho^\perp\subset W_\tau \hookrightarrow W_\sigma$ restricts from the $\Cl_{-1}$-action on~$W_\tau^\perp\subset W_\sigma$.
\end{defn}

\begin{cor}\label{eq:cormaptoKO}
Given a $\Cl_n$-bundle on $M$, there exists a continuous functor
\beq\label{eq:constructtoKO}
 \check{{\rm C}}(U_i)\to \mathcal{M}_n
\eeq
that is unique up to a contractible space of choices. After realization, this determines a map $M\to B\mathcal{M}_n\simeq \KO^n$ into the $n$th space of $\KO$ spectrum.
\end{cor}
\bp The functor~\eqref{eq:constructtoKO} comes from choosing a $\Cl_n$-equivariant embedding $W_\sigma\subset \mathcal{H}_n\times U_\sigma$ over $U_\sigma$ for each $\sigma\subset I$. By Kuiper's theorem~\cite{Kuiper}, the space of such choices is contractible. Then the corollary follows from Proposition~\ref{prop:ABScat} and Lemma~\ref{lem:Cechrealize}. 
\ep

\begin{lem}[{\cite[\S6]{DBEBis}}]\label{prop:CechdeRham}
Given a $\Cl_n$ bundle on $M$, choices of self-adjoint $\Cl_n$-linear superconnections $\A_\sigma$ on each $W_\sigma$ uniquely determine a \v{C}ech--de~Rham representative of the Pontryagin character $\Ch(W)\in \H^\bullet(M;\R[u,u^{-1}])$ of $[W]=[\{U_\sigma,W_\sigma,g_{\sigma\tau},e_{\sigma\tau}\}]\in \KO^n(M)$. 
\end{lem}
\begin{proof}[Sketch of proof.] Over each $U_i$, a choice of superconnection on $W_i$ gives a differential form $\Ch_n(\A_i)\in \Omega^\bullet(U_i;\R[u,u^{-1}])$ of total degree~$n$. The data $g_{\sigma\tau},e_{\sigma\tau}$ allow one to extend this to a \v{C}ech--de~Rham cocycle. 
\ep

\begin{defn} \label{defn:superconngeneral}
A \emph{differential refinement} of a $\Cl_n$-bundle $V=\{V_\sigma,g_{\sigma\tau},e_{\sigma\tau}\}$ defined relative to an open cover $\{U_i\}_{i\in I}$ is data
\begin{enumerate}
\item[i.] for each $\sigma\subset I$, a self-adjoint $\Cl_n$-linear superconnection $\A_\sigma$ on $V_\sigma\to U_\sigma$, and
\item[ii.] a closed global form $ \Ch(V,\A,\phi)\in \Omega^\bullet(M;\R[\alpha^{\pm 1}])$ of total degree~$n$;
\item[iii.] a coboundary $\{\phi_\sigma\in \Omega^\bullet(U_\sigma;\R[\alpha^{\pm 1}])\}$ between $\Ch(V,\A,\phi)$ and the \v{C}ech--de~Rham cocycle $\Ch_n(\{U_i\},\A_i)\in \Omega^\bullet(U_*;\R[u^{\pm 1}])$ in Lemma~\ref{prop:CechdeRham}.
\end{enumerate}
A differential refinement of a $\Cl_n$-bundle is a \emph{differential cocycle} denoted by $(V,\A,\phi)=\{(V_\sigma,\A_\sigma,\phi_\sigma),(g_{\sigma\tau},e_{\sigma\tau})\}$.
\end{defn}

%\bp
%The fact that the form is closed follows from \cite[Theorem~9.17]{BGV}, and the degree restriction comes from properties of the Clifford super trace.
%\ep

%{Chern--Simons forms. Have the version for a 1-parameter family (page 46 of BGV), and then want a version for a $k$-parameter family. }
\subsection{Index bundles}\label{sec:dKO} 
Given a self-adjoint, $\Cl_{n}$-linear superconnection~$\A$ on a vector bundle $V\to M$, define the subset $U_\lambda\subset M$ for $\lambda\in \R_{>0}$ as follows. A point $m\in M$ belongs to $U_\lambda$ if $\lambda$ is not an eigenvalue of the (non-negative) operator $(\A^{[0]})^2_m$ (using the notation~\eqref{eq:superconnectionZ}):
\beq
U_\lambda=:\{m\in M\mid \lambda\notin {\rm Spec}((\A^{[0]})^2_m)\}\subset M.\label{eq:Ulambda2}
\eeq
Let $\mathcal{H}\to M$ denote the fiberwise Hilbert completion of the metrized vector bundle $V\to M$; typically $\mathcal{H}\to M$ is only a continuous (rather than smooth) vector bundle. Over each subset $U_\lambda\subset M$, define a subset $\mathcal{H}^{<\lambda}\subset \mathcal{H}|_{U_\lambda}$ as the direct sum of $\nu$-eigenspaces $\mathcal{H}_{m}^\nu$ of $(\A^{[0]})^2_m$ with eigenvalue $\nu<\lambda$
\beq
 \mathcal{H}^{<\lambda}_m:=\bigoplus_{\nu<\lambda} \mathcal{H}_{m}^\nu \subset \mathcal{H}_m\label{eq:bundleofClnmodules}
 \eeq
Letting $m$ vary, $\mathcal{H}^{<\lambda}$ denotes the $\Cl_{-n}$-invariant subset of~$\mathcal{H}|_{U_\lambda}$ whose fibers are $ \mathcal{H}^{<\lambda}_m$. Define the maps (of sets) over $U_\lambda$
\beq\label{eq:fiberwiseincludeproject}
p^\lambda\colon \mathcal{H}|_{U_\lambda}\to \mathcal{H}^{<\lambda},\qquad i_\lambda\colon \mathcal{H}^{<\lambda}\to \mathcal{H}|_{U_\lambda}
\eeq
where $p^\lambda$ is the fiberwise orthogonal projection, and $i_\lambda$ is the fiberwise inclusion. 

\begin{defn} \label{defn:superconncutoff}
A self-adjoint Clifford linear superconnection $\A$ \emph{admits smooth index bundles} if 
\begin{enumerate}
\item the collection of subsets $\{U_\lambda\}_{\lambda\in \R_{>0}}$ defined in~\eqref{eq:Ulambda2} form a smooth open cover;
\item the subsets $\mathcal{H}^{<\lambda}\subset \mathcal{H}|_{U_\lambda}\to U_\lambda$ are smooth vector bundles;
\item the composition defines a smooth connection:
$$
\nabla^{<\lambda}:=p^\lambda\circ \A^{[1]}|_{U_\lambda}\circ p^\lambda\colon \Gamma(\mathcal{H}^{<\lambda})\to \Omega^1(U_\lambda;\mathcal{H}^{<\lambda}),
$$
\item the limit of Chern--Simons forms
\beq
\eta_{\lambda}(t):=\lim_{r\to \infty} \int_1^r \sTr_{\Cl_n}(\frac{d\A_{\lambda}(s)}{ds}\exp(\A_{\lambda}(s)))ds \in \Omega^\bullet(U_\lambda ;C^\infty(\R_{>0})[u^{\pm 1}])\label{eq:improperintegral}
\eeq
exists and satisfies
\beq\label{eq:satisfies}
d\eta_\lambda(t)=\Ch(\A(t))|_{U_\lambda}-\Ch(\nabla^{<\lambda})
\eeq
where in~\eqref{eq:improperintegral}, $\A_{\lambda}(s)$ is the 1-parameter family of superconnnections on $V|_{U_\lambda} \oplus \Pi \mathcal{H}^{<\lambda}$,
\beq\label{eq:nullconcordFT}
&&\A_{\lambda}(s)=\A_{\lambda}(s)^{[0]}+\A^{[1]} \oplus \nabla^{<\lambda} +\sum_{i>1} s^{-i} \A^{[i]} \oplus 0_{\Pi \mathcal{H}^{<\lambda}},\quad \A_{\lambda}(s)^{[0]}=\left[\begin{array}{cc} \A^{[0]} & s\cdot i_\lambda \\ s \cdot p^\lambda & 0 \end{array}\right]
\eeq
where $i_\lambda$ and $p^\lambda$ are the maps~\eqref{eq:fiberwiseincludeproject}.
%Can change the $s$-family to depend on some $\chi(s)$-cutoff where we choose $\chi$ to be invertible. Then to view the $\eta$-form as a function of $t$, the integral needs to go from zero to the value of $s$ where $\chi(s)=t$. 
\end{enumerate}
\end{defn}

We sketch the proof of the equality~\eqref{eq:satisfies}. At $s=0$, the superconnection is a direct sum of the Bismut superconnection and the (parity reversal of the) ordinary connection $\nabla^{<\lambda}$, giving the right hand side of~\eqref{eq:satisfies}. For $s>0$, the operator $\A_{\lambda}(s)^{[0]}$ is invertible \cite[Lemma 7.18]{FreedLott}, and so its Chern character is zero; the limit~\eqref{eq:improperintegral} gives an explicit coboundary for the Chern form. 

For a superconnection admitting cutoffs, we review the index bundle construction from~\cite[\S8.2-8.3]{DBEBis}. We note that to construct an index map to KO-theory, only (1) and (2) are necessary in Definition~\ref{defn:superconncutoff}; conditions (3) and (4) allow a refinement to differential KO-theory. 

\begin{prop} \label{prop:appenindexbundle}
Given a $\Cl_n$-linear superconnection admitting cutoffs, a choice of locally finite open cover $\{U_\lambda\}_{\lambda \in \Lambda}$ associated with a discrete subset $\Lambda\subset \R_{>0}$ determines a map 
\beq\label{eq:KOindexspec}
\Ind_\Lambda(\A)\colon M\to \KO^n,
\eeq
to the $n$th space in the $\KO$-spectrum that is unique up to contractible choice. Given choices $\Lambda\subset \Lambda' \subset \R_{>0}$ related by a refinement of open covers, there is a unique homotopy between the associated maps~\eqref{eq:KOindexspec}.
\end{prop}
\bp
The $\Cl_n$-modules $\mathcal{H}^{<\lambda}\to U_\lambda$ determine class $[\mathcal{H}^{<\lambda}]\in \KO^n(U_\lambda)$ for each $\lambda\in \R_{>0}$. On overlaps $U_\lambda\bigcap U_{\mu}$ for $\lambda<\mu$ there are canonical inclusions
\beq
g_{\lambda\mu}\colon \mathcal{H}^{<\lambda}\hookrightarrow \mathcal{H}^{<\mu}\label{eq:compatibilitydata}
\eeq
and the orthogonal complement to the image of $g_{\lambda\mu}$ has a $\Cl_{n}\otimes \Cl_{-1}$-action extending the $\Cl_n$-action, determined by the invertible odd endomorphisms
\beq\label{eq:compatibilitydata2}
e^{\lambda\mu}|_{\mathcal{H}_b^{\nu}}= \frac{1}{\sqrt{\nu}} \A^{[0]}|_{\mathcal{H}_b^{\nu}},\quad \mathcal{H}_b^{\nu}\subset (\mathcal{H}^{<\lambda})^\perp\subset \mathcal{H}^{<\mu}
\eeq
that define the action of the generator $e\in \Cl_{-1}$. The inclusions $g_{\lambda\mu}$ and odd endomorphisms $e^{\lambda\mu}$ satisfy a compatibility property on triple intersections. By Corollary~\ref{eq:cormaptoKO}, we obtain a map to the $n$th space in the $\KO$-spectrum~\eqref{eq:KOindexspec}.
A refinement of open covers associated with $\Lambda\subset \Lambda'\subset \R_{>0}$, determines a continuous functor between the associated \v{C}ech categories that after realization is a homotopy equivalence. This leads to a pair of functors~\eqref{eq:cormaptoKO} and restriction of the operators~\eqref{eq:compatibilitydata2} furnish a natural isomorphism between these functors. By \cite[Proposition~2.1]{Segalclassifying}, the resulting maps $\Ind_\Lambda(\A),\Ind_{\Lambda'}(\A)\colon M\to \KO^{-n}$ are homotopic. 
\ep

The connections $\nabla^{<\lambda}$ on $\mathcal{H}^{<\lambda}$ determine a \v{C}ech--de~Rham cocycle representing the Pontryagin character of the class associated with~\eqref{eq:KOindexspec} that we denote by
\beq\label{eq:appencdR}
&&\Ch(\Ind_\Lambda(\A))=[U_\lambda,\Ch(\nabla^{<\lambda})]\in C^*(\{U_\lambda\},\Omega^\bullet).
\eeq
%(see Bott and Tu page 90 or so for more notation). On an $l$-fold overlap for $\sigma\subset \Lambda$, there are $l$ different bundles with connection on $U_\sigma$. Choosing the largest $\lambda$, and some additional connection datum on the orthogonal complements (that is unique up to contractible choice) we get an $l$-parameter family of connections on $\mathcal{H}^{<\lambda}$, and then this gives a connection on $\mathcal{H}^{<\lambda}$ over $U_\sigma\times \R^l$. Integrals over cubes in $\R^l$ give all the (higher) Chern--Simons form compatibilities we need.

\begin{lem}\label{lem:anotherdamnhomotopy}
Let $\A$ be a self-adjoint Clifford linear superconnection that admits smooth index bundles. The differential forms~\eqref{eq:improperintegral} witness a coboundary in the \v{C}ech--de~Rham complex between~\eqref{eq:appencdR} and the (globally defined) cocycle $\Ch(\A)$. Hence, 
\beq\label{eq:superconnindex}
\widehat{\Ind}_\Lambda(\A)=[\{\mathcal{H}^{<\lambda}\to U_\lambda,\nabla^{<\lambda},\eta_\lambda,g_{\lambda\mu},e^{\lambda\mu}\}]\in \dKO^n(M)
\eeq
exacts a differential cocycle from the superconnection $\A$.
\end{lem}
\bp
The \v{C}ech--de~Rham coboundary follows from~\eqref{eq:satisfies}. The differential cocycle~\eqref{eq:superconnindex} uses the description of differential $\KO$-theory from~\cite[\S7]{DBEBis}.
\ep

The Bismut superconnection is the main example of a superconnection that admits smooth index bundles, as we review below. 

\begin{defn} A \emph{spin structure} for a rank $n$ real, oriented, metrized vector bundle $V\to M$ is a $\Spin(n)$-principal bundle $P_V\to M$ equipped with a double cover $P_V\to {\rm Fr}(V)$ of the oriented frame bundle of $V$ that is equivariant relative to $\Spin(n)\to \SO(n)$. 
\end{defn}

The collection of spin structures on $V$ forms a category (possibly empty) whose morphisms are isomorphisms of principal bundles compatible with the map of the frame bundle of $V$. For a family $\pi\colon X\to M$ of Riemannian manifolds, a spin structure on the family is a spin structure on the vertical tangent bundle $T(X/M)=\ker(d\pi)\subset TM$. We will require two types of bundles of Clifford modules built out of a spin structure. 

\begin{defn}
For $\pi\colon X\to M$ a bundle of Riemannian manifolds for which $T(X/M)$ has a spin structure, the \emph{spinor bundle} is 
\beq\label{eq:spinorsdefn}\label{eq:familiespsinor}
\bS_{X/M}=P_{T(X/M)} \times_{\Spin(n)} \Cl_{n}.
\eeq
The sections of $\bS_{X/M}$ support an $M$-family of Dirac operators $\slashed{D}$ that graded commute with the left $\Cl_{-n}=\Cl_n^\op$ action. 
\end{defn}

Using the work of Quillen~\cite{Quillensuperconn} and Bismut~\cite{Bismutindex}, the index bundle data~\eqref{eq:bundleofClnmodules},~\eqref{eq:compatibilitydata}, and~\eqref{eq:compatibilitydata2} may be refined to a differential KO-class. There is a unitary $\Cl_{-n}$-linear connection $\nabla^{\pi_*(\bS\otimes V)}$ on $\pi_*(\bS\otimes V)$, constructed from the Levi-Civita connection on $\bS$, the connection on $V$, and the mean curvature of the fibers of $\pi\colon X\to M$~\cite[Proposition~9.13]{BGV}. The \emph{Bismut superconnection} is defined as~\cite[Proposition~10.15]{BGV}
\beq
\A=\slashed{D}\otimes V+\nabla^{\pi_*(\bS\otimes V)}-\frac{c(T)}{4}\label{eq:Bismutsuper}
\eeq
where $c(T)$ is Clifford multiplication by the curvature 2-form $T$ of the bundle $\pi\colon X\to M$. This incarnation of the Bismut superconnection is real and $\Cl_{-n}$-linear: each of the three terms on the right of~\eqref{eq:Bismutsuper} independently commutes with the $\Cl_{-n}$-action (the operator $c(T)$ uses the $\Cl(T(X/M))$-action, which commutes with the $\Cl_{-n}$-action). 

\begin{prop}\label{prop:Bismutcutoffs}
Any finite-rank self-adjoint Clifford linear superconnnection admits cutoffs. Self-adjoint Clifford linear superconnections adapted to a family of Clifford linear Dirac operators admit cutoffs. 
\end{prop}
\bp
The finite rank case is~\cite[Theorem~9.7]{BGV}. For superconnections adapted to Dirac operators (without Clifford actions), this follows from~\cite[Proposition~9.10 and Theorem~9.26]{BGV}, see also \cite[\S7.2]{FreedLott}. The straightforward extension of these results to Clifford linear Dirac operators is~\cite[\S8.2-8.3]{DBEBis}.
\ep

\bibliographystyle{amsalpha}
\bibliography{references}

\end{document}